\documentclass{m2an_RN}
\usepackage{amssymb,amsmath,amsthm,epsfig,verbatim,pstricks,url}
\usepackage{tikz,pgfplots,tikz-3dplot}
\usetikzlibrary{fit}
\newcommand\addvmargin[1]{
\node[fit=(current bounding box),inner ysep=#1,inner xsep=0]{};
}
\usepackage{ifpdf}  
\newtheorem{ass}[thrm]{Assumption}
\renewcommand{\theequation}{\arabic{section}.\arabic{equation}}
\newcommand{\bRplus}{{\mathbb R}_{>0}}
\newcommand{\bRgeq}{{\mathbb R}_{\geq 0}}
\newcommand{\RZ}{{\mathbb R} \slash {\mathbb Z}}
\newcommand{\bS}{{\mathbb S}}
\newcommand{\bR}{{\mathbb R}}

\newcommand{\bZ}{\mathbb{Z}}
\newcommand{\bH}{{\mathbb H}}
\newcommand{\spa}{\operatorname{span}}
\newcommand{\torsion}{\mathfrak{t}}
\newcommand{\doctorkappa}{\mathfrak{K}}
\newcommand{\doctorZ}{\mathfrak{Z}}
\newcommand{\dH}[1]{\;{\rm d}{\mathcal{H}}^{#1}} 

\newcommand{\drho}{\;{\rm d}\rho}
\newcommand{\spont}{{\overline\varkappa}}

\newcommand{\Vh}{\underline{V}^h}
\newcommand{\Vhzero}{\underline{V}^h_0}

\newcommand{\Vhpartialzero}{\underline{V}^h_{\partial_0}}

\newcommand{\Vpartialzero}{\underline{V}_{\partial_0}}
\newcommand{\Whpartialzero}{W^h_{\partial_0}}
\newcommand{\Whpartialzerob}{W^h_{(\partial_0)}}

\newcommand{\nabS}{\nabla_{\!\mathcal{S}}}
\newcommand{\Id}{\rm Id}
\newcommand{\id}{\rm id}
\newcommand{\deldel}[1]{\frac{\delta}{{\delta}#1}}
\newcommand{\dd}[1]{\frac{\rm d}{{\rm d}#1}}
\newcommand{\ddt}{\dd{t}}
\newcommand{\xspace}{\mathbb {X}}
\newcommand{\xspaceh}{\mathbb {X}^h}
\newcommand{\yspace}{\mathbb {Y}}
\newcommand{\yspaceh}{\mathbb {Y}^h}

\newcommand{\unitn}{\vec{\rm n}}

\newcommand{\ek}{e}
\newcommand{\jf}{j_1}
\newcommand{\BGNpwf}{\mathcal{P}}
\newcommand{\BGNpwfwf}{\mathcal{P}_{\mathcal{S}}}
\newcommand{\Gauss}{{\mathcal{K}}}
\newcommand{\Domain}{{\mathcal{D}}}
\newcommand{\AADE}{{\mathcal{A}}}
\newcommand{\ttau}{\Delta t}
\def\epsilon{\varepsilon} 

\newcommand{\mat}[1]{\underline{\underline{#1}}\rule{0pt}{0pt}}

\newcommand{\errorXx}{\|\Gamma - \Gamma^h\|_{L^\infty}}
\newcommand{\errorXxLL}{\|\Gamma - \Gamma^h\|_{L^\infty(L^2)}}
\allowdisplaybreaks

\begin{document}
\title{
Stable approximations for axisymmetric
{W}illmore flow for closed and open surfaces
}
\thanks{John died on 30 June 2019, when this manuscript was nearly completed.
We dedicate this article to his memory.}

\author{John W. Barrett} 
\address{Department of Mathematics, Imperial College, London, SW7 2AZ, U.K.}
\author{Harald Garcke}
\address{Fakult{\"a}t f{\"u}r Mathematik, Universit{\"a}t Regensburg, 
93040 Regensburg, Germany}
\author{Robert N\"urnberg} 
\address{Department of Mathematics, University of Trento, Trento, Italy}

\begin{abstract}
For a hypersurface in $\bR^3$, Willmore flow is defined as the $L^2$--gradient
flow of the classical Willmore energy: the integral of the squared
mean curvature. This geometric evolution law is of interest in differential
geometry, image reconstruction and mathematical biology. In this paper, we
propose novel numerical approximations for the Willmore flow of axisymmetric
hypersurfaces. For the semidiscrete continuous-in-time variants we prove a
stability result. 
We consider both closed surfaces, and surfaces with a boundary. In the latter
case, we carefully derive 
weak formulations of
suitable boundary conditions. Furthermore, we
consider many generalizations of the classical Willmore energy, particularly
those that play a role in the study of biomembranes. In the generalized
models we include
spontaneous curvature and area difference elasticity (ADE) effects, 
Gaussian curvature and line energy contributions.
Several numerical experiments demonstrate the efficiency and robustness of our
developed numerical methods.
\end{abstract} 



\maketitle

\renewcommand{\thefootnote}{\arabic{footnote}}

\setcounter{equation}{0}
\section{Introduction} \label{sec:1}

Geometric functionals involving the principal curvatures of a
two-dimensional surface play an important role in mechanics, geometry,
imaging and biology. In plate and shell theories such functionals go
back to the work of \cite{Poisson1814}, \cite{Germain1821} and
\cite{Kirchhoff1850}. In geometry, an energy given by the integrated
square of the mean curvature has been studied intensively since the
pioneering work of \cite{Willmore93}. Especially variational problems
are of interest and the famous Willmore conjecture, which states that
the minimizer among genus $1$ surfaces is given by the Clifford torus,
was only solved recently by \cite{MarquesN14}. In imaging, 
boundary value problems involving
the Willmore functional have been used in problems related to 
image inpainting and surface restoration, see \cite{ClarenzDDRR04} and
\cite{BobenkoS05}. In the theory of biological membranes and vesicles, 
side constraints on surface area and enclosed volume,
and more general curvature functionals play a role. 
In a work of \cite{Canham70} a
possible explanation of the shape of the human red blood was given
using a curvature functional together with volume and area
constraints. In this approach the membrane is modeled as a
two-dimensional surface. Later \cite{Helfrich73}, in a seminal paper, 
introduced the energy
\begin{equation} \label{eq:helfrich}
\tfrac{\alpha}{2}\, \int_{\mathcal{S}} (k_m - \spont)^2 \dH{2} 
+\alpha_G \int_{\mathcal{S}} k_g \dH{2}\,,
\end{equation}
for a surface $\mathcal{S}$ in $\bR^3$,
where $k_m$ is the mean curvature, $k_g$ is the Gaussian curvature,
$\dH{2}$ stands for integration with respect to the two-dimensional
surface measure and $\alpha,\alpha_G$ are so-called bending
rigidities. The important new ingredient is the term $\spont$, the
so-called
spontaneous curvature,   
which reflects a possible  asymmetry in the
membrane. In biological applications, membranes in equilibrium minimize
(\ref{eq:helfrich}) under volume and area constraints on the surface
$\mathcal{S}$.

The simplest evolution law which decreases the energy (\ref{eq:helfrich}),
and which can be used to obtain minimizers, is the $L^2$--gradient flow
\begin{equation}\label{eq:L2grad}
\mathcal{V}_{\mathcal{S}}  = -\alpha\,\Delta_{\mathcal{S}} \,k_m
+2\,\alpha (k_m - \spont)\,k_g-\tfrac{\alpha}{2}\,(k_m^2 -
\spont^2)\,k_m\,,
\end{equation}
where $\mathcal{V}_{\mathcal{S}}$ is the normal velocity of an
evolving surface $(\mathcal{S}(t))_{t\in[0,T]}$. The above formula
shows that the Gauss curvature term $\alpha_G\int_{\mathcal{S}}
k_g\,\dH{2}$ does not give a contribution to the flow, which is due to
the fact that for closed surfaces $\int_{\mathcal{S}} k_g\,\dH{2}$ is
a topological invariant. For surfaces with boundary, however, the term
$\int_{\mathcal{S}} k_g\,\dH{2}$ enters the $L^2$--gradient flow via
boundary conditions.
The equation (\ref{eq:L2grad}) also shows that the $L^2$--gradient flow
is highly nonlinear and for open surfaces also highly nonlinear
boundary conditions have to be considered.

It is the goal of this paper to introduce variational discretizations
for an axisymmetric formulation of (\ref{eq:L2grad}), and to show
stability estimates in a semi-discrete setting. The main contributions
of this paper are as follows. 
\begin{itemize}
\item Using a Lagrangian calculus we derive two mixed formulations of
  (\ref{eq:L2grad}), which can be used for all boundary conditions which appear
  in geometry and applications.
\item The derivation of continuous-in-time, discrete-in-space formulations
  for which stability bounds can be shown.
\item A proof of an equidistribution property for one of the schemes,
  which relies on an implicit tangential motion of vertices
and leads to a uniform distribution of vertices on the polygonal curve
everywhere where the curve is not locally flat.
  We refer
  to our review article \cite{bgnreview} for more information
  on the background of this tangential motion.
\item For fully discrete variants, existence and uniqueness results are
  shown under mild assumptions.
\item Numerical computations show the efficiency of the approach.
\end{itemize}
To our knowledge, this is the first time that weak formulations 
involving  general boundary conditions are derived in the axisymmetric setting.

To describe earlier literature in more detail let us discuss the geometry
under consideration. We consider the case that $\mathcal{S}(t)$ is an
axisymmetric surface, which is rotationally symmetric with respect to
the $x_2$--axis, see Figure~\ref{fig:sketch}.
Besides the geometry of a closed surface, we allow for
open surfaces, i.e.\ the boundary of the rotationally symmetric
surface can consist of either one or two circles. Altogether four
different topologies can be considered: surfaces of spherical
topology, surfaces of toroidal topology, surfaces suspended between
two rings and surfaces suspended at one ring. We have to compute the
evolution of a curve which then has to be rotated around the
$x_2$--axis. For boundary points on the $x_2$--axis, singular and
degenerate behaviour in the resulting equations appear, which makes the
analytical and numerical treatment difficult. 
At other boundary points, which correspond to a boundary ring of the surface
$\mathcal{S}(t)$, one has to describe further conditions, 
which can be of the following form:
\begin{itemize}
\item Clamped boundary conditions (position and angle fixed at the
  boundary).
\item Navier boundary conditions (position fixed and a natural
  boundary condition involving the mean curvature).
\item Semifree boundary conditions, i.e.\ the boundary is free to move
  on a plane.
\item Free boundary conditions, for which several natural boundary
  conditions have to hold.
\end{itemize}
The motivation behind considering axisymmetric geometries is clear: a
vastly more efficient numerical method, compared to truly
three-dimensional computations. On the other hand, many situations of
interest in practice do have rotational symmetry. Moreover, some
qualitative aspects of the considered evolution equations, or the
impact of certain physical parameters, can often be studied in the
axisymmetric setting. For example, in
\cite[Fig.\ 24]{axisd}, we numerically studied the onset of a
singularity for Willmore flow 
for a surface of genus 0, and we perform an analogous investigation
for genus-1 surfaces in Appendix~\ref{sec:B} of this paper.
The axisymmetric setting is hence very popular
in the (bio-)physics literature and we refer to 
\cite{JulicherS94,CapovillaGS02,TuO-Y03} for a derivation of the equilibrium
equations for axisymmetric open membranes.  

Earlier results on the numerical approximation of geometric
evolution problems in the axisymmetric setting can be found in
\cite{aximcf,axisd}, 
as well as on the surprisingly closely related problem of curve evolutions 
in Riemannian manifolds, see \cite{hypbol,hypbolpwf}.
There appears to be little numerical
analysis for such evolution problems in the literature. In the case of
Willmore flow, we refer to \cite{MayerS02,DeckelnickS10,CoxL15,axisd}
for existing numerical approaches.
The present paper fills the gap left by
\cite{axisd}, where two schemes for axisymmetric Willmore flow of
closed surfaces were considered, for which no stability proofs appear
to exist. The numerical analysis will share certain features with our
earlier works \cite{pwf,pwftj,pwfade,pwfopen,hypbolpwf}, work that has
been critically influenced and inspired by the seminal works
\cite{Dziuk08,DeckelnickD09}. However, as mentioned already above,
the axisymmetry introduces additional difficulties due to degenerate or 
singular coefficients
that have to be taken care of in the analysis and in the numerical treatment.
For numerical approaches of the flow (\ref{eq:L2grad}) for open
membranes without the restriction of axisymmetry, we refer to
\cite{pwfopen} and \cite{WangD08}, where the latter authors use
a phase field approach.

The rest of the paper is organized as follows. In Section~\ref{sec:maths} we
postulate the mathematical problems we would like to consider, and in
Section~\ref{sec:weak} we state suitable weak formulations for the evolution
problems. Based on these weak formulations we introduce two types of
semidiscrete schemes in Section~\ref{sec:sd} and prove their stability.
We also consider approximations for the area and volume
preserving variants of Willmore flow.
The corresponding fully discrete schemes are presented in Section~\ref{sec:fd}, 
and numerical results are shown in Section~\ref{sec:nr}. 
Finally, in Appendix~\ref{sec:A} we prove the consistency
of the weak formulations introduced in Section~\ref{sec:weak}, including
the considered boundary conditions, while in
Appendix~\ref{sec:B} we present numerical evidence for the onset of
a singularity for Willmore flow of genus-1 surfaces.

\setcounter{equation}{0}
\section{Mathematical formulation} \label{sec:maths}
\subsection{Generating curve}
\begin{figure}
\center
\newcommand{\AxisRotator}[1][rotate=0]{%
    \tikz [x=0.25cm,y=0.60cm,line width=.2ex,-stealth,#1] \draw (0,0) arc (-150:150:1 and 1);%
}
\begin{tikzpicture}[every plot/.append style={very thick}, scale = 1]
\begin{axis}[axis equal,axis line style=thick,axis lines=center, xtick style ={draw=none}, 
ytick style ={draw=none}, xticklabels = {}, 
yticklabels = {}, 
xmin=-0.2, xmax = 0.8, ymin = -0.4, ymax = 2.55]
after end axis/.code={  
   \node at (axis cs:0.0,1.5) {\AxisRotator[rotate=-90]};
   \draw[blue,->,line width=2pt] (axis cs:0,0) -- (axis cs:0.5,0);
   \draw[blue,->,line width=2pt] (axis cs:0,0) -- (axis cs:0,0.5);
   \node[blue] at (axis cs:0.5,-0.2){$\vec\ek_1$};
   \node[blue] at (axis cs:-0.2,0.5){$\vec\ek_2$};
   \draw[red,very thick] (axis cs: 0,0.7) arc[radius = 70, start angle= -90, end angle= 90];
   \node[red] at (axis cs:0.7,1.9){$\Gamma$};
}
\end{axis}
\end{tikzpicture} \qquad \qquad
\tdplotsetmaincoords{120}{50}
\begin{tikzpicture}[scale=2, tdplot_main_coords,axis/.style={->},thick]
\draw[axis] (-1, 0, 0) -- (1, 0, 0);
\draw[axis] (0, -1, 0) -- (0, 1, 0);
\draw[axis] (0, 0, -0.2) -- (0, 0, 2.7);
\draw[blue,->,line width=2pt] (0,0,0) -- (0,0.5,0) node [below] {$\vec\ek_1$};
\draw[blue,->,line width=2pt] (0,0,0) -- (0,0.0,0.5);
\draw[blue,->,line width=2pt] (0,0,0) -- (0.5,0.0,0);
\node[blue] at (0.2,0.4,0.1){$\vec\ek_3$};
\node[blue] at (0,-0.2,0.3){$\vec\ek_2$};
\node[red] at (0.7,0,1.9){$\mathcal{S}$};
\node at (0.0,0.0,2.4) {\AxisRotator[rotate=-90]};

\tdplottransformmainscreen{0}{0}{1.4}
\shade[tdplot_screen_coords, ball color = red] (\tdplotresx,\tdplotresy) circle (0.7);
\end{tikzpicture}
\caption{Sketch of $\Gamma$ and $\mathcal{S}$, as well as 
the unit vectors $\vec\ek_1$, $\vec\ek_2$ and $\vec\ek_3$.}
\label{fig:sketch}
\end{figure}
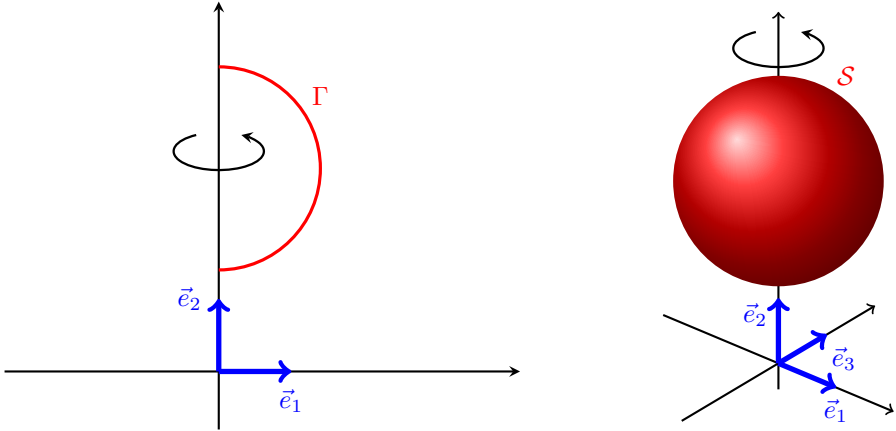

Let $\RZ$ be the periodic interval $[0,1]$, and set
\[
I = \RZ\,, \text{ with } \partial I = \emptyset\,,\quad \text{or}\quad
I = (0,1)\,, \text{ with } \partial I = \{0,1\}\,.
\]
We consider the axisymmetric situation, where 
$\vec x(\cdot,t) : \overline{I} \to \bRgeq \times \bR$ 
is a parameterization of $\Gamma(t)$. 
Throughout $\Gamma(t)$ represents the generating curve of a
surface $\mathcal{S}(t)$ 
that is axisymmetric with respect to the $x_2$--axis, see
Figure~\ref{fig:sketch}. In particular, on defining
$\vec\Pi_3^3(r, z, \theta) = 
(r\,\cos\theta, z, r\,\sin\theta)^T$ for $r\in \bRgeq$, $z \in \bR$, 
$\theta \in [0,2\,\pi]$ and
$\Pi_2^3(r, z) = \{\vec\Pi_3^3(r, z, \theta) : \theta \in [0,2\,\pi)\}$,
we have that
\begin{equation} \label{eq:SGamma}
\mathcal{S}(t) = 
\bigcup_{(r,z)^T \in \Gamma(t)} \Pi_2^3(r, z)
= \bigcup_{\rho \in \overline{I}} \Pi_2^3(\vec x(\rho,t))\,.
\end{equation}
Here we allow $\Gamma(t)$ to be either a closed curve, parameterized over
$\RZ$, which corresponds to $\mathcal{S}(t)$ being a genus-1 surface
without boundary.
Or $\Gamma(t)$ may be an open curve, parameterized over $[0,1]$.
If both ends of $\Gamma(t)$ are attached to the $x_2$--axis, 
then $\mathcal{S}(t)$ is a genus-0 surface without boundary.
If only one end of $\Gamma(t)$ is attached to the $x_2$--axis, 
then $\mathcal{S}(t)$ is an open surface with boundary, where the boundary
consists of a single connected component.
If no endpoint of $\Gamma(t)$ is attached to the $x_2$--axis, 
then $\mathcal{S}(t)$ is an open surface with boundary, where the boundary
consists of two connected components. 
On the boundary we either prescribe clamped boundary conditions, or Navier
boundary conditions, or semifree boundary conditions, or
free boundary conditions. For clamped and Navier boundary conditions, the
boundary point is fixed in space, while for the semifree boundary conditions
the boundary point is allowed to move on a line parallel to one of the two 
axes. As the name for the free boundary condition suggests, the endpoint is
free to move in space. In order to define the different boundary conditions, we
let $\partial_0 I \cup \partial_C I \cup \partial_N I \cup \partial_{1} I
\cup \partial_{2} I \cup \partial_F I$ 
be a disjoint partitioning of $\partial I$, 
with $\partial_0 I$ denoting the subset of boundary points 
of $I$ that correspond to endpoints of $\Gamma(t)$ attached to the $x_2$--axis.
Moreover, $\partial_C I$, $\partial_N I$, 
$\partial_{SF} I = \partial_{1} I \cup \partial_{2} I$ and 
$\partial_F I$ correspond to clamped, Navier, semifree and free boundary
conditions, respectively. 
See Table~\ref{tab:diagram} for a visualization of the different types of 
boundary nodes.
\begin{table}
\center
\caption{The different types of boundary nodes enforced by 
(\ref{eq:axibc})--(\ref{eq:freeslipbc}), and their effect on the possible
movement of the boundary circles $\partial\mathcal{S}$. 
Here the boundary circles in $\bR^3$ 
are shown with the help of an oblique projection.
}
\begin{tabular}{ccc}
\hline
$\partial I$ & $\partial \Gamma$ & $\partial\mathcal{S}$ \\ \hline
$\partial_0 I$ &
\begin{tikzpicture}[scale=0.5,baseline=40]
\begin{axis}[axis equal,axis line style=thick,axis lines=center, 
xtick style ={draw=none}, ytick style ={draw=none}, xticklabels = {}, 
yticklabels = {}, xmin=-0.1, xmax = 2, ymin = -2, ymax = 2]
\addplot[mark=*,color=blue,mark size=6pt] coordinates {(0,1)};
\draw[<->,line width=3pt,color=red] (axis cs:0.3,0.5) -- (axis cs:0.3,1.5);
\node at (axis cs:0.5,-0.3){\Large$\vec\ek_1$};
\node at (axis cs:-0.3,0.5){\Large$\vec\ek_2$};
\end{axis}
\addvmargin{1mm}
\end{tikzpicture} 
& N/A \\ 
$\partial_C I \cup \partial_N I$ &
\begin{tikzpicture}[scale=0.5,baseline=40]
\begin{axis}[axis equal,axis line style=thick,axis lines=center, 
xtick style ={draw=none}, ytick style ={draw=none}, xticklabels = {}, 
yticklabels = {}, xmin=-0.1, xmax = 2, ymin = -2, ymax = 2]
\addplot[mark=*,color=blue,mark size=6pt] coordinates {(2,1)};
\node at (axis cs:0.5,-0.3){\Large$\vec\ek_1$};
\node at (axis cs:-0.3,0.5){\Large$\vec\ek_2$};
\end{axis}
\addvmargin{1mm}
\end{tikzpicture} 
& 
\begin{tikzpicture}[baseline=0]
\draw[color=blue,thick] (0,0) circle [x radius=2cm, y radius=1cm];
\addvmargin{1mm}
\end{tikzpicture} 
\\ 
$\partial_1 I$ &
\begin{tikzpicture}[scale=0.5,baseline=40]
\begin{axis}[axis equal,axis line style=thick,axis lines=center, 
xtick style ={draw=none}, ytick style ={draw=none}, xticklabels = {}, 
yticklabels = {}, xmin=-0.1, xmax = 2, ymin = -2, ymax = 2]
\addplot[mark=*,color=blue,mark size=6pt] coordinates {(2,1)};
\draw[<->,line width=3pt,color=red] (axis cs:2.3,0.5) -- (axis cs:2.3,1.5);
\draw[thick,color=blue] (axis cs:2,-2) -- (axis cs:2,2);
\node at (axis cs:0.5,-0.3){\Large$\vec\ek_1$};
\node at (axis cs:-0.3,0.5){\Large$\vec\ek_2$};
\end{axis}
\addvmargin{1mm}
\end{tikzpicture} 
& 
\begin{tikzpicture}[baseline=0]
\draw[color=blue,thick] (0,0) circle [x radius=2cm, y radius=1cm];
\draw[<->,color=red,line width=2pt] (2.2,-0.5) -- (2.2,0.5);
\draw[color=blue,thin] (2,-1) -- (2,1);
\draw[color=blue,thin] (-2,-1) -- (-2,1);
\addvmargin{1mm}
\end{tikzpicture} 
\\ 
$\partial_2 I$ &
\begin{tikzpicture}[scale=0.5,baseline=40]
\begin{axis}[axis equal,axis line style=thick,axis lines=center, 
xtick style ={draw=none}, ytick style ={draw=none}, xticklabels = {}, 
yticklabels = {}, xmin=-0.1, xmax = 2, ymin = -2, ymax = 2]
\addplot[mark=*,color=blue,mark size=6pt] coordinates {(2,1)};
\draw[<->,line width=3pt,color=red] (axis cs:1.5,0.7) -- (axis cs:2.5,0.7);
\draw[thick,color=blue] (axis cs:0,1) -- (axis cs:4,1);
\node at (axis cs:0.5,-0.3){\Large$\vec\ek_1$};
\node at (axis cs:-0.3,0.5){\Large$\vec\ek_2$};
\end{axis}
\addvmargin{1mm}
\end{tikzpicture} 
& 
\begin{tikzpicture}[baseline=0]
\draw[color=blue,thick] (0,0) circle [x radius=2cm, y radius=1cm];
\draw[color=blue,thin] (0,0) circle [x radius=1.5cm, y radius=0.75cm];
\draw[color=blue,thin] (0,0) circle [x radius=2.5cm, y radius=1.25cm];
\draw[<->,color=red,line width=2pt] (1.5,0) -- (2.5,0);
\draw[<->,color=red,line width=2pt] (-2.5,0) -- (-1.5,0);
\addvmargin{1mm}
\end{tikzpicture} \\ \hline
\end{tabular}
\label{tab:diagram}
\end{table}%

Hence, we always assume that, for all $t \in [0,T]$,
\begin{subequations} 
\begin{align} 
\vec x(\rho,t) \,.\,\vec\ek_1 > 0 \quad &
\forall\ \rho \in \overline{I}\setminus \partial_0 I\,,
\label{eq:xpos} \\
\vec x(\rho,t) \,.\,\vec\ek_1 = 0 \quad &
\forall\ \rho \in \partial_0 I\,,
\label{eq:axibc} \\
\vec x_t(\rho,t) = \vec 0 \quad &
\forall\ \rho \in \partial_C I \cup \partial_N I\,, \label{eq:noslipbc} \\
\vec x_t(\rho,t) \,.\,\vec\ek_i = 0 \quad &
\forall\ \rho \in \partial_i I \,, \ i =1,2\,.  \label{eq:freeslipbc} 
\end{align}
\end{subequations}
We will discuss the additional boundary conditions to 
(\ref{eq:axibc})--(\ref{eq:freeslipbc}) later in this section.

On assuming that $|\vec x_\rho| \geq c_0 > 0$ in $\overline I \times [0,T]$,
we introduce the arclength $s$ of the curve, i.e.\ $\partial_s =
|\vec x_\rho|^{-1}\,\partial_\rho$, and set
\begin{equation} \label{eq:tau}
\vec\tau(\rho,t) = \vec x_s(\rho,t) = 
\frac{\vec x_\rho(\rho,t)}{|\vec x_\rho(\rho,t)|} \qquad \mbox{and}
\qquad \vec\nu(\rho,t) = -[\vec\tau(\rho,t)]^\perp
\qquad \mbox{in } \overline{I} \,,
\end{equation}
where $(\cdot)^\perp$ denotes a clockwise rotation by $\frac{\pi}{2}$.

On recalling (\ref{eq:SGamma}), we observe that the normal
$\unitn_{\mathcal{S}}$ on $\mathcal{S}(t)$ is given by
\begin{equation} \label{eq:nuS}
\unitn_{\mathcal{S}}(\vec\Pi_3^3(\vec x(\rho,t),\theta)) = 
\begin{pmatrix}
(\vec\nu(\rho,t)\,.\,\vec\ek_1)\,\cos\theta \\
\vec\nu(\rho,t)\,.\,\vec\ek_2 \\
(\vec\nu(\rho,t)\,.\,\vec\ek_1)\,\sin\theta 
\end{pmatrix}
 \quad\text{for}\quad
\rho \in \overline{I}\,,\ t \in [0,T]\,,\ \theta \in [0,2\,\pi)\,.
\end{equation}
Similarly, the normal velocity $\mathcal{V}_{\mathcal{S}}$ of $\mathcal{S}(t)$ 
in the direction $\unitn_{\mathcal{S}}$ is given by
\begin{equation} \label{eq:calVS}
\mathcal{V}_{\mathcal{S}} = \vec x_t(\rho,t)\,.\,\vec\nu(\rho,t) \quad\text{on }
\Pi_2^3(\vec x(\rho,t)) \subset \mathcal{S}(t)\,,
\quad \forall\ \rho \in \overline I\,,\ t \in [0,T]\,.
\end{equation}

For the curvature $\varkappa$ of $\Gamma(t)$ it holds that
\begin{equation} \label{eq:varkappa}
\varkappa\,\vec\nu = \vec\varkappa = \vec\tau_s =
\frac1{|\vec x_\rho|} \left[ \frac{\vec x_\rho}{|\vec x_\rho|} \right]_\rho
\quad \mbox{in }\ \overline I\,.
\end{equation}
We recall that the so-called
mean curvature, i.e.\ the sum of the principal curvatures,
and Gaussian curvature of $\mathcal{S}(t)$ are given by
\begin{equation} \label{eq:meanGaussS}
\varkappa_{\mathcal{S}} = 
\varkappa - 
\frac{\vec\nu\,.\,\vec\ek_1}{\vec x\,.\,\vec\ek_1}
\quad\text{and}\quad
\Gauss_{\mathcal{S}} = -
\varkappa\,\frac{\vec\nu\,.\,\vec\ek_1}{\vec x\,.\,\vec\ek_1}
= \varkappa\,(\varkappa_{\mathcal{S}}-\varkappa)
\quad\text{in }\ \overline{I}\,,
\end{equation}
respectively; see e.g.\ \cite[(2.11)]{axisd}. 
More precisely, if $k_m$ and $k_g$
denote the mean and Gaussian curvatures of $\mathcal{S}(t)$, then
\begin{equation} \label{eq:kmkg}
k_m = \varkappa_{\mathcal{S}}(\rho,t) 
\ \text{ and }\
k_g = \Gauss_{\mathcal{S}}(\rho,t) 
\quad\text{on }
\Pi_2^3(\vec x(\rho,t)) \subset \mathcal{S}(t)\,,
\quad \forall\ \rho \in \overline I\,,\ t \in [0,T]\,.
\end{equation}
In the literature, the two terms making up $\varkappa_{\mathcal{S}}$
in (\ref{eq:meanGaussS}) 
are often referred to as in-plane and azimuthal curvatures,
respectively, with their sum being equal to the mean curvature.
We note that combining (\ref{eq:meanGaussS}) and (\ref{eq:varkappa}) yields 
that
\begin{equation} \label{eq:kappaS}
\varkappa_{\mathcal{S}}\,\vec\nu = 
\frac1{|\vec x_\rho|} \left[ \frac{\vec x_\rho}{|\vec x_\rho|} \right]_\rho
- \frac{\vec\nu\,.\,\vec\ek_1}{\vec x\,.\,\vec\ek_1}\,\vec\nu
\quad \mbox{in }\ \overline I\,.
\end{equation}
Weak formulations of (\ref{eq:varkappa}) and (\ref{eq:kappaS}) will form the basis of our 
approximations for Willmore flow.
Clearly, for a smooth surface with bounded curvatures it follows from
(\ref{eq:meanGaussS}) that 
\begin{equation} \label{eq:bcnu}
\vec\nu(\rho,t) \,.\,\vec\ek_1 = 0
\qquad \forall\ \rho \in \partial_0 I\,,\ t\in[0,T]\,,
\end{equation}
which is clearly equivalent to
\begin{equation} \label{eq:bc}
\vec x_\rho(\rho,t) \,.\,\vec\ek_2 = 0
\qquad \forall\ \rho \in \partial_0 I\,,\ t\in[0,T]\,.
\end{equation}
A precise derivation of (\ref{eq:bc}) in the context of a weak formulation
of (\ref{eq:kappaS}) can be found in \cite[Appendix~A]{aximcf}.
We note that for the singular fraction in \eqref{eq:meanGaussS} 
it follows from \eqref{eq:bc} and (\ref{eq:bcnu}), on recalling 
(\ref{eq:varkappa}), that
\begin{equation}
\lim_{\rho\to \rho_0} 
\frac{\vec\nu(\rho,t)\,.\,\vec\ek_1}{\vec x(\rho,t)\,.\,\vec\ek_1}
= \lim_{\rho\to \rho_0} 
\frac{\vec\nu_\rho(\rho,t)\,.\,\vec\ek_1}{\vec x_\rho(\rho,t)\,.\,\vec\ek_1}
= \vec\nu_s(\rho_0,t)\,.\,\vec\tau(\rho_0,t) 
= -\varkappa(\rho_0,t)
\qquad \forall\ \rho_0\in\partial_0 I\,,\ t \in [0,T]\,.
\label{eq:bclimit}
\end{equation}

\subsection{Willmore flow}

We now define the generalized Willmore energy of the surface $\mathcal{S}(t)$ 
as
\begin{equation} \label{eq:E}
E(t) = \tfrac12\,\alpha \int_{\mathcal{S}(t)} (k_m - \spont)^2 \dH{2} 
= \pi\,\alpha \int_I \vec x\,.\,\vec\ek_1\,
(\varkappa_{\mathcal{S}} -\spont)^2 \,|\vec x_\rho| \drho\,,
\end{equation}
where we have recalled (\ref{eq:kmkg}); 
see also \cite[(6),(7)]{CoxL15}. Here $\alpha\in\bRplus$ and
$\spont\in\bR$ are given constants, with $\spont$ denoting the
so-called spontaneous curvature.
On $\mathcal{S}(t)$, Willmore flow, i.e.\ the 
$L^2$--gradient flow for (\ref{eq:E}), is given by
\begin{equation}\label{eq:Willmore_flow}
\frac1\alpha\,\mathcal{V}_{\mathcal{S}} 
= -\Delta_{\mathcal{S}} \,k_m + 2\, (k_m- \spont)\,k_g
-\tfrac{1}{2}\,(k_m^2 - \spont^2)\,k_m \qquad\text{on }\ \mathcal{S}(t)\,,
\end{equation}
recall (\ref{eq:calVS}) for the definition of $\mathcal{V}_{\mathcal{S}}$, 
see e.g.\ \cite{willmore}.
Here $\Delta_{\mathcal{S}} = \nabS\,.\,\nabS$ 
is the Laplace--Beltrami operator on $\mathcal{S}(t)$.
Associated with (\ref{eq:Willmore_flow}) 
are boundary conditions,
but we will discuss these once we have generalized 
the energy \eqref{eq:E}.

For applications to biomembranes, and for surfaces with boundary, 
the considered energy can be more general than (\ref{eq:E}). 
In particular, surface area constraints, 
area difference elasticity (ADE) effects, Gaussian curvature contributions and line energy
now also play a role. Hence, we adapt (\ref{eq:E}) to
\begin{equation} \label{eq:EEa}
E(t) = \tfrac12\,\alpha
\int_{\mathcal{S}(t)} (k_m - \spont)^2 \dH{2} 
+ \lambda\,\mathcal{H}^2(\mathcal{S}(t))
+ \tfrac{\beta}{2} \,\AADE^2_{\mathcal{S}}(t)
+ \alpha_G\int_{\mathcal{S}(t)} k_g \dH{2}
+ \varsigma\,\mathcal{H}^1(\partial\mathcal{S}(t)) 
\end{equation}
with $\AADE_{\mathcal{S}}(t) = \int_{\mathcal{S}(t)} k_m \dH{2} -M_0$,
and given constants $\beta\in\bRgeq$, $M_0, \alpha_G, \varsigma 
\in \bR$, see e.g.\
\cite{Canham70, Helfrich73,Nitsche93, JulicherL96,%
Seifert97,CapovillaGS02,TuO-Y03,pwfade,pwfopen} for more details. 
In addition, $\lambda \in \bR$ is a constant that can penalise or
encourage surface area growth. If chosen time-dependent, it can act as
a Lagrange multiplier for a surface area constraint. We make this more
explicit later on, see \S\ref{subsec:cons} below.
We remark that the contributions 
$\tfrac12\,\alpha \int_{\mathcal{S}(t)} k_m^2 \dH{2} 
+ \alpha_G \int_{\mathcal{S}(t)} k_g \dH{2}$
to the energy $E(t)$ are positive semidefinite with respect to the principal 
curvatures only if $\alpha_G \in [-2\,\alpha,0]$.
We note that this constraint is likely to have implications for the existence 
and regularity of the corresponding $L^2$--gradient flow.  

Noting the Gauss--Bonnet theorem, see \cite{Kuhnel15}, 
\begin{equation} \label{eq:GB} 
\int_{\mathcal{S}} k_g \dH{2} = 2\,\pi\,m(\mathcal{S}) 
+ \int_{\partial\mathcal{S}} k_{\partial \mathcal{S},\mu} \dH{1} \,,
\end{equation}
where $m(\mathcal{S}) \in \bZ$ denotes the Euler characteristic of 
$\mathcal{S}$ and $k_{\partial \mathcal{S},\mu}$ is the geodesic curvature of
$\partial\mathcal{S}$, we can rewrite (\ref{eq:EEa}) 
as
\begin{equation} \label{eq:EE}
E(t) = \tfrac12\,\alpha \int_{\mathcal{S}(t)} (k_m - \spont)^2 \dH{2} 
+ \lambda\,\mathcal{H}^2(\mathcal{S}(t))
+ \tfrac{\beta}{2} \,\AADE_{\mathcal{S}}^2(t)
+ \alpha_G\left[\int_{\partial\mathcal{S}(t)} 
k_{\partial\mathcal{S},\mu} \dH{1} 
+ 2\,\pi\,m(\mathcal{S}(t)) \right]
+ \varsigma\,\mathcal{H}^1(\partial\mathcal{S}(t))\,.
\end{equation}
In order to define $k_{\partial \mathcal{S},\mu}$, we first define 
the conormal, $\vec\mu_{\partial\mathcal{S}}$, to $\mathcal{S}(t)$ on 
$\partial\mathcal{S}(t)$ to be 
\begin{equation} \label{eq:mupS}
\vec\mu_{\partial\mathcal{S}}= \pm\, \unitn_{\mathcal{S}} \times \vec\id_s
\quad \text{on }\ \partial{\mathcal{S}}(t)\,,
\end{equation}
where $\vec\id$ denotes the identity in $\bR^3$
and $s$ denotes arclength on the curve
$\partial\mathcal{S}(t) \subset \bR^3$, 
so that $\vec\id_s$ is its unit tangent vector. The
sign in (\ref{eq:mupS}) is chosen so that $\vec\mu_{\partial\mathcal{S}}$
points out of $\mathcal{S}(t)$.
It holds that
$\vec\id_{ss} = \vec k_{\partial\mathcal{S}} =
k_{\partial \mathcal{S},\rm n}\,\unitn_{\mathcal{S}} +
k_{\partial \mathcal{S},\mu}\,\vec\mu_{\partial\mathcal{S}}$ 
on $\partial\mathcal{S}(t)$,
where $\vec k_{\partial\mathcal{S}}$ is the curvature vector on
$\partial\mathcal{S}(t)$, and where
$k_{\partial \mathcal{S},\rm n}$ is the normal curvature 
and $k_{\partial \mathcal{S},\mu}$ is the geodesic curvature
of
$\partial\mathcal{S}(t)$. 

Similarly to (\ref{eq:kmkg}), it is easily seen that
\begin{equation} \label{eq:kdSmu}
k_{\partial\mathcal{S},\rm n} = - 
\frac{\vec\nu(\rho,t)\,.\,\vec\ek_1}{\vec x(\rho,t)\,.\,\vec\ek_1}
\text{ and }
k_{\partial\mathcal{S},\mu} = - 
\frac{\vec\mu(\rho,t)\,.\,\vec\ek_1}{\vec x(\rho,t)\,.\,\vec\ek_1}
\quad\text{on }
\Pi_2^3(\vec x(\rho,t)) \subset \partial\mathcal{S}(t)\,,
\quad \forall\ \rho \in \partial I\setminus\partial_0 I\,,\ t \in [0,T]\,,
\end{equation}
where $\vec\nu(\cdot,t)$ is the unit normal on $\Gamma(t)$ 
as defined in (\ref{eq:tau}) and (\ref{eq:nuS}) and 
\begin{equation} \label{eq:mu}
\vec\mu(p,t) = (-1)^{p+1}\,\vec\tau(p,t)
\quad \forall\ p \in \partial I\,,\ t \in [0,T]\,,
\end{equation}
denotes the corresponding conormal 
of $\Gamma(t)$ at the endpoint $\vec x(p,t)$, for $p \in \partial I$.
Here, we have recalled that the conormal $\vec\mu_{\partial\mathcal{S}}$ 
points out of $\mathcal{S}(t)$.

Hence, an energy equivalent to (\ref{eq:EE}), 
for flows of axisymmetric surfaces without topological changes, 
can be written as
\begin{align}
\widetilde E(t) & = E(t) - 2\,\pi\,\alpha_G\,m(\mathcal{S}(t)) 
\nonumber \\ & =
\pi\,\alpha \int_I \vec x\,.\,\vec\ek_1
\left[ \varkappa_{\mathcal{S}}-\spont \right]^2
|\vec x_\rho| \drho
+ { 2\,\pi\,\lambda \int_I (\vec x\,.\,\vec\ek_1)\,|\vec x_\rho| \drho}
+ { \tfrac\beta2\,\AADE_{\mathcal{S}}^2(t)}
\nonumber \\ & \qquad
- 2\,\pi\,\alpha_G\sum_{p\in \partial I \setminus \partial_0 I} \left[
\vec x\,.\,\vec\ek_1\,
\frac{\vec\mu\,.\,\vec\ek_1}{\vec x\,.\,\vec\ek_1}\right](p)
+ 2\,\pi\,\varsigma\sum_{p\in \partial I \setminus \partial_0 I}
\vec x(p)\,.\,\vec\ek_1 \nonumber \\ &
= \pi\,\alpha \int_I \vec x\,.\,\vec\ek_1
\left[ \varkappa_{\mathcal{S}} -\spont \right]^2
|\vec x_\rho| \drho
+ { 2\,\pi\,\lambda \int_I (\vec x\,.\,\vec\ek_1)\,|\vec x_\rho| \drho}
+ { \tfrac\beta2\, \AADE_{\mathcal{S}}^2(t)}
 \nonumber \\ & \qquad
- 2\,\pi\,\alpha_G\sum_{p\in \partial I \setminus \partial_0 I}
\vec\mu(p)\,.\,\vec\ek_1
+ 2\,\pi\,\varsigma\sum_{p\in \partial I \setminus \partial_0 I}
\vec x(p)\,.\,\vec\ek_1\,, 
\label{eq:EEax}
\end{align}
where
\begin{equation} \label{eq:AADE}
\AADE_{\mathcal{S}}(t)
= 2\,\pi \int_I \vec x\,.\,\vec\ek_1\, 
\varkappa_{\mathcal{S}}\,|\vec x_\rho| \drho - M_0\,.
\end{equation}

In this general situation, (\ref{eq:Willmore_flow}) is replaced by
\begin{equation}
\mathcal{V}_{\mathcal{S}} 
= -\alpha\,\Delta_{\mathcal{S}} \,k_m + 2\left[ 
\alpha\,(k_m- \spont) + \beta\, \AADE_{\mathcal{S}}\right] k_g
-\left[\tfrac12\,\alpha\,(k_m^2 - \spont^2) -\lambda \right] k_m
\quad\text{on }\ \mathcal{S}(t)\,,
\label{eq:gradflowlambda}
\end{equation}
see e.g.\ \cite[(1.21)]{pwfopen}. 
A strong formulation for the flow (\ref{eq:gradflowlambda}) on $I$,
on recalling \cite[(B.3)]{axisd}, as well as (\ref{eq:calVS}) and
(\ref{eq:kmkg}), is given by
\begin{equation} \label{eq:xtbgnlambda}
(\vec x\,.\,\vec\ek_1)\,\vec x_t\,.\,\vec\nu = 
- \alpha\,[\vec x\,.\,\vec\ek_1\,(\varkappa_{\mathcal{S}})_{s}]_s
+ 2\,\vec x\,.\,\vec\ek_1 \left[ \alpha\,(\varkappa_{\mathcal{S}} - \spont
) + \beta \, \AADE_{\mathcal{S}}\right]\Gauss_{\mathcal{S}}
- \vec x\,.\,\vec\ek_1
\left[\tfrac12\,\alpha\,(\varkappa_{\mathcal{S}}^2 -\spont^2)
- \lambda \right] \varkappa_{\mathcal{S}}
\quad\text{in }\ I\,.
\end{equation}
Next we discuss the boundary conditions associated with 
(\ref{eq:gradflowlambda}) and (\ref{eq:xtbgnlambda}).

\subsection{Boundary conditions}

We recall from \cite{pwfopen} the following boundary conditions one can
consider for $\mathcal{S}(t)$ on $\partial\mathcal{S}(t)$.
A connected component of the boundary can either move freely, 
or move along the boundary of a fixed domain $\Domain$, or it will be fixed. 
For the latter case two types of boundary conditions arise: 
clamped and Navier. 
Corresponding to (\ref{eq:axibc})--(\ref{eq:freeslipbc}), we now partition
$\partial\mathcal{S}$ into
$\partial_C \mathcal{S} \cup \partial_N \mathcal{S} \cup \partial_1\mathcal{S}
\cup \partial_2\mathcal{S} \cup \partial_F\mathcal{S}$, and we also set
$\partial_{SF} \mathcal{S} = \partial_1\mathcal{S}
\cup \partial_2\mathcal{S}$.
In the free boundary case, the three natural boundary conditions are, for $t \in (0,T]$, given by
\begin{subequations}
\label{eq:free}
\begin{align}
\alpha\,(\nabS\,k_m)\,.\,\vec\mu_{\partial\mathcal{S}}
+\varsigma\,k_{\partial\mathcal{S},\rm n}
& = 0 \quad \text{on }\ \partial_F\mathcal{S}(t)\,,
\label{eq:freea} \\
-\tfrac12\,\alpha\,(k_m - \spont)^2 
-  \beta\,\AADE_{\mathcal{S}}\,k_m +\varsigma\,k_{\partial\mathcal{S},\mu}
- \alpha_G\,k_g & = { \lambda}
\quad \text{on }\ \partial_F\mathcal{S}(t)\,, \label{eq:freeb} \\
\alpha\,(k_m - \spont) + \beta\,\AADE_{\mathcal{S}}
+ \alpha_G\,k_{\partial\mathcal{S},\rm n} & = 0 
\quad \text{on }\ \partial_F\mathcal{S}(t)\,. \label{eq:freec}
\end{align}
\end{subequations}
In general, the term $\alpha_G\,\torsion_s$ features on the right hand side of
(\ref{eq:freea}), where 
$\torsion=-(\unitn_{\mathcal{S}})_s\,.\,\vec\mu_{\partial\mathcal{S}}$
denotes the torsion of 
$\partial\mathcal{S}(t)$, see \cite[(1.15), (1.22)]{pwfopen}. However, in the
axisymmetric case $\partial\mathcal{S}(t)$ is made up of circles, and so the
torsion is zero.
For the semifree case, when 
$\partial_{SF} \mathcal{S}(t) \subset \partial\Domain$ for all $t\in[0,T]$, 
where $\partial\Domain$ is the
boundary of a fixed domain $\Domain \subset \bR^3$,
we let $\partial\Domain$ be given by a function $F\in C^1(\bR^3)$ 
such that
\begin{equation*} 
\partial\Domain = \{ \vec z \in \bR^3 : F(\vec z) = 0 \} \qquad \mbox{and}
\qquad |\nabla F(\vec z)| = 1 \quad \forall\ \vec z\in \partial \Domain\,,
\end{equation*}
and we denote the normal to $\Domain$ on $\partial\Domain$ by
$\unitn_{\Domain} = \nabla F$. For the special axisymmetric setting considered
here, recall (\ref{eq:freeslipbc}), we restrict ourselves to
the two cases 
\begin{equation} \label{eq:normal12}
\unitn_{\Domain} = \unitn_1 = 
\frac{\vec\id - (\vec\id\,.\,\vec\ek_2)\,\vec\ek_2}
{|\vec\id - (\vec\id\,.\,\vec\ek_2)\,\vec\ek_2|}
\quad\mbox{on } \ \partial_1\mathcal{S}(t)\,, \qquad
\unitn_{\Domain} = \unitn_2 = \vec\ek_2 
\quad\mbox{on } \ \partial_2\mathcal{S}(t)\,.
\end{equation}
The semifree boundary conditions are, for $t \in (0,T]$, then
\begin{subequations}
\label{eq:sfree}
\begin{align}
& \partial_{SF}\mathcal{S}(t) \subset \partial\Domain \label{eq:sfreea} \\
&- [-\tfrac12\,\alpha\,(k_m - \spont)^2 
   - \beta\, \AADE_{\mathcal{S}}\,k_m +\varsigma\,k_{\partial\mathcal{S},\mu}
- \alpha_G\,k_g-\lambda]\,(\unitn_{\mathcal{S}}\,.\,\unitn_i) 
+ [\alpha\,(\nabS\,k_m)\,.\,\vec\mu_{\partial\mathcal{S}} 
+\varsigma\,k_{\partial\mathcal{S},\rm n}  
]\,(\vec\mu_{\partial\mathcal{S}}\,.\,\unitn_i) 
 = 0 
\nonumber \\ & \hspace{9cm}
\text{on }\ \partial_i\mathcal{S}(t)\,,\ i=1,2\,,
\label{eq:sfreeb} \\ &
\alpha\,(k_m - \spont) + \beta\,\AADE_{\mathcal{S}}
+ \alpha_G\,k_{\partial\mathcal{S},\rm n}  = 0 
\qquad \qquad \qquad \text{on }\ \partial_{SF}\mathcal{S}(t) =
\partial_1\mathcal{S}(t) \cup \partial_2\mathcal{S}(t) \,. \label{eq:sfreec}
\end{align}
\end{subequations}
Note that compared to \cite[(1.17), (1.22)]{pwfopen}, we have once again omitted the
vanishing torsion term.
Clamped boundary conditions are, for $t \in (0,T]$, given by
\begin{equation} \label{eq:clamped}
\partial_C \mathcal{S}(t) = \partial_C \mathcal{S}(0) \quad\text{and} \quad
\vec\mu_{\partial\mathcal{S}} (t) 
= \vec\zeta_{\mathcal{S}} \quad \text{on }\ \partial_C \mathcal{S}(0)\,,
\end{equation}
where $\vec\zeta_{\mathcal{S}} \in C^{0}(\partial_C\mathcal{S}(0), \bS^2)$,
with $\bS^2 := \{\vec z \in \bR^{3} : |\vec z| = 1\}$,
needs to be axisymmetric.
Similarly, Navier boundary conditions are, for $t \in (0,T]$, given by
\begin{equation} \label{eq:Navier}
\partial_N\mathcal{S}(t) = \partial_N\mathcal{S}(0) \quad\text{and} \quad
\alpha\,(k_m -\spont) +  \beta\, \AADE_{\mathcal{S}} + \alpha_G\,
k_{\partial\mathcal{S},\rm n} =0 
\quad \text{on }\ \partial_N\mathcal{S}(0)\,.
\end{equation}

We now translate the above boundary conditions to the axisymmetric case. 
On noting (\ref{eq:free}), (\ref{eq:mu}), (\ref{eq:kdSmu}) and (\ref{eq:kmkg}),
we obtain, for $t \in (0,T]$,
the free boundary conditions
\begin{subequations}
\label{eq:axifree}
\begin{align}
(-1)^{p+1}\, \alpha\, (\varkappa_{\mathcal{S}})_s
- \varsigma\,\frac{\vec\nu\,.\,\vec\ek_1}{\vec x\,.\,\vec\ek_1} &= 0
\quad\text{on }\  \partial_F I\,, \label{eq:axifreea} \\
 -\tfrac12\,\alpha\,(\varkappa_{\mathcal{S}} - \spont)^2 
- \beta\,\AADE_{\mathcal{S}}\,\varkappa_{\mathcal{S}} -
\varsigma\,\frac{\vec\mu\,.\,\vec\ek_1}{\vec x\,.\,\vec\ek_1}
- \alpha_G\,\Gauss_{\mathcal{S}} &= \lambda
\quad\text{on } \ \partial_F I\,, \label{eq:axifreeb} \\
 \alpha\,(\varkappa_{\mathcal{S}} - \spont) + \beta\,\AADE_{\mathcal{S}}
- \alpha_G\,\frac{\vec\nu\,.\,\vec\ek_1}{\vec x\,.\,\vec\ek_1} &= 0
\quad\text{on } \ \partial_F I\,.
\label{eq:axifreec}
\end{align}
\end{subequations}
Similarly, (\ref{eq:sfree}), (\ref{eq:normal12}),
(\ref{eq:mu}), (\ref{eq:kdSmu}) and (\ref{eq:kmkg})
yield, for $t \in (0,T]$, the semifree boundary conditions 
\begin{subequations} \label{eq:axisfree}
\begin{align}
&\vec x(\cdot,t)\,.\,\vec\ek_i = \vec x(\cdot,0)\,.\,\vec\ek_i 
\quad\text{on } \ \partial_i I\,, \; i = 1,2\,,\label{eq:axisfreea} \\
&- \left[-\tfrac12\,\alpha\,(\varkappa_{\mathcal{S}} - \spont)^2 
- \beta\,\AADE_{\mathcal{S}}\,\varkappa_{\mathcal{S}} -
\varsigma\,\frac{\vec\mu\,.\,\vec\ek_1}{\vec x\,.\,\vec\ek_1}
- \alpha_G\,\Gauss_{\mathcal{S}}
-\lambda \right]\,(\vec\nu\,.\,\vec\ek_i) 
+ \left[(-1)^{p+1}\, \alpha\, (\varkappa_{\mathcal{S}})_s
- \varsigma\,\frac{\vec\nu\,.\,\vec\ek_1}{\vec x\,.\,\vec\ek_1}
\right]\,(\vec\mu\,.\,\vec\ek_i) 
 = 0 
\nonumber \\ & \hspace{9cm} 
\text{on }\ \partial_i I\,,\ i=1,2\,,
\label{eq:axisfreeb} 
 \\
& \alpha\,(\varkappa_{\mathcal{S}} - \spont) + \beta\,\AADE_{\mathcal{S}}
- \alpha_G\,\frac{\vec\nu\,.\,\vec\ek_1}{\vec x\,.\,\vec\ek_1} = 0 
\quad \text{on } \ \partial_{SF} I =\partial_1 I \cup \partial_2 I\,.
\label{eq:axisfreec}
\end{align}
\end{subequations}
Taking into account the clamped boundary conditions (\ref{eq:clamped}), 
we define, similarly to (\ref{eq:nuS}),
\begin{equation} \label{eq:veczetaS}
\vec\zeta_{\mathcal{S}}(\vec\Pi_3^3(\vec x(p,0),\theta)) = 
\begin{pmatrix}
\vec\zeta(p)\,.\,\vec\ek_1\,\cos\theta \\
\vec\zeta(p)\,.\,\vec\ek_2 \\
\vec\zeta(p)\,.\,\vec\ek_1\,\sin\theta 
\end{pmatrix}
 \quad\text{for}\quad 
p \in \partial_C I\,,
\quad \theta \in [0,2\,\pi)\,,
\end{equation}
to be the conormal of $\mathcal{S}(t)$ on
$\partial_C\mathcal{S}(0) = \partial_C\mathcal{S}(t)$. 
Here $\vec\zeta(p)$, for $p \in \partial_C I$,
are given unit vectors that prescribe the clamping direction for the conormals 
$\vec\mu(p,t)$ of $\Gamma(t)$ at the endpoints $\vec x(p,t) = \vec x(p,0)$, 
for $p \in \partial_C I$. In particular, for clamped boundary conditions
we will enforce, for $t \in (0,T]$, on recalling (\ref{eq:mu}),
\begin{subequations}
\label{eq:CL}
\begin{alignat}{2} 
\vec x(p,t) &= \vec x(p,0) &&\quad \text{for}\quad p \in \partial_C I\,,
\label{eq:fixedC} \\
\label{eq:clampedI}
(-1)^{p+1}\,\vec\tau(p,t) = 
\vec\mu(p,t) &= \vec\zeta(p) &&\quad \text{for}\quad p \in \partial_C I\,.
\end{alignat}
\end{subequations}
Similarly, for Navier boundary conditions
we will enforce, for $t \in (0,T]$, on recalling (\ref{eq:Navier}), (\ref{eq:kmkg}) and 
(\ref{eq:kdSmu}),
\begin{subequations}
\label{eq:Nav}
\begin{alignat}{2} 
\vec x(p,t) &= \vec x(p,0) &&\quad \text{for}\quad p \in \partial_N I\,,
\label{eq:fixedN} \\
\label{eq:NavierI} 
\alpha\,(\varkappa_{\mathcal{S}} - \spont) + \beta\,\AADE_{\mathcal{S}}
- \alpha_G\,\frac{\vec\nu\,.\,\vec\ek_1}{\vec x\,.\,\vec\ek_1} &= 0  
&&\quad \text{on }\ \partial_N I\,.
\end{alignat}
\end{subequations}
Finally, we impose the following boundary conditions on $\partial_0 I$  
\begin{subequations}
\label{eq:part0I}
\begin{alignat}{2}
\vec x\,.\,\vec\ek_1 &= 0 \quad &&\text{on }\ 
 \partial_0 I, \quad \mbox{for } t \in [0,T]\,,
\label{eq:fixed0} \\
\vec x_\rho\,.\,\vec\ek_2 & = 0 \quad &&\text{on }\ 
 \partial_0 I, \quad \mbox{for } t \in [0,T]\,, \label{eq:bcbc} \\
(\varkappa_{\mathcal{S}})_\rho &= 0 \quad &&\text{on }\ 
 \partial_0 I, \quad \mbox{for } t \in (0,T]  \,.
\label{eq:sdbca} 
\end{alignat}
\end{subequations}
Here (\ref{eq:sdbca}) ensures that the radially
symmetric function $k_m$ on $\mathcal{S}(t)$ induced by 
$\varkappa_{\mathcal S}$, recall \eqref{eq:kmkg},
is differentiable, while \eqref{eq:bcbc} is the same as \eqref{eq:bc}. 
It is natural to ask for differentiability of $k_m$ due to the regularisation
property of parabolic equations. We also remark that 
$\Delta_{\mathcal{S}} \,k_m$ in \eqref{eq:L2grad}, for radially symmetric 
solutions, is only defined even in the weak sense if (\ref{eq:sdbca}) holds.
We note that the boundary conditions (\ref{eq:axibc})--(\ref{eq:freeslipbc})
are incorporated in (\ref{eq:fixed0}), (\ref{eq:fixedC}), (\ref{eq:fixedN})
and (\ref{eq:axisfreea}), respectively.  

We now introduce an energy equivalent to (\ref{eq:EEax}), which takes into
account the clamped,  (\ref{eq:CL}), Navier, (\ref{eq:fixedN}), 
and semifree, (\ref{eq:axisfreea}), boundary conditions, 
\begin{align}
\widehat E(t) & = \widetilde E(t) + 2\,\pi\,\alpha_G\sum_{p\in \partial_C I}
\vec\zeta(p)\,.\,\vec\ek_1
- 2\,\pi\,\varsigma\sum_{p\in \partial_C I \cup \partial_N I \cup \partial_1 I}
\vec x(p)\,.\,\vec\ek_1 
\nonumber \\ & 
= \pi\,\alpha \int_I \vec x\,.\,\vec\ek_1
\left[ \varkappa_{\mathcal{S}} -\spont \right]^2
|\vec x_\rho| \drho
+ { 2\,\pi\,\lambda \int_I (\vec x\,.\,\vec\ek_1)\,|\vec x_\rho| \drho}
+ { \tfrac\beta2\, \AADE_{\mathcal{S}}^2(t)}
 \nonumber \\ & \qquad
- 2\,\pi\,\alpha_G\sum_{p\in \partial_M I }
\vec\mu(p)\,.\,\vec\ek_1
+ 2\,\pi\,\varsigma\sum_{p\in \partial_{2} I \cup \partial_F I}
\vec x(p)\,.\,\vec\ek_1\,,
\label{eq:EEaxCN}
\end{align}
where $\AADE_{\mathcal{S}}$ is defined in (\ref{eq:AADE}) and
\begin{equation} \label{eq:partialM}
\partial_M I = \partial_N I \cup \partial_{SF} I \cup \partial_F I\,.
\end{equation}

\subsection{Conserved flows} \label{subsec:cons}

In a number applications, such as biomembranes,
the $L^2$--gradient flow of (\ref{eq:EEa}) is considered under conservation of 
the total surface area and, in the case of a closed surface, 
conservation of the enclosed volume.
Before we state these variants, we recall the following useful results.
We have, similarly to \eqref{eq:EEax}, that
\begin{equation} \label{eq:dAdt}
\ddt\, \mathcal{H}^2(\mathcal{S}(t))
= \ddt\,2\,\pi \int_I \vec x\,.\,\vec\ek_1\,|\vec x_\rho| \drho
= 2\,\pi \int_I \left[ \vec x_t\,.\,\vec\ek_1\,|\vec x_\rho| +
 (\vec x\,.\,\vec\ek_1) \,(\vec x_t)_\rho\,.\,\vec\tau \right] \drho \,.
\end{equation}
In the case of a closed surface $\mathcal{S}(t)$, we have from (\ref{eq:calVS})
that
\begin{equation} \label{eq:dVdt}
\ddt\,\mathcal{L}^3(\Omega(t)) = 
\int_{\mathcal{S}(t)} \mathcal{V}_{\mathcal{S}} \dH{2}
= 2\,\pi\,\int_I (\vec x\,.\,\vec\ek_1)
\,\vec x_t\,.\,\vec\nu\,|\vec x_\rho|\drho \,,
\end{equation}
where $\mathcal{L}^3$ denotes the Lebesgue measure in $\bR^3$, 
$\mathcal{S}(t) = \partial\Omega(t)$, 
and where 
we assume from now on that
$\unitn_{\mathcal{S}}$ is the outer or inner normal
to $\Omega(t)$ on $\mathcal{S}(t)$, recall (\ref{eq:nuS}) and (\ref{eq:tau}). 

Generalized Helfrich flow is the surface area and
volume conserving variant of \eqref{eq:gradflowlambda}, 
and its strong form can be stated as
\begin{equation}
\mathcal{V}_{\mathcal{S}} = -\alpha\,\Delta_{\mathcal{S}} \,k_m + 2\left[ 
\alpha\,(k_m- \spont) + \beta\, \AADE_{\mathcal{S}}\right] k_g
-\left[\tfrac12\,\alpha\,(k_m^2 - \spont^2) -\lambda \right] k_m
+ \lambda_A\,k_m - \lambda_V 
\quad\text{on }\ \mathcal{S}(t)\,,
\label{eq:Helfrich}
\end{equation}
where $(\lambda_A(t),\lambda_V(t))^T \in \bR^2$ are chosen such that 
\begin{equation} \label{eq:sidedtSAV}
\ddt\,\mathcal{H}^2(\mathcal{S}(t)) = 0 
\,,\qquad
\ddt\,\mathcal{L}^3(\Omega(t)) = 0\,. 
\end{equation}

For axisymmetric surfaces the flow \eqref{eq:Helfrich} with 
\eqref{eq:sidedtSAV} can be equivalently formulated as
\begin{align} \label{eq:xHelfrich}
(\vec x\,.\,\vec\ek_1)\,\vec x_t\,.\,\vec\nu & = 
- \alpha\,[\vec x\,.\,\vec\ek_1\,(\varkappa_{\mathcal{S}})_{s}]_s
+ 2\,\vec x\,.\,\vec\ek_1 \left[ \alpha\,(\varkappa_{\mathcal{S}} - \spont
) + \beta \, \AADE_{\mathcal{S}}\right]\Gauss_{\mathcal{S}}
- \vec x\,.\,\vec\ek_1
\left[\tfrac12\,\alpha\,(\varkappa_{\mathcal{S}}^2 -\spont^2)
- \lambda \right] \varkappa_{\mathcal{S}}
\nonumber \\ & \qquad
+ \lambda_A\,\vec x\,.\,\vec\ek_1\,\varkappa_{\mathcal{S}} 
 -  \lambda_V \,\vec x\,.\,\vec\ek_1
\quad\text{in }\ I\,,
\end{align}
where $(\lambda_A(t),\lambda_V(t))^T \in \bR^2$ are chosen such that 
\begin{equation} \label{eq:xsidedtSAV}
\int_I \left[ \vec x_t\,.\,\vec\ek_1\,|\vec x_\rho| +
 (\vec x\,.\,\vec\ek_1) \,(\vec x_t)_\rho\,.\,\vec\tau \right] \drho = 0\,,
\qquad
\int_I (\vec x\,.\,\vec\ek_1)
\,\vec x_t\,.\,\vec\nu\,|\vec x_\rho|\drho = 0\,,
\end{equation}
where we recall \eqref{eq:xtbgnlambda}, \eqref{eq:dAdt} and \eqref{eq:dVdt}. 

\setcounter{equation}{0}
\section{Weak formulations} \label{sec:weak}

On recalling (\ref{eq:tau}), we have for all 
$\vec a,\,\vec b \in \bR^2$ that
\begin{subequations} 
\begin{align}
&\vec a \,.\,\vec b^\perp= -\vec a^\perp .\,\vec b, 
\label{eq:abperp}\\
&\vec a^\perp 
= (\vec a^\perp\,.\,\vec\tau)\,\vec\tau +(\vec a^\perp\,.\,\vec\nu)\,\vec\nu
= (\vec a^\perp\,.\,\vec\nu^\perp)\,\vec\tau  
  - (\vec a^\perp\,.\,\vec\tau^\perp)\,\vec\nu
= (\vec a\,.\,\vec\nu)\,\vec\tau  
  - (\vec a\,.\,\vec\tau)\,\vec\nu \,.
\label{eq:aperp}
\end{align}
\end{subequations}

We define the first variation of a differentiable quantity $B(\vec x)$, in the
direction $\vec\chi$ as
\begin{equation} \label{eq:A}
\left[\deldel{\vec x}\,B(\vec x)\right](\vec\chi)
= \lim_{\varepsilon \rightarrow 0}
\frac{B(\vec x + \varepsilon\, \vec\chi)-B(\vec x)}{\varepsilon}\,.
\end{equation}
For later use, on noting (\ref{eq:A}) and (\ref{eq:tau}), we observe that 
\begin{subequations} \label{eq:24}
\begin{align}
\left[\deldel{\vec x}\, |\vec x_\rho|\right](\vec\chi)
& = \frac{\vec x_\rho\,.\,\vec\chi_\rho}{|\vec x_\rho|}
= \vec\tau\,.\,\vec\chi_\rho = 
\vec\tau\,.\,\vec\chi_s\, |\vec x_\rho|\,, \label{eq:ddxrho} \\
\left[\deldel{\vec x}\, \vec\tau\right](\vec\chi)
& = \left[\deldel{\vec x}\,\frac{\vec x_\rho}{|\vec x_\rho|}
\right](\vec\chi) = \frac{\vec\chi_\rho}{|\vec x_\rho|}
- \frac{\vec x_\rho}{|\vec x_\rho|^2}\,
 \frac{\vec x_\rho\,.\,\vec\chi_\rho}{|\vec x_\rho|}
= \vec\chi_s - \vec\tau\,(\vec\chi_s\,.\,\vec\tau)
= (\vec\chi_s\,.\,\vec\nu)\,\vec\nu \,, \label{eq:ddtau} \\
\left[\deldel{\vec x}\, \vec\nu\right](\vec\chi)
& = - \left[\deldel{\vec x}\, \vec\tau^\perp\right](\vec\chi) = 
- (\vec\chi_s\,.\,\vec\nu)\,\vec\nu^\perp =
- (\vec\chi_s\,.\,\vec\nu)\,\vec\tau \,, \label{eq:ddnu} \\
\left[\deldel{\vec x}\,\vec\nu\,|\vec x_\rho|\right](\vec\chi) &=
- \left[\deldel{\vec x}\,\vec x_{\rho}^\perp\right](\vec\chi)
= - \vec\chi_\rho^\perp = - \vec\chi_s^\perp \,|\vec x_\rho|\,, 
\label{eq:ddnuxrho}
\end{align}
\end{subequations}
where we always assume that $\vec\chi$ is sufficiently smooth so that all the
quantities are defined almost everywhere; e.g.\ $\vec\chi \in 
[W^{1,\infty}(I)]^2$.
In addition, we note that
\begin{equation} \label{eq:ddA}
\left[\deldel{\vec x}\, B(\vec x) \right] (\vec x_t) 
= \ddt\, B(\vec x)\,.
\end{equation}

Let 
\begin{subequations} 
\begin{align} \label{eq:Vpartialzero}
\Vpartialzero &= \{ \vec\eta \in [H^1(I)]^2 : \vec\eta(\rho)\,.\,\vec\ek_1 = 0
\quad \forall\ \rho \in \partial_0 I\}\,, \\
\xspace & = \left\{ \vec\eta \in \Vpartialzero : 
\vec\eta(\rho) = \vec 0\quad \forall\ \rho \in \partial_C I \cup \partial_N I\,,
\ \vec\eta(\rho)\,.\,\vec\ek_i = 0\quad \forall\ \rho \in \partial_i I\,, 
i = 1,2
\right\}\,. \label{eq:xspace}
\end{align}
\end{subequations}
Here $\xspace$ is the space of all directions in which we can vary a
given curve.
For a given $\vec z \in \bR^2$, on recalling (\ref{eq:partialM}), we define
\begin{equation} \label{eq:yspace}
\yspace(\vec z) = \left\{ \vec\eta \in \Vpartialzero : 
\vec\eta(\rho) = \vec z \quad \forall\ \rho \in \partial_M I \right\}\,.
\end{equation}
On recalling \eqref{eq:partialM}, we note that
\begin{equation} \label{eq:partialCempty}
\text{if } \partial_C I = \emptyset 
\quad\text{then } \yspace(\vec 0) \subset \xspace\,.
\end{equation}

Let $(\cdot,\cdot)$ denote the $L^2$--inner product on $I$.
We now consider the following weak formulation of (\ref{eq:varkappa})
with $\vec x \in \Vpartialzero$ and $\varkappa \in L^2(I)$ such that 
\begin{equation} \label{eq:varkappaweak}
\left( \varkappa\,\vec\nu,\vec\eta\, |\vec x_\rho| \right)
+ \left(\vec\tau,\vec\eta_\rho \right) =
\sum_{p \in \partial_C I} \left[\vec\zeta\,.\,\vec\eta\right](p)
+ \sum_{p \in \partial_M I} \left[\vec{\rm m}\,.\,\vec\eta\right](p)
\qquad \forall\ \vec\eta \in \Vpartialzero\,,
\end{equation}
where we recall (\ref{eq:tau}).
We note that (\ref{eq:varkappaweak}) weakly imposes (\ref{eq:bcbc})
and (\ref{eq:clampedI}), where 
$\vec\zeta(p) \in \bS^1$, $p\in \partial_C I$, 
are given data. 
However, (\ref{eq:varkappaweak}) also yields that  
$\vec{\rm m}(p) = \vec\mu(p) \in \bR^2$, $p\in \partial_M I$.
This will not be the case under discretization, where 
$\vec{\rm m}(p) \in \bR^2$, $p\in \partial_M I$,
are approximations to the conormals $\vec\mu(p)$, $p\in \partial_M I$.

Similarly, we consider the following weak formulation of (\ref{eq:kappaS})
with $\vec x \in \Vpartialzero$ and $\varkappa_{\mathcal{S}} \in L^2(I)$
such that 
\begin{equation} 
 \left(\vec x\,.\,\vec\ek_1\,
\varkappa_{\mathcal{S}}\,\vec\nu + \vec\ek_1,\vec\eta \,|\vec x_\rho| \right)
+ \left( (\vec x\,.\,\vec\ek_1)\,\vec\tau, \vec\eta_\rho \right) 
= \sum_{p \in \partial_C I} \left[
(\vec x\,.\,\vec\ek_1)\,\vec\zeta\,.\,\vec\eta\right](p)
+ \sum_{p \in \partial_M I} \left[(\vec x\,.\,\vec\ek_1)\,
\vec{\rm m}\,.\,\vec\eta\right](p)
\qquad \forall\ \vec\eta \in \Vpartialzero\,.
\label{eq:varkappaSweak}
\end{equation}
It is shown in Appendix~A of \cite{aximcf} that, despite the degenerate weight,
(\ref{eq:varkappaSweak}) weakly imposes (\ref{eq:bc}). 
In addition, (\ref{eq:varkappaSweak}) weakly imposes (\ref{eq:clampedI}).

\subsection{Based on $\varkappa$} \label{sec:Lag}

We begin with a weak formulation based on \eqref{eq:varkappaweak}.
Finite element approximations based on this weak formulation
will exhibit an equidistribution property.

On recalling (\ref{eq:EEaxCN}), (\ref{eq:AADE}), (\ref{eq:meanGaussS}), 
(\ref{eq:varkappaweak}) and that $\vec\mu =\vec {\rm m}$ on $\partial_M I$,
we define the Lagrangian
\begin{align} \label{eq:Lag}
&\mathcal{L}(\vec x, \varkappa^\star, \vec{\rm m}, \vec y) =
\pi\left( \alpha
\left[ \varkappa^\star - \frac{\vec\nu\,.\,\vec\ek_1}{\vec x\,.\,\vec\ek_1}
-\spont \right]^2 
+ 2\,\lambda, \vec x\,.\,\vec\ek_1\,|\vec x_\rho| \right)
+ \tfrac\beta2 \left[
2\,\pi \left( 
\varkappa^\star - \frac{\vec\nu\,.\,\vec\ek_1}{\vec x\,.\,\vec\ek_1},
\vec x\,.\,\vec\ek_1\,|\vec x_\rho| \right)  - M_0 \right]^2
\nonumber \\& \quad 
- \left(\varkappa^\star\,\vec\nu,\vec y\, |\vec x_\rho|\right) 
- \left(\vec\tau,\vec y_\rho\right) 
+ 2\,\pi\,\varsigma \sum_{p \in \partial_2 I \cup \partial_F I}
\vec x(p)\,.\,\vec\ek_1 
+ \sum_{p \in \partial_C I} \left[\vec\zeta\,.\,\vec y\right](p)
+\sum_{p \in \partial_M I} \left[\vec{\rm m}\,.\,(
\vec y - 2\,\pi\,\alpha_G\,\vec\ek_1)\right](p)\,,
\end{align}
for $\vec x \in \Vpartialzero$, $\varkappa^\star \in L^2(I)$,
$\vec{\rm m} : \partial_M I \to  \bR^2$
and $\vec y \in \Vpartialzero$. 
Here, we recall from (\ref{eq:axibc}) and (\ref{eq:bclimit}) 
that on the continuous level the Lagrangian (\ref{eq:Lag}) 
is well-defined for curves generating a smooth surface also in the case $\partial_0 I \not= \emptyset$. 

Taking variations $\vec\eta \in \Vpartialzero$ in $\vec y$, and setting
$\left[\deldel{\vec y}\, \mathcal{L}\right](\vec\eta) = 0$
we obtain 
\begin{equation} \label{eq:kappaeta}
(\varkappa^\star\,\vec\nu, \vec\eta\,|\vec x_\rho|)
+ (\vec\tau,\vec\eta_\rho) = 
\sum_{p \in \partial_C I} \left[ \vec\zeta\,.\,\vec\eta\right](p)
+ \sum_{p \in \partial_M I} \left[\vec{\rm m}\,.\,\vec\eta\right](p)
\qquad \forall\ \vec\eta \in \Vpartialzero\,,
\end{equation}
and so combining with (\ref{eq:varkappaweak}) yields that
$\varkappa^\star=\varkappa$. We are going to use this identity from now on.
Taking variations $\chi \in L^2(I)$ in $\varkappa^\star$ and setting
$\left[\deldel{\varkappa^\star}\, \mathcal{L}\right](\chi) = 0$
we obtain, on using $\varkappa^\star=\varkappa$, that
\begin{equation} \label{eq:kappay}
2\,\pi\left(\alpha
\left[\varkappa - \frac{\vec\nu\,.\,\vec\ek_1}{\vec x\,.\,\vec\ek_1}
-\spont\right] + \beta\,\AADE,
\vec x\,.\,\vec\ek_1\,\chi\,|\vec x_\rho|\right)
- \left(\vec\nu\,.\,\vec y,\chi\,|\vec x_\rho|\right) = 0
\qquad \forall\ \chi \in L^2(I)\,,
\end{equation}
where
\begin{equation} \label{eq:Axkappa}
\AADE(t) = 2\,\pi\left( 
\varkappa - \frac{\vec\nu\,.\,\vec\ek_1}{\vec x\,.\,\vec\ek_1},
\vec x\,.\,\vec\ek_1\,|\vec x_\rho| \right) - M_0 
= 2\,\pi\left(\vec x\,.\,\vec\ek_1\,\varkappa - \vec\nu\,.\,\vec\ek_1,
|\vec x_\rho| \right) - M_0 \,.
\end{equation}
We note that
\begin{equation} \label{eq:deldelA}
\left[\deldel{\vec x}\, \AADE(t)\right](\vec\chi)
= 2\,\pi\left(\varkappa, 
\left[\deldel{\vec x}\, \vec x\,.\,\vec\ek_1\,|\vec x_\rho|\right](\vec\chi)
 \right) 
- 2\,\pi\left(\vec\ek_1, 
\left[\deldel{\vec x}\, 
\vec\nu \,|\vec x_\rho|\right](\vec\chi) \right) 
\qquad \forall\ \vec\chi\in \xspace\,.
\end{equation}
Taking variations in $\vec{\rm m}$, and setting them to zero, yields
that
\begin{equation} \label{eq:my}
\vec y = 
2\,\pi\,\alpha_G\,\vec\ek_1 \qquad\text{on } \ \partial_M I\,.
\end{equation}
Taking variations $\vec\chi \in \xspace$ in $\vec x$, and 
setting $2\,\pi\,((\vec x\,.\,\vec\ek_1)\,\vec x_t\,.\,\vec\nu, 
\vec\chi\,.\,\vec\nu\,|\vec x_\rho|)= 
- \left[\deldel{\vec x}\, \mathcal{L}\right](\vec\chi)$
we obtain, on noting (\ref{eq:deldelA}) and (\ref{eq:xspace}), that 
\begin{align}
& 2\,\pi\left((\vec x\,.\,\vec\ek_1)\,\vec x_t\,.\,\vec\nu, 
\vec\chi\,.\,\vec\nu\,|\vec x_\rho|\right) 
= -\pi \left( \alpha
\left[\varkappa - \frac{\vec\nu\,.\,\vec e_1}{\vec x\,.\,\vec e_1} 
-\spont \right]^2 + 2\,\lambda + 2\,\beta\,\AADE\,\varkappa, 
\left[\deldel{\vec x}\,(\vec x\,.\,\vec\ek_1)\,|\vec x_\rho|\right]
(\vec\chi) \right) 
\nonumber \\ & \quad 
+2\,\pi\,\alpha \left(
\varkappa - \frac{\vec\nu\,.\,\vec e_1}{\vec x\,.\,\vec e_1} 
-\spont, \vec x\,.\,\vec\ek_1
\left[\deldel{\vec x}\,\frac{\vec\nu\,.\,\vec e_1}{\vec x\,.\,\vec e_1} 
\right]
(\vec\chi) \,|\vec x_\rho| \right) 
+ \left(\varkappa\,\vec y + 2\,\pi\,\beta\,\AADE\,\vec e_1,
\left[\deldel{\vec x}\,\vec\nu\,|\vec x_\rho|\right](\vec\chi)\right)
+ \left(\vec y_\rho , 
\left[\deldel{\vec x}\,\vec\tau\right](\vec\chi)\right)
\nonumber \\ & \quad
- 2\,\pi\,\varsigma \sum_{p \in \partial_{2} I \cup \partial_F I}
\vec\chi(p)\,.\,\vec\ek_1
\qquad \forall\ \vec\chi \in \xspace\,.
\label{eq:Lagdx}
\end{align}
Choosing $\vec\chi = \vec x_t \in \xspace$ in (\ref{eq:Lagdx}) yields, 
on noting (\ref{eq:ddA}), that
\begin{align}
& 2\,\pi\left(\vec x\,.\,\vec\ek_1\,(\vec x_t\,.\,\vec\nu)^2, 
|\vec x_\rho|\right)
= -\pi \left( \alpha
\left[\varkappa - \frac{\vec\nu\,.\,\vec e_1}{\vec x\,.\,\vec e_1} 
-\spont \right]^2 + 2\,\lambda + 2\,\beta\,\AADE\,\varkappa, 
\left[(\vec x\,.\,\vec\ek_1)\,|\vec x_\rho|\right]_t \right) 
\nonumber \\ & \
+2\,\pi\,\alpha
 \left( \varkappa - \frac{\vec\nu\,.\,\vec e_1}{\vec x\,.\,\vec e_1} 
-\spont,  \vec x\,.\,\vec\ek_1
\left[\frac{\vec\nu\,.\,\vec e_1}{\vec x\,.\,\vec e_1} \right]_t
|\vec x_\rho| \right) 
+ \left(\varkappa\,\vec y + 2\,\pi\,\beta\,\AADE\,\vec e_1,
\left[\vec\nu\,|\vec x_\rho|\right]_t\right)
+ \left(\vec y_\rho , \vec\tau_t\right)
- 2\,\pi\,\varsigma \sum_{p \in \partial_{2} I \cup \partial_F I}
\vec x_t(p)\,.\,\vec\ek_1\, .
\label{eq:Pxtxt}
\end{align}

Differentiating (\ref{eq:kappaeta}) with respect to time, and then
choosing $\vec\eta = \vec y$, on recalling $\varkappa^\star = \varkappa$
and $\vec\zeta$ is independent of $t$,
yields that
\begin{equation} \label{eq:dtside}
\left( \varkappa_t, \vec y\,.\,\vec\nu\,|\vec x_\rho| \right) 
+ \left( \varkappa\,\vec y, \left[\vec\nu\,|\vec x_\rho|\right]_t \right) 
+ \left(\vec\tau_t,\vec y_\rho \right) = 
\sum_{p \in \partial_M I} \left[ \vec{\rm m}_t\,.\,\vec y\right](p) \,.
\end{equation}
It follows from (\ref{eq:dtside}), (\ref{eq:my}) and (\ref{eq:kappay}) with $\chi=\varkappa_t 
\in L^2(I)$ that
\begin{align}
\left( \varkappa\,\vec y, \left[\vec\nu\,|\vec x_\rho|\right]_t \right) 
+ \left(\vec\tau_t,\vec y_\rho \right)
& = - \left( \varkappa_t, \vec y\,.\,\vec\nu\,|\vec x_\rho| \right) 
+ 2\,\pi\,\alpha_G \sum_{p \in \partial_M I} \vec{\rm m}_t(p)\,.\,\vec\ek_1
\nonumber \\ & 
= - 2\,\pi \left(  \alpha
\left[\varkappa - \frac{\vec\nu\,.\,\vec\ek_1}{\vec x\,.\,\vec\ek_1}
-\spont\right] + \beta\,\AADE, 
\vec x\,.\,\vec\ek_1\, \varkappa_t\,|\vec x_\rho| \right) 
+ 2\,\pi\,\alpha_G \sum_{p \in \partial_M I} \vec{\rm m}_t(p)\,.\,\vec\ek_1\,.
\label{eq:mykappay}
\end{align}
Combining (\ref{eq:Pxtxt}) and (\ref{eq:mykappay}) yields that
\begin{align}
& 2\,\pi\left(\vec x\,.\,\vec\ek_1\,(\vec x_t\,.\,\vec\nu)^2,
|\vec x_\rho|\right) 
 = -\pi \left( \alpha
\left[\varkappa - \frac{\vec\nu\,.\,\vec e_1}{\vec x\,.\,\vec e_1} 
-\spont \right]^2 + 2\,\lambda + 2\,\beta\,\AADE\,\varkappa, 
\left[(\vec x\,.\,\vec\ek_1)\,|\vec x_\rho|\right]_t \right) 
\nonumber \\ & \quad
+2\,\pi\,\alpha \left( 
\varkappa - \frac{\vec\nu\,.\,\vec e_1}{\vec x\,.\,\vec e_1} -\spont , 
\vec x\,.\,\vec\ek_1 
\left[\frac{\vec\nu\,.\,\vec e_1}{\vec x\,.\,\vec e_1} \right]_t
|\vec x_\rho| \right) 
+ 2\,\pi\,\beta\,\AADE \left(\vec\ek_1, 
\left[ \vec\nu\, |\vec x_\rho| \right]_t \right) 
\nonumber \\ & \quad 
- 2\,\pi
\left(\alpha\left[\varkappa - \frac{\vec\nu\,.\,\vec\ek_1}{\vec x\,.\,\vec\ek_1}
-\spont\right] + \beta\,\AADE,
\vec x\,.\,\vec\ek_1\,\varkappa_t\,|\vec x_\rho|\right) 
+ 2\,\pi\,\alpha_G \sum_{p \in \partial_M I} \vec{\rm m}_t(p)\,.\,\vec\ek_1
- 2\,\pi\,\varsigma \sum_{p \in \partial_{2} I \cup \partial_F I}
\vec x_t(p)\,.\,\vec\ek_1 
\nonumber \\ & \
= - \pi\,\ddt 
\left[\alpha\left
[\varkappa - \frac{\vec\nu\,.\,\vec\ek_1}{\vec x\,.\,\vec\ek_1}-\spont
\right]^2 + 2\,\lambda,\vec x\,.\,\vec\ek_1\,|\vec x_\rho|\right]
- \tfrac\beta2\,\ddt\,\AADE^2 
- 2\,\pi\,\varsigma \sum_{p \in \partial_{2} I \cup \partial_F I}
\vec x_t(p)\,.\,\vec\ek_1 
+2\,\pi\,\alpha_G \sum_{p \in \partial_M I} \vec{\rm m}_t(p)\,.\,
 \vec\ek_1 \nonumber \\ & \
= - \ddt \, \widehat E(t)\,,
\label{eq:Pxtstab}
\end{align}
where we have recalled (\ref{eq:EEaxCN}), (\ref{eq:AADE}), (\ref{eq:meanGaussS}), 
(\ref{eq:Axkappa}) and that $\vec\mu =\vec {\rm m}$ on $\partial_M I$. 

\begin{rmrk} \label{rem:Pxtstab}
The property $\ddt \, \widehat E(t) + 
2\,\pi\,(\vec x\,.\,\vec\ek_1\,(\vec x_t\,.\,\vec\nu)^2,
|\vec x_\rho|) = 0$ shown in \eqref{eq:Pxtstab} demonstrates the gradient
flow property of the derived weak formulation.
\end{rmrk}

We now return to (\ref{eq:Lagdx}), which, on recalling (\ref{eq:24}), can be rewritten as
\begin{align}
& 2\,\pi\left((\vec x\,.\,\vec\ek_1)\,\vec x_t\,.\,\vec\nu,
\vec\chi\,.\,\vec\nu \,|\vec x_\rho|\right) 
=- \pi \left( \alpha
\left[\varkappa - \frac{\vec\nu\,.\,\vec\ek_1}{\vec x\,.\,\vec e_1} 
-\spont \right]^2 + 2\,\lambda + 2\,\beta\,\AADE\,\varkappa, 
\vec\chi\,.\,\vec\ek_1\,|\vec x_\rho| + (\vec x\,.\,\vec\ek_1)\,\vec\tau\,.\,
\vec\chi_\rho \right) 
\nonumber \\ & \qquad 
- 2\,\pi\,\alpha\left( 
\varkappa - \frac{\vec\nu\,.\,\vec\ek_1}{\vec x\,.\,\vec e_1} - \spont, 
\frac{\vec\nu\,.\,\vec e_1}{\vec x\,.\,\vec e_1}\,\vec\chi\,.\,\vec\ek_1
\,|\vec x_\rho| \right) 
- 2\,\pi\,\alpha\left(
\varkappa - \frac{\vec\nu\,.\,\vec\ek_1}{\vec x\,.\,\vec e_1} - \spont ,
(\vec\chi_\rho\,.\,\vec\nu)\,\vec\tau\,.\,\vec\ek_1 \right) 
- \left( \varkappa\, \vec y + 2\,\pi\,\beta\,\AADE\,\vec e_1 ,
\vec\chi_\rho^\perp \right)
\nonumber \\ & \qquad 
+ \left(\vec y_\rho\,.\,\vec\nu , \vec\chi_\rho\,.\,\vec\nu\,|\vec x_\rho|^{-1}
\right) 
- 2\,\pi\,\varsigma \sum_{p \in \partial_{2} I \cup \partial_F I}
\vec\chi(p)\,.\,\vec\ek_1
\qquad \forall\ \vec\chi \in \xspace\,.
\label{eq:Pxt2}
\end{align}

Our finite element approximation is going to be based on the following
formulation, on combining (\ref{eq:Pxt2}), (\ref{eq:kappay}),
(\ref{eq:kappaeta}), (\ref{eq:my}) and (\ref{eq:Axkappa}), 
and on recalling $\varkappa^\star=\varkappa$, (\ref{eq:abperp}) and (\ref{eq:tau}). 

\noindent
$(\BGNpwf)$
Let $\vec x(\cdot,0) \in \Vpartialzero$ and $\alpha \in \bRplus$,
$\spont,M_0,\alpha_G,\lambda,\varsigma \in \bR$, $\beta \in \bRgeq$, 
$\vec\zeta : \partial_C I \to \bS^1$ be given.
For $t \in (0,T]$, find $\vec x(\cdot,t) \in \Vpartialzero$, 
with $\vec x_t(\cdot,t) \in \xspace$, $\varkappa(\cdot,t) \in L^2(I)$, 
$\vec y(\cdot,t) \in \yspace(2\,\pi\,\alpha_G\,\vec\ek_1)$ 
and $\vec{\rm m}(\cdot,t) : \partial_M I \to \bR^2$
such that
\begin{subequations} \label{eq:weak3}
\begin{align}
& 2\,\pi
\left((\vec x\,.\,\vec\ek_1)\,\vec x_t\,.\,\vec\nu, \vec\chi\,.\,\vec\nu\,
|\vec x_\rho|\right) - \left(\vec y_\rho\,.\,\vec\nu, \vec\chi_\rho\,.\,\vec\nu
\, |\vec x_\rho|^{-1} \right) 
\nonumber \\ & 
= - \pi \left( \alpha
\left[\varkappa - \frac{\vec\nu\,.\,\vec\ek_1}{\vec x\,.\,\vec e_1} -\spont 
\right]^2 + 2\,\lambda+ 2\,\beta\,\AADE\,\varkappa, 
\vec\chi\,.\,\vec\ek_1\,|\vec x_\rho| + (\vec x\,.\,\vec\ek_1)\,
\vec\tau\,.\,\vec\chi_\rho \right) 
- 2\,\pi\,\alpha
\left(\varkappa - \frac{\vec\nu\,.\,\vec\ek_1}{\vec x\,.\,\vec e_1} - \spont,
\frac{\vec\nu\,.\,\vec e_1}{\vec x\,.\,\vec e_1}\,\vec\chi\,.\,\vec\ek_1
\,|\vec x_\rho| \right) 
\nonumber \\ & \
- 2\,\pi\,\alpha \left(
\varkappa - \frac{\vec\nu\,.\,\vec\ek_1}{\vec x\,.\,\vec e_1} - \spont,
(\vec\tau\,.\,\vec\ek_1)\,\vec\chi_\rho\,.\,\vec\nu \right) 
+ \left( \varkappa\, \vec y^\perp - 2\,\pi\,\beta\,\AADE \,\vec\ek_2,
\vec\chi_\rho \right)
- 2\,\pi\,\varsigma \sum_{p \in  \partial_{2} I \cup \partial_F I}
\vec\chi(p)\,.\,\vec\ek_1
\quad \forall\ \vec\chi \in \xspace \,,
\label{eq:weak3a} \\
&2\,\pi\left(\alpha
\left[\varkappa - \frac{\vec\nu\,.\,\vec\ek_1}{\vec x\,.\,\vec\ek_1}
-\spont\right] + \beta\,\AADE,
\vec x\,.\,\vec\ek_1\,\chi\,|\vec x_\rho|\right)
- \left(\vec\nu\,.\,\vec y,\chi\,|\vec x_\rho|\right) = 0
\qquad \forall\ \chi \in L^2(I)\,, \label{eq:weak3b} \\
& (\varkappa\,\vec\nu, \vec\eta\,|\vec x_\rho|)
+ (\vec x_\rho ,\vec\eta_\rho \,|\vec x_\rho|^{-1}) = 
\sum_{p \in \partial_C I} \left[\vec\zeta\,.\,\vec\eta\right](p)
+ \sum_{p \in \partial_M I} \left[\vec{\rm m}\,.\,\vec\eta\right](p)
 \qquad \forall\ \vec\eta \in \Vpartialzero\,, \label{eq:weak3c} 
\end{align}
\end{subequations}
where $\AADE(t)$ is given by (\ref{eq:Axkappa}).
We note that the number of unknowns fixed via $\vec y \in
\yspace(2\,\pi\,\alpha_G\,\vec\ek_1)$ on $\partial_M I$ is matched by the new
degrees of freedom arising from $\{\vec{\rm m}(p)\}_{p\in\partial_M I}$.

We remark that (\ref{eq:weak3}) is independent of the tangential part
$\vec y\,.\,\vec\tau$ of $\vec y$. To see this, we note that it follows from
(\ref{eq:tau}) and (\ref{eq:varkappa}) that
\begin{equation} \label{eq:ysnu}
(\vec y_s\,.\,\vec\nu)\,\vec\nu
= (\vec y\,.\,\vec\nu)_s\,\vec\nu - (\vec y\,.\,\vec\nu_s)\,\vec\nu
= (\vec y\,.\,\vec\nu)_s\,\vec\nu + \varkappa\,(\vec
y\,.\,\vec\tau)\,\vec\nu \qquad \mbox{in }\ \overline I \,.
\end{equation}
Hence the only terms involving $\vec y$ in (\ref{eq:weak3a}) are
\begin{equation} \label{eq:kyperp}
\varkappa\,\vec y^\perp + (\vec y_s\,.\,\vec\nu)\,\vec\nu
= \varkappa\,(\vec y^\perp + (\vec y\,.\,\vec\tau)\,\vec\nu)
+ (\vec y\,.\,\vec\nu)_s\,\vec\nu
= \varkappa\,(\vec y\,.\,\vec\nu)\,\vec\tau + (\vec y\,.\,\vec\nu)_s\,\vec\nu
\qquad \mbox{in }\ \overline I \,,
\end{equation}
where we have recalled (\ref{eq:aperp}).
This shows that (\ref{eq:weak3}) only depends on $\vec y\,.\,\vec\nu$,
and not on $\vec y\,.\,\vec\tau$.
We refer to Appendix~\ref{sec:A1}, where we show that 
(\ref{eq:weak3}) for a sufficiently smooth solution gives rise to the strong 
form (\ref{eq:xtbgnlambda}) and (\ref{eq:varkappa}).

\subsubsection{Conserved flows}\label{sec:LagC}

In this subsection we present a weak formulation for the conserving flow
\eqref{eq:Helfrich}. To this end, we assume that 
the hypersurface $\mathcal{S}(t)$ has no boundary, and so
\begin{equation} \label{eq:ppI}
\partial I = \partial_0 I \qquad\implies\qquad \xspace = \Vpartialzero\,.
\end{equation}
Then, on writing (\ref{eq:weak3a}) as
\begin{equation*} 
2\,\pi \left((\vec x\,.\,\vec\ek_1)\,\vec x_t\,.\,\vec\nu,
\vec\chi\,.\,\vec\nu\,|\vec x_\rho|\right)
- \left(\vec y_\rho\,.\,\vec\nu, \vec\chi_\rho \,.\,\vec\nu \,
|\vec x_\rho|^{-1}\right) =
\left( \vec f, \vec\chi\,|\vec x_\rho|\right)
\quad \forall\ \vec\chi \in \xspace
\,,
\end{equation*}
a weak formulation of \eqref{eq:xHelfrich} and \eqref{eq:xsidedtSAV} 
is given by \eqref{eq:weak3}, with \eqref{eq:weak3a} replaced by
\begin{align} \label{eq:weak4LMa}
& 2\,\pi
\left((\vec x\,.\,\vec\ek_1)\,\vec x_t\,.\,\vec\nu,
\vec\chi\,.\,\vec\nu\,|\vec x_\rho|\right)
- \left(\vec y_\rho\,.\,\vec\nu, \vec\chi_\rho \,.\,\vec\nu 
\,|\vec x_\rho|^{-1}\right) = 
\left( \vec f, \vec\chi\,|\vec x_\rho|\right)
\nonumber \\ & \qquad
- 2\,\pi\, \lambda_A 
\left[
\left( \vec\ek_1, \vec\chi\,|\vec x_\rho|\right) 
+ \left( (\vec x\,.\,\vec\ek_1)\,\vec\tau
, \vec\chi_\rho\right) \right]
- 2\,\pi\, \lambda_V\left((\vec x\,.\,\vec\ek_1)\,\vec\nu, 
\vec\chi\,|\vec x_\rho|\right)\quad \forall\ \vec\chi \in \xspace
\,,
\end{align}
where $(\lambda_A(t),\lambda_V(t))^T \in \bR^2$ are chosen such that 
\eqref{eq:xsidedtSAV} holds.

Choosing $\vec\eta = (\vec x\,.\,\vec\ek_1)\,\vec x_t
\in \xspace = \Vpartialzero$ in (\ref{eq:weak3c}) and noting (\ref{eq:tau})
yields that
\begin{align}
\left( \vec\ek_1, \vec x_t\,|\vec x_\rho|\right) 
+ \left( (\vec x\,.\,\vec\ek_1)\,\vec\tau, (\vec x_t)_\rho\right)
& 
= - \left( (\vec x\,.\,\vec\ek_1)\,\varkappa,\vec x_t \,.\,\vec\nu\,
|\vec x_\rho| \right) + \left(\vec x_t\,.\,
\left[\vec\ek_1 - (\vec\ek_1\,.\,\vec\tau)\,\vec\tau \right], |\vec x_\rho| \right) 
 \nonumber \\ & 
=- \left( \vec x\,.\,\vec\ek_1\left(
\varkappa - \frac{ \vec\nu\,.\, \vec\ek_1}
{\vec x\,.\,\vec\ek_1} \right) ,\vec x_t \,.\,\vec\nu\,
|\vec x_\rho| \right).
\label{eq:dAdtk}
\end{align} 
Now choosing 
$\vec\chi = \left(\varkappa - \frac{\vec\nu\,.\,\vec\ek_1}{\vec x\,.\,\vec\ek_1}
 \right)\vec\nu \in \xspace = \Vpartialzero$
and $\vec\chi = \vec\nu \in \xspace = \Vpartialzero$ in (\ref{eq:weak4LMa}), 
recall (\ref{eq:bcnu}), (\ref{eq:ppI}) and (\ref{eq:Vpartialzero}),
we see that the two side constraints in \eqref{eq:xsidedtSAV}
will be satisfied if $(\lambda_A(t), \lambda_V(t))^T \in \bR^2$ solve 
the symmetric system
\begin{align}
& 2\,\pi
\begin{pmatrix}
\left( \vec x\,.\,\vec\ek_1,
\left(\varkappa - \frac{\vec\nu\,.\,\vec\ek_1}{\vec x\,.\,\vec\ek_1}
 \right)^2 |\vec x_\rho|\right)
& \left( \vec x\,.\,\vec\ek_1,
\left(\varkappa - \frac{\vec\nu\,.\,\vec\ek_1}{\vec x\,.\,\vec\ek_1}
 \right) |\vec x_\rho|\right) \\
\left( \vec x\,.\,\vec\ek_1,
\left(\varkappa - \frac{\vec\nu\,.\,\vec\ek_1}{\vec x\,.\,\vec\ek_1}
 \right) |\vec x_\rho|\right) &
\left( \vec x\,.\,\vec\ek_1,|\vec x_\rho|\right)
\end{pmatrix}
\begin{pmatrix} \lambda_A \\ \lambda_V \end{pmatrix}
\nonumber \\ & \qquad\qquad
= 
\begin{pmatrix} 
 \left(\vec y_\rho\,.\,\vec\nu, 
\left[\left(\varkappa - \frac{\vec\nu\,.\,\vec\ek_1}{\vec x\,.\,\vec\ek_1}
 \right)\vec\nu\right]_\rho .\, \vec\nu
\,|\vec x_\rho|^{-1} \right) 
+
\left( \vec f, \left(\varkappa - \frac{\vec\nu\,.\,\vec\ek_1}{\vec x\,.\,\vec\ek_1}
 \right)\vec\nu
\,|\vec x_\rho|\right) \\
\left( \vec f, \vec\nu \,|\vec x_\rho|\right)
\end{pmatrix}\,.
\label{eq:lambdamu}
\end{align}
The matrix in (\ref{eq:lambdamu}) is symmetric positive
semidefinite, and it is singular if and only if 
$\varkappa_{\mathcal{S}} = 
\varkappa - \frac{\vec\nu\,.\,\vec\ek_1}{\vec x\,.\,\vec\ek_1}$ 
is a constant. It can be shown that this is equivalent to 
$\mathcal{S}(t)$ being a sphere and hence to $\Gamma(t)$ being an open 
halfcircle. 
Choosing 
$\vec\chi = \vec x_t \in \xspace=\Vpartialzero$ in (\ref{eq:weak4LMa}) 
and noting \eqref{eq:xsidedtSAV}
proves the stability result (\ref{eq:Pxtstab})
for the weak formulation of the conserved flow, 
on recalling (\ref{eq:Pxtxt}), (\ref{eq:mykappay}) and (\ref{eq:weak3a}).

An alternative formulation of the two conservation side constraints can
be obtained by observing that \eqref{eq:sidedtSAV} is equivalent to
\begin{equation*} 
\mathcal{H}^2(\mathcal{S}(t)) = \mathcal{H}^2(\mathcal{S}(0))
\,,\qquad
\mathcal{L}^3(\Omega(t)) = \mathcal{L}^3(\Omega(0))\,.
\end{equation*}
In order to formulate \eqref{eq:sidedtSAV} in terms of $\vec x$, we define,
on recalling \eqref{eq:dAdt}, 
\begin{equation} \label{eq:Area}
A(\vec x(t))
= 2\,\pi \left(\vec x\,.\,\vec\ek_1 , |\vec x_\rho| \right)
= \mathcal{H}^2(\mathcal{S}(t))
\end{equation}
and, see e.g.\ \cite[(3.10)]{axisd},
\begin{equation} \label{eq:Volume}
V(\vec x(t)) = \pi\left( (\vec x\,.\,\vec\ek_1)^2 ,\vec\nu\,.\,\vec\ek_1\,
|\vec x_\rho| \right) =
\mathcal{L}^3(\Omega(t)) \,.
\end{equation}
Hence $(\lambda_A(t), \lambda_V(t))^T \in \bR^2$ in \eqref{eq:dAdtk} will
be such that \eqref{eq:xsidedtSAV} holds, which is equivalent to
\begin{equation} \label{eq:xsideSAV}
A(\vec x(t)) = A(\vec x(0))\,,\qquad V(\vec x(t)) = V(\vec x(0))\,.
\end{equation}
Our discretization in \S\ref{sec:LagC} will be based on \eqref{eq:dAdtk} 
and \eqref{eq:xsideSAV}.

\subsection{Based on $\varkappa_{\mathcal{S}}$} \label{sec:Lag2}

A drawback of the formulation used in \S\ref{sec:Lag} is that 
in the case $\partial_0 I \not= \emptyset$
the integrals featuring the singular fraction 
$\frac{\vec\nu\,.\,\vec\ek_1}{\vec x\,.\,\vec\ek_1}$ in
the Lagrangian \eqref{eq:Lag} are not well-defined on the discrete
level, and so an appropriate interpretation of these terms is needed.
An alternative is
to use the mean curvature of the surface as a variable in the weak formulation,
i.e.\ to use a formulation that features \eqref{eq:varkappaSweak}. Then the
discretization follows naturally, and a similar approach has been followed by
the present authors for flows in Riemannian manifolds in \cite{hypbolpwf}.

On recalling (\ref{eq:EEaxCN}), (\ref{eq:AADE}), (\ref{eq:varkappaSweak})
and that $\vec\mu = \vec {\rm m}$ on $\partial_M I$,  
we define the Lagrangian
\begin{align} \label{eq:Lag2}
& \mathcal{L}_{\mathcal{S}}(\vec x, \varkappa^\star_{\mathcal{S}}, 
\vec{\rm m}, \vec y_{\mathcal{S}})\nonumber 
= \pi\left( \alpha \left[ \varkappa^\star_{\mathcal{S}} -\spont \right]^2
+ 2\,\lambda
, \vec x\,.\,\vec\ek_1 \, |\vec x_\rho| \right)
+ \tfrac\beta2 \left[
2\,\pi \left( \vec x\,.\,\vec\ek_1\,\varkappa_{\mathcal{S}}^\star, 
|\vec x_\rho| \right)  - M_0 \right]^2
- \left(\vec x\,.\,\vec\ek_1\,\varkappa^\star_{\mathcal{S}}\,\vec\nu
 + \vec\ek_1,\vec y_{\mathcal{S}}\, |\vec x_\rho|\right)
\nonumber \\& \
 - \left( (\vec x\,.\,\vec\ek_1)\,\vec\tau,(\vec y_{\mathcal{S}})_\rho\right)
+ 2\,\pi\,\varsigma\hspace{-4mm}\sum_{p \in \partial_2 I \cup \partial_F I}
\vec x(p)\,.\,\vec\ek_1 
+ \sum_{p \in \partial_C I} \left[
(\vec x\,.\,\vec\ek_1)\,\vec\zeta\,.\,\vec y_{\mathcal{S}}\right](p)
+\sum_{p \in \partial_M I} \left[\vec{\rm m}\,.\,(
(\vec x\,.\,\vec\ek_1)\,\vec y_{\mathcal{S}} 
- 2\,\pi\,\alpha_G\,\vec\ek_1)\right](p)\,,
\end{align}
for $\vec x \in \Vpartialzero$,
$\varkappa^\star_{\mathcal{S}} \in L^2(I)$,
$\vec{\rm m} : \partial_M I \to \bR^2$
and $\vec y_{\mathcal{S}} \in \Vpartialzero$.

Taking variations $\vec\eta \in \Vpartialzero$ in 
$\vec y_{\mathcal{S}}$, and setting
$\left[\deldel{\vec y_{\mathcal{S}}}\, 
\mathcal{L}_{\mathcal{S}}\right](\vec\eta) = 0$,
we obtain 
\begin{equation} \label{eq:kappaSeta}
\left(
\vec x\,.\,\vec\ek_1\,\varkappa^\star_{\mathcal{S}}\,\vec\nu + \vec\ek_1, 
\vec\eta\,|\vec x_\rho|\right)
+ \left((\vec x\,.\,\vec\ek_1)\,\vec\tau,\vec\eta_\rho\right) 
= \sum_{p \in \partial_C I} \left[
(\vec x\,.\,\vec\ek_1)\,\vec\zeta\,.\,\vec\eta\right](p)
+ \sum_{p \in \partial_M I} \left[(\vec x\,.\,\vec\ek_1)\,
\vec{\rm m}\,.\,\vec\eta\right](p) 
\qquad \forall\ \vec\eta \in \Vpartialzero\,,
\end{equation}
and so combining with (\ref{eq:varkappaSweak}) yields that
$\varkappa_{\mathcal{S}}^\star=\varkappa_{\mathcal{S}}$. 
We are going to use this identity from now on.
Taking variations $\chi \in L^2(I)$ in $\varkappa_{\mathcal{S}}^\star$, 
and setting
$\left[\deldel{\varkappa_{\mathcal{S}}^\star}\, 
\mathcal{L}_{\mathcal{S}}\right](\chi) = 0$, we obtain,
on using $\varkappa_{\mathcal{S}}^\star=\varkappa_{\mathcal{S}}$, that
\begin{equation} \label{eq:kappaSy}
2\,\pi
\left( \alpha\,(\varkappa_{\mathcal{S}} - \spont)
+ \beta\,\AADE_{\mathcal{S}}
, \vec x\,.\,\vec\ek_1\,\chi\,|\vec x_\rho|\right)
- \left((\vec x\,.\,\vec\ek_1)\,\vec y_{\mathcal{S}}, \chi\,\vec\nu
\,|\vec x_\rho|\right) = 0 \qquad \forall\ \chi \in L^2(I)\,,
\end{equation}
where
\begin{equation} \label{eq:ASxkappa}
\AADE_{\mathcal{S}}(t) = 2\,\pi\left(\varkappa_{\mathcal{S}},
\vec x\,.\,\vec\ek_1\,|\vec x_\rho| \right) - M_0 \,.
\end{equation}
We note that
\begin{equation} \label{eq:deldelAS}
\left[\deldel{\vec x}\, \AADE_{\mathcal{S}}\right](\vec\chi)
= 2\,\pi\left(\varkappa_{\mathcal{S}}, 
\left[\deldel{\vec x}\, \vec x\,.\,\vec\ek_1\,|\vec x_\rho|\right](\vec\chi)
 \right) \qquad \forall\ \vec\chi\in \xspace\,.
\end{equation}
Taking variations in $\vec{\rm m}$, and setting them to zero, yields that
\begin{equation} \label{eq:mSy}
(\vec x\,.\,\vec\ek_1)\,\vec y_{\mathcal{S}} = 
2\,\pi\,\alpha_G\,\vec\ek_1 \qquad\text{ on } \quad \partial_M I\,.
\end{equation}
Taking variations $\vec\chi \in \xspace$ in $\vec x$,
and setting
$2\,\pi\,((\vec x\,.\,\vec\ek_1)\,\vec x_t\,.\,\vec\nu, 
\vec\chi\,.\,\vec\nu\,|\vec x_\rho|)= 
- \left[\deldel{\vec x}\, \mathcal{L}_{\mathcal{S}}\right](\vec\chi)$
we obtain, on noting (\ref{eq:deldelAS}), (\ref{eq:xspace}) and (\ref{eq:partialM}), that 
\begin{align}
&2\,\pi\left((\vec x\,.\,\vec\ek_1)\,\vec x_t\,.\,\vec\nu, 
\vec\chi\,.\,\vec\nu\,|\vec x_\rho|\right)
= -\pi \left(\alpha \left[\varkappa_{\mathcal{S}} -\spont \right]^2 
+  2\,\lambda + 2\,\beta\,\AADE_{\mathcal{S}}\,\varkappa_{\mathcal{S}}, 
\left[\deldel{\vec x}\,\vec x\,.\,\vec\ek_1\,|\vec x_\rho|\right](\vec\chi) 
\right) 
\nonumber \\ & \qquad 
+ \left(\varkappa_{\mathcal{S}}\,\vec y_{\mathcal{S}},
\left[\deldel{\vec x}\,(\vec x\,.\,\vec\ek_1)\,\vec\nu\,|\vec x_\rho|\right](\vec\chi)\right)
 + \left(\vec\ek_1\,.\,\vec y_{\mathcal{S}}, 
 \left[\deldel{\vec x}\,|\vec x_\rho|\right](\vec\chi)\right)
\nonumber \\ & \qquad 
+ \left((\vec y_{\mathcal{S}})_\rho , 
\left[\deldel{\vec x}\,(\vec x\,.\,\vec\ek_1)\,\vec\tau\right](\vec\chi)
\right) 
-\sum_{p \in \partial_{2} I \cup \partial_F I}
\left[
\left[2\,\pi\,\varsigma + \vec{\rm m}\,.\,\vec y_{\mathcal{S}}
\right] \vec\chi\,.\,\vec\ek_1 \right](p) 
\qquad \forall\ \vec\chi\in\xspace\,.
\label{eq:Lag2dx}
\end{align}
Choosing $\vec\chi = \vec x_t \in \xspace$ in (\ref{eq:Lag2dx}), and recalling
(\ref{eq:ddA}), yields
\begin{align}
&2\,\pi\left(\vec x\,.\,\vec\ek_1\,(\vec x_t\,.\,\vec\nu)^2, 
|\vec x_\rho|\right) 
= -\pi \left(\alpha \left[\varkappa_{\mathcal{S}} -\spont \right]^2 
+ 2\,\lambda + 2\,\beta\,\AADE_{\mathcal{S}}\,\varkappa_{\mathcal{S}}, 
\left[\vec x\,.\,\vec\ek_1\,|\vec x_\rho|\right]_t\right) 
+ \left(\varkappa_{\mathcal{S}}\,\vec y_{\mathcal{S}},
\left[(\vec x\,.\,\vec\ek_1)\,\vec\nu\,|\vec x_\rho|\right]_t\right)
\nonumber \\ & \qquad 
 + \left(\vec\ek_1\,.\,\vec y_{\mathcal{S}}, 
 \left[|\vec x_\rho|\right]_t\right)
+ \left((\vec y_{\mathcal{S}})_\rho , 
\left[(\vec x\,.\,\vec\ek_1)\,\vec\tau\right]_t \right)
-\sum_{p \in \partial_{2} I \cup \partial_F I}
\left[
\left[2\,\pi\,\varsigma + \vec{\rm m}\,.\,\vec y_{\mathcal{S}}
\right] \vec x_t\,.\,\vec\ek_1 \right](p) 
\,.
\label{eq:Lag2dxxt}
\end{align}

Differentiating (\ref{eq:kappaSeta}) with respect to $t$ and then choosing
$\vec\eta = \vec y_{\mathcal{S}} \in \Vpartialzero$, on recalling 
$\varkappa_{\mathcal{S}}^\star=\varkappa_{\mathcal{S}}$, $\vec x_t \in \xspace$ 
and $\vec\zeta$ is independent of $t$
and so the term on $\partial_C I$ vanishes, yields that
\begin{align}
& \left( (\varkappa_{\mathcal{S}})_t , (\vec x\,.\,\vec\ek_1)\,\,
\vec y_{\mathcal{S}}\,.\,\vec\nu\,|\vec x_\rho|\right)
+ \left( \varkappa_{\mathcal{S}}\,\vec y_{\mathcal{S}}, 
\left[\vec x\,.\,\vec\ek_1\,|\vec x_\rho|\,\vec\nu\right]_t \right)
+ \left( \vec\ek_1, \vec y_{\mathcal{S}}\left[|\vec x_\rho|\right]_t\right)
+ \left(\left[(\vec x\,.\,\vec\ek_1)\,\vec\tau\right]_t,
(\vec y_{\mathcal{S}})_\rho\right) 
\nonumber \\ & \qquad
= \sum_{p \in \partial_M I} \left[(\vec x_t\,.\,\vec\ek_1)\,
\vec{\rm m}\,.\,\vec y_{\mathcal{S}} + (\vec x\,.\,\vec\ek_1)\,
\vec{\rm m}_t\,.\,\vec y_{\mathcal{S}} \right](p) \,.
\label{eq:new5}
\end{align}
It follows from (\ref{eq:new5}), (\ref{eq:mSy}) and (\ref{eq:kappaSy}) with $\chi
= [\varkappa_{\mathcal{S}}]_t$ that
\begin{align}
&
\left(\varkappa_{\mathcal{S}}\,\vec y_{\mathcal{S}},
\left[(\vec x\,.\,\vec\ek_1)\,\vec\nu\,|\vec x_\rho|\right]_t\right)
 + \left(\vec\ek_1\,.\,\vec y_{\mathcal{S}}, 
 \left[|\vec x_\rho|\right]_t\right)
+ \left((\vec y_{\mathcal{S}})_\rho , 
\left[(\vec x\,.\,\vec\ek_1)\,\vec\tau\right]_t \right)
\nonumber \\ & \quad 
= - \left( [\varkappa_{\mathcal{S}}]_t , (\vec x\,.\,\vec\ek_1)\,\,
\vec y_{\mathcal{S}}\,.\,\vec\nu\,|\vec x_\rho|\right)
+ \sum_{p \in \partial_M I} \left[(\vec x_t\,.\,\vec\ek_1)\,
\vec{\rm m}\,.\,\vec y_{\mathcal{S}} + 2\,\pi\,\alpha_G\,\,
\vec{\rm m}_t\,.\,\vec\ek_1 \right](p)
\nonumber \\ & \quad 
= - 2\,\pi
\left(\alpha\, (\varkappa_{\mathcal{S}} - \spont) + \beta\,\AADE_{\mathcal{S}}, 
\vec x\,.\,\vec\ek_1\,[\varkappa_{\mathcal{S}}]_t\,|\vec x_\rho|\right)
+ \sum_{p \in \partial_M I} \left[(\vec x_t\,.\,\vec\ek_1)\,
\vec{\rm m}\,.\,\vec y_{\mathcal{S}} + 2\,\pi\,\alpha_G\,\,
\vec{\rm m}_t\,.\,\vec\ek_1 \right](p)\,.
\label{eq:new6}
\end{align}
Combining (\ref{eq:Lag2dxxt}) and (\ref{eq:new6}) yields that
\begin{align}
& 2\,\pi\left(\vec x\,.\,\vec\ek_1\,(\vec x_t\,.\,\vec\nu)^2, 
|\vec x_\rho|\right) 
 = -\pi \left(\alpha \left[\varkappa_{\mathcal{S}} -\spont \right]^2 
+ 2\,\lambda + 2\,\beta\,\AADE_{\mathcal{S}}\,\varkappa_{\mathcal{S}}, 
\left[\vec x\,.\,\vec\ek_1\,|\vec x_\rho|\right]_t\right) 
\nonumber \\ & \quad
- 2\,\pi
\left(\alpha\, (\varkappa_{\mathcal{S}} - \spont) + \beta\,\AADE_{\mathcal{S}}, 
\vec x\,.\,\vec\ek_1\,[\varkappa_{\mathcal{S}}]_t\,|\vec x_\rho|\right)
- 2\,\pi\,\varsigma\,\sum_{p \in \partial_{2} I \cup \partial_F I}
\vec x_t(p)\,.\,\vec\ek_1 
+2\,\pi\,\alpha_G\,\sum_{p \in \partial_M I} \vec{\rm m}_t(p)\,.\,
 \vec\ek_1
\nonumber \\ & \
= - \pi\,\ddt \left(\alpha\left
[\varkappa_{\mathcal{S}} -\spont \right]^2 + { 2\,\lambda},
\vec x\,.\,\vec\ek_1\,|\vec x_\rho|\right) 
- \tfrac\beta2\,\ddt\,\AADE_{\mathcal{S}}^2 
- 2\,\pi\,\varsigma\,\sum_{p \in \partial_{2} I \cup \partial_F I}
\vec x_t(p)\,.\,\vec\ek_1 
+2\,\pi\,\alpha_G\,\sum_{p \in \partial_M I} \vec{\rm m}_t(p)\,.\,
 \vec\ek_1
\nonumber \\ & \
= - \ddt \, \widehat E(t)\,,
\label{eq:PxtSstab}
\end{align}
where we have recalled the definition (\ref{eq:EEaxCN}). 
Of course, Remark~\ref{rem:Pxtstab} also applies to \eqref{eq:PxtSstab}.
   
In order to derive a suitable weak formulation, we now return to
(\ref{eq:Lag2dx}). Using (\ref{eq:24}) and noting (\ref{eq:tau}), 
(\ref{eq:Lag2dx}) can be rewritten as 
\begin{align}
& 
2\,\pi\left((\vec x\,.\,\vec\ek_1)\,\vec x_t\,.\,\vec\nu, \vec\chi\,.\,\vec\nu
\,|\vec x_\rho|\right) 
=-\left(\pi\,\vec x\,.\,\vec\ek_1\left[\alpha\,
\left(\varkappa_{\mathcal{S}}-\spont \right)^2 
+ 2\,\lambda + 2\,\beta\,\AADE_{\mathcal{S}}\,\varkappa_{\mathcal{S}}\right]
-\vec y_{\mathcal{S}}\,.\,\vec\ek_1, \vec\chi_\rho\,.\,\vec\tau \right) 
\nonumber \\ & \
- \left( \left[ \pi\,[\alpha\left(\varkappa_{\mathcal{S}}-\spont \right)^2
+ 2\,\lambda + 2\,\beta\,\AADE_{\mathcal{S}}\,\varkappa_{\mathcal{S}}]
- \varkappa_{\mathcal{S}}\,\vec y_{\mathcal{S}}\,.\,\vec\nu\right]\, |\vec x_\rho|
-(\vec y_{\mathcal{S}})_\rho\,.\,\vec\tau, 
\vec\chi\,.\,\vec\ek_1 \right)
\nonumber \\ & \ 
- \left(
\vec x\,.\,\vec\ek_1\,\varkappa_{\mathcal{S}}\,\vec y_{\mathcal{S}},
\vec\chi^{\perp}_\rho \right)
+ \left((\vec x\,.\,\vec\ek_1)\,(\vec y_{\mathcal{S}})_\rho \,.\,\vec\nu , 
 \vec\chi_\rho\,.\,\vec\nu \,|\vec x_{\rho}|^{-1}\right)
-\sum_{p \in \partial_{2} I \cup \partial_F I}
\left[
\left[2\,\pi\,\varsigma + \vec{\rm m}\,.\,\vec y_{\mathcal{S}}
\right] \vec\chi\,.\,\vec\ek_1 \right](p) 
\quad \forall \ \vec\chi \in \xspace \,.
\label{eq:Lag2dx2}
\end{align}   

Overall, we obtain the following weak formulation from (\ref{eq:Lag2dx2}),
(\ref{eq:kappaSy}), (\ref{eq:kappaSeta}), (\ref{eq:mSy}) and (\ref{eq:ASxkappa}), 
on recalling $\varkappa_{\mathcal{S}}^\star = \varkappa_{\mathcal{S}}$,
(\ref{eq:abperp}) and (\ref{eq:tau}). 

\noindent
$(\BGNpwfwf)$
Let $\vec x(\cdot,0) \in \Vpartialzero$ and $\alpha \in \bRplus$,
$\spont,M_0,\alpha_G,\lambda,\varsigma \in \bR$, $\beta \in \bRgeq$, 
$\vec\zeta : \partial_C I \to \bS^1$ be given.
For $t \in (0,T]$, find $\vec x(\cdot,t) \in \Vpartialzero$, 
with $\vec x_t(\cdot,t) \in
\xspace$, $\varkappa_{\mathcal{S}}(\cdot,t) \in L^2(I)$, 
$\vec y_{\mathcal{S}}(\cdot,t) \in \Vpartialzero$, with
$[(\vec x\,.\,\vec\ek_1)\,\vec y_{\mathcal{S}}](\cdot,t) \in 
\yspace(2\,\pi\,\alpha_G\,\vec\ek_1)$, 
and $\vec{\rm m}(\cdot,t) : \partial_M I \to \bR^2$ such that
\begin{subequations} \label{eq:PS}
\begin{align}
& 
2\,\pi\left((\vec x\,.\,\vec\ek_1)\,\vec x_t\,.\,\vec\nu, \vec\chi\,.\,\vec\nu\,
|\vec x_\rho|\right) 
- \left((\vec x\,.\,\vec\ek_1)\,(\vec y_{\mathcal{S}})_\rho \,.\,\vec\nu , 
 \vec\chi_\rho\,.\,\vec\nu \,|\vec x_{\rho}|^{-1}\right)
\nonumber \\ & \quad
=-\left(\pi\,\vec x\,.\,\vec\ek_1\left[\alpha\,
\left(\varkappa_{\mathcal{S}}-\spont \right)^2 
+ 2\,\lambda + 2\,\beta\,\AADE_{\mathcal{S}}\,\varkappa_{\mathcal{S}}\right]
-\vec y_{\mathcal{S}}\,.\,\vec\ek_1, \vec\chi_\rho\,.\,\vec\tau \right) 
\nonumber \\ & \qquad 
- \left( \left[ \pi\,[\alpha\left(\varkappa_{\mathcal{S}}-\spont \right)^2
+ 2\,\lambda + 2\,\beta\,\AADE_{\mathcal{S}}\,\varkappa_{\mathcal{S}}]
- \varkappa_{\mathcal{S}}\,\vec y_{\mathcal{S}}\,.\,\vec\nu\right]\, |\vec x_\rho|
-(\vec y_{\mathcal{S}})_\rho\,.\,\vec\tau, 
\vec\chi\,.\,\vec\ek_1 \right)
\nonumber \\
& \qquad + \left(
\vec x\,.\,\vec\ek_1\,\varkappa_{\mathcal{S}}\,\vec y_{\mathcal{S}}^\perp,
\vec\chi_\rho \right) 
-\sum_{p \in \partial_{2} I \cup \partial_F I}
\left[
\left[2\,\pi\,\varsigma + \vec{\rm m}\,.\,\vec y_{\mathcal{S}}
\right] \vec\chi\,.\,\vec\ek_1 \right](p) 
\qquad \forall \ \vec\chi \in \xspace \,,
\label{eq:PSa} \\
&
2\,\pi
\left( \vec x\,.\,\vec\ek_1\left[\alpha\,(\varkappa_{\mathcal{S}} - \spont)
+ \beta\,\AADE_{\mathcal{S}}\right]
, \chi\,|\vec x_\rho|\right)
- \left((\vec x\,.\,\vec\ek_1)\,\vec y_{\mathcal{S}},
\chi\,\vec\nu\,|\vec x_\rho|\right) = 0 \qquad \forall\ \chi \in L^2(I)\,,
\label{eq:PSb} \\
& \left(
\vec x\,.\,\vec\ek_1\,\varkappa_{\mathcal{S}}\,\vec\nu + \vec\ek_1, 
\vec\eta\,|\vec x_\rho|\right)
+ \left((\vec x\,.\,\vec\ek_1)\,\vec x_\rho ,\vec\eta_\rho \,|\vec x_\rho|^{-1}\right) 
= 
\sum_{p \in \partial_C I} \left[(\vec x\,.\,\vec\ek_1)\,
\vec\zeta\,.\,\vec\eta\right](p)
+ \sum_{p \in \partial_M I} \left[(\vec x\,.\,\vec\ek_1)\,
\vec{\rm m}\,.\,\vec\eta\right](p) 
\quad \forall\ \vec\eta \in \Vpartialzero\,,
\label{eq:PSc}
\end{align}
\end{subequations}
where $\AADE_{\mathcal{S}}(t)$ is given by (\ref{eq:ASxkappa}).
Similarly to (\ref{eq:weak3}),  
we note that the number of unknowns fixed via 
$(\vec x\,.\,\vec\ek_1)\,\vec y_{\mathcal{S}} 
\in \yspace(2\,\pi\,\alpha_G\,\vec\ek_1)$ on $\partial_M I$ 
is matched by the new
degrees of freedom arising from $\{\vec{\rm m}(p)\}_{p\in\partial_M I}$.

Similarly to (\ref{eq:ysnu}) and (\ref{eq:kyperp}), 
one can show that (\ref{eq:PS}) is independent of the tangential part
$\vec y_{\mathcal{S}} \,.\,\vec\tau$ of $\vec y_{\mathcal{S}}$. 
We refer to Appendix~\ref{sec:A2}, where we show that 
(\ref{eq:PS}) for a sufficiently smooth solution gives rise to the strong 
form (\ref{eq:xtbgnlambda}) and (\ref{eq:kappaS}).

\begin{rmrk} \label{rem:Lag2C}
Similarly to the procedure in {\rm \S\ref{sec:LagC}}, we can state a weak
formulation for the conserving flow \eqref{eq:xHelfrich}. In particular,
on writing \eqref{eq:PSa} as
\begin{equation*} 
2\,\pi \left((\vec x\,.\,\vec\ek_1)\,\vec x_t\,.\,\vec\nu,
\vec\chi\,.\,\vec\nu\,|\vec x_\rho|\right)
- \left((\vec x\,.\,\vec\ek_1) \,(\vec y_{\mathcal{S}})_\rho \,.\,\vec\nu, 
\vec\chi_\rho \,.\,\vec\nu \,|\vec x_\rho|^{-1} \right) =
\left( \vec f_{\mathcal{S}}, \vec\chi\,|\vec x_\rho| \right)
\qquad \forall\ \vec\chi \in \xspace
\,,
\end{equation*}
we can formulate the conserving flow as
\begin{align*} 
& 2\,\pi
\left((\vec x\,.\,\vec\ek_1)\,\vec x_t\,.\,\vec\nu,
\vec\chi\,.\,\vec\nu\,|\vec x_\rho|\right)
- \left((\vec x\,.\,\vec\ek_1) \,(\vec y_{\mathcal{S}})_\rho \,.\,\vec\nu, 
\vec\chi_\rho \,.\,\vec\nu \,|\vec x_\rho|^{-1} \right) 
= \left( \vec f_{\mathcal{S}}, \vec\chi\,|\vec x_\rho| \right)
\nonumber \\ &\qquad 
- 2\,\pi\, \lambda_A \left[
\left( \vec\ek_1, \vec\chi\,|\vec x_\rho|\right) 
+ \left( (\vec x\,.\,\vec\ek_1)\,\vec\tau
, \vec\chi_\rho\right) \right]
- 2\,\pi\, \lambda_V \left( (\vec x\,.\,\vec\ek_1)\,\vec\nu,
\vec\chi\,|\vec x_\rho|\right)
\quad \forall\ \vec\chi \in \xspace
\,,
\end{align*}
where $(\lambda_A(t),\lambda_V(t))^T \in \bR^2$ are chosen such that 
\eqref{eq:xsidedtSAV} holds.
\end{rmrk}

\setcounter{equation}{0}
\section{Semidiscrete schemes} \label{sec:sd}

Let $[0,1]=\bigcup_{j=1}^J I_j$, $J\geq3$, be a
decomposition of $[0,1]$ into intervals given by the nodes $q_j$,
$I_j=[q_{j-1},q_j]$. 
For simplicity, and without loss of generality,
we assume that the subintervals form an equipartitioning of $[0,1]$,
i.e.\ that 
\begin{equation} \label{eq:Jequi}
q_j = j\,h\,,\quad \mbox{with}\quad h = J^{-1}\,,\qquad j=0,\ldots, J\,.
\end{equation}
Clearly, if $I=\RZ$ we identify $0=q_0 = q_J=1$. In addition, we let $q_{J+1}=q_1$.

The necessary finite element spaces are defined as follows:
\begin{align*}
V^h = \{\chi \in C^{0}(\overline I) : \chi\!\mid_{I_j} 
\mbox{ is linear, $j=1,\ldots,J$}\}\,,\ V^h_0 = H^1_0(I)\cap V^h 
\quad\text{and}\quad
\Vh = [V^h]^2\,,\quad
\Vhzero = [V^h_0]^2\,.
\end{align*}
In addition, we define 
$\Vhpartialzero  = \Vh \cap \Vpartialzero$,
$W^h = V^h$, $\Whpartialzero = \{ \eta \in V^h : \eta(\rho) = 0 \ 
\forall\ \rho \in \partial_0 I\}$,
$\xspaceh = \xspace \cap \Vh$ and,
for a given $\vec z \in \bR^2$, $\yspaceh(\vec z) = \yspace(\vec z) \cap \Vh$.
Let $\{\chi_j\}_{j=j_0}^J$ denote the standard basis of $V^h$,
where $j_0 = 0$ if $I = (0,1)$ and $j_0 = 1$ if $I=\RZ$.
We also set $\jf = J-1$ if $I = (0,1)$
and $\jf = J$ if $I=\RZ$.
For later use, we let $\pi^h:C^{0}(\overline I)\to V^h$ 
be the standard interpolation operator at the nodes $\{q_j\}_{j=0}^J$,
and similarly $\pi^h_{\partial_0}:C^{0}(\overline I)\to \Whpartialzero$,
as well as $\vec\pi^h:[C^{0}(\overline I)]^2 \to \Vh$.

Let $(\cdot,\cdot)$ denote the $L^2$--inner product on $I$, and 
define the mass lumped $L^2$--inner product $(f,g)^h$,
for two piecewise continuous functions, with possible jumps at the 
nodes $\{q_j\}_{j=1}^J$, via
\begin{equation}
( f, g )^h = \tfrac12\,h\,\sum_{j=1}^J
\left[(f\,g)(q_j^-) + (f\,g)(q_{j-1}^+)\right], \label{eq:ip0}
\end{equation}
where we define
$f(q_j^\pm)=\underset{\delta\searrow 0}{\lim}\ f(q_j\pm\delta)$.
The definition (\ref{eq:ip0}) naturally extends to vector valued functions.

Let $(\vec X^h(t))_{t\in[0,T]}$, 
with $\vec X^h(\cdot,t)\in \Vhpartialzero$,
be an approximation to $(\vec x(t))_{t\in[0,T]}$ and define
$\Gamma^h(t) = \vec X^h(\overline I,t)$. 
A natural discrete analogue of the well-posedness assumptions for 
the continuous solution $\vec x$ is given as follows.

\begin{ass} \label{ass:Ah1}
Let
\begin{equation} \label{eq:Xhpos}
\vec X^h(\rho,t) \,.\,\vec\ek_1 > 0 \quad 
\forall\ \rho \in \overline I\setminus \partial_0 I\,,\ t \in [0,T]\,.
\end{equation}
In addition, let $\vec X^h(q_j,t) \ne \vec X^h(q_{j+1},t)$,
$j=0,\ldots,J-1$, for all $t\in [0,T]$.
\end{ass}

Then, similarly to (\ref{eq:tau}), we set
\begin{equation} \label{eq:tauh}
\vec\tau^h = \vec X^h_s = \frac{\vec X^h_\rho}{|\vec X^h_\rho|} 
\qquad \mbox{and} \qquad \vec\nu^h = -(\vec\tau^h)^\perp
\qquad \text{in } \overline I\,,
\end{equation}
which is well-defined if Assumption~\ref{ass:Ah1} holds.
We note that \eqref{eq:Xhpos} implies $\vec\tau^h\,.\,\vec\ek_1 \not=0$
on elements touching the $x_2$--axis, and so
\begin{equation*} 
\vec\nu^h\,.\,\vec\ek_2 \not=0 \qquad\text{ on } \partial_0 I\,,
\end{equation*}
compare also with \eqref{eq:bcnu} and \eqref{eq:bc}.

\begin{ass} \label{ass:Ah2}
Let {\rm Assumption~\ref{ass:Ah1}} hold and let
$\vec X^h(q_{j-1},t) \ne \vec X^h(q_{j+1},t)$,
$j=1,\ldots,\jf$, for all $t\in [0,T]$.
\end{ass}
We note that Assumption~\ref{ass:Ah2} is only violated if two
neighbouring elements of $\Gamma^h(t)$ lie identically on top of each 
other. For our fully discrete schemes, this never happens in practice.
For later use, we let 
$\vec\omega^h \in \Vh$ be the mass-lumped 
$L^2$--projection of $\vec\nu^h$ onto $\Vh$, i.e.\
\begin{equation} \label{eq:omegah}
\left(\vec\omega^h, \vec\varphi \, |\vec X^h_\rho| \right)^h 
= \left( \vec\nu^h, \vec\varphi \, |\vec X^h_\rho| \right)
= \left( \vec\nu^h, \vec\varphi \, |\vec X^h_\rho| \right)^h
\qquad \forall\ \vec\varphi\in \Vh\,.
\end{equation}
Combining \eqref{eq:omegah}, \eqref{eq:ip0} and \eqref{eq:tauh} yields that
\begin{equation*} 
\vec\omega^h(q_j) 
 = \begin{cases}
- \dfrac{\left(\vec X^h(q_{j+1}) - \vec X^h(q_{j-1})\right)^\perp}
{|\vec X^h(q_{j+1}) - \vec X^h(q_j)|
+|\vec X^h(q_j) - \vec X^h(q_{j-1})|}
& q_j \in \overline I \setminus \partial I\,, \\
\vec\nu^h(q_j) & q_j \in \partial I\,.
\end{cases}
\end{equation*}
It follows that $\vec v^h \in \Vh$, defined by
\begin{equation} \label{eq:vh}
\vec v^h = \vec\pi^h\left[\frac{\vec\omega^h}{|\vec\omega^h|}\right],
\end{equation}
is well-defined if Assumption~\ref{ass:Ah2} holds.
We also define $\mat Q^h \in [V^h]^{2 \times 2}$ defined by
\begin{equation} \label{eq:Qh}
\mat Q^h(q_j) = \begin{cases}
\mat\Id & q_j \in \partial I \setminus \partial_0 I\,,\\
\vec v^h \otimes \vec v^h & q_j \in I \cup \partial_0 I\,.
\end{cases}
\end{equation}
Later on we will describe the evolution of $\Gamma^h(t)$ through 
$\vec\pi^h[\mat Q^h\,\vec X^h_t]$, for $\vec X^h_t \in \xspace^h$. This will
allow tangential motion for interior nodes, because we will only let a discrete
normal component of $\vec X^h_t$ be specified through an appropriate variation
of the discrete energy. But crucially, we will specify the full velocity
$\vec X^h_t$ through this energy variation at boundary nodes
$q_j \in \partial I \setminus \partial_0 I$. This is because at these boundary
nodes we cannot allow an arbitrary tangential motion, as this would affect the
evolution of $\Gamma^h(t)$ itself, and not just the evolution of its
parameterization $\vec X^h$. A similar strategy has been pursued by the authors
in \cite[(3.19)]{pwfopen}.

Similarly to (\ref{eq:24}), 
on noting (\ref{eq:A}) and (\ref{eq:tauh}), we have for all 
$\vec\chi \in \xspaceh$ on $I_j$, $j=1,\ldots,J,$ that 
\begin{subequations} \label{eq:57}
\begin{align}
\left[\deldel{\vec X^h}\, |\vec X^h_\rho|\right](\vec\chi)
& = \frac{\vec X^h_\rho\,.\,\vec\chi_\rho}{|\vec X^h_\rho|}
= \vec\tau^h\,.\,\vec\chi_\rho 
\,, \label{eq:ddxrhoh} \\
\left[\deldel{\vec X^h}\, \vec\tau^h\right](\vec\chi)
& = \left[\deldel{\vec X^h}\,\frac{\vec X^h_\rho}{|\vec X^h_\rho|}
\right](\vec\chi) = \frac{\vec\chi_\rho}{|\vec X^h_\rho|}
- \frac{\vec X^h_\rho}{|\vec X^h_\rho|^2}\,
 \frac{\vec X^h_\rho\,.\,\vec\chi_\rho}{|\vec X^h_\rho|}
= \vec\chi_s - \vec\tau^h\,(\vec\chi_s\,.\,\vec\tau^h)
= (\vec\chi_s\,.\,\vec\nu^h)\,\vec\nu^h \,, \label{eq:ddtauh} \\
\left[\deldel{\vec X^h}\, \vec\nu^h\right](\vec\chi)
& = - \left[\deldel{\vec X^h}\, (\vec\tau^h)^\perp\right](\vec\chi) = 
- (\vec\chi_s\,.\,\vec\nu^h)\,(\vec\nu^h)^\perp =
- (\vec\chi_s\,.\,\vec\nu^h)\,\vec\tau^h \,, \label{eq:ddnuh} \\
\left[\deldel{\vec X^h}\,\vec\nu^h\,|\vec X^h_\rho|\right](\vec\chi) &=
- \left[\deldel{\vec X^h}\,(\vec X^h_{\rho})^\perp\right](\vec\chi)
= - \vec\chi_\rho^\perp 
\,. \label{eq:ddnuxrhoh} 
\end{align}
\end{subequations}
In addition to (\ref{eq:57}), we will require 
$\left[\deldel{\vec X^h}\, \vec\omega^h\right](\vec\chi)$.
It follows from (\ref{eq:omegah}), (\ref{eq:57}) and (\ref{eq:tauh}) that
\begin{align}
\left(\left[\deldel{\vec X^h}\, \vec\omega^h\right](\vec\chi), 
\vec\varphi \, |\vec X^h_\rho| \right)^h 
& = \left(\left[\deldel{\vec X^h}\, \vec\nu^h\right](\vec\chi), 
\vec\varphi \, |\vec X^h_\rho| \right)^h
- \left( \vec\omega^h - \vec\nu^h,\vec\varphi 
\left[\deldel{\vec X^h}\, |\vec X^h_\rho|\right](\vec\chi)
\right)^h  \nonumber \\
& = -\left( (\vec\nu^h \,.\,\vec\chi_\rho)\,\vec\tau^h , 
\vec\varphi \right)^h
- \left( 
\vec\tau^h\,.\,\vec\chi_\rho\,
(\vec\omega^h - \vec\nu^h)
,\vec\varphi 
\right)^h 
\qquad \forall\ \vec\varphi\in\Vhpartialzero\,.
\label{eq:ddomegah}
\end{align}

\subsection{Based on $\kappa^h$} \label{sec:sdkappa}

As the discrete analogue of (\ref{eq:varkappaweak}), 
we let $\vec X^h \in \Vhpartialzero$,
$\kappa^h \in V^h$ and $\vec{\rm m}^h : \partial_M I \to \bR^2$
be such that
\begin{equation} \label{eq:sideh}
\left( \kappa^h\,\vec\nu^h, \vec\eta \,|\vec X^h_\rho|\right)^h
+ \left( \vec\tau^h, \vec\eta_\rho \right) = 
\sum_{p \in \partial_C I} \left[\vec\zeta\,.\,\vec\eta\right](p)
+ \sum_{p \in \partial_M I} \left[\vec{\rm m}^h\,.\,\vec\eta\right](p)
\qquad \forall\ \vec\eta \in \Vhpartialzero\,,
\end{equation}
where we recall (\ref{eq:tauh}).

We would like to mimic on the discrete level the procedure in
Section~\ref{sec:Lag}. However, a naive discretization of (\ref{eq:Lag}) will
not give a well-defined Lagrangian, since a discrete variant of
(\ref{eq:bclimit}) will in general not hold. To overcome the arising
singularity in a discretization of (\ref{eq:Lag}), we now introduce
the following discrete approximation of $\varkappa_{\mathcal{S}}$, which will
be based on $\kappa^h$.
In particular, on
recalling (\ref{eq:bclimit}) and (\ref{eq:omegah}), 
we introduce,
given $\vec X^h \in \Vhpartialzero$ and $\kappa^h \in V^h$, the function 
$\doctorkappa^h(\vec X^h, \kappa^h) \in V^h$ such that
\begin{equation} \label{eq:calKh}
[\doctorkappa^h (\vec X^h,\kappa^h)](q_j) = \begin{cases}
\kappa^h(q_j) - 
\dfrac{\vec\omega^h(q_j)\,.\,\vec\ek_1}{\vec X^h(q_j)\,.\,\vec\ek_1}
& q_j \in \overline I \setminus \partial_0 I\,, \\
2\, \kappa^h(q_j) & q_j \in \partial_0 I\,.
\end{cases}
\end{equation}
Clearly, using $\vec\nu^h$ in place of $\vec\omega^h$ in 
\eqref{eq:calKh} would not be well-defined for interior nodes.
For later use we also define $\doctorZ^h \in V^h$ such that
\begin{equation} \label{eq:doctorZh}
\doctorZ^h (q_j) = \begin{cases}
1 & q_j \in \overline I \setminus \partial_0 I\,, \\
2 & q_j \in \partial_0 I\,.
\end{cases}
\end{equation}

On noting (\ref{eq:calKh}), we define the discrete analogue of the energy (\ref{eq:EEaxCN})
\begin{align} 
\widehat E^h(t) &= \pi\left(\alpha 
\left[ \doctorkappa^h(\vec X^h, \kappa^h) - \spont \right]^2 +2 \,\lambda,  
\vec X^h\,.\,\vec\ek_1 \,|\vec X^h_\rho|\right)^h + \tfrac{\beta}{2}
\left(\AADE^h(t)\right)^2 \nonumber \\ & \qquad
- 2\,\pi\,\alpha_G\sum_{p\in \partial_M I}
\vec{\rm m}^h(p)\,.\,\vec\ek_1
+ 2\,\pi\,\varsigma\sum_{p\in \partial_2 I \cup \partial_F I}
\vec X^h(p)\,.\,\vec\ek_1\,,
\label{eq:Eh}
\end{align}
where   
\begin{equation} \label{eq:Ah}
\AADE^h(t) = 
2\,\pi\left((\vec X^h\,.\,\vec\ek_1)\,\kappa^h - \vec\nu^h.\,\vec\ek_1, 
|\vec X^h_\rho| \right)^h - M_0\,.
\end{equation}

\begin{rmrk} \label{rem:kappahpartial0}
We observe that the energy $\widehat E^h(t)$ does not depend
on the values of $\kappa^h$ on $\partial_0 I$. We will thus fix these values to
be zero from now on. A welcome side effect of this procedure is that on assuming
that e.g.\ $0 \in \partial_0 I$, then choosing
$\vec\eta = \chi_0\,\vec\ek_2$ in \eqref{eq:sideh} yields that
$(\vec X^h (q_1) - \vec X^h (q_0) ) \,.\,\vec\ek_2 = 0$.
Similarly we get $(\vec X^h (q_J) - \vec X^h (q_{J-1}) ) \,.\,\vec\ek_2 = 0$
if $1 \in \partial_0 I$.

Without fixing $\kappa^h$ to be zero on $\partial_0 I$, we observe numerical
difficulties in practice for fully discrete variants of the semidiscrete
approximation that we are going to derive.
\end{rmrk}

Similarly to (\ref{eq:Lag}), we define the discrete Lagrangian
\begin{align*} 
& \mathcal{L}^h(\vec X^h, \kappa^h, \vec{\rm m}^h, \vec Y^h) =
\pi\left(\alpha 
\left[ \doctorkappa^h(\vec X^h, \kappa^h) - \spont \right]^2 +2 \,\lambda,  
\vec X^h\,.\,\vec\ek_1 \,|\vec X^h_\rho|\right)^h 
\nonumber \\ & \
+ \tfrac{\beta}{2}
\left[
2\,\pi\left((\vec X^h\,.\,\vec\ek_1)\,\kappa^h - \vec\nu^h.\,\vec\ek_1, 
|\vec X^h_\rho| \right)^h - M_0
\right]^2 
- \left( \kappa^h\,\vec\nu^h, \vec Y^h \,|\vec X^h_\rho|\right)^h
- \left( \vec\tau^h, \vec Y^h_\rho \right)
+ 2\,\pi\,\varsigma\sum_{p\in \partial_2 I \cup \partial_F I}
\vec X^h(p)\,.\,\vec\ek_1
\nonumber \\ & \
+\sum_{p\in \partial_C I} \left[ \vec\zeta \,.\,
\vec Y^h \right](p)
+\sum_{p\in \partial_M I} \left[ \vec{\rm m}^h \,.
\left( \vec Y^h - 2\,\pi\,\alpha_G \,\vec\ek_1\right) \right](p)\,,
\end{align*}
for the minimization of the energy \eqref{eq:Eh} subject to the side constraint
\eqref{eq:sideh}, 
where $\vec X^h \in \Vhpartialzero$, $\kappa^h \in \Whpartialzero$,
$\vec{\rm m}^h : \partial_M I \to \bR^2$
and $\vec Y^h \in \Vhpartialzero$.

Taking variations $\vec\eta \in \Vhpartialzero$ in $\vec Y^h$, and setting
$\left[\deldel{\vec Y^h}\, \mathcal{L}^h\right](\vec\eta) = 0$
we obtain (\ref{eq:sideh}), similarly to (\ref{eq:kappaeta}). 
Taking variations $\chi \in \Whpartialzero$ in $\kappa^h$ and setting
$\left[\deldel{\kappa^h}\, \mathcal{L}^h\right](\chi) = 0$ we obtain 
\begin{equation} 
2\,\pi\left( \vec X^h\,.\,\vec\ek_1 \left(\alpha\,[\doctorkappa^h(
\vec X^h,\kappa^h) - \spont ] +\beta\,\AADE^h\right),\chi 
\,|\vec X^h_\rho|\right)^h
- \left( \vec Y^h, \chi\,\vec\nu^h \,|\vec X^h_\rho|\right)^h = 0
\qquad \forall\ \chi \in \Whpartialzero\,, \label{eq:kappahYh}
\end{equation}
where we have recalled (\ref{eq:calKh}). 
Taking variations in $\vec{\rm m}^h$, and setting them to zero, yields,
similarly to (\ref{eq:my}), that
\begin{equation} \label{eq:mhy}
\vec Y^h = 2\,\pi\,\alpha_G\,\vec\ek_1 \qquad\text{on }\  \partial_M I\,.
\end{equation}
Taking variations $\vec\chi \in \xspaceh$ in $\vec X^h$, and 
setting $2\,\pi\,((\vec X^h\,.\,\vec\ek_1)\,\mat Q^h\,\vec X^h_t,\vec\chi
\,|\vec X^h_\rho|)^h = 
- \left[\deldel{\vec X^h}\, \mathcal{L}^h\right](\vec\chi)$
we obtain  
\begin{align}
& 2\,\pi\left((\vec X^h\,.\,\vec\ek_1)\,\mat Q^h\,\vec X^h_t,
\vec\chi\,|\vec X^h_\rho|\right)^h 
= - \pi \left( \alpha
\left[ \doctorkappa^h(\vec X^h,\kappa^h) - \spont \right]^2
+ 2\,\lambda + 2\,\beta\,\AADE^h\,\kappa^h,
\left[\deldel{\vec X^h}\,(\vec X^h\,.\,\vec\ek_1)\,|\vec X^h_\rho|\right]
(\vec\chi) \right)^h 
\nonumber \\ 
& \
- 2\,\pi\,\alpha \left( 
\left[ \doctorkappa^h(\vec X^h,\kappa^h) - \spont \right], 
\left[\deldel{\vec X^h}\,\doctorkappa^h(\vec X^h,\kappa^h)\right]
(\vec\chi)\,(\vec X^h\,.\,\vec\ek_1)\,|\vec X^h_\rho| \right)^h 
\nonumber \\ & \
+ \left( \kappa^h\, \vec Y^h + 2\,\pi\,\beta\,\AADE^h\,\vec\ek_1,
\left[\deldel{\vec X^h}\,\vec\nu^h\,|\vec X^h_\rho|\right](\vec\chi)\right)^h
+ \left(\vec Y^h_\rho , 
\left[\deldel{\vec X^h}\,\vec\tau^h\right](\vec\chi)\right) 
- 2\,\pi\,\varsigma\sum_{p\in \partial_{2} I \cup \partial_F I}
\vec\chi(p)\,.\,\vec\ek_1
\quad \forall\ \vec\chi \in \xspaceh\,.
\label{eq:Pxth}
\end{align}
Choosing $\vec\chi = \vec X^h_t$ in (\ref{eq:Pxth}) yields, on noting 
(\ref{eq:ddA}), the discrete analogue of (\ref{eq:Pxtxt})
\begin{align}
& 2\,\pi\left(\vec X^h\,.\,\vec\ek_1\,|\mat Q^h\,\vec X^h_t|^2,
|\vec X^h_\rho|\right)^h 
 = - \pi \left( 
\alpha \left[ \doctorkappa^h(\vec X^h,\kappa^h) - \spont \right]^2
+ 2\,\lambda + 2\,\beta\,\AADE^h\,\kappa^h, 
\left[(\vec X^h\,.\,\vec\ek_1)\,|\vec X^h_\rho|\right]_t \right)^h 
\nonumber \\ & \qquad \qquad 
- 2\,\pi\,\alpha \left( 
 \doctorkappa^h(\vec X^h,\kappa^h) - \spont , 
\left[\dfrac{\vec\omega^h\,.\,\vec\ek_1}{\vec X^h\,.\,\vec\ek_1}
\right]_t (\doctorZ^h - 2)
\,(\vec X^h\,.\,\vec\ek_1)\,|\vec X^h_\rho| \right)^h 
\nonumber \\ & \qquad \qquad 
+ \left( \kappa^h\, \vec Y^h + 2\,\pi\,\beta\,\AADE^h\,\vec\ek_1,
\left[\vec\nu^h\,|\vec X^h_\rho|\right]_t \right)^h
+ \left(\vec Y^h_\rho , \vec\tau^h_t \right) 
- 2\,\pi\,\varsigma\sum_{p\in \partial_{2} I \cup \partial_F I}
\vec X^h_t(p)\,.\,\vec\ek_1\, .
\label{eq:Pxtxth}
\end{align}

Differentiating (\ref{eq:sideh}) with respect to $t$, and then choosing
$\vec\eta = \vec Y^h \in \Vhpartialzero$ yields, 
similarly to (\ref{eq:dtside}), that 
\begin{equation}
 \left( \kappa^h_t, 
\vec Y^h\,.\,\vec\nu^h\,|\vec X^h_\rho|\right)^h
+ \left( \kappa^h\,\vec Y^h, \left[\vec\nu^h\,|\vec X^h_\rho|\right]_t\right)^h
+ \left(\vec\tau^h_t,\vec Y^h_\rho\right)
= \sum_{p \in \partial_M I} \left[ \vec{\rm m}^h_t\,.\,\vec Y^h\right](p)\,.
\label{eq:dtyyh}
\end{equation}
It follows from (\ref{eq:dtyyh}), (\ref{eq:mhy}) and (\ref{eq:kappahYh}) 
with $\chi = \kappa^h_t \in \Whpartialzero$ that
\begin{align}
&\left( \kappa^h\,\vec Y^h, \left[\vec\nu^h\,|\vec X^h_\rho|\right]_t\right)^h
+ \left(\vec\tau^h_t,\vec Y^h_\rho\right)
 = -\left( \kappa^h_t, 
\vec Y^h\,.\,\vec\nu^h\,|\vec X^h_\rho|\right)^h
+ 2\,\pi\, \alpha_G \sum_{p \in \partial_M I} \left[ \vec{\rm m}^h_t\,.\,\vec\ek_1 \right](p)
\nonumber \\
& \qquad = -2\,\pi\left( \vec X^h\,.\,\vec\ek_1 \left( \alpha\,[\doctorkappa^h(
\vec X^h,\kappa^h) - \spont ] + \beta\, \AADE^h \right), 
\kappa^h_t \,|\vec X^h_\rho| 
\right)^h 
 + 2\,\pi\, \alpha_G 
\sum_{p \in \partial_M I} \left[ \vec{\rm m}^h_t(p)\,.\,\vec\ek_1 \right]\,.
\label{eq:mykappah}
\end{align}
Combining (\ref{eq:Pxtxth}) 
and (\ref{eq:mykappah}) 
yields, on recalling (\ref{eq:Eh}), the discrete analogue of (\ref{eq:Pxtstab})
\begin{align}
& 2\,\pi\left(\vec X^h\,.\,\vec\ek_1\,|\mat Q^h\,\vec X^h_t|^2,
|\vec X^h_\rho|\right)^h 
 = - \pi \left( \alpha
\left[ \doctorkappa^h(\vec X^h,\kappa^h) - \spont \right]^2
+ 2\,\lambda + 2\,\beta\,\AADE^h\,\kappa^h, 
\left[(\vec X^h\,.\,\vec\ek_1)\,|\vec X^h_\rho|\right]_t \right)^h 
\nonumber \\ & \quad
- 2\,\pi\,\alpha \left( 
 \doctorkappa^h(\vec X^h,\kappa^h) - \spont , 
\left[\dfrac{\vec\omega^h\,.\,\vec\ek_1}{\vec X^h\,.\,\vec\ek_1}
\right]_t (\doctorZ^h - 2)
\,(\vec X^h\,.\,\vec\ek_1)\,|\vec X^h_\rho| \right)^h 
\nonumber \\ & \quad
-2\,\pi\left( \vec X^h\,.\,\vec\ek_1 \left(\alpha\,[\doctorkappa^h(
\vec X^h,\kappa^h) - \spont ] + \beta\, \AADE^h \right), 
\kappa^h_t \,|\vec X^h_\rho| 
\right)^h
+ 2\,\pi\,\beta\,\AADE^h \left( \vec\ek_1, \left[ \vec\nu^h\,|\vec X^h_\rho|\right]_t
\right)^h
\nonumber \\ & \quad
+ 2\,\pi\, \alpha_G 
\sum_{p \in \partial_M I} \left[ \vec{\rm m}^h_t(p)\,.\,\vec\ek_1 \right]
- 2\,\pi\,\varsigma\sum_{p\in \partial_{2} I \cup \partial_F I}
\vec X^h_t(p)\,.\,\vec\ek_1
\nonumber \\ & \
= - \ddt \,\widehat E^h(t)\, .
\label{eq:Pxtstabh}
\end{align}

We now return to (\ref{eq:Pxth}), which, on recalling (\ref{eq:57}),
(\ref{eq:calKh}), (\ref{eq:doctorZh})
and (\ref{eq:ddomegah}), can be rewritten as
\begin{align}
& 2\,\pi\left((\vec X^h\,.\,\vec\ek_1)\,\mat Q^h\,\vec X^h_t,
\vec\chi\,|\vec X^h_\rho|\right)^h = 
- \pi \left( \alpha
\left[ \doctorkappa^h(\vec X^h,\kappa^h) - \spont \right]^2 
+ 2\,\lambda + 2\,\beta\,\AADE^h\,\kappa^h,
\vec\chi\,.\,\vec\ek_1\,|\vec X^h_\rho| + (\vec X^h\,.\,\vec\ek_1)\,\vec\tau^h\,.\,
\vec\chi_\rho \right)^h
\nonumber \\ & \quad 
+ 2\,\pi\,\alpha\left( \left[\doctorkappa^h(\vec X^h,\kappa^h) - \spont\right] 
 (\doctorZ^h - 2), 
\frac{\vec\omega^h\,.\,\vec\ek_1}{\vec X^h\,.\,\vec e_1}\,\vec\chi\,.\,\vec\ek_1
\,|\vec X^h_\rho| \right)^h 
\nonumber \\ & \quad 
+ 2\,\pi\,\alpha \left( \left[\doctorkappa^h(\vec X^h,\kappa^h) - \spont\right]
(\doctorZ^h - 2)\,\vec\ek_1, 
(\vec\nu^h\,.\,\vec\chi_\rho) \,\vec\tau^h + (\vec\tau^h\,.\,\vec\chi_\rho) \,
(\vec\omega^h-\vec\nu^h) \right)^h 
\nonumber \\ & \quad 
- \left( \kappa^h\,\vec Y^h + 2\,\pi\,\beta\,\AADE^h\,\vec e_1,
\vec\chi_\rho^\perp \right)^h
+ \left(\vec Y^h_\rho\,.\,\vec\nu^h , \vec\chi_\rho\,.\,\vec\nu^h\,|\vec X^h_\rho|^{-1}
\right) 
- 2\,\pi\,\varsigma\sum_{p\in \partial_{2} I \cup \partial_F I}
\vec\chi(p)\,.\,\vec\ek_1
\qquad \forall\ \vec\chi \in \xspaceh\,.
\label{eq:Pxt2h}
\end{align}
Combining (\ref{eq:Pxt2h}), (\ref{eq:kappahYh}), (\ref{eq:sideh}),
(\ref{eq:mhy}) and (\ref{eq:Ah}), our semidiscrete approximation based on 
$\kappa^h$ is given, on noting (\ref{eq:abperp}) and (\ref{eq:tauh}), 
as follows. 

\noindent
$(\BGNpwf^h)^{h}$
Let $\vec X^h(\cdot,0) \in \Vhpartialzero$ and $\alpha \in \bRplus$,
$\spont,M_0,\alpha_G,\lambda,\varsigma \in \bR$, $\beta \in \bRgeq$, 
$\vec\zeta : \partial_C I \to \bS^1$ be given.
Then, for $t \in (0,T]$ find
$\vec X^h(\cdot,t) \in \Vhpartialzero$, 
with $\vec X^h_t(\cdot,t) \in \xspace^h$,
$\kappa^h(\cdot,t) \in \Whpartialzero$, 
$\vec Y^h(\cdot,t) \in \yspaceh(2 \,\pi \, \alpha_G\,\vec\ek_1)$ 
and $\vec{\rm m}^h(\cdot,t) : \partial_M I \to \mathbb R^2$ such that
\begin{subequations} \label{eq:sd}
\begin{align}
& 2\,\pi\left((\vec X^h\,.\,\vec\ek_1)\,\mat Q^h\,\vec X^h_t,
\vec\chi\,|\vec X^h_\rho|\right)^h 
- \left(\vec Y^h_\rho\,.\,\vec\nu^h , \vec\chi_\rho\,.\,\vec\nu^h\,
|\vec X^h_\rho|^{-1} \right) 
\nonumber \\ & \quad 
= - \pi \left( 
\alpha \left[ \doctorkappa^h(\vec X^h,\kappa^h) - \spont \right]^2
+ 2\,\lambda + 2\,\beta\,\AADE^h\,\kappa^h, 
\vec\chi\,.\,\vec\ek_1\,|\vec X^h_\rho| + (\vec X^h\,.\,\vec\ek_1)\,\vec\tau^h\,.\,
\vec\chi_\rho \right)^h
\nonumber \\ & \qquad 
+ 2\,\pi\,\alpha\left( \left[\doctorkappa^h(\vec X^h,\kappa^h) - \spont\right] 
 (\doctorZ^h - 2), 
\frac{\vec\omega^h\,.\,\vec\ek_1}{\vec X^h\,.\,\vec\ek_1}\,\vec\chi\,.\,\vec\ek_1
\,|\vec X^h_\rho| \right)^h 
\nonumber \\ & \qquad 
+ 2\,\pi\,\alpha \left( \left[\doctorkappa^h(\vec X^h,\kappa^h) - \spont\right]
(\doctorZ^h - 2)\,\vec\ek_1, 
(\vec\nu^h\,.\,\vec\chi_\rho) \,\vec\tau^h + (\vec\tau^h\,.\,\vec\chi_\rho) \,
(\vec\omega^h-\vec\nu^h) \right)^h 
\nonumber \\ & \qquad 
+ \left( \kappa^h\,(\vec Y^h)^\perp - 2\,\pi\,\beta\,\AADE^h\,\vec\ek_2 ,
\vec\chi_\rho \right)^h - 2\,\pi\,\varsigma\sum_{p\in \partial_{2} I \cup \partial_F I}
\vec\chi(p)\,.\,\vec\ek_1
\qquad \forall\ \vec\chi \in \xspaceh\,, \label{eq:sda} \\
& 
2\,\pi\left( \vec X^h\,.\,\vec\ek_1 \left( \alpha\,[\doctorkappa^h(
\vec X^h,\kappa^h) - \spont ] + \beta\,\AADE^h \right)
,\chi \,|\vec X^h_\rho|\right)^h
- \left( \vec Y^h, \chi\,\vec\nu^h \,|\vec X^h_\rho|\right)^h = 0
 \qquad \forall\ \chi \in \Whpartialzero\,,\label{eq:sdb} \\
&\left( \kappa^h\,\vec\nu^h, \vec\eta \,|\vec X^h_\rho|\right)^h
+ \left( \vec X^h_\rho, \vec\eta_\rho\,|\vec X^h_\rho|^{-1} \right) = 
\sum_{p \in \partial_C I} \left[\vec\zeta\,.\,\vec\eta\right](p)
+ \sum_{p \in \partial_M I} \left[\vec{\rm m}^h\,.\,\vec\eta\right](p)
 \qquad
\forall\ \vec\eta \in \Vhpartialzero\,, \label{eq:sdc}
\end{align}
\end{subequations}
where $\AADE^h(t)$ is given by (\ref{eq:Ah}). 

\begin{thrm} \label{thm:stab}
Let {\rm Assumption~\ref{ass:Ah2}} be satisfied and 
let $(\vec X^h(t),\kappa^h(t),\vec Y^h(t),\vec{\rm m}^h(t))_{t\in (0,T]}$ 
be a solution to \mbox{\rm (\ref{eq:sd})}. 
Then the solution satisfies the stability bound
\[
\ddt \,\widehat E^h(t)
+ 2\,\pi\left(\vec X^h\,.\,\vec\ek_1\,|\mat Q^h\,\vec X^h_t|^2,
|\vec X^h_\rho|\right)^h=0\,. 
\]
\end{thrm}
\begin{proof}
The desired result follows as (\ref{eq:sd}) is just a rewrite of (\ref{eq:Pxth}),
(\ref{eq:kappahYh}), (\ref{eq:sideh}) and (\ref{eq:Ah}), and then noting 
(\ref{eq:Pxtxth})--(\ref{eq:Pxtstabh}). 
\end{proof}

\begin{rmrk} \label{rem:equid}
We note that on choosing 
$\vec\eta = \chi_j\,[\vec\omega^h(q_j)]^\perp$, for 
$j \in \{1,\ldots,\jf\}$ so $\vec\eta \in \Vhpartialzero$, 
in {\rm (\ref{eq:sdc})} we obtain that
\[
|\vec X^h(q_j) - \vec X^h(q_{j-1})| = 
|\vec X^h(q_{j+1}) - \vec X^h(q_{j})| 
\ \text{or} \
\vec X^h(q_j) - \vec X^h(q_{j-1}) \parallel
\vec X^h(q_{j+1}) - \vec X^h(q_{j})
\]
for $j=1,\ldots,\jf$. 
See {\rm \cite[Remark 2.4]{triplej}} 
for details. Hence the curve $\Gamma^h(t)$
will be equidistributed where-ever two neighbouring elements are not parallel.
This aspect of the solution for the 
semidiscrete problem \eqref{eq:sd} means that it can be viewed as a 
highly nonlinear and degenerate system of differential-algebraic equations,
see also \cite[Remark~81]{bgnreview} for a related discussion.
In particular, at present we are unable to prove the existence of
solutions to \eqref{eq:sd}. In addition, an error analysis for our
semidiscrete approximations also appears to be out of reach.
\end{rmrk}

\begin{rmrk} \label{rem:kappahpartial02}
We now revisit the discussion in {\rm Remark~\ref{rem:kappahpartial0}}.
Assuming that e.g.\ $0 = q_0 \in \partial_0 I$, then choosing
$\vec\eta = \chi_0\,\vec\ek_2$ in \eqref{eq:sdc} yields that
$\kappa^h(q_0)\,\vec\nu^h(q_0)\,.\,\vec\ek_2
= 2 \frac{(\vec X^h (q_1) - \vec X^h (q_0) )\,.\,\vec\ek_2} 
{|\vec X^h (q_1) - \vec X^h (q_0)|^2}
$,
and so
\begin{equation} \label{eq:tauhe2}
\vec\tau^h(q_0)\,.\,\vec\ek_2 = \tfrac12\,|\vec X^h (q_1) - \vec X^h (q_0)|\,
\kappa^h(q_0)\,\vec\nu^h(q_0)\,.\,\vec\ek_2\,,
\end{equation}
where we have noted \eqref{eq:tauh}. 
Clearly, \eqref{eq:tauhe2} 
is a discrete approximation of \eqref{eq:bc}, which stipulates a $90^\circ$
degree contact angle between $\Gamma(t)$ and the $x_2$--axis.

For the scheme \eqref{eq:sd} we fix $\kappa^h \in \Whpartialzero$, and so the
right hand side in \eqref{eq:tauhe2} is zero. For a more general scheme,
where we allow $\kappa^h \in V^h$ and let $\chi \in V^h$ in \eqref{eq:sdb},
the right hand side in \eqref{eq:tauhe2} is still of order $O(J^{-1})$ on
assuming that $\kappa^h(q_0)$ is bounded. However, we observe that
$\kappa^h(q_0)$, if $0 \in \partial_0 I$, only appears in \eqref{eq:sdc} 
for this more general scheme. Hence there is no reason to assume that
$\kappa^h(q_0)$ should remain bounded. In fact, $\kappa^h(q_0)$ simply acts
as a register for the value
$\frac{\vec\tau^h(q_0)\,.\,\vec\ek_2}{\vec\nu^h(q_0)\,.\,\vec\ek_2}\,
\frac{2\,J}{\mathcal{H}^1(\Gamma^h(t))}$.
Moreover, a fully discrete variant of the discussed more general version
of \eqref{eq:sd} leads to incorrect contact angles at the $x_2$--axis and to
the numerical breakdown of the scheme. This is the main reason why we consider
the scheme \eqref{eq:sd} as it is.
\end{rmrk}

\subsubsection{Conserved flows} \label{sec:new51}

We rewrite \eqref{eq:sda} as
\begin{equation*} 
2\,\pi \left((\vec X^h\,.\,\vec\ek_1)\,\mat Q^h\,\vec X^h_t, 
\vec\chi\,|\vec X^h_\rho|\right)^h 
- \left( \vec Y^h_\rho\,.\,\vec\nu^h, \vec\chi_\rho\,.\,\vec\nu^h
  \,|\vec X^h_\rho|^{-1}\right)
= \left( \vec f^h, \vec\chi \,|\vec X^h_\rho|\right)^h 
\quad  \forall\ \vec\chi \in \xspaceh\,.
\end{equation*}
Then the natural generalization of $(\BGNpwf^h)^{h}$, \eqref{eq:sd}, that
approximates the weak formulation \eqref{eq:weak4LMa}, \eqref{eq:weak3b},
\eqref{eq:weak3c} and \eqref{eq:xsideSAV} is given by \eqref{eq:sd}, with
\eqref{eq:sda} replaced by
\begin{align}
& 2\,\pi \left((\vec X^h\,.\,\vec\ek_1)\,\mat Q^h\,\vec X^h_t, 
\vec\chi\,|\vec X^h_\rho|\right)^h 
- \left( \vec Y^h_\rho\,.\,\vec\nu^h, \vec\chi_\rho\,.\,\vec\nu^h
  \,|\vec X^h_\rho|^{-1}\right) = \left( \vec f^h, \vec\chi \,|\vec X^h_\rho|\right)^h
 \nonumber \\ & \quad
- 2\,\pi\,\lambda_A^h
\left[
\left( \vec\ek_1, \vec\chi\,|\vec X^h_\rho|\right) 
+ \left( (\vec X^h\,.\,\vec\ek_1)\,\vec\tau^h
, \vec\chi_\rho\right) \right]
- 2\,\pi\, \lambda_V^h
\left((\vec X^h\,.\,\vec\ek_1)\, 
\vec\nu^h, \vec\chi \,|\vec X^h_\rho|\right)
 \quad \forall\ \vec\chi \in \xspaceh\,, \label{eq:Q0hXtLMa}
\end{align}
where $(\lambda_A^h(t), \lambda_V^h(t))^T \in \bR^2$ are such that
\begin{equation} \label{eq:sdwfb}
A(\vec X^h(t)) = A(\vec X^h(0))\quad\text{and}\quad
V(\vec X^h(t)) = V(\vec X^h(0))\,.
\end{equation}
Here, on recalling \eqref{eq:Area} and \eqref{eq:Volume}, we note that
$A(\vec X^h(t))$ denotes the surface area of $\mathcal{S}^h(t)$,
where, similarly to (\ref{eq:SGamma}), we set
\begin{equation*} 
\mathcal{S}^h(t) = \bigcup_{\rho \in \overline{I}} \Pi_2^3(\vec X^h(\rho,t))\,.
\end{equation*}
Moreover, $V(\vec X^h(t))$
is the volume of the domain $\Omega^h(t)$ with
$\partial\Omega^h(t) = \mathcal{S}^h(t)$ in the case that $\mathcal{S}^h(t)$ 
has no boundary. We remark that
\begin{equation*} 
A(\vec Z^h) = 2\,\pi\left(\vec Z^h\,.\,\vec\ek_1 ,
|\vec Z^h_\rho|\right) \quad \vec Z^h \in \Vhpartialzero
\end{equation*}
and
\begin{equation*} 
V(\vec Z^h) = -\pi \left( (\vec Z^h\,.\,\vec\ek_1)^2, 
[\vec Z^h_\rho]^\perp\,.\,\vec\ek_1\right) \quad \vec Z^h \in \Vhpartialzero\,,
\end{equation*}
recall \eqref{eq:Area}, \eqref{eq:Volume} and \eqref{eq:tauh}.
Moreover, we recall from \cite[(3.7), (3.11)]{axisd} that, similarly to
\eqref{eq:dAdt} and \eqref{eq:dVdt}, it holds that
\begin{equation} \label{eq:dAhdt}
\ddt\, A(\vec X^h(t))
= 2\,\pi \left[ \left( \vec\ek_1, \vec X^h_t \,|\vec X^h_\rho| \right) +
 \left( (\vec X^h\,.\,\vec\ek_1) \,\vec\tau^h , (\vec X^h_t)_\rho \right) 
\right]
\end{equation}
and
\begin{equation} \label{eq:dVhdt}
\ddt\,V(\vec X^h(t))
= 2\,\pi\left( (\vec X^h\,.\,\vec\ek_1)\,\vec\nu^h, \vec X^h_t
\,|\vec X^h_\rho| \right) .
\end{equation}

\begin{thrm} \label{thm:stabC}
Let {\rm Assumption~\ref{ass:Ah2}} be satisfied and 
let $(\vec X^h(t),\kappa^h(t),\vec Y^h(t),\vec{\rm m}^h(t),
\lambda^h_A(t), \lambda^h_V(t))_{t\in (0,T]}$ 
be a solution to \eqref{eq:Q0hXtLMa}, \eqref{eq:sdb}, \eqref{eq:sdc}, 
\eqref{eq:sdwfb}.
Then the solution satisfies the stability bound
\[
\ddt \,\widehat E^h(t)
+ 2\,\pi\left(\vec X^h\,.\,\vec\ek_1\,|\mat Q^h\,\vec X^h_t|^2,
|\vec X^h_\rho|\right)^h=0\,.
\]
\end{thrm}
\begin{proof}
Differentiating the two equations in \eqref{eq:sdwfb} with respect to $t$,
noting \eqref{eq:dAhdt}, \eqref{eq:dVhdt}, and choosing
$\vec\chi = \vec X^h_t$ in \eqref{eq:Q0hXtLMa} yields
\[
2\,\pi\left(\vec X^h\,.\,\vec\ek_1\,|\mat Q^h\,\vec X^h_t|^2,
|\vec X^h_\rho|\right)^h
- \left( \vec Y^h_\rho\,.\,\vec\nu^h, (\vec X^h_t)_\rho\,.\,\vec\nu^h
  \,|\vec X^h_\rho|^{-1}\right)
= \left( \vec f^h, \vec X^h_t \,|\vec X^h_\rho|\right)^h ,
\]
which is equivalent to \eqref{eq:Pxtxth}. Hence the stability result follows as
in the proof of Theorem~\ref{thm:stab}.
\end{proof}

\subsection{Based on $\kappa^h_{\mathcal{S}}$} \label{sec:sdkappaS}

The approach in \S\ref{sec:sdkappa} required the introduction of the
auxiliary finite element function \eqref{eq:calKh} in order to be able to
obtain a well-defined Lagrangian. An alternative is to define the Lagrangian in
terms of $\kappa^h_{\mathcal{S}}$, a direct discrete analogue of 
$\varkappa_{\mathcal{S}}$. Then it is possible to repeat the proof of 
(\ref{eq:PxtSstab}) on the discrete level.

As the discrete analogue of (\ref{eq:varkappaSweak}), we let
$\vec X^h \in \Vhpartialzero$, 
$\kappa_{\mathcal{S}}^h \in W^h_{(\partial_0)}$ and $\vec{\rm m}^h : \partial_M I \to \bR^2$
be such that
\begin{align} \label{eq:sideSh}
& \left( \vec X^h\,.\,\vec\ek_1\,\kappa_{\mathcal{S}}^h\,\vec\nu^h,
\vec\eta \,|\vec X^h_\rho|\right)^{(h)}
+ \left(\vec\ek_1, \vec\eta \,|\vec X^h_\rho|\right)
+ \left( (\vec X^h\,.\,\vec\ek_1)\,
\vec\tau^h, \vec\eta_\rho \right)
\nonumber \\ & \qquad
= \sum_{p \in \partial_C I} \left[(\vec X^h\,.\,\vec\ek_1)\,
\vec\zeta\,.\,\vec\eta\right](p)
+ \sum_{p \in \partial_M I} \left[(\vec X^h\,.\,\vec\ek_1)\,
\vec{\rm m}^h\,.\,\vec\eta\right](p) \qquad
\forall\ \vec\eta \in \Vhpartialzero\,.
\end{align}
Here, and throughout, we use the notation $\cdot^{(h)}$ to denote an 
expression with or without the superscript $h$, and simultaneously
for the notation with subscripts $\cdot_{(\partial_0)}$.
I.e.\ we consider two separate situations: either mass lumping is used
on the first integral, and then we let 
$\kappa^h_{\mathcal{S}} \in W^h_{\partial_0}$, or true integration
is employed and we let $\kappa^h_{\mathcal{S}} \in W^h$.
We define the discrete analogue of (\ref{eq:EEaxCN}) 
\begin{equation} \label{eq:ESh}
\widehat E^h_{\mathcal{S}}(t) = \pi\left(\alpha\, 
[ \kappa_{\mathcal{S}}^h -\spont ]^2 + 2\,\lambda ,
\vec X^h\,.\,\vec\ek_1 |\vec X^h_\rho|\right)^{(h)}
+ \tfrac\beta2\, (\AADE^h_{\mathcal{S}}(t))^2
- 2\,\pi\,\alpha_G\sum_{p\in \partial_M I}
\vec{\rm m}^h(p)\,.\,\vec\ek_1
+ 2\,\pi\,\varsigma \hspace{-3mm}\sum_{p\in \partial_2 I \cup \partial_F I}
\vec X^h(p)\,.\,\vec\ek_1\,,
\end{equation}
where 
\begin{equation} \label{eq:AhS}
\AADE^h_{\mathcal{S}}(t) = 
2\,\pi\left(\vec X^h\,.\,\vec\ek_1\,\kappa^h_{\mathcal{S}}, 
|\vec X^h_\rho| \right)^{(h)} - M_0\,.
\end{equation}
Similarly to (\ref{eq:Lag2}), we define the discrete Lagrangian
\begin{align} \label{eq:Lag2h}
& \mathcal{L}^h_{\mathcal{S}}(\vec X^h, \kappa^h_{\mathcal{S}}, 
\vec{\rm m}^h, \vec Y^h_{\mathcal{S}})
 = \pi\left(\alpha\, 
[ \kappa_{\mathcal{S}}^h -\spont ]^2 + 2\,\lambda ,
\vec X^h\,.\,\vec\ek_1 |\vec X^h_\rho|\right)^{(h)}
+ \tfrac\beta2\left[ 
2\,\pi\left(\vec X^h\,.\,\vec\ek_1\,\kappa^h_{\mathcal{S}}, 
|\vec X^h_\rho| \right)^{(h)} - M_0 \right]^2
\nonumber \\ & \qquad
- \left( \vec X^h\,.\,\vec\ek_1\,\kappa_{\mathcal{S}}^h\,\vec\nu^h +\vec\ek_1, 
\vec Y^h_{\mathcal{S}} \,|\vec X^h_\rho|\right)^{(h)}
- \left( (\vec X^h\,.\,\vec\ek_1)\,
\vec\tau^h, (\vec Y^h_{\mathcal{S}})_\rho \right) 
+ 2\,\pi\,\varsigma\,\sum_{p \in \partial_2 I \cup \partial_F I}
\vec X^h(p)\,.\,\vec\ek_1 
\nonumber \\ & \qquad
+ \sum_{p \in \partial_C I} \left[
(\vec X^h\,.\,\vec\ek_1)\,\vec\zeta\,.\,\vec Y^h_{\mathcal{S}}\right](p)
+\sum_{p \in \partial_M I} \left[\vec{\rm m}^h\,.\,(
(\vec X^h\,.\,\vec\ek_1)\,\vec Y^h_{\mathcal{S}} 
- 2\,\pi\,\alpha_G\,\vec\ek_1)\right](p)\,,
\end{align}
for $\vec X^h \in \Vhpartialzero$, 
$\kappa_{\mathcal{S}}^h \in \Whpartialzerob$,
$\vec{\rm m}^h : \partial_M I \to \bR^2$
and $\vec Y^h_{\mathcal{S}} \in \Vhpartialzero$.
Here we observe that in the case of numerical integration, we fix
$\kappa_{\mathcal{S}}^h$ to be zero on the boundary $\partial_0 I$, as the
Lagrangian (\ref{eq:Lag2h}) does not depend on these boundary values at all.

Taking variations $\vec\eta \in \Vhpartialzero$ in $\vec Y^h_{\mathcal{S}}$, and setting
$\left[\deldel{\vec Y^h_{\mathcal{S}}}\, \mathcal{L}_{\mathcal{S}}^h\right](\vec\eta) = 0$
we obtain (\ref{eq:sideSh}), similarly to (\ref{eq:kappaSeta}).
Taking variations $\chi \in \Whpartialzerob$ in $\kappa_{\mathcal{S}}^h$ and 
setting
$\left[\deldel{\kappa_{\mathcal{S}}^h}\, \mathcal{L}_{\mathcal{S}}^h\right](\chi) = 0$
we obtain, similarly to (\ref{eq:kappaSy}), 
\begin{equation} 
2\,\pi\left( \alpha\,(\kappa_{\mathcal{S}}^h - \spont)
+ \beta\,\AADE_{\mathcal{S}}^h,
\vec X^h\,.\,\vec\ek_1\,\chi\,|\vec X^h_\rho|\right)^{(h)}
- \left( (\vec X^h\,.\,\vec\ek_1)\,\vec Y^h_{\mathcal{S}}, 
\chi\,\vec\nu^h \,|\vec X^h_\rho|\right)^{(h)} = 0
\qquad \forall\ \chi \in \Whpartialzerob\,. \label{eq:kappaShYh}
\end{equation}
Taking variations in $\vec{\rm m}^h$, and setting them to zero, yields,
similarly to (\ref{eq:mSy}), that
\begin{equation} \label{eq:mhSy}
(\vec X^h\,.\,\vec\ek_1)\,\vec Y^h_{\mathcal{S}} = 
2\,\pi\,\alpha_G\,\vec\ek_1 \qquad\text{on }\ \partial_M I\,.
\end{equation}
Taking variations $\vec\chi \in \xspaceh$ in $\vec X^h$ 
and setting
$2\,\pi\left((\vec X^h\,.\,\vec\ek_1)\,\mat Q^h\,\vec X^h_t, 
\vec\chi\,|\vec X^h_\rho|\right)^h=
- \left[\deldel{\vec X^h}\, \mathcal{L}^h_{\mathcal{S}}\right](\vec\chi)$
we obtain, similarly to (\ref{eq:Lag2dx}), that 
\begin{align}
&2\,\pi\left((\vec X^h\,.\,\vec\ek_1)\,\mat Q^h\,\vec X^h_t, 
\vec\chi\,|\vec X^h_\rho|\right)^h  =
-\pi \left(\alpha\left[\kappa_{\mathcal{S}}^h -\spont \right]^2 
+ 2\,\lambda + 2\,\beta\, \AADE_{\mathcal{S}}^h\,\kappa_{\mathcal{S}}^h, 
\left[\deldel{\vec X^h}\,(\vec X^h\,.\,\vec\ek_1)\,|\vec X^h_\rho|\,\right](\vec\chi) 
\right)^{(h)} \nonumber \\& \
+ \left(\kappa_{\mathcal{S}}^h\,\vec Y_{\mathcal{S}}^h,
\left[\deldel{\vec X^h}\,(\vec X^h\,.\,\vec\ek_1)\,\vec\nu^h\,|\vec X^h_\rho|\right]
(\vec\chi)\right)^{(h)}
 + \left(\vec\ek_1\,.\,\vec Y_{\mathcal{S}}^h, 
 \left[\deldel{\vec X^h}\,|\vec X^h_\rho|\right](\vec\chi)\right)
+ \left((\vec Y_{\mathcal{S}}^h)_\rho , 
\left[\deldel{\vec X^h}\,(\vec X^h\,.\,\vec\ek_1)\,\vec\tau^h\right]
(\vec\chi)\right)
\nonumber \\& \
- \sum_{p \in \partial_{2} I \cup \partial_F I}
\left[ \left[2\,\pi\,\varsigma + \vec{\rm m}^h\,.\,\vec Y^h_{\mathcal{S}}
\right]
(\vec\chi\,.\,\vec\ek_1)\,\right](p)  
\qquad \forall \ \vec\chi \in \xspaceh \,.
\label{eq:Lag2dx1h}
\end{align}
Choosing $\vec\chi = \vec X^h_t \in \xspaceh$ in (\ref{eq:Lag2dx1h}),
on noting (\ref{eq:ddA}), yields the discrete analogue of (\ref{eq:Lag2dxxt})
\begin{align}
&2\,\pi\left(\vec X^h\,.\,\vec\ek_1\,|\mat Q^h\,\vec X^h_t|^2, 
|\vec X^h_\rho|\right)^h  
 = -\pi \left( \alpha\left[\kappa_{\mathcal{S}}^h -\spont \right]^2
+ 2\,\lambda+ 2\,\beta\, \AADE_{\mathcal{S}}^h\,\kappa_{\mathcal{S}}^h, 
[(\vec X^h\,.\,\vec\ek_1)\,|\vec X^h_\rho|]_t \right)^{(h)} 
\nonumber \\& \qquad 
+ \left(\kappa_{\mathcal{S}}^h\,\vec Y_{\mathcal{S}}^h,
[(\vec X^h\,.\,\vec\ek_1)\,\vec\nu^h\,|\vec X^h_\rho|]_t \right)^{(h)}
 + \left(\vec\ek_1\,.\,\vec Y_{\mathcal{S}}^h, 
 (|\vec X^h_\rho|)_t \right)
+ \left((\vec Y_{\mathcal{S}}^h)_\rho , 
[(\vec X^h\,.\,\vec\ek_1)\,\vec\tau^h]_t\right)
\nonumber \\& \qquad 
- \sum_{p \in \partial_{2} I \cup \partial_F I}
\left[ \left[2\,\pi\,\varsigma + \vec{\rm m}^h\,.\,\vec Y^h_{\mathcal{S}}
\right]
(\vec X^h_t\,.\,\vec\ek_1)\,\right](p)\,.
\label{eq:newetaxt}
\end{align}

Differentiating (\ref{eq:sideSh}) with respect to $t$, and then 
choosing $\vec\eta = \vec Y^h_{\mathcal{S}} \in \Vhpartialzero$, we obtain, 
on recalling that $\vec X^h_t \in \xspace^h$,
\begin{align}
& \left( (\kappa^h_{\mathcal{S}})_t\,\vec Y^h_{\mathcal{S}},
(\vec X^h\,.\,\vec\ek_1)\,\vec\nu^h\,|\vec X^h_\rho| \right)^{(h)}
+ \left( \kappa^h_{\mathcal{S}}\,\vec Y^h_{\mathcal{S}},
\left[(\vec X^h\,.\,\vec\ek_1)\,\vec\nu^h\,|\vec X^h_\rho|\right]_t \right)^{(h)}
+ \left( \vec\ek_1\,.\,\vec Y^h_{\mathcal{S}},
\left[|\vec X^h_\rho|\right]_t \right)
\nonumber \\ & \qquad
+ \left( (\vec Y^h_{\mathcal{S}})_\rho, 
\left[(\vec X^h\,.\,\vec\ek_1)\,\vec\tau^h\right]_t
\right)
= \sum_{p \in \partial_M I} \left[(\vec X^h_t\,.\,\vec\ek_1)\,
\vec{\rm m}^h\,.\,\vec Y^h_{\mathcal{S}} + (\vec X^h\,.\,\vec\ek_1)\,
(\vec{\rm m}^h)_t\,.\,\vec Y^h_{\mathcal{S}} \right](p) \,.
\label{eq:dtsideSh}
\end{align}
It follows from (\ref{eq:dtsideSh}), (\ref{eq:mhSy}) and
(\ref{eq:kappaShYh}) with $\chi = (\kappa^h_{\mathcal{S}})_t \in 
\Whpartialzerob$ that
\begin{align}
&
\left(\kappa^h_{\mathcal{S}}\,\vec Y^h_{\mathcal{S}},
\left[(\vec X^h\,.\,\vec\ek_1)\,\vec\nu^h\,|\vec X^h_\rho|\right]_t
\right)^{(h)}
 + \left(\vec\ek_1\,.\,\vec Y^h_{\mathcal{S}}, 
 \left[|\vec X^h_\rho|\right]_t\right)
+ \left((\vec Y^h_{\mathcal{S}})_\rho , 
\left[(\vec X^h\,.\,\vec\ek_1)\,\vec\tau^h\right]_t \right)
\nonumber \\ & \
= - \left( (\kappa^h_{\mathcal{S}})_t , (\vec X^h\,.\,\vec\ek_1)\,\,
\vec Y^h_{\mathcal{S}}\,.\,\vec\nu^h\,|\vec X^h_\rho|\right)^{(h)}
+ \sum_{p \in \partial_M I} \left[(\vec X^h_t\,.\,\vec\ek_1)\,
\vec{\rm m}^h\,.\,\vec Y^h_{\mathcal{S}} + 2\,\pi\,\alpha_G\,\,
\vec{\rm m}^h_t\,.\,\vec\ek_1 \right](p)
\nonumber \\ & \
= - 2\,\pi
\left(\alpha\, (\kappa^h_{\mathcal{S}} - \spont) 
+ \beta\,\AADE^h_{\mathcal{S}}, 
\vec X^h\,.\,\vec\ek_1\,(\kappa^h_{\mathcal{S}})_t\,|\vec X^h_\rho|
\right)^{(h)}
+ \sum_{p \in \partial_M I} \left[(\vec X^h_t\,.\,\vec\ek_1)\,
\vec{\rm m}^h\,.\,\vec Y^h_{\mathcal{S}} + 2\,\pi\,\alpha_G\,\,
\vec{\rm m}^h_t\,.\,\vec\ek_1 \right](p)\,.
\label{eq:new6h}
\end{align}
Combining (\ref{eq:newetaxt}) and (\ref{eq:new6h}) yields the discrete 
analogue of (\ref{eq:PxtSstab})
\begin{align}
& 2\,\pi\left(\vec X^h\,.\,\vec\ek_1\,|\mat Q^h\,\vec X^h_t|^2, 
|\vec X^h_\rho|\right)^h 
= -\pi \left(\alpha \left[\kappa^h_{\mathcal{S}} -\spont \right]^2 
+ 2\,\lambda + 2\,\beta\,\AADE^h_{\mathcal{S}}\,\kappa^h_{\mathcal{S}}, 
\left[\vec X^h\,.\,\vec\ek_1\,|\vec X^h_\rho|\right]_t\right)^{(h)}  
\nonumber \\ & \quad
- 2\,\pi
\left(\alpha\, (\kappa^h_{\mathcal{S}} - \spont) 
+ \beta\,\AADE^h_{\mathcal{S}}, 
\vec X^h\,.\,\vec\ek_1\,(\kappa^h_{\mathcal{S}})_t\,|\vec X^h_\rho|
\right)^{(h)} 
- 2\,\pi\,\varsigma\,\sum_{p \in \partial_{2} I \cup \partial_F I}
\vec X^h_t(p)\,.\,\vec\ek_1 
+2\,\pi\,\alpha_G\,\sum_{p \in \partial_M I} \vec{\rm m}^h_t(p)\,.\,
 \vec\ek_1
\nonumber \\ & = - \pi\,\ddt 
\left(\alpha\left
[\kappa^h_{\mathcal{S}} -\spont \right]^2 + { 2\,\lambda},
\vec X^h\,.\,\vec\ek_1\,|\vec X^h_\rho|\right)^{(h)}  
- \tfrac\beta2\,\ddt\,(\AADE^h_{\mathcal{S}})^2 
- 2\,\pi\,\varsigma\hspace{-4mm}\sum_{p \in \partial_{2} I \cup \partial_F I}
\vec X^h_t(p)\,.\,\vec\ek_1 
+2\,\pi\,\alpha_G\,\sum_{p \in \partial_M I} \vec{\rm m}^h_t(p)\,.\,
 \vec\ek_1
\nonumber \\ & 
= - \ddt \, \widehat E^h_{\mathcal{S}}(t)\,,
\label{eq:Shstab}
\end{align}
where we have recalled the definition (\ref{eq:ESh}). 

In order to derive a suitable semidiscrete approximation, we now return to
(\ref{eq:Lag2dx1h}). Substituting (\ref{eq:57}) into 
(\ref{eq:Lag2dx1h}) yields that 
\begin{align}
& 2\,\pi\left( (\vec X^h\,.\,\vec\ek_1)\,\mat Q^h\,\vec X^h_t, 
\vec\chi\,|\vec X^h_\rho|\right)^h 
 = - \left( \pi\,\vec X^h\,.\,\vec\ek_1\,
\left[ \alpha \,[\kappa^h_{\mathcal{S}} - \spont]^2 
+ 2\,\lambda + 2\,\beta\, \AADE_{\mathcal{S}}^h\,\kappa_{\mathcal{S}}^h
\right]
- \vec Y^h_{\mathcal{S}}\,.\,\vec\ek_1  ,
\vec\chi_\rho\,.\,\vec\tau^h \right)^{(h)}
\nonumber \\ & \ 
- \left( \left[ \pi \left[ 
\alpha\left[\kappa_{\mathcal{S}}^h -\spont \right]^2 
+ 2\,\lambda + 2\,\beta\, \AADE_{\mathcal{S}}^h\,\kappa_{\mathcal{S}}^h \right]
- \kappa^h_{\mathcal{S}}\,\vec Y^h_{\mathcal{S}}\,.\,\vec\nu^h \right] 
|\vec X^h_\rho| - (\vec Y^h_{\mathcal{S}})_\rho\,.\,\vec\tau^h,
\vec\chi\,.\,\vec\ek_1 \right)^{(h)}
- \left( \vec X^h\,.\,\vec\ek_1\,\kappa^h_{\mathcal{S}}\, 
\vec Y^h_{\mathcal{S}} , \vec\chi_\rho^\perp \right)^{(h)}
\nonumber \\ & \
+ \left( (\vec X^h\,.\,\vec\ek_1)\,(\vec Y^h_{\mathcal{S}})_\rho\,.\,\vec\nu^h,
\vec\chi_\rho\,.\,\vec\nu^h |\vec X^h_\rho|^{-1} \right)
- \sum_{p \in \partial_{2} I \cup \partial_F I}
\left[ \left[2\,\pi\,\varsigma + \vec{\rm m}^h\,.\,\vec Y^h_{\mathcal{S}}
\right]
(\vec\chi\,.\,\vec\ek_1)\,\right](p) 
\qquad \forall \ \vec\chi \in \xspaceh \,.
\label{eq:Lag2dx2h}
\end{align}

On combining (\ref{eq:Lag2dx2h}), (\ref{eq:kappaShYh}),
(\ref{eq:sideSh}), (\ref{eq:mhSy}) and (\ref{eq:AhS}),  
on noting (\ref{eq:abperp}) and (\ref{eq:tauh}),
our semidiscrete approximation based on $\kappa^h_{\mathcal{S}}$
is given as follows.

\noindent
$(\BGNpwfwf^h)^{(h)}$
Let $\vec X^h(\cdot,0) \in \Vhpartialzero$ and $\alpha \in \bRplus$,
$\spont,M_0,\alpha_G\,\lambda,\varsigma \in \bR$, $\beta \in \bRgeq$, 
$\vec\zeta : \partial_C I \to \bS^1$ be given.
For $t \in (0,T]$ find $\vec X^h(\cdot,t)$, 
with $\vec X^h_t(\cdot,t) \in \xspaceh$,
$\kappa^h_{\mathcal{S}}(\cdot,t) \in \Whpartialzerob$,
$\vec Y^h_{\mathcal{S}}(\cdot,t) \in \Vhpartialzero$,
with 
$\vec\pi^h\left[(\vec X^h\,.\,\vec\ek_1)\, \vec Y^h_{\mathcal{S}}\right]
(\cdot,t) \in \yspaceh(2\,\pi\,\alpha_G \,\vec\ek_1)$,
and $\vec{\rm m}^h(\cdot,t) : \partial_M I \to \bR^2$ such that
\begin{subequations} \label{eq:sdS}
\begin{align}
& 2\,\pi\left( (\vec X^h\,.\,\vec\ek_1)\,\mat Q^h\,\vec X^h_t, \vec\chi
\,|\vec X^h_\rho|\right)^h
- \left( (\vec X^h\,.\,\vec\ek_1)\,(\vec Y^h_{\mathcal{S}})_\rho\,.\,\vec\nu^h,
\vec\chi_\rho\,.\,\vec\nu^h |\vec X^h_\rho|^{-1} \right)
\nonumber \\ 
& \quad =
-  \left( 
\pi\,\vec X^h\,.\,\vec\ek_1\,
\left[ \alpha \,[\kappa^h_{\mathcal{S}} - \spont]^2 
+ 2\,\lambda + 2\,\beta\, \AADE_{\mathcal{S}}^h\,\kappa_{\mathcal{S}}^h
\right]
- \vec Y^h_{\mathcal{S}}\,.\,\vec\ek_1  ,
\vec\chi_\rho\,.\,\vec\tau^h \right)^{(h)}
\nonumber \\ & \qquad  -
\left( \left[ \pi \left[ 
\alpha\left[\kappa_{\mathcal{S}}^h -\spont \right]^2 
+ 2\,\lambda + 2\,\beta\, \AADE_{\mathcal{S}}^h\,\kappa_{\mathcal{S}}^h \right]
- \kappa^h_{\mathcal{S}}\,\vec Y^h_{\mathcal{S}}\,.\,\vec\nu^h \right] 
|\vec X^h_\rho| - (\vec Y^h_{\mathcal{S}})_\rho\,.\,\vec\tau^h,
\vec\chi\,.\,\vec\ek_1 \right)^{(h)}
\nonumber \\ & \qquad 
+ \left( \vec X^h\,.\,\vec\ek_1\,\kappa^h_{\mathcal{S}}\, 
(\vec Y^h_{\mathcal{S}})^\perp , \vec\chi_\rho \right)^{(h)}
- \sum_{p \in \partial_{2} I \cup \partial_F I}
\left[ \left[2\,\pi\,\varsigma + \vec{\rm m}^h\,.\,\vec Y^h_{\mathcal{S}}
\right]
(\vec\chi\,.\,\vec\ek_1)\,\right](p)
 \qquad \forall \ \vec\chi \in \xspaceh \,,
\label{eq:sdSa}\\
& 2\,\pi\left( \vec X^h\,.\,\vec\ek_1
\left[\alpha\,(\kappa_{\mathcal{S}}^h - \spont) +\beta\,\AADE^h_{\mathcal{S}}\right],
\chi\,|\vec X^h_\rho|\right)^{(h)}
- \left( (\vec X^h\,.\,\vec\ek_1)\,\vec Y^h_{\mathcal{S}}, 
\chi\,\vec\nu^h \,|\vec X^h_\rho|\right)^{(h)} = 0
\qquad \forall\ \chi \in \Whpartialzerob\,, \label{eq:sdSb} \\
& \left( \vec X^h\,.\,\vec\ek_1\,\kappa_{\mathcal{S}}^h\,\vec\nu^h,
\vec\eta \,|\vec X^h_\rho|\right)^{(h)}
+ \left(\vec\ek_1, \vec\eta \,|\vec X^h_\rho|\right)
+ \left( (\vec X^h\,.\,\vec\ek_1)\,
\vec X^h_\rho, \vec\eta_\rho \, |\vec X^h_\rho|^{-1}\right)
\nonumber \\ & \hspace{4cm}
= \sum_{p \in \partial_C I} \left[(\vec X^h\,.\,\vec\ek_1)\,
\vec\zeta\,.\,\vec\eta\right](p)
+ \sum_{p \in \partial_M I} \left[(\vec X^h\,.\,\vec\ek_1)\,
\vec{\rm m}^h\,.\,\vec\eta\right](p) \qquad
\forall\ \vec\eta \in \Vhpartialzero\,, \label{eq:sdSc}
\end{align}
\end{subequations}
where $\AADE^h_{\mathcal{S}}(t)$ is given by (\ref{eq:AhS}).

\begin{thrm} \label{thm:stabS}
Let {\rm Assumption~\ref{ass:Ah2}} be satisfied and 
let $(\vec X^h(t),\kappa^h_{\mathcal{S}}(t),\vec Y_{\mathcal{S}}^h(t),
\vec{\rm m}^h(t))_{t\in(0,T]}$ 
be a solution to \mbox{\rm (\ref{eq:sdS})}. 
Then the solution satisfies the stability bound
\[
\ddt \,\widehat E^h_{\mathcal{S}}(t)
+ 2\,\pi\left(\vec X^h\,.\,\vec\ek_1\,|\mat Q^h\,\vec X^h_t|^2,
|\vec X^h_\rho|\right)^h=0\,. 
\]
\end{thrm}
\begin{proof}
The desired result follows as (\ref{eq:sdS}) is just a rewrite of 
(\ref{eq:Lag2dx1h}), (\ref{eq:kappaShYh}), (\ref{eq:sideSh}) and 
(\ref{eq:AhS}), and then noting (\ref{eq:newetaxt})--(\ref{eq:Shstab}). 
\end{proof}

\begin{rmrk} \label{rem:fdSkappa0h}
For the scheme $(\BGNpwfwf^h)^{h}$ we note that if
$0 \in \partial_0 I$, then choosing
$\vec\eta = \chi_0\,\vec\ek_2$ in \eqref{eq:sdSc} yields that
$(\vec X^h (q_1) - \vec X^h (q_0) ) \,.\,\vec\ek_2 = 0$,
recall also {\rm Remark~\ref{rem:kappahpartial02}}.
While the same is not true for the scheme $(\BGNpwfwf^h)$, a weaker form of
these constraints is still enforced via \eqref{eq:sdSc}, leading to a nearly
$90^\circ$ degree contact angle on the $x_2$--axis in practice
for fully discrete variants based on $(\BGNpwfwf^h)$.
\end{rmrk}

\setcounter{equation}{0}
\section{Fully discrete scheme}  \label{sec:fd}

Let $0= t_0 < t_1 < \ldots < t_{M-1} < t_M = T$ be a
partitioning of $[0,T]$ into possibly variable time steps 
$\ttau_m = t_{m+1} - t_{m}$, $m=0,\ldots, M-1$.
For $\vec X^m \in \Vhpartialzero$, we let 
$\vec\tau^m$ and $\vec\nu^m$ be the natural fully discrete analogues
of $\vec\tau^h$ and $\vec\nu^h$ on $\Gamma^m = \vec X^m(I)$, 
recall (\ref{eq:tauh}).  
In addition, let 
$\vec\omega^m \in \underline{V}^h$ be the natural 
fully discrete analogue of
$\vec\omega^h \in \underline{V}^h$,
recall (\ref{eq:omegah}); 
and similarly for $\vec v^m \in \Vh$, recall (\ref{eq:vh}). 
Finally, let $\vec Q^m \in [V^h]^{2\times 2}$ be the natural fully discrete 
analogue of $\vec Q^h$, recall \eqref{eq:Qh}.

\subsection{Based on $\kappa^{m+1}$}

We propose the following fully discrete approximation of $(\BGNpwf^h)^h$.

\noindent
$(\BGNpwf^m)^{h}$
Let $\vec X^0 \in \Vhpartialzero$, $\kappa^0 \in V^h$, 
$\vec Y^0 \in \Vhpartialzero$ and $\alpha \in \bRplus$,
$\spont,M_0,\alpha_G,\lambda,\varsigma \in \bR$, $\beta \in \bRgeq$, 
$\vec\zeta : \partial_C I \to \bS^1$ be given.
For $m=0,\ldots,M-1$, find $\delta\vec X^{m+1} \in \xspace^h$, 
with $\vec X^{m+1} = \vec X^m + \delta\vec X^{m+1}$, 
$\kappa^{m+1} \in \Whpartialzero$,
$\vec Y^{m+1} \in \yspaceh( 2\,\pi\,\alpha_G\,\vec\ek_1)$ 
and $\vec{\rm m}^{m+1} : \partial_M I \to \bR^2$ such that
\begin{subequations} \label{eq:fd}
\begin{align}
& 2\,\pi\left(\vec X^m\,.\,\vec\ek_1\,\mat Q^m\,
\frac{\vec X^{m+1} - \vec X^m}{\ttau_m}, \vec\chi\,|\vec X^m_\rho|\right)^h 
- \left(\vec Y^{m+1}_\rho , \vec\chi_\rho\,|\vec X^m_\rho|^{-1} \right) 
+ \left( \vec Y^{m}_\rho\,.\,\vec\tau^m, 
\vec\chi_\rho\,.\,\vec\tau^m\,|\vec X^m_\rho|^{-1} \right)
\nonumber \\ & \quad 
= - \pi \left( 
\alpha \left[ \doctorkappa^h(\vec X^m,\kappa^m) - \spont \right]^2
+ 2\,\lambda + 2 \,\beta\,\AADE^m\,\kappa^m, 
\vec\chi\,.\,\vec\ek_1\,|\vec X^m_\rho| + (\vec X^m\,.\,\vec\ek_1)\,\vec\tau^m\,.\,
\vec\chi_\rho \right)^h
\nonumber \\ & \qquad \quad
+ 2\,\pi\,\alpha\left( \left[\doctorkappa^h(\vec X^m,\kappa^m) - \spont\right] 
 (\doctorZ^h - 2), 
\frac{\vec\omega^m\,.\,\vec\ek_1}{\vec X^m\,.\,\vec e_1}\,
\vec\chi\,.\,\vec\ek_1  \,|\vec X^m_\rho| \right)^h 
\nonumber \\ & \qquad \quad
+ 2\,\pi\,\alpha \left( \left[\doctorkappa^h(\vec X^m,\kappa^m) - \spont\right]
(\doctorZ^h - 2)\,\vec\ek_1, 
(\vec\nu^m\,.\,\vec\chi_\rho) \,\vec\tau^m + (\vec\tau^m\,.\,\vec\chi_\rho) \,
(\vec\omega^m-\vec\nu^m) \right)^h 
\nonumber \\ & \qquad \quad
+ \left( \kappa^m\,(\vec Y^m)^\perp - 2\,\pi\,\beta\,\AADE^m\,\vec\ek_2,
\vec\chi_\rho \right)^h
- 2\,\pi\,\varsigma\sum_{p\in \partial_{2} I \cup \partial_F I}
\vec\chi(p)\,.\,\vec\ek_1
\qquad \forall\ \vec\chi \in \xspaceh\,, \label{eq:fda} \\
&2\,\pi \left( \vec X^m\,.\,\vec\ek_1 \left( \alpha\,[\doctorkappa^h(
\vec X^m,\kappa^{m+1}) - \spont ] + \beta\, \AADE^{m}\right)
,\chi \,|\vec X^m_\rho|\right)^h
- \left( \vec Y^{m+1}, \chi\,\vec\nu^m \,|\vec X^m_\rho|\right)^h = 0
\qquad \forall\ \chi \in \Whpartialzero\,,\label{eq:fdb} \\
&\left( \kappa^{m+1}\,\vec\nu^m, \vec\eta \,|\vec X^m_\rho|\right)^h
+ \left( \vec X^{m+1}_\rho, \vec\eta_\rho\,|\vec X^m_\rho|^{-1} \right) = 
\sum_{p \in \partial_C I} \left[\vec\zeta\,.\,\vec\eta\right](p)
+ \sum_{p \in \partial_M I} \left[\vec{\rm m}^{m+1}\,.\,\vec\eta\right](p)
 \qquad \forall\ \vec\eta \in \Vhpartialzero\,, \label{eq:fdc}
\end{align}
where 
\begin{equation}
\AADE^m = 
2\,\pi\left(\vec X^m\,.\,\vec\ek_1\,\kappa^m - \vec\nu^m.\,\vec\ek_1, 
|\vec X^m_\rho| \right)^h - M_0\,. 
\label{eq:fdd}
\end{equation}
\end{subequations}

\begin{ass} \label{ass:spanvm0}
Let $\vec X^m$ satisfy {\rm Assumption~\ref{ass:Ah2}} with $\vec X^h$ 
replaced by $\vec X^m$, and let 
$\dim\spa\{\vec v^m(q_j)\}_{j = 1}^{J-1} = 2$.
\end{ass}

\begin{ass} \label{ass:C}
Let $\vec X^m$ satisfy {\rm Assumption~\ref{ass:Ah1}}
with $\vec X^h$ replaced by $\vec X^m$,
and be such that the
following holds. If $\vec U \in \yspace^h(\vec 0)$ with
\begin{equation} \label{eq:assC}
 \left(\vec U_\rho , \vec\chi_\rho\,|\vec X^m_\rho|^{-1} \right) = 0
\quad \forall\ \vec\chi \in \xspaceh
\quad\text{and}\quad
\left( \vec U, \chi\,\vec\nu^m \,|\vec X^m_\rho|\right)^h = 0
\quad \forall\ \chi \in \Whpartialzero\,,
\end{equation}
then $\vec U = \vec 0$.
\end{ass}

We note that Assumption~\ref{ass:C} is only violated in very rare cases.
For example, if 
$\emptyset \not= \partial_C I = \partial I$ and
if $\vec X^m$ parameterizes a straight line, then
$\vec U = \vec\tau^m$ constant in $I$ satisfies \eqref{eq:assC}.
However, the following lemma shows that in most cases the assumption holds.

\begin{lmm} \label{lem:assC}
Let $\vec X^m$ satisfy {\rm Assumption~\ref{ass:spanvm0}}. Then if 
$\partial_C I = \emptyset$, or if 
$\partial_C I \not=\partial I \setminus \partial_0 I$, then
{\rm Assumption~\ref{ass:C}} holds.
\end{lmm}
\begin{proof}
Let $\vec U \in \yspace^h(\vec 0)$ satisfy \eqref{eq:assC}. 
If $\partial_C I = \emptyset$ then
we can choose $\vec\chi = \vec U \in \yspace^h(\vec 0) \subset \xspace^h$
in \eqref{eq:assC}, recall \eqref{eq:partialCempty}, to obtain that 
$\vec U$ is constant in $\overline I$.
The second property in \eqref{eq:assC}, on recalling \eqref{eq:omegah},
then yields that 
$\vec U\,.\,\vec\omega^m(q_j) = \vec U\,.\,\vec v^m(q_j) = 0$, 
$j=1,\ldots,J-1$.
Hence Assumption~\ref{ass:spanvm0} gives that $\vec U = \vec 0$.

We now consider the case $\partial_C I \not= \emptyset$. As we assume
$\partial_C I \not=\partial I \setminus \partial_0 I$, it holds that also
$\partial_M I \not= \emptyset$, recall \eqref{eq:partialM}. 
For ease of exposition, let $\partial_M I = \{0\}$ and $\partial_C I = \{1\}$.
It follows from the first condition in \eqref{eq:assC} that
\[
\frac{\vec U(q_{j+1}) - \vec U(q_j)}{|\vec X^m(q_{j+1}) - \vec X^m(q_j)|}
= \frac{\vec U(q_{j}) - \vec U(q_{j-1})}{|\vec X^m(q_{j}) - \vec X^m(q_{j-1})|}\,,
\qquad j = 1,\ldots,J-1\,,
\]
and so, on recalling Assumption \ref{ass:spanvm0}, 
there exist positive numbers $\alpha_j$ such that
\begin{equation} \label{eq:Ualphaj}
\vec U(q_{j+1}) = (1+ \alpha_j)\,\vec U(q_j) - \alpha_j\,\vec U(q_{j-1})\,,
\qquad j = 1,\ldots,J-1\,.
\end{equation}
Combining \eqref{eq:Ualphaj} and the fact that $\vec U \in \yspace^h(\vec 0)$,
i.e.\ $\vec U(q_0) = \vec 0$, we obtain, via induction, that there exist  
numbers $\zeta_{J-1} \geq  \zeta_{J-2}\geq \ldots \geq \zeta_1 >0$ such that
\[
\vec U(q_{j+1}) = (1+\zeta_j)\,\vec U(q_1)\,, \qquad j = 1,\ldots,J-1\,.
\]
Hence it follows from the second property in \eqref{eq:assC} and the 
Assumption~\ref{ass:spanvm0}, recall \eqref{eq:omegah}, that
$\vec U(q_1) = \vec 0$, and so $\vec U =  \vec 0$ in $I$. 
This completes the proof.
\end{proof}

\begin{lmm} \label{lem:ex}
Let {\rm Assumption~\ref{ass:spanvm0}} hold.
Moreover, if $\emptyset \not= \partial_C I$
and if $\partial_C I= \partial I \setminus \partial_0 I$
then let also {\rm Assumption~\ref{ass:C}} hold.
Let $\vec X^m, \vec Y^m \in \Vhpartialzero$, 
$\kappa^m \in V^h$, $\vec{\rm m}^m \in \bR^2$ 
and $\alpha \in \bRplus$,
$\spont,M_0,\alpha_G \in \bR$, $\beta \in \bRgeq$, 
$\vec\zeta : \partial_C I \to \bR^2$ be given.
Then there exists a unique solution to $(\BGNpwf^m)^h$, \eqref{eq:fd}.
\end{lmm}
\begin{proof}
As we have a linear system of equations, with the same number of equations as
unknowns, existence follows from uniqueness. Hence we consider a solution
to the homogeneous equivalent of \eqref{eq:fd}, and need to show that this
solution is in fact zero. In particular, let
$\delta\vec X \in \xspaceh$, 
$\kappa \in \Whpartialzero$, $\vec Y \in
\yspaceh(\vec0)$, and $\vec{\rm{m}} : \partial_M I \to \bR^2$
be such that
\begin{subequations} 
\begin{align}
& 2\,\pi\left( (\vec X^m\,.\,\vec\ek_1)\,\mat Q^m\,
\delta\vec X, \vec\chi\,|\vec X^m_\rho|\right)^h 
- \ttau_m \left(\vec Y_\rho , \vec\chi_\rho\,|\vec X^m_\rho|^{-1} \right) 
= 0 \qquad \forall\ \vec\chi \in \xspaceh\,, \label{eq:fdproofa} \\
&2\,\pi\,\alpha \left( 
\vec X^m\,.\,\vec\ek_1\,\kappa ,
\chi \,|\vec X^m_\rho|\right)^h
- \left( \vec Y, \chi\,\vec\nu^m \,|\vec X^m_\rho|\right)^h = 0
\qquad \forall\ \chi \in \Whpartialzero\,,\label{eq:fdproofb} \\
&\left( \kappa\,\vec\nu^m, \vec\eta \,|\vec X^m_\rho|\right)^h
+ \left( (\delta\vec X)_\rho, \vec\eta_\rho\,|\vec X^m_\rho|^{-1} \right) = 
\sum_{p \in \partial_M I} \left[\vec{\rm m}\,.\,\vec\eta\right](p)
\qquad \forall\ \vec\eta \in \Vhpartialzero\,. \label{eq:fdproofc}
\end{align}
\end{subequations}
Choosing $\vec\chi=\delta\vec X$ in \eqref{eq:fdproofa}, 
$\chi = \kappa$ in \eqref{eq:fdproofb} and
$\vec\eta = \vec Y$ in \eqref{eq:fdproofc} yields
that
\begin{equation} \label{eq:fdproof1}
2\,\pi\left( \vec X^m\,.\,\vec\ek_1\,|\mat Q^m\,\delta\vec X|^2,
|\vec X^m_\rho|\right)^h
+ 2\,\pi\,\alpha\,\ttau_m \left( \vec X^m\,.\,\vec\ek_1\,
\kappa^2 , |\vec X^m_\rho|\right)^h = 0\,.
\end{equation}
It follows from \eqref{eq:fdproof1} and $\kappa \in \Whpartialzero$
that $\kappa=0$.
Similarly, it follows from \eqref{eq:fdproof1}, 
$\delta\vec X \in \Vhpartialzero$,
\eqref{eq:Vpartialzero}, \eqref{eq:vh} and \eqref{eq:Qh} 
that $\pi^h\,[\delta\vec X\,.\,\vec v^m] = 0$
with $\delta\vec X = \vec 0$ on $\partial I \setminus \partial_0 I$.
Then choosing $\vec\eta =\delta\vec X \in \xspace^h \subset \Vhpartialzero$
in \eqref{eq:fdproofc} yields, on recalling \eqref{eq:partialM}
and $\delta\vec X = \vec 0$ on $\partial I \setminus \partial_0 I$, that
\begin{equation} \label{eq:fdproof2}
\left( |(\delta\vec X)_\rho|^2, |\vec X^m_\rho|^{-1}\right)
= \sum_{p \in \partial_M I} \left[\vec{\rm m}\,.\,\delta\vec X\right](p)
= 0\,.
\end{equation}
It follows from \eqref{eq:fdproof2} that $\delta\vec X$ is constant in 
$\overline I$. Together with $\pi^h\,[\delta\vec X\,.\,\vec v^m] = 0$
and Assumption~\ref{ass:spanvm0} we obtain that $\delta\vec X = \vec 0$.
Now \eqref{eq:fdproofc} implies that $\vec{\rm m} = \vec 0$.

It remains to show that $\vec Y = \vec 0$. We have from \eqref{eq:fdproofa} 
and \eqref{eq:fdproofb} that
\[
 \left(\vec Y_\rho , \vec\chi_\rho\,|\vec X^m_\rho|^{-1} \right) = 0
\quad \forall\ \vec\chi \in \xspaceh
\quad\text{and}\quad
\left( \vec Y, \chi\,\vec\nu^m \,|\vec X^m_\rho|\right)^h = 0
\quad \forall\ \chi \in \Whpartialzero\,.
\]
Moreover, our assumptions and Lemma~\ref{lem:assC} yield that
Assumption~\ref{ass:C} holds. Hence we have that $\vec Y = \vec 0$,
and thus we have shown the existence of a unique
solution to $(\BGNpwf^m)^h$.

\end{proof}

\begin{rmrk} \label{rem:clamp}
We note that it is not possible to prove the existence of a unique solution to
\eqref{eq:fd} in the case of clamped boundary conditions, when
$\partial_C I = \partial I = \{0,1\}$.
The authors faced a similar issue in the context of the approximation of
Willmore flow for general open surfaces in {\rm \cite{pwfopen}}. There
the problem could be overcome by a suitable tweak to the discretization,
see {\rm (3.19)} and {\rm Theorem~4.1} there. In particular, the approach 
applied in {\rm \cite{pwfopen}} relied on the discretization of the side
constraint {\rm (2.22)} there, which is formulated in terms of the mean
curvature vector $\vec k_m = k_m\,\vec{\rm n}_{\mathcal{S}}$ of $\mathcal{S}$, 
rather than in terms of the scalar curvature variables $\varkappa$ and
$\varkappa_{\mathcal{S}}$ that we consider here, recall 
\eqref{eq:varkappaweak} and \eqref{eq:varkappaSweak}. 
Hence the approach from {\rm \cite{pwfopen}}
is not applicable to the situation considered in this paper.
\end{rmrk}

\begin{rmrk} \label{rem:rmm}
We note that in practice it is easiest to find the solution to 
\eqref{eq:fd} by first eliminating $\vec{\rm m}^{m+1}$ from \eqref{eq:fd}
via replacing the test space $\Vhpartialzero$ in \eqref{eq:fdc} 
with $\yspace^h(\vec 0)$.
Having computed $(\delta\vec X^{m+1}, \kappa^{m+1}, \vec Y^{m+1})$ in this way,
for example with the help of a sparse factorization package like UMFPACK,
see {\rm \cite{Davis04}}, the values $\vec{\rm m}^{m+1}$ can be obtained
from \eqref{eq:fdc}. 
For example, if $q_0 \in \partial_M I$ then we have, on recalling
\eqref{eq:omegah}, that
\[ 
\vec{\rm m}^{m+1}(q_0) = 
\left( 1, \chi_0\,|\vec X^m_\rho|\right) \kappa^{m+1}(q_0)\,
\vec\omega^m(q_0)
+ \left( 1 , (\chi_0)_\rho \,|\vec X^m_\rho|^{-1}\right) 
\vec X^{m+1}_\rho(q_0) \,.
\]
\end{rmrk}

\begin{rmrk} \label{rem:stab}
In light of our stability result {\rm Theorem~\ref{thm:stab}} for the
semidiscrete scheme \eqref{eq:sd}, it would be desirable to also prove
(conditional) stability for the fully discrete scheme \eqref{eq:fd}.
However, at present this remains an open problem, 
in line with other fully
discrete schemes for Willmore flow and elastic flow in the literature, see
e.g.\ \cite{willmore,Dziuk08,DeckelnickD09,pwf,pwfade}.
Note that the only stability result for a fully discrete scheme
for Willmore flow, that we are aware of, is given in \cite{Rusu05}, 
for a scheme where the tangential velocity is zero. 
\end{rmrk}

\subsubsection{Conserved flows} \label{sec:fdC}

Here, following the approach in \cite[\S4.3.1]{axisd}, we consider fully
discrete variants of the semidiscrete conserving approximations in
\S\ref{sec:new51}. In particular, on rewriting \eqref{eq:fda} as
\[
 2\,\pi\left(\vec X^m\,.\,\vec\ek_1\,\mat Q^m\,
\frac{\vec X^{m+1} - \vec X^m}{\ttau_m}, \vec\chi\,|\vec X^m_\rho|\right)^h 
- \left(\vec Y^{m+1}_\rho , \vec\chi_\rho\,|\vec X^m_\rho|^{-1} \right) 
 = \left( \vec f^m, \vec\chi \, |\vec X^m_\rho| \right)^h ,
\]
we can formulate our surface area and volume conserving variant for
$(\BGNpwf^m)^h$ as follows. Here, for ease of presentation, we assume
that $\partial_M I =  \emptyset$, so that we do not need to consider
the discrete conormals $\vec{\rm m}^{m+1}$.
Moreover, on recalling (\ref{eq:yspace}), we have that
$\yspace^h(2\,\pi\,\alpha_G\,\vec\ek_1) = \Vhpartialzero$.

\noindent
$(\BGNpwf^m_{A,V})^h$:
Let $\vec X^0 \in \Vhpartialzero$, $\kappa^0 \in V^h$, 
$\vec Y^0 \in \Vhpartialzero$ and $\alpha \in \bRplus$,
$\spont,M_0,\alpha_G,\lambda,\varsigma \in \bR$, $\beta \in \bRgeq$, 
$\vec\zeta : \partial_C I \to \bS^1$ be given.
For $m=0,\ldots,M-1$, find $\delta\vec X^{m+1} \in \xspace^h$, 
with $\vec X^{m+1} = \vec X^m + \delta\vec X^{m+1}$, 
$\kappa^{m+1} \in \Whpartialzero$,
$\vec Y^{m+1} \in \Vhpartialzero$,
and $\lambda_A^{m+1},\lambda_V^{m+1} \in \bR$ such that 
\eqref{eq:fdb}, \eqref{eq:fdc} and
\begin{subequations} \label{eq:fdwf}
\begin{align}
& 2\,\pi\left(\vec X^m\,.\,\vec\ek_1\,\mat Q^m\,
\frac{\vec X^{m+1} - \vec X^m}{\ttau_m}, \vec\chi\,|\vec X^m_\rho|\right)^h 
- \left(\vec Y^{m+1}_\rho , \vec\chi_\rho\,|\vec X^m_\rho|^{-1} \right) 
= \left( \vec f^m, \vec\chi \, |\vec X^m_\rho| \right)^h
\nonumber \\ & \qquad
- 2\,\pi\,\lambda_A^{m+1} \left[
\left( \vec\ek_1, \vec\chi\,|\vec X^m_\rho|\right) 
+ \left( (\vec X^m\,.\,\vec\ek_1)\,\vec\tau^m
, \vec\chi_\rho\right) \right]
- 2\,\pi\,\lambda_V^{m+1} \left((\vec X^m\,.\,\vec\ek_1)\,\vec\nu^m,
\vec\chi \, |\vec X^m_\rho| \right)
\qquad \forall\ \vec\chi \in \xspaceh\,, \label{eq:fdwfa} \\
&{\rm (i)}\ A(\vec X^{m+1}) = A(\vec X^0)\,,\quad
{\rm (ii)}\ V(\vec X^{m+1}) = V(\vec X^0) \label{eq:fdwfb} 
\end{align}
\end{subequations}
hold.
The nonlinear system of equations arising at each time level of
$(\BGNpwf^m_{A,V})^h$ can be solved with a suitable iterative solution method, 
see below.
In the simpler case of surface area conserving flow, we need to find
$(\delta\vec X^{m+1}, \kappa^{m+1}, \vec Y^{m+1},$ $ 
\lambda_A^{m+1}, \lambda_V^{m+1}) 
\in \xspace^h \times \Whpartialzero \times \Vhpartialzero
\times \bR \times \{0\}$
such that \eqref{eq:fdb}, \eqref{eq:fdc}, \eqref{eq:fdwfa} and 
\eqref{eq:fdwfb}(i) hold. 
Similarly, for volume conserving flow, we need to find
$(\delta\vec X^{m+1}, \kappa^{m+1}, \vec Y^{m+1}, 
\lambda_A^{m+1}, \lambda_V^{m+1}) 
\in \xspace^h \times \Whpartialzero \times \Vhpartialzero
\times \{0\} \times \bR$
such that \eqref{eq:fdb}, \eqref{eq:fdc}, \eqref{eq:fdwfa} and 
\eqref{eq:fdwfb}(ii) hold. 

Adapting the strategy in \cite[\S4.3.1]{axisd}, 
we now describe a Newton method for solving the nonlinear system 
(\ref{eq:fdwf}), \eqref{eq:fdb} and \eqref{eq:fdc}.
The linear system (\ref{eq:fdwfa}), \eqref{eq:fdb} and \eqref{eq:fdc}, 
with $(\lambda_A^{m+1}, \lambda_V^{m+1})$ in
(\ref{eq:fdwfa}) replaced by $(\lambda_A, \lambda_V)$, can be written as:
Find $(\delta\vec X^{m+1}(\lambda_A,\lambda_V), 
\kappa^{m+1}(\lambda_A,\lambda_V), \vec Y^{m+1}(\lambda_A,\lambda_V)) 
\in \xspace^h\times \Whpartialzero \times \Vhpartialzero$
such that
\begin{equation*} 
\mathbb{T}^m\,\begin{pmatrix}
\vec Y^{m+1}(\lambda_A,\lambda_V)\\[1mm]
\delta\vec X^{m+1}(\lambda_A,\lambda_V)\\[1mm]
\kappa^{m+1}(\lambda_A,\lambda_V)
\end{pmatrix}
= 
\begin{pmatrix} \vec{\underline{\mathfrak g}}^m_1 \\[1mm] 
\vec{\underline{\mathfrak g}}^m_2 \\[1mm] 
\vec{\underline{\mathfrak g}}^m_3
\end{pmatrix}
+ \lambda_A\, \begin{pmatrix} \vec{\underline{\mathcal K}}^m 
\\[1mm] 0 \\[1mm]\vec 0
\end{pmatrix}
+ \lambda_V\, \begin{pmatrix} \vec{\underline{\mathcal N}}^m 
\\[1mm] 0 \\[1mm]\vec 0
\end{pmatrix}.
\end{equation*}
Assuming the linear operator $\mathbb{T}^m$ is invertible, we obtain that
\begin{align}
\begin{pmatrix}
\vec Y^{m+1}(\lambda_A,\lambda_V)\\[1mm]
\delta\vec X^{m+1}(\lambda_A,\lambda_V)\\[1mm]
\kappa^{m+1}(\lambda_A,\lambda_V)
\end{pmatrix} &
= (\mathbb{T}^m)^{-1}
\left[
\begin{pmatrix} \vec{\underline{\mathfrak g}}^m_1 \\[1mm] 
\vec{\underline{\mathfrak g}}^m_2 \\[1mm] 
\vec{\underline{\mathfrak g}}^m_3
\end{pmatrix}
+ \lambda_V\, 
\begin{pmatrix} \vec{\underline{\mathcal N}}^m \\[1mm] 0 \\[1mm] \vec 0
\end{pmatrix}\right] 
=: (\mathbb{T}^m)^{-1} 
\begin{pmatrix} \vec{\underline{\mathfrak g}}^m_1 \\[1mm] 
\vec{\underline{\mathfrak g}}^m_2 \\[1mm] 
\vec{\underline{\mathfrak g}}^m_3
\end{pmatrix}
+ \lambda_A \begin{pmatrix} \vec{\underline s}^m_1 \\[1mm] 
\vec{\underline s}^m_2 \\[1mm] {\underline s}^m_3
\end{pmatrix}
+ \lambda_V\, \begin{pmatrix} \vec{\underline q}^m_1 \\[1mm] 
\vec{\underline q}^m_2 \\[1mm] {\underline q}^m_3
\end{pmatrix} .
\label{eq:lmsysinverse}
\end{align}
It immediately follows from (\ref{eq:lmsysinverse}) that
\begin{equation*}
\partial_{\lambda_A} \vec X^{m+1}(\lambda_A,\lambda_V) 
= \vec{\underline s}^m_2\,,\quad
\partial_{\lambda_V} \vec X^{m+1}(\lambda_A,\lambda_V)
= \vec{\underline q}^m_2\,,
\end{equation*}
where $\vec X^{m+1}(\lambda_A,\lambda_V) = \vec X^m + 
\delta\vec X^{m+1}(\lambda_A,\lambda_V)$.
Hence we can proceed as in \cite[(4.13)]{axisd} to define a Newton iteration
for finding a solution to the nonlinear system $(\BGNpwf^m_{A,V})^h$.
In practice this Newton iteration always 
converged within a couple of iterations.

\subsection{Based on $\kappa^{m+1}_{\mathcal{S}}$} \label{sec:fdkappaS}

We can consider the following fully discrete approximation of (\ref{eq:sdS}).

\noindent
$(\BGNpwfwf^m)^{(h)}$
Let $\vec X^0, \vec Y^0_{\mathcal{S}} \in \Vhpartialzero$, 
$\kappa^0_{\mathcal{S}} \in
V^h$, $\vec{\rm m}^0 \in \bR^2$ and $\alpha \in \bRplus$,
$\spont,M_0,\alpha_G,\lambda,\varsigma \in \bR$, $\beta \in \bRgeq$, 
$\vec\zeta : \partial_C I \to \bS^1$ be given.
For $m=0,\ldots,M-1$, find $\delta\vec X^{m+1} \in \xspaceh$, with
$\vec X^{m+1} = \vec X^m + \delta\vec X^{m+1}$, 
$\kappa^{m+1}_{\mathcal{S}} \in W^h_{(\partial_0)}$, 
$\vec Y^{m+1}_{\mathcal{S}} \in \Vhpartialzero$,
with 
$\vec\pi^h\left[(\vec X^m\,.\,\vec\ek_1)\, \vec Y^{m+1}_{\mathcal{S}}\right] 
\in \yspaceh(2\,\pi\,\alpha_G \,\vec\ek_1)$, 
and $\vec{\rm m}^{m+1} : \partial_M I \to \bR^2$ such that
\begin{subequations} \label{eq:fdS}
\begin{align}
& 2\,\pi\left( (\vec X^m\,.\,\vec\ek_1)\,\mat Q^m\,
\frac{\vec X^{m+1} - \vec X^m}{\ttau_m}, \vec\chi\,|\vec X^m_\rho|\right)^{h} 
- \left( (\vec X^m\,.\,\vec\ek_1)\,
(\vec Y^{m+1}_{\mathcal{S}})_\rho, \vec\chi_\rho\,|\vec X^m_\rho|^{-1}\right)
 \nonumber \\ & \qquad
+ \left( (\vec X^m\,.\,\vec\ek_1)\,
(\vec Y^{m}_{\mathcal{S}})_\rho\,.\,\vec\tau^m, 
\vec\chi_\rho\,.\,\vec\tau^m\,|\vec X^m_\rho|^{-1} \right) 
\nonumber \\  
& \quad =
-  \left( 
\pi\,\vec X^m\,.\,\vec\ek_1\,
\left[ \alpha \,[\kappa^m_{\mathcal{S}} - \spont]^2 
+ 2\,\lambda + 2\,\beta\, \AADE_{\mathcal{S}}^m\,\kappa_{\mathcal{S}}^m
\right]
- \vec Y^m_{\mathcal{S}}\,.\,\vec\ek_1  ,
\vec\chi_\rho\,.\,\vec\tau^m \right)^{(h)}
\nonumber \\ & \qquad  -
\left( \left[ \pi \left[ 
\alpha\left[\kappa_{\mathcal{S}}^m -\spont \right]^2 
+ 2\,\lambda + 2\,\beta\, \AADE_{\mathcal{S}}^m\,\kappa_{\mathcal{S}}^m \right]
- \kappa^m_{\mathcal{S}}\,\vec Y^m_{\mathcal{S}}\,.\,\vec\nu^m \right] 
|\vec X^m_\rho| - (\vec Y^m_{\mathcal{S}})_\rho\,.\,\vec\tau^m,
\vec\chi\,.\,\vec\ek_1 \right)^{(h)}
\nonumber \\ & \qquad 
+ \left( \vec X^m\,.\,\vec\ek_1\,\kappa^m_{\mathcal{S}}\, 
(\vec Y^m_{\mathcal{S}})^\perp , \vec\chi_\rho \right)^{(h)}
- \sum_{p \in \partial_{2} I \cup \partial_F I}
\left[ \left[2\,\pi\,\varsigma + \vec{\rm m}^m\,.\,\vec Y^m_{\mathcal{S}}
\right]
(\vec\chi\,.\,\vec\ek_1)\,\right](p)
 \qquad \forall \ \vec\chi \in \xspaceh \,,
\label{eq:fdSa}\\
& 2\,\pi\left( \vec X^m\,.\,\vec\ek_1\left(
\alpha\,(\kappa_{\mathcal{S}}^{m+1} - \spont) +\beta\,\AADE^m_{\mathcal{S}}
\right) , \chi\,|\vec X^m_\rho|\right)^{(h)}\!
- \left( (\vec X^m\,.\,\vec\ek_1)\,\vec Y^{m+1}_{\mathcal{S}}, 
\chi\,\vec\nu^m \,|\vec X^m_\rho|\right)^{(h)} = 0 
\quad \forall\ \chi \in \Whpartialzerob\,, \label{eq:fdSb} \\
& \left( \vec X^m\,.\,\vec\ek_1\,\kappa_{\mathcal{S}}^{m+1}\,\vec\nu^m,
\vec\eta \,|\vec X^m_\rho|\right)^{(h)}
+ \left(\vec\ek_1, \vec\eta \,|\vec X^m_\rho|\right)
+ \left( (\vec X^m\,.\,\vec\ek_1)\,
\vec X^{m+1}_\rho, \vec\eta_\rho \,|\vec X^m_\rho|^{-1}\right)
\nonumber \\ & \hspace{4cm}
= \sum_{p \in \partial_C I} \left[(\vec X^m\,.\,\vec\ek_1)\,
\vec\zeta\,.\,\vec\eta\right](p)
+ \sum_{p \in \partial_M I} \left[(\vec X^m\,.\,\vec\ek_1)\,
\vec{\rm m}^{m+1}\,.\,\vec\eta\right](p) \qquad
\forall\ \vec\eta \in \Vhpartialzero\,, \label{eq:fdSc}
\end{align}
where 
\begin{equation}
\AADE_{\mathcal{S}}^m =
 2\,\pi\left(\vec X^m\,.\,\vec\ek_1\,\kappa^m_{\mathcal{S}}, 
|\vec X^m_\rho| \right)^{(h)} - M_0\,.
\label{eq:fdSd}
\end{equation}
\end{subequations}

We now state the analogues of Assumptions \ref{ass:spanvm0} and \ref{ass:C}.  
\begin{ass} \label{ass:spanX1num}
Let $\dim\spa\left\{ 
\left( (\vec X^m\,.\,\vec\ek_1)\, \vec\nu^m,
\chi \,|\vec X^m_\rho|\right)^{(h)} : \chi \in \Whpartialzerob \right\} = 2$.
\end{ass}

\begin{ass} \label{ass:CS}
Let $\vec X^m$ satisfy {\rm Assumption~\ref{ass:Ah1}}
with $\vec X^h$ replaced by $\vec X^m$ and be such that the
following holds. If $\vec U \in \yspace^h(\vec 0)$ with
\[
 \left( (\vec X^m\,.\,\vec \ek_1)\,\vec U_\rho , \vec\chi_\rho\,|\vec X^m_\rho|^{-1} \right) = 0
\quad \forall\ \vec\chi \in \xspaceh 
\qquad\text{and}\qquad
\left( (\vec X^m\,.\,\vec\ek_1)\,\vec U, \chi\,\vec\nu^m \,|\vec X^m_\rho|
\right)^{(h)} = 0 \quad \forall\ \chi \in \Whpartialzerob\,, 
\]
then $\vec U = \vec 0$.
\end{ass}

\begin{lmm} \label{lem:exS}
Let {\rm Assumptions~\ref{ass:spanvm0}} and {\rm \ref{ass:spanX1num}} hold.
Moreover, if $\partial_C I \not= \emptyset$ then let
{\rm Assumptions~\ref{ass:CS}} hold.
Let $\vec X^m, \vec Y^m_{\mathcal{S}} \in \Vhpartialzero$, 
$\kappa^m_{\mathcal{S}} \in V^h$, $\vec{\rm m}^m \in \bR^2$ 
and $\alpha \in \bRplus$,
$\spont,\lambda,M_0,\alpha_G \in \bR$, $\beta \in \bRgeq$, 
$\vec\zeta : \partial_C I \to \bS^1$ be given.
Then there exists a unique solution to $(\BGNpwfwf^m)^{(h)}$, (\ref{eq:fdS}).
\end{lmm}
\begin{proof}
As we have a linear system of equations, with the same number of equations as
unknowns, existence follows from uniqueness. Hence we consider a solution
to the homogeneous equivalent of \eqref{eq:fdS}, and need to show that this
solution is in fact zero. In particular, let
$\delta\vec X \in \xspaceh$, 
$\kappa_{\mathcal{S}} \in W^h_{(\partial_0)}$, 
$\vec Y_{\mathcal{S}} \in \Vhpartialzero$,
with $\pi^h\left[(\vec X^m\,.\,\vec\ek_1)\, \vec Y_{\mathcal{S}}\right] \in 
\yspaceh(\vec0)$, and $\vec{\rm{m}} \in \bR^2$ be such that
\begin{subequations} 
\begin{align}
& 2\,\pi\left( (\vec X^m\,.\,\vec\ek_1)\,\mat Q^m\,
\frac{\delta\vec X}{\ttau_m}, \vec\chi \,|\vec X^m_\rho|\right)^h
- \left( (\vec X^m\,.\,\vec\ek_1)\,
(\vec Y_{\mathcal{S}})_\rho, \vec\chi_\rho\,|\vec X^m_\rho|^{-1}\right)
= 0 \qquad \forall \ \vec\chi \in \xspaceh \,,
\label{eq:fdSproofa}\\
& 2\,\pi\,\alpha\left( \vec X^m\,.\,\vec\ek_1\,
\kappa_{\mathcal{S}} , \chi\,|\vec X^m_\rho|\right)^{(h)}
- \left( (\vec X^m\,.\,\vec\ek_1)\,\vec Y_{\mathcal{S}}, 
\chi\,\vec\nu^m \,|\vec X^m_\rho|\right)^{(h)} = 0 
\qquad \forall\ \chi \in \Whpartialzerob\,, \label{eq:fdSproofb} \\
& \left( \vec X^m\,.\,\vec\ek_1\,\kappa_{\mathcal{S}}\,\vec\nu^m,
\vec\eta \,|\vec X^m_\rho|\right)^{(h)}
+ \left( (\vec X^m\,.\,\vec\ek_1)\,
(\delta\vec X)_\rho, \vec\eta_\rho \,|\vec X^m_\rho|^{-1}\right)
= \sum_{p \in \partial_M I} \left[(\vec X^m\,.\,\vec\ek_1)\,
\vec{\rm m}\,.\,\vec\eta\right](p) \quad
\forall\ \vec\eta \in \Vhpartialzero\,. \label{eq:fdSproofc}
\end{align}
\end{subequations}
Choosing $\vec\chi=\delta\vec X$ in \eqref{eq:fdSproofa}, 
$\chi = \kappa_{\mathcal{S}}$ in \eqref{eq:fdSproofb} and
$\vec\eta = \vec Y_{\mathcal{S}}$ in \eqref{eq:fdSproofc} yields,
on noting 
$\pi^h\left[(\vec X^m\,.\,\vec\ek_1)\, \vec Y_{\mathcal{S}}\right] \in 
\yspaceh(\vec0)$, that
\begin{equation} \label{eq:fdSproof1}
2\,\pi\left( \vec X^m\,.\,\vec\ek_1\,|\mat Q^m\,\delta\vec X|^2,
|\vec X^m_\rho|\right)^h
+ 2\,\pi\,\alpha\,\ttau_m \left( \vec X^m\,.\,\vec\ek_1\,
(\kappa_{\mathcal{S}})^2 , |\vec X^m_\rho|\right)^{(h)} = 0\,.
\end{equation}
It follows from \eqref{eq:fdSproof1} and 
$\kappa_{\mathcal{S}} \in W^h_{(\partial_0)}$ that $\kappa_{\mathcal{S}}=0$.
Similarly, it follows from \eqref{eq:fdSproof1}
that $[\delta\vec X\,.\,\vec v^m](q_j) = 0$, $j=1,\ldots,J-1$, and, 
on recalling \eqref{eq:Qh}, that $\delta\vec X = \vec 0$
on $\partial\,I \setminus \partial_0 I$.
Then choosing $\vec\eta =\delta\vec X \in \xspace^h \subset \Vhpartialzero$
in \eqref{eq:fdSproofc} yields, on recalling \eqref{eq:partialM}
and $\delta\vec X = \vec 0$ on $\partial\,I \setminus \partial_0 I$, that
\begin{equation} \label{eq:fdSproof2}
\left( \vec X^m\,.\,\vec\ek_1\,|(\delta\vec X)_\rho|^2, 
|\vec X^m_\rho|^{-1}\right) = \sum_{p \in \partial_M I} 
\left[(\vec X^m\,.\,\vec\ek_1)\,\vec{\rm m}\,.\,\delta\vec X\right](p)
=0 \,.
\end{equation}
It follows from \eqref{eq:fdSproof2} that $\delta\vec X$ is constant in 
$\overline I$. Together with 
$[\delta\vec X\,.\,\vec v^m](q_j) = 0$, $j=1,\ldots,J-1$, and
Assumption~\ref{ass:spanvm0} we obtain that $\delta\vec X = \vec 0$.
Then \eqref{eq:fdSproofc} implies that $\vec{\rm m} = \vec 0$.

If $\partial_C I = \emptyset$, then
we can choose $\vec\chi = \vec Y_{\mathcal{S}} \in \yspace^h(\vec 0) \subset 
\xspace^h$ in 
\eqref{eq:fdSproofa}, recall \eqref{eq:partialCempty},  
to obtain that $\vec Y_{\mathcal{S}}$ is constant in $\overline I$.
Then \eqref{eq:fdSproofb}, together with Assumption~\ref{ass:spanX1num},
gives that $\vec Y_{\mathcal{S}} = \vec 0$. 
If $\partial_C I \not= \emptyset$, on the other hand,
then Assumption~\ref{ass:CS} directly gives that $\vec Y_{\mathcal{S}}=\vec 0$.
Hence we have shown the existence of a unique solution to $(\BGNpwfwf^m)^h$.
\end{proof}

\begin{rmrk}
Similarly to {\rm Remark~\ref{rem:rmm}}, in practice the system
\eqref{eq:fdS} is easiest solved by first eliminating 
$\vec{\rm m}^{m+1}$. This can be achieved by replacing the test space 
$\Vhpartialzero$ in \eqref{eq:fdSc} with $\yspace^h(\vec 0)$. 
Having computed
$(\delta\vec X^{m+1}, \kappa^{m+1}_{\mathcal{S}}, \vec Y^{m+1}_{\mathcal{S}})$,
the values $\vec{\rm m}^{m+1}$ can then be obtained from \eqref{eq:fdSc}. 
For example, if $q_0 \in \partial_M I$ then we have that
\begin{align*} 
(\vec X^m(q_0)\,.\,\vec\ek_1)\,\vec{\rm m}^{m+1}(q_0) & = 
\left( \vec X^m\,.\,\vec\ek_1\,\kappa_{\mathcal{S}}^{m+1},
\chi_0\,|\vec X^m_\rho|\right)^{(h)}\vec\nu^m(q_0)
+ \left(1, \chi_0 \,|\vec X^m_\rho|\right) \vec\ek_1 \nonumber \\ & \qquad
+ \left( \vec X^m\,.\,\vec\ek_1, (\chi_0)_\rho \,|\vec X^m_\rho|^{-1}\right) 
\vec X^{m+1}_\rho(q_0) \,.
\end{align*}
\end{rmrk}

\setcounter{equation}{0}
\section{Numerical results} \label{sec:nr}
On recalling (\ref{eq:Eh}), (\ref{eq:kappahYh}) and (\ref{eq:ESh}),
as the fully discrete energy for the two schemes $(\BGNpwf^m)^{h}$
and $(\BGNpwfwf^{m})^{(h)}$ we consider, respectively,
\begin{subequations}
\begin{align} \label{eq:Em}
\widehat E^{m+1} & = \pi\left(\alpha 
\left[ \doctorkappa^h(\vec X^m, \kappa^{m+1}) - \spont \right]^2 +2 \,\lambda,  
\vec X^m\,.\,\vec\ek_1 \,|\vec X^m_\rho|\right)^h 
+ \tfrac{\beta}{2} \left[
2\,\pi\left(\vec X^m\,.\,\vec\ek_1\,\kappa^{m+1} - \vec\nu^m.\,\vec\ek_1, 
|\vec X^m_\rho| \right)^h - M_0 \right]^2 \nonumber \\ & \qquad
- 2\,\pi\,\alpha_G\sum_{p\in \partial_M I}
\vec{\rm m}^{m+1}(p)\,.\,\vec\ek_1
+ 2\,\pi\,\varsigma\sum_{p\in \partial_2 I \cup \partial_F I}
\vec X^{m+1}(p)\,.\,\vec\ek_1\,, \\
\label{eq:EmS}
\widehat E^{m+1}_{\mathcal{S}} & = \pi\left(\alpha\, 
[ \kappa_{\mathcal{S}}^{m+1} -\spont ]^2 + 2\,\lambda ,
\vec X^m\,.\,\vec\ek_1 |\vec X^m_\rho|\right)^{(h)}
+ \tfrac\beta2\left[ 
2\,\pi\left(\vec X^m\,.\,\vec\ek_1\,\kappa^{m+1}_{\mathcal{S}}, 
|\vec X^m_\rho| \right)^{(h)} - M_0 \right]^2
 \nonumber \\ & \quad
- 2\,\pi\,\alpha_G\sum_{p\in \partial_M I}
\vec{\rm m}^{m+1}(p)\,.\,\vec\ek_1
+ 2\,\pi\,\varsigma\sum_{p\in \partial_2 I \cup \partial_F I}
\vec X^{m+1}(p)\,.\,\vec\ek_1\,.
\end{align}
\end{subequations}

Given $\Gamma^0 = \vec X^0(I)$, we define the following initial data.
First, we let $\vec\mu^0$ be the true conormal to $\Gamma^0$, i.e.\
$\vec\mu^0 = (-1)^{p+1}\,\vec\tau^0(p)$ for $p \in \partial I\setminus
\partial_0 I$, recall \eqref{eq:mu}, and then set 
$\vec{\rm m}^0 = \vec\mu^0$ on $\partial_M I$.
Next, on recalling (\ref{eq:sideh}), we let $\vec\kappa^0\in \Vhpartialzero$ 
be such that
\begin{equation*} 
\left( \vec\kappa^{0},\vec\eta \,|\vec X^0_\rho|\right)^h
+ \left( \vec\tau^{0} , \vec\eta_\rho \right)
= \sum_{p \in \partial I\setminus \partial_0 I} 
\left[\vec\mu^0\,.\,\vec\eta\right](p)
\quad \forall\ \vec\eta \in \Vhpartialzero\,,
\end{equation*}
and then define $\kappa^0 = \pi^h_{\partial_0}
[\vec\kappa^{0}\,.\,\vec v^0]\in \Whpartialzero$.
Moreover, we let $\vec Y^0_\star \in \Vh$ be such that
\begin{equation*} 
\vec Y^0_\star = 2\,\pi\,\vec\pi^h\left[|\vec\omega^0|^{-1}
\vec X^0\,.\,\vec\ek_1\left( \alpha\left[ \doctorkappa^h(\vec X^0, \kappa^0) -
\spont\right] + \beta\,\AADE^0 \right) \vec v^0 \right],
\end{equation*}
recall (\ref{eq:kappahYh}), \eqref{eq:omegah}, \eqref{eq:vh} and \eqref{eq:fdd}, and then define
$\vec Y^0 \in \yspace^h(2\,\pi\,\alpha_G\,\vec\ek_1)$ via
\begin{equation*} 
\vec Y^0(q_j) = \begin{cases}
(\vec Y^0_\star(q_j)\,.\,\vec\ek_2)\,\vec\ek_2 & q_j \in \partial_0 I\,,\\
2\,\pi\,\alpha_G\,\vec\ek_1 & q_j \in \partial_M I\,,\\
\vec Y^0_\star(q_j) & q_j \in \overline I \setminus (\partial_0 I \cup 
\partial_M I)\,.
\end{cases}
\end{equation*}
In addition, we set $\kappa^0_{\mathcal{S}} = \doctorkappa^h(\vec X^0, \kappa^0)
\in V^h$, and let $\vec Y^0_{\mathcal{S},\star} \in \Vh$ be such 
that
\begin{equation*} 
\vec Y^0_{\mathcal{S},\star} = 2\,\pi\,\vec\pi^h\left[
|\vec\omega^0|^{-1}\left(\alpha\,[\kappa^0_{\mathcal{S}} - \spont]
+ \beta\,\AADE^0_{\mathcal{S}}\right) \vec v^0\right],
\end{equation*}
recall (\ref{eq:sdSb}), \eqref{eq:vh} and \eqref{eq:fdSd}. Then we define
$\vec Y^0_{\mathcal{S}} \in \Vhpartialzero$, with
$\vec\pi^h\left[(\vec X^0\,.\,\vec\ek_1)\,\vec Y^0_{\mathcal{S}}\right]
\in \yspace^h(2\,\pi\,\alpha_G\,\vec\ek_1)$, via
\begin{equation*} 
\vec Y^0_{\mathcal{S}}(q_j) = \begin{cases}
(\vec Y^0_{\mathcal{S},\star}(q_j)\,.\,\vec\ek_2)\,\vec\ek_2 
& q_j \in \partial_0 I\,,\\
\dfrac{2\,\pi\,\alpha_G\,\vec\ek_1}{\vec X^0(q_j)\,.\,\vec\ek_1} 
& q_j \in \partial_M I\,,\\
\vec Y^0_{\mathcal{S},\star}(q_j) 
& q_j \in \overline I \setminus (\partial_0 I \cup \partial_M I)\,.
\end{cases}
\end{equation*}

Unless otherwise stated we use $\alpha = 1$, 
$\spont=\lambda=\beta=\alpha_G=\varsigma=0$
and employ uniform time steps, $\ttau_m = \ttau$,
$m=0,\ldots,M-1$.

\subsection{Surfaces without boundary}

\subsubsection{Genus 1 surfaces}
In this subsection, we consider genus 1 surfaces without boundary, so that
$\partial I = \emptyset$.

\begin{figure}
\center
\includegraphics[angle=-90,width=0.3\textwidth]{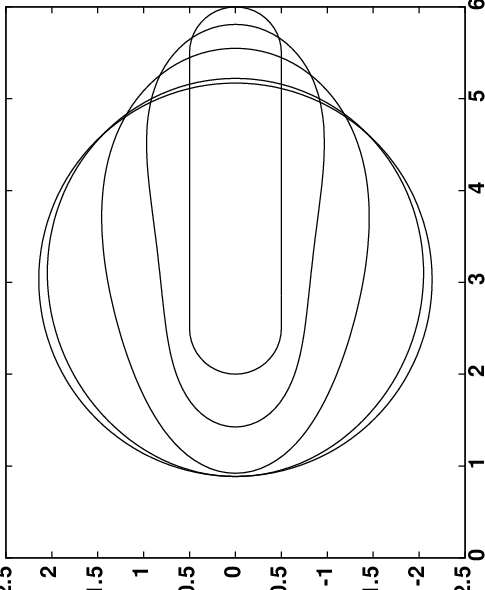} 
\includegraphics[angle=-90,width=0.34\textwidth]{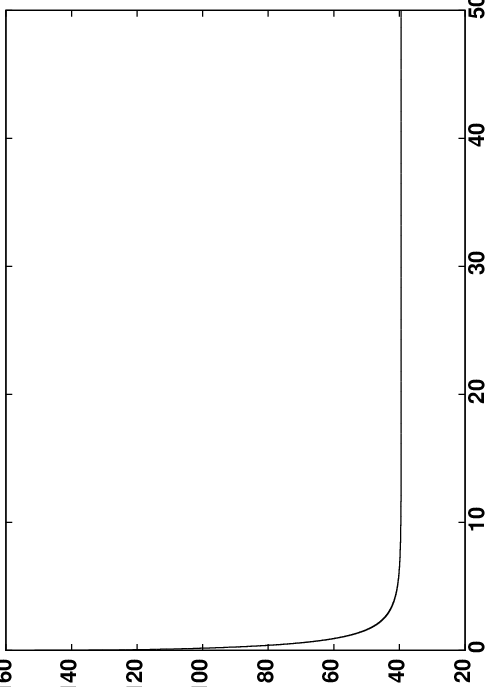}
\includegraphics[angle=-90,width=0.34\textwidth]{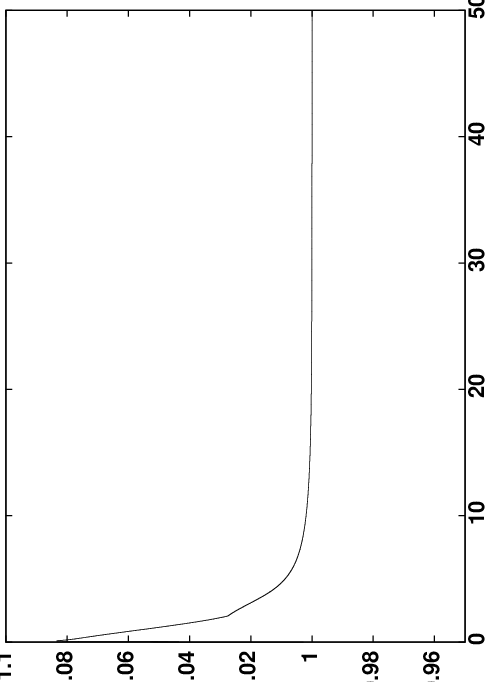}
\caption{$(\BGNpwf^m)^{h}$
Willmore flow for a torus.
Solution at times $t=0,0.5,2, 10, 50$. 
Below a plot of the discrete energy (\ref{eq:Em}) 
and of the ratio (\ref{eq:ratio}).}
\label{fig:pwfcigar41}
\end{figure}

Starting with an elongated cigar-like shape for $\Gamma^0$, we observe
the evolution shown in Figure~\ref{fig:pwfcigar41} for the scheme 
$(\BGNpwf^m)^{h}$, \eqref{eq:fd}.
The discretization parameters are $J=128$ and $\ttau = 10^{-4}$.
The observed final radius is $2.14$, with the centre at $(3.03,0)$. Hence
$R/r = 3.03/2.14 = 1.414 \approx \sqrt{2}$.
Here we recall that the ratio $\sqrt{2}$ characterizes the 
Clifford torus, the known minimizer of the Willmore energy (\ref{eq:E}), with
$\spont=0$ and $\alpha=1$, 
among all genus $1$ surfaces, see \cite{MarquesN14}, with Willmore
energy equal to $4\,\pi^2 = 39.478$.
We note that, as expected, the energy \eqref{eq:Em} is monotonically
decreasing, while the ratio
\begin{equation} \label{eq:ratio}
r^m = \dfrac{\max_{j=1,\ldots, J} |\vec X^m(q_j) - \vec X^m(q_{j-1})|}
{\min_{j=1,\ldots,J} |\vec X^m(q_j) - \vec X^m(q_{j-1})|}
\end{equation}
approaches one as time increases.

\begin{figure}
\center
\includegraphics[angle=-90,width=0.3\textwidth]{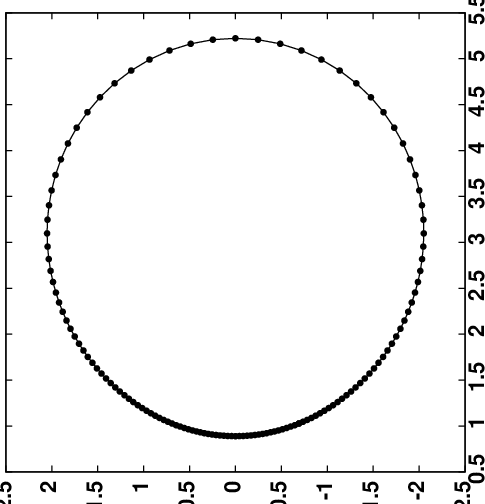}
\includegraphics[angle=-90,width=0.3\textwidth]{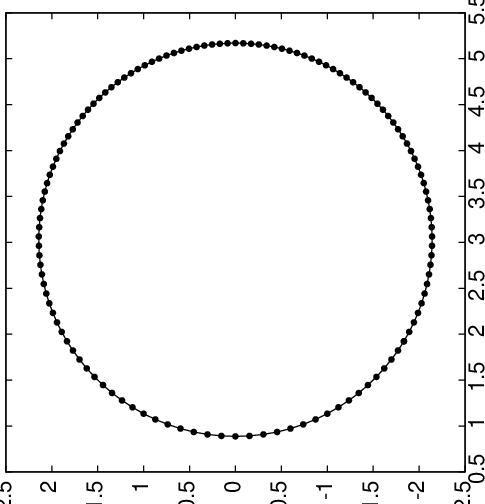} 
\includegraphics[angle=-90,width=0.35\textwidth]{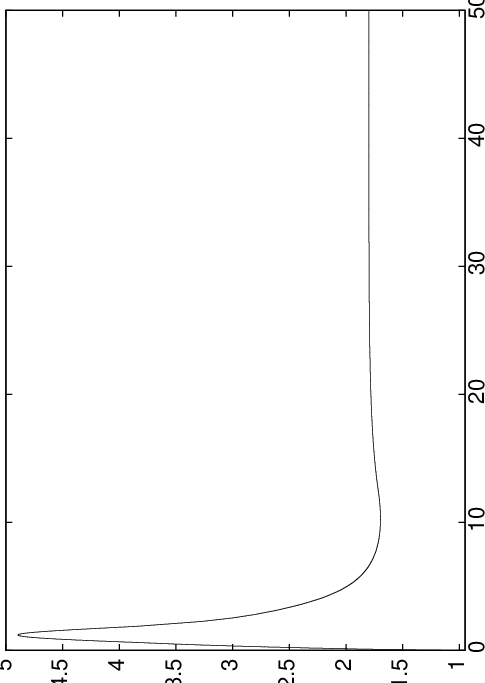}
\caption{$(\BGNpwfwf^m)^{(h)}$
Same evolution as in Figure~\ref{fig:pwfcigar41}.
Left: $\Gamma^m$ at time $t=10$ for $(\BGNpwfwf^m)^{h}$,
middle: $\Gamma^m$ at time $t=50$ for $(\BGNpwfwf^m)$,
right: the ratio \eqref{eq:ratio} over time for $(\BGNpwfwf^m)$.}
\label{fig:pwfcigar41new}
\end{figure}

When we repeat the simulation for the two schemes $(\BGNpwfwf^m)^{(h)}$,
we note markedly different tangential motions. For the scheme with mass lumping
throughout, $(\BGNpwfwf^m)^{h}$, the vertices coalesce on the left side of the
curve and eventually the algorithm breaks down. For the scheme 
$(\BGNpwfwf^m)$, on the other hand, the density of vertices is higher on the
right side of the circular shape, with the ratio \eqref{eq:ratio} smaller than
2. We demonstrate this in Figure~\ref{fig:pwfcigar41new}, where we also show
an evolution of the ratio \eqref{eq:ratio} for the scheme $(\BGNpwfwf^m)$
over time.

Because of the bad behaviour of the scheme $(\BGNpwfwf^m)^{h}$ in practice,
we will discard that scheme from now on. We will mainly concentrate on the
scheme $(\BGNpwf^m)^{h}$, which in practice leads to nearly equidistributed
polygonal curves, and at times compare it to the scheme $(\BGNpwfwf^m)$.

\subsubsection{Genus 0 surfaces}
In this subsection, we consider genus 0 surfaces without boundary, so that
$\partial_0 I = \{0,1\}$. We will parameterize $\Gamma$ clockwise, so that
$\vec\nu$ induces the outer normal 
$\vec{\rm n}_{\mathcal{S}}$ on $\mathcal{S}$,
recall \eqref{eq:tau} and \eqref{eq:nuS}.  

\begin{figure}
\center
\includegraphics[angle=-90,width=0.35\textwidth]{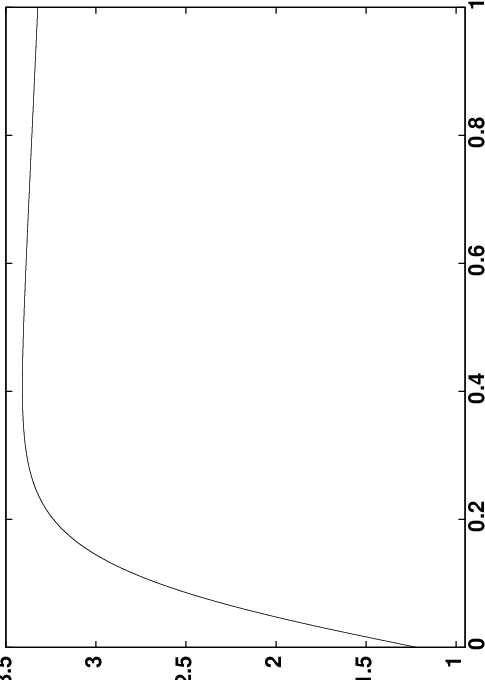} 
\quad
\includegraphics[angle=-90,width=0.15\textwidth]{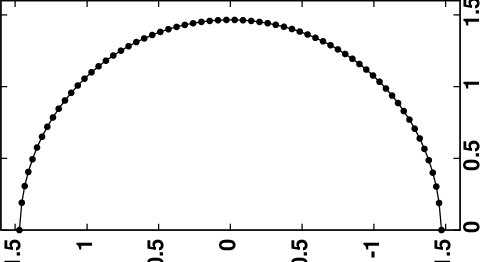}
\caption{$(\BGNpwfwf^m)$
Generalized Willmore flow with $\spont=-1$ for a sphere.
A plot of (\ref{eq:ratio}) over time and a plot of the solution at
time $t=1$ for the run with $J=64$ from Table~\ref{tab:kspwfspont-1}.
}
\label{fig:kspwfspont-1}
\end{figure}

We begin with a convergence experiment. To this end, we note that a sphere of 
radius $R(t)$, where $R(t)$ satisfies
\begin{equation} \label{eq:ODE}
R'(t) = - \tfrac\spont{R(t)}\,(\tfrac2{R(t)} + \spont)\,,
\quad R(0) = R_0 \in \bRplus\,,
\end{equation}
is a solution to (\ref{eq:Willmore_flow}).
The nonlinear ODE (\ref{eq:ODE}), in the case $\spont\not=0$, is solved by 
$R(t) = z(t) - \tfrac2\spont$, where $z(t)$ is such that 
$\tfrac12\,(z^2(t) - z_0^2) - \tfrac4\spont\,(z(t)-z_0) + \tfrac4{\spont^2}\,
\ln \tfrac{z(t)}{z_0} + \spont^2\,t = 0$, with $z_0 = R_0 + \tfrac2\spont$.
We use the solution to (\ref{eq:ODE}), with $\spont = -1$, and 
a sequence of approximations for the unit sphere ($R_0=1$) to compute the 
error 
$\errorXx = \max_{m=1,\ldots,M} \max_{j=0,\ldots,J} \left| 
|\vec X^m(q_j)| - R(t_m) \right|$
over the time interval $[0,1]$ between
the true solution and the discrete solutions for the schemes 
$(\BGNpwf^m)^h$ and $(\BGNpwfwf^m)$.
In particular, we choose $\vec X^0 \in \Vhpartialzero$ with
\begin{equation*} 
\vec X^0(q_j) = \begin{pmatrix} 
\cos[(\tfrac12 - q_j)\,\pi + 0.1\,\cos((\tfrac12 - q_j)\,\pi)] \\
\sin[(\tfrac12 - q_j)\,\pi + 0.1\,\cos((\tfrac12 - q_j)\,\pi)]
\end{pmatrix}, \quad j = 0,\ldots,J\,,
\end{equation*}
recall (\ref{eq:Jequi}), to ensure an initially non-uniform
distribution of vertices.
We also define the error
\begin{equation*} 
\errorXxLL = \max_{m=1,\ldots,M} \left[ \left( \left|
 |\vec X^m| - R(t_m)|\right|^2, |\vec X^m_\rho| \right)^h
\right]^\frac12 .
\end{equation*}
Here we used the time step size $\ttau=0.1\,h^2_{\Gamma^0}$,
where $h_{\Gamma^0}$ is the maximal edge length of $\Gamma^0$.
The computed errors, together with their experimental order of 
convergence (EOC), are reported in Tables~\ref{tab:kpwfspont-1} and
\ref{tab:kspwfspont-1}. 
We remark that repeating these simulations for the scheme $(\BGNpwfwf^m)^h$
failed, because of tangential motion leading to vertices in the left halfplane.
\begin{table}
\center
\caption{Errors for the convergence test for the scheme $(\BGNpwf^m)^h$
with $\spont = -1$. We also show the ratio \eqref{eq:ratio} at time $t=1$.}
\begin{tabular}{rrccccc}
\hline
$J$ & $h_{\Gamma^0}$ & $\errorXx$ & EOC & $\errorXxLL$ & EOC & $r^M$ 
\\ \hline
32  & 1.0792e-01 & 1.3951e-02 & ---  & 7.3622e-03 & ---  & 1.06 \\ 
64  & 5.3988e-02 & 4.1092e-03 & 1.76 & 1.8232e-03 & 2.02 & 1.06 \\
128 & 2.6997e-02 & 1.1867e-03 & 1.79 & 4.5434e-04 & 2.00 & 1.06 \\
256 & 1.3499e-02 & 3.3690e-04 & 1.82 & 1.1347e-04 & 2.00 & 1.06 \\
512 & 6.7495e-03 & 9.4318e-05 & 1.84 & 2.8358e-05 & 2.00 & 1.06 \\
\hline
\end{tabular}
\label{tab:kpwfspont-1}
\end{table}%
\begin{table}
\center
\caption{Errors for the convergence test for the scheme $(\BGNpwfwf^m)$
with $\spont = -1$. We also show the ratio \eqref{eq:ratio} at time $t=1$.}
\begin{tabular}{rrccccc}
\hline
$J$ & $h_{\Gamma^0}$ & $\errorXx$ & EOC & $\errorXxLL$ & EOC & $r^M$ 
\\ \hline
32  & 1.0792e-01 & 1.5595e-02 & ---  & 9.8711e-03 & ---  & 2.79 \\ 
64  & 5.3988e-02 & 5.5381e-03 & 1.49 & 2.7237e-03 & 1.86 & 3.32 \\
128 & 2.6997e-02 & 1.9703e-03 & 1.49 & 7.5278e-04 & 1.86 & 3.96 \\
256 & 1.3499e-02 & 7.0062e-04 & 1.49 & 2.0783e-04 & 1.86 & 4.72 \\
512 & 6.7495e-03 & 2.4888e-04 & 1.49 & 5.7280e-05 & 1.86 & 5.61 \\
\hline
\end{tabular}
\label{tab:kspwfspont-1}
\end{table}%
We remark that for the experiments in Table~\ref{tab:kpwfspont-1}, the
ratio \eqref{eq:ratio}, which at time $t=0$ starts off at about $r^0=1.22$,
always decreases monotonically and approaches the value 1, 
so that the final semicircle is nearly equidistributed.
For the experiments in Table~\ref{tab:kspwfspont-1}, however, this is not the
case, and the distribution of vertices remains very nonuniform. In each case,
the longest elements are the two elements touching the $\vec\ek_2$--axis,
while the shortest elements are found far away from the axis.
As an example, we show the plot
of (\ref{eq:ratio}) over time, as well as the distribution of vertices at time
$t=1$ for the run with $J=64$ in Figure~\ref{fig:kspwfspont-1}, and these
are generic for the behaviour for every value of $J$.
We conjecture that these long elements at the boundary are the reason for the
suboptimal convergence rates seen in Table~\ref{tab:kspwfspont-1}. We note that
in $L^\infty(L^1)$ the convergence order is quadratic for the scheme 
$(\BGNpwfwf^m)$, although we do not display the precise numbers here.

\begin{figure}
\center
\includegraphics[angle=-90,width=0.3\textwidth]{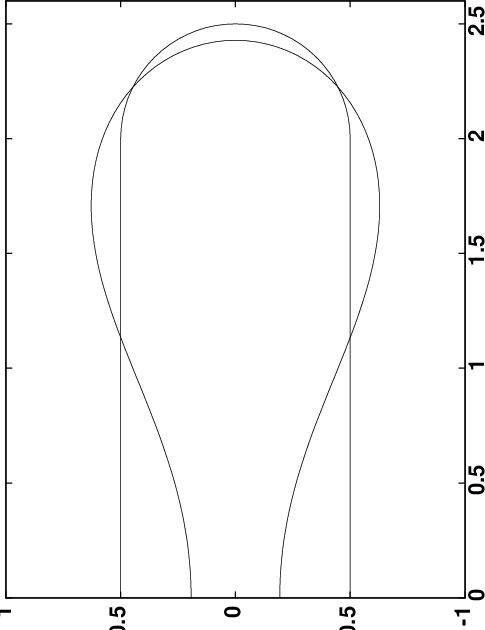} 
\
\includegraphics[angle=-90,width=0.3\textwidth]{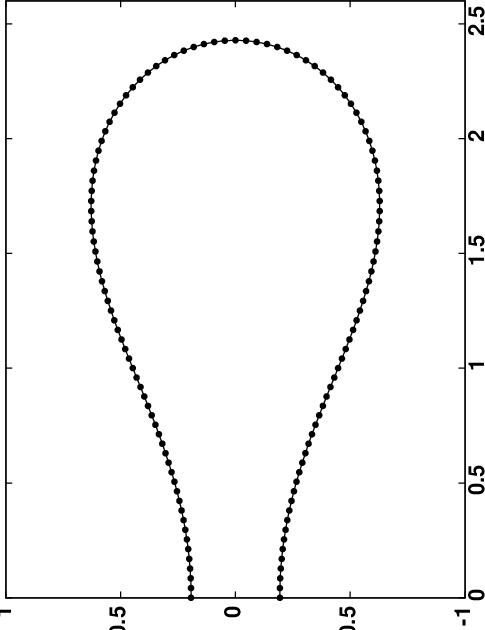}
\qquad
\includegraphics[angle=-90,width=0.3\textwidth]{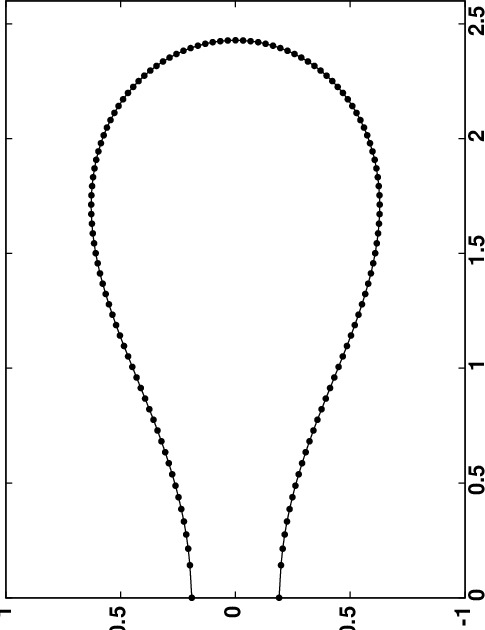}
\includegraphics[angle=-90,width=0.3\textwidth]{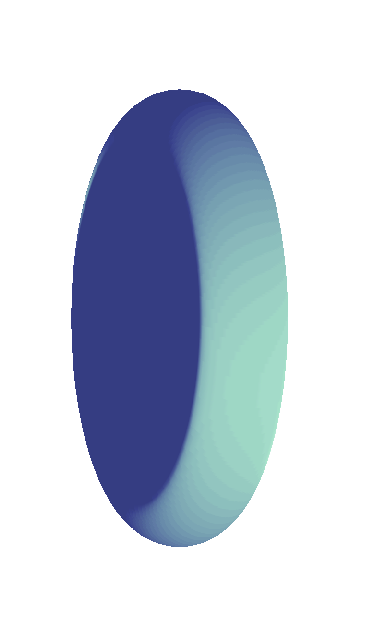} 
\includegraphics[angle=-90,width=0.3\textwidth]{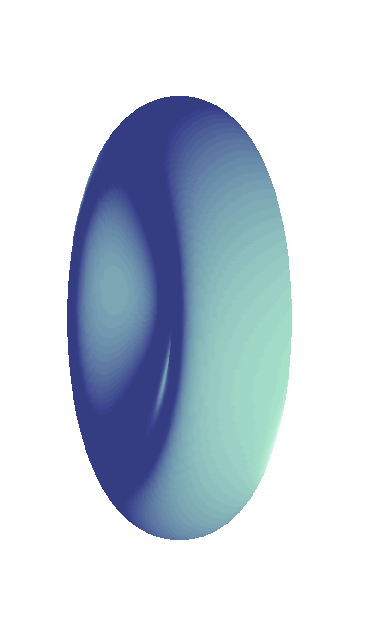} 
\caption{$(\BGNpwf^m)^h$
Helfrich flow for a flat disc of dimension $5\times1\times5$.
Solution at times $t=0,0.2$, and separately at time $t=0.2$.
On the right we show the final distribution of vertices for the
scheme $(\BGNpwfwf^m)$.
Below we visualize the axisymmetric surfaces generated by
$\Gamma^m$ at times $t=0$ and $t=0.2$.}
\label{fig:kpwf_flatcigar_helfrich}
\end{figure}

Next we consider a numerical experiment for Helfrich flow, i.e.\
surface area and volume preserving Willmore flow, for a flat disc.
The dimensions of the initial disc are chosen as 
$5\times1\times5$, so that the
flow evolves towards a sharp that resembles a human red blood cell, 
see Figure~\ref{fig:kpwf_flatcigar_helfrich}.
The discretization parameters for the two schemes are $J=128$, $\ttau=10^{-4}$
and $T=0.2$.

\subsection{Surfaces with boundary}

In the case of clamped boundary conditions, recall \eqref{eq:veczetaS}, 
we define $\vec\zeta(p)$, for $p\in \partial_C I$, via
$\vec\zeta(p) = (\sin \vartheta(p), \cos \vartheta(p))^T$,
where $\vartheta(p) \in \bR$, $p\in \partial_C I$, 
denote the prescribed contact angle with the $x_1$--axis.

We recall from Lemmas~\ref{lem:ex} and \ref{lem:exS} that we can prove
existence of a unique solution to the linear systems arising at each time 
level in the presence of clamped boundary conditions only if 
Assumption~\ref{ass:C} holds. This assumption is violated if $\Gamma^m$ is a
straight line, and so the linear systems are indeed singular. In all other
cases, however, the linear systems for the simulations we present in the
following are nonsingular, and so we can find their unique solutions.

\subsubsection{Surfaces with one connected boundary component}
In this subsection, we consider the situation
$\partial_0 I = \{0\}$ and $\partial I = \{0,1\}$. On recalling 
\eqref{eq:tau} this means that the normal $\vec\nu(0,t)$ will point upwards.

\begin{figure}
\center
\includegraphics[angle=-90,width=0.17\textwidth]{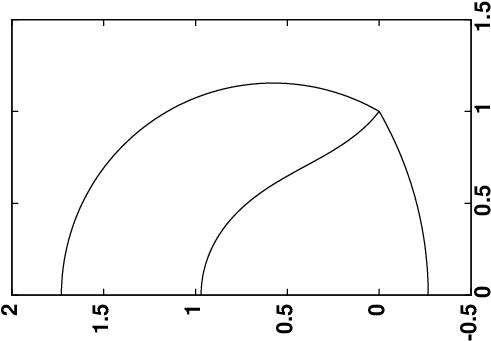} 
\includegraphics[angle=-90,width=0.35\textwidth]{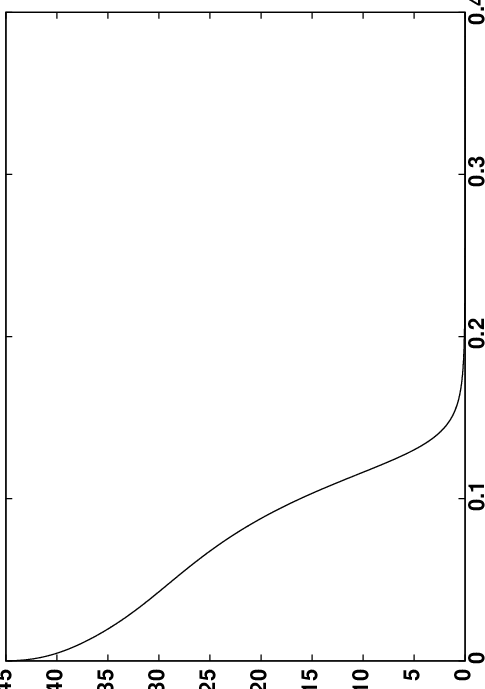} 
\includegraphics[angle=-90,width=0.2\textwidth]{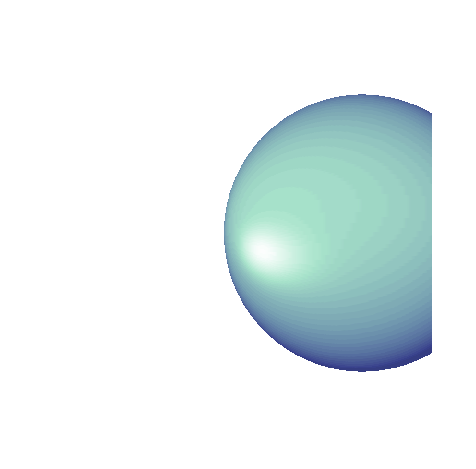} 
\includegraphics[angle=-90,width=0.2\textwidth]{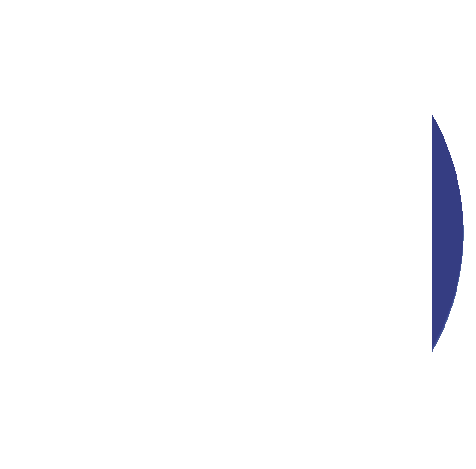} 
\quad
\includegraphics[angle=-90,width=0.14\textwidth]{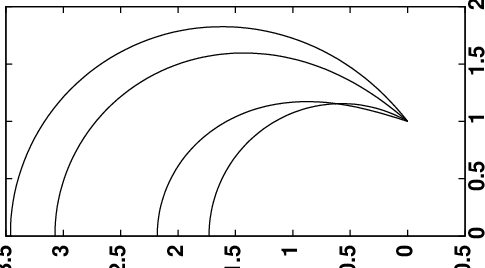} 
\includegraphics[angle=-90,width=0.35\textwidth]{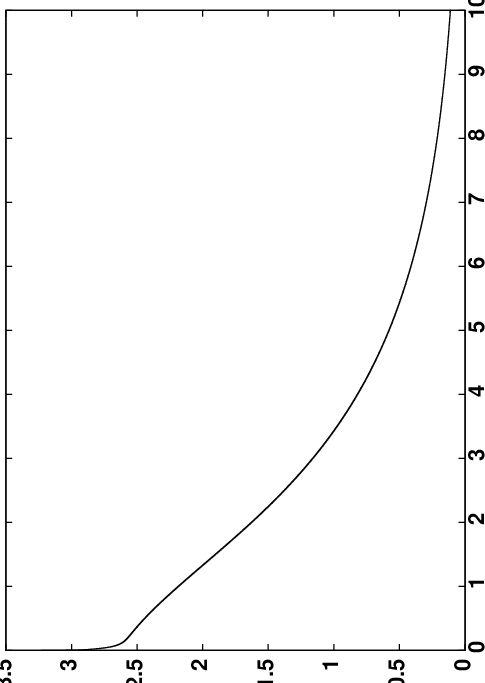} 
\includegraphics[angle=-90,width=0.2\textwidth]{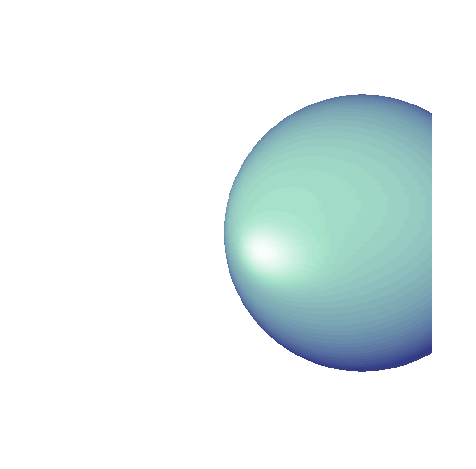} 
\includegraphics[angle=-90,width=0.2\textwidth]{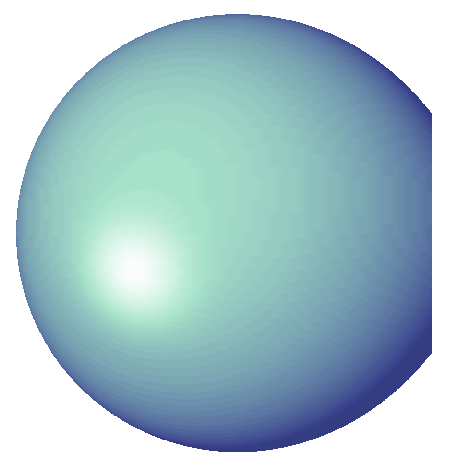} 
\caption{$(\BGNpwf^m)^h$
Generalized Willmore flow with Navier boundary conditions for a sphere-like 
cap with $\spont=1$, top, and $\spont=-1$, bottom.
Solution at times $t=0,0.1,0.4$ (top) and $t=0,1,5,10$ (bottom),
as well as a plot of the discrete energy \eqref{eq:Em} over time. 
We also show the axisymmetric surfaces generated by
$\Gamma^0$ and $\Gamma^M$.
}
\label{fig:kpwf_Navier}
\end{figure}

\begin{figure}
\center
\includegraphics[angle=-90,width=0.17\textwidth]{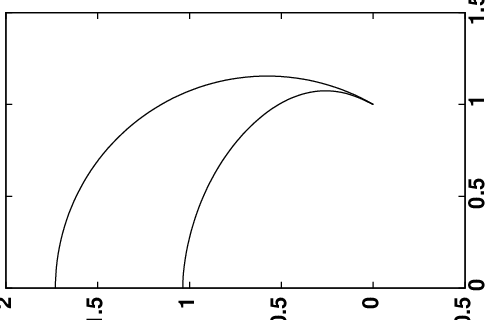} 
\includegraphics[angle=-90,width=0.25\textwidth]{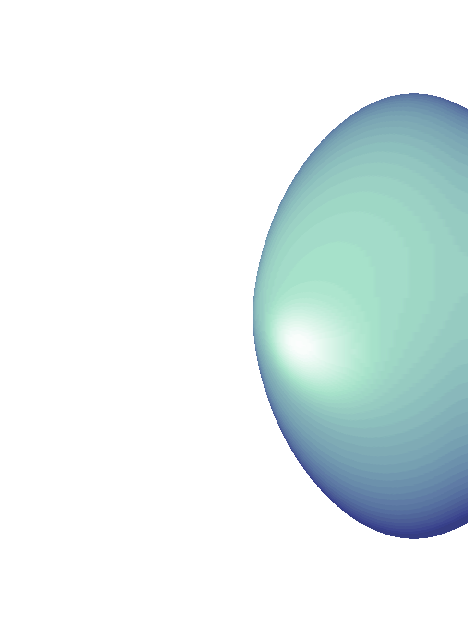} 
\quad
\includegraphics[angle=-90,width=0.15\textwidth]{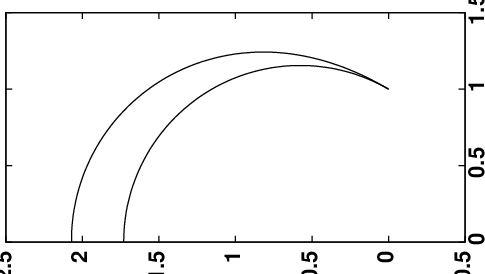} 
\includegraphics[angle=-90,width=0.25\textwidth]{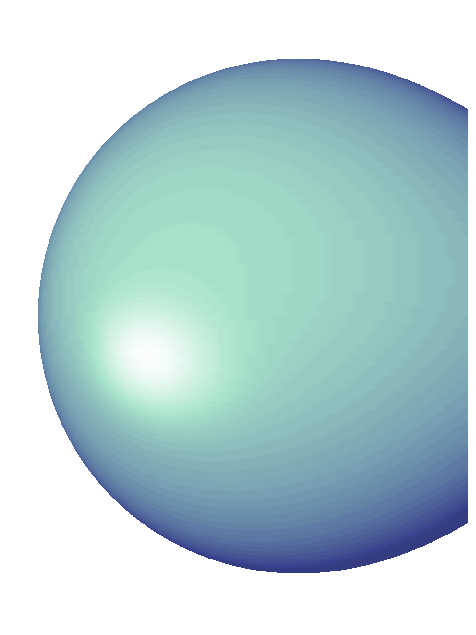} 
\caption{$(\BGNpwf^m)^h$
Generalized Willmore flow with clamped boundary conditions with conormal angle
$\vartheta(1) = \frac76\,\pi = 210^\circ$ for a sphere-like cap
with $\spont=1$, left, and $\spont=-1$, right.
The solutions are shown at times $t=0,0.5$, and we also visualize
the axisymmetric surfaces generated by $\Gamma^M$.
}
\label{fig:kpwf_clamped}
\end{figure}

We show two experiments for Navier boundary conditions in 
Figure~\ref{fig:kpwf_Navier}, where we observe that depending on the sign of
the spontaneous curvature, the sphere-like cap is either convex or concave.
Our numerical results are in agreement with Figures~9 and 10 in \cite{pwfopen}.
As the discretization parameters we choose $J=64$ and $\ttau=10^{-4}$.
The same simulations as in Figures~\ref{fig:kpwf_Navier}, but for clamped
boundary conditions are shown in Figure~\ref{fig:kpwf_clamped}.

\subsubsection{Surfaces with two connected boundary components}
In this subsection, we consider the case 
$\partial_0 I = \emptyset$ and $\partial I = \{0,1\}$. We will always
parameterize $\Gamma(0)$ from left to right, so that the normal $\vec\nu$ for a
straight line points upwards.
When the curve $\Gamma(0)$ is vertical, we parameterize it from top to bottom
so that the normal $\vec\nu$ points to the right.
Throughout this subsection we choose
the discretization parameters $J=64$ and $\ttau=10^{-4}$.

\begin{figure}
\center
\includegraphics[angle=-90,width=0.15\textwidth]{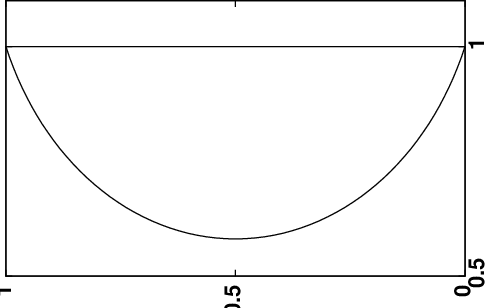} 
\includegraphics[angle=-90,width=0.25\textwidth]{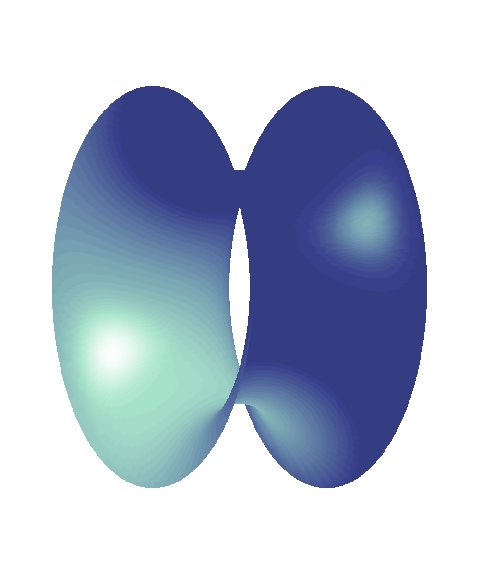} 
\
\includegraphics[angle=-90,width=0.15\textwidth]{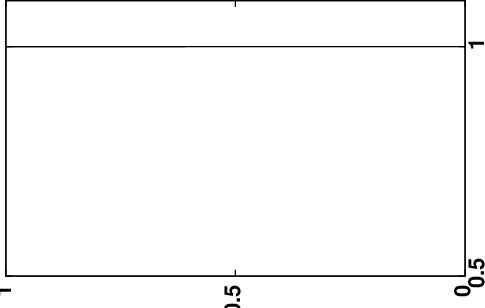} 
\includegraphics[angle=-90,width=0.25\textwidth]{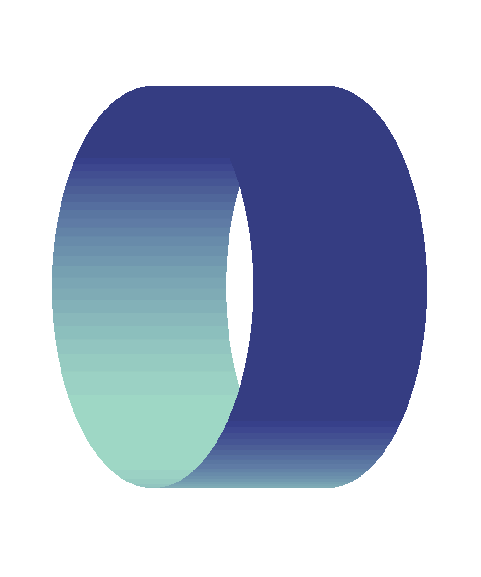} 
\caption{$(\BGNpwf^m)^h$
Generalized Willmore flow with Navier boundary conditions for an open cylinder.
Left with $\spont=1$, right with $\spont=-1$. 
The solutions are shown at times $t=0,0.5$, and we also visualize
the axisymmetric surfaces generated by $\Gamma^M$.
}
\label{fig:kpwf_Navier_cylinder}
\end{figure}

\begin{figure}
\center
\includegraphics[angle=-90,width=0.15\textwidth]{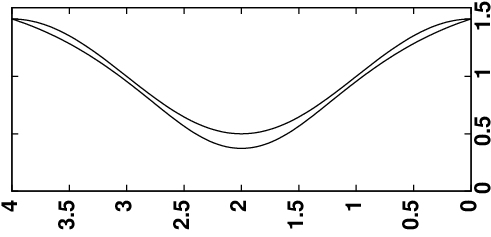} 
\includegraphics[angle=-90,width=0.25\textwidth]{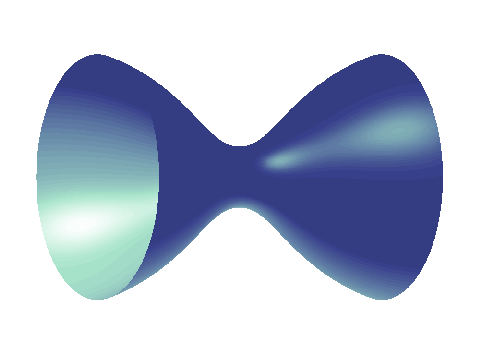} 
\includegraphics[angle=-90,width=0.15\textwidth]{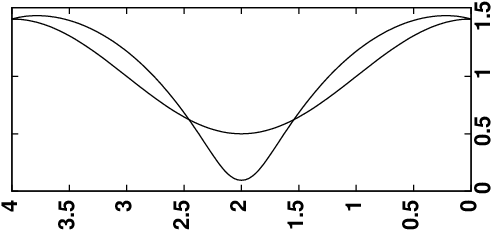} 
\includegraphics[angle=-90,width=0.25\textwidth]{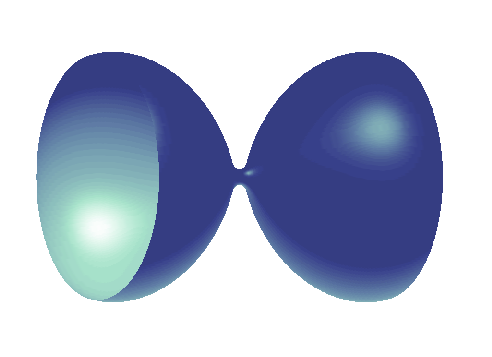} 
\caption{$(\BGNpwf^m)^h$
Generalized Willmore flow with Navier boundary conditions for a dumbbell-like
open cylinder shape. Left with $\spont=-1$, right with $\spont=-2$. 
The solutions are shown at times $t=0,1$, and we also visualize
the axisymmetric surfaces generated by $\Gamma^M$.
}
\label{fig:kpwf_Navier_dumbbell}
\end{figure}

Navier boundary conditions for an open cylinder, when $\Gamma(0)$ is a straight
vertical line, are shown in Figure~\ref{fig:kpwf_Navier_cylinder}. 
Here we note that for $\spont=-1$ the evolution is stationary, which is in
contrast to Figure~\ref{fig:kpwf_Navier_cyl_Gauss}, below.
In order to see the influence of a negative spontaneous curvature on the
evolution, we start two simulations for Navier conditions for a cut cylinder 
with a dumbbell shape, see Figure~\ref{fig:kpwf_Navier_dumbbell}. We observe
that the dumbbell becomes more and more pronounced, the smaller we choose
$\spont$. For $\spont = -2$ the evolution nearly leads to a pinch-off.
Choosing $\spont=-3$ does indeed lead to pinch-off, which we do not show 
in Figure~\ref{fig:kpwf_Navier_dumbbell}.

\begin{figure}
\center
\includegraphics[angle=-90,width=0.15\textwidth]{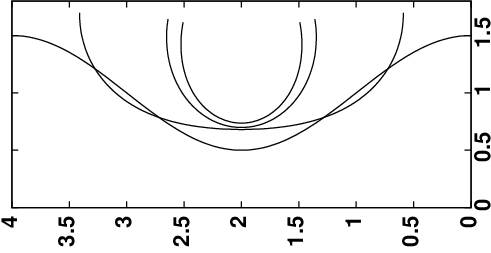} 
\includegraphics[angle=-90,width=0.25\textwidth]{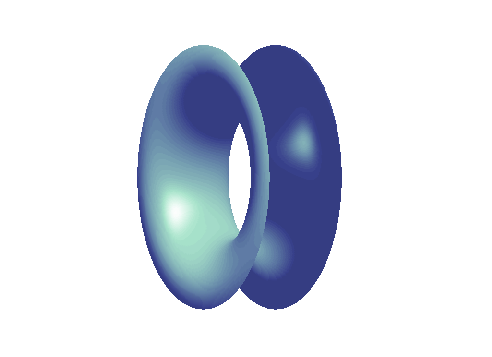} 
\
\includegraphics[angle=-90,width=0.15\textwidth]{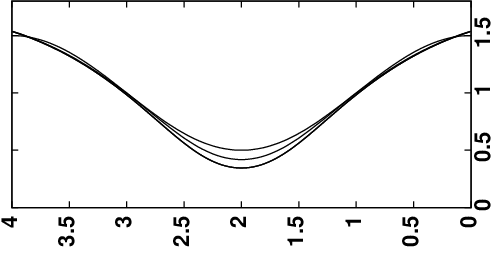} 
\includegraphics[angle=-90,width=0.25\textwidth]{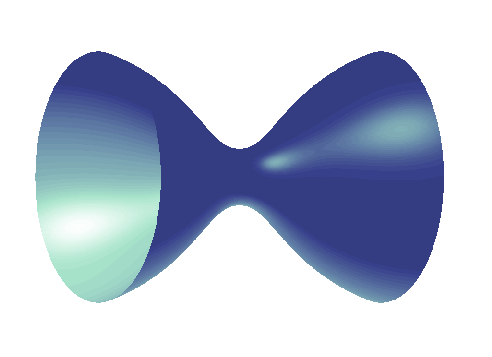} 
\caption{$(\BGNpwf^m)^h$
Generalized Willmore flow with free boundary conditions for a dumbbell-like
open cylinder shape. Left with $\spont=1$, right with $\spont=-1$. 
The solutions are shown at times $t=0,0.1,0.5,1$, and we also visualize
the axisymmetric surfaces generated by $\Gamma^M$.
}
\label{fig:kpwf_free_dumbbell}
\end{figure}

\begin{figure}
\center
\includegraphics[angle=-90,width=0.2\textwidth]{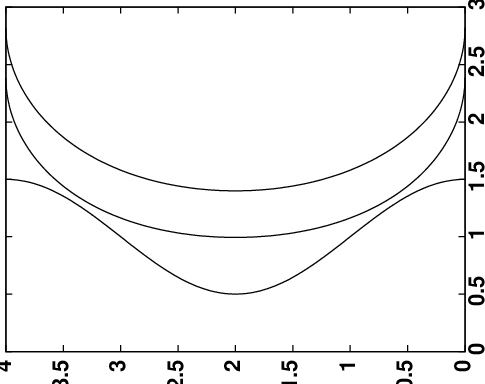} 
\includegraphics[angle=-90,width=0.25\textwidth]{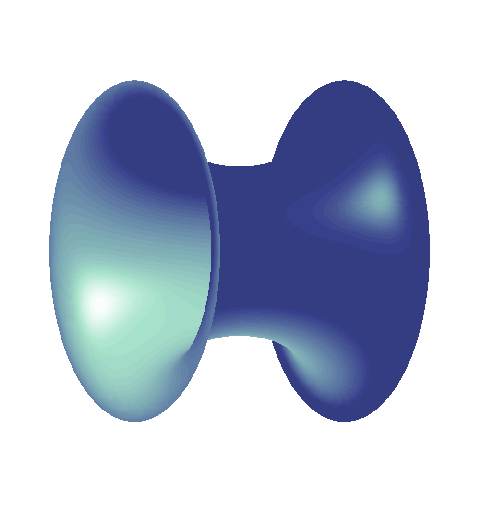} 
\
\includegraphics[angle=-90,width=0.2\textwidth]{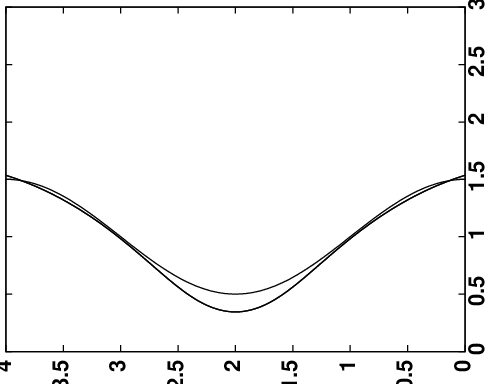} 
\includegraphics[angle=-90,width=0.25\textwidth]{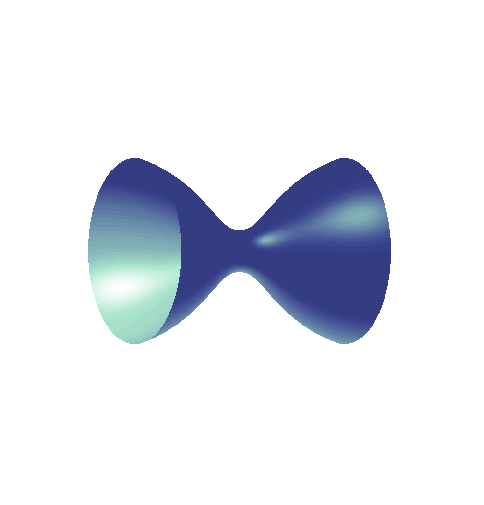} 
\caption{%
$(\BGNpwf^m)^h$
Generalized Willmore flow with semifree boundary conditions for a 
dumbbell-like open cylinder shape. 
Left with $\spont=1$, right with $\spont=-1$. 
The solutions}
are shown at times $t=0,1,10$, and we also visualize
the axisymmetric surfaces generated by $\Gamma^M$.
\label{fig:kpwf_sfree_dumbbell}
\end{figure}

Two simulations for the same initial data, 
but now with $\partial_F I = \partial I$, 
are shown in Figure~\ref{fig:kpwf_free_dumbbell}. Here the evolution for
$\spont=1$ appears to approach the inner half of a torus, while for $\spont=-1$
the limiting shape is similar to an hourglass. 
In both cases the final discrete energy is approximately zero.
In the former case we can prevent the energy converging to zero by requiring
the endpoints to remain on lines parallel to the $x_1$-axis. I.e.\ we use the
same initial data, but now let $\partial_{SF} I = \partial_{2} I = \partial I$.
The evolution is shown in Figure~\ref{fig:kpwf_sfree_dumbbell}, where 
a numerical steady state is reached with discrete energy 
$\widehat E^{m+1} > 25$. For completeness we also display the corresponding
simulation with $\spont = -1$, which is indistinguishable from the
one in Figure~\ref{fig:kpwf_free_dumbbell}.

\begin{figure}
\center
\mbox{
\includegraphics[angle=-90,width=0.25\textwidth]{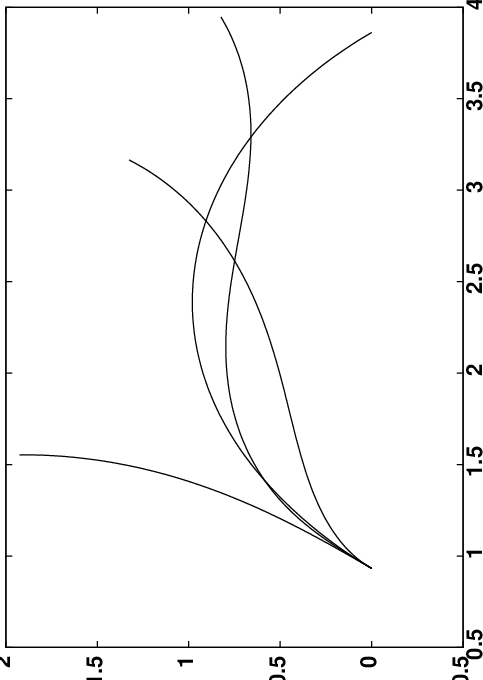} 
\includegraphics[angle=-90,width=0.2\textwidth]{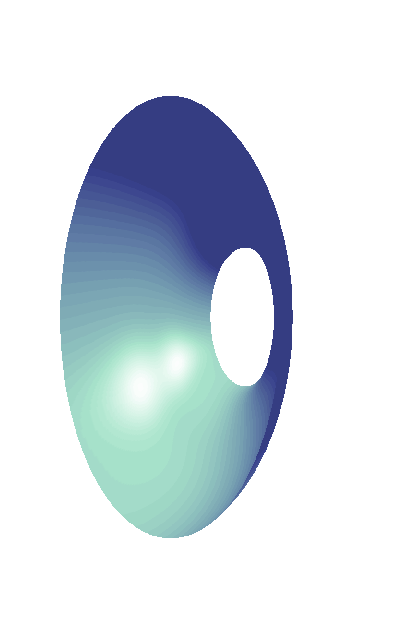} 
\includegraphics[angle=-90,width=0.33\textwidth]{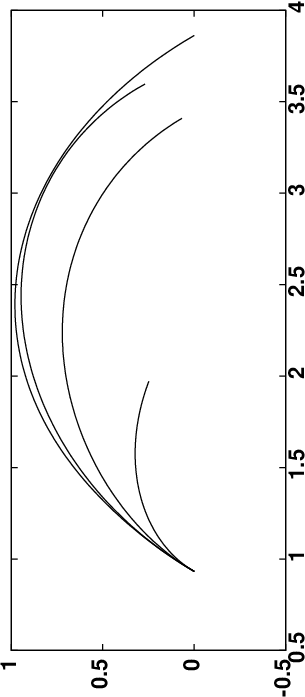} 
\includegraphics[angle=-90,width=0.2\textwidth]{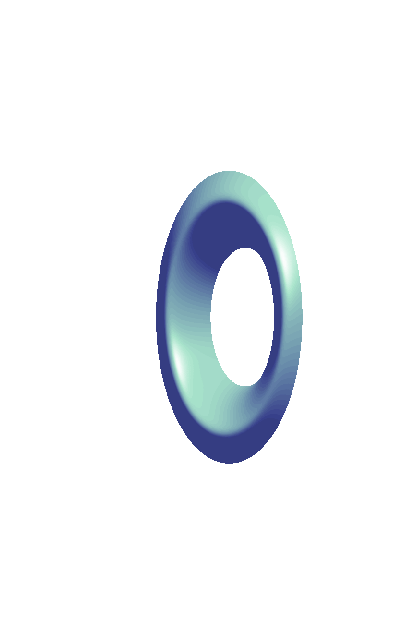} 
}
\caption{$(\BGNpwf^m)^h$
Generalized Willmore flow with mixed clamped and free boundary conditions 
for a torus-like cap. Left with $\spont=1$, right with $\spont=-1$.
The contact angle for the clamped node is chosen as 
$\vartheta(0) = \tfrac76\,\pi=210^\circ$.
The solutions are shown at times $t=0,0.1,1,5$, and we also visualize
the axisymmetric surfaces generated by $\Gamma^M$.
}
\label{fig:kpwf_clampedfree}
\end{figure}

Different simulations for the top part of a torus, with clamped boundary
conditions at the inner ring, and free boundary conditions at the outer ring,
are shown in Figure~\ref{fig:kpwf_clampedfree}.

\begin{figure}
\center
\mbox{
\includegraphics[angle=-90,width=0.25\textwidth]{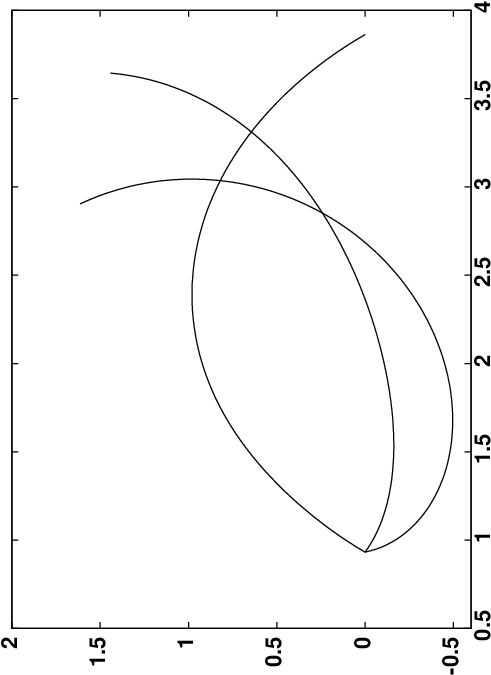} 
\includegraphics[angle=-90,width=0.25\textwidth]{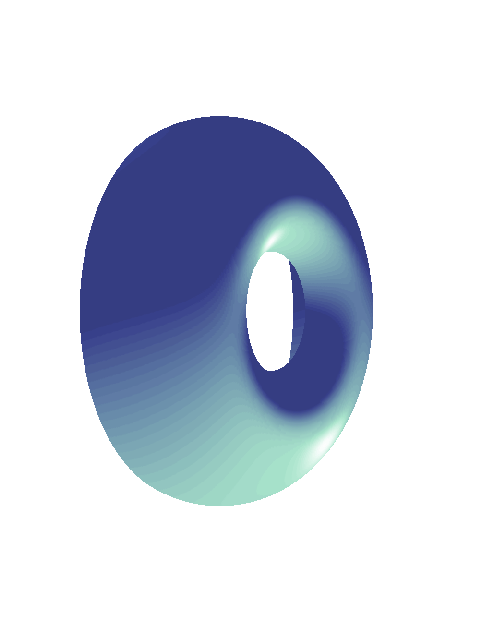} 
\includegraphics[angle=-90,width=0.25\textwidth]{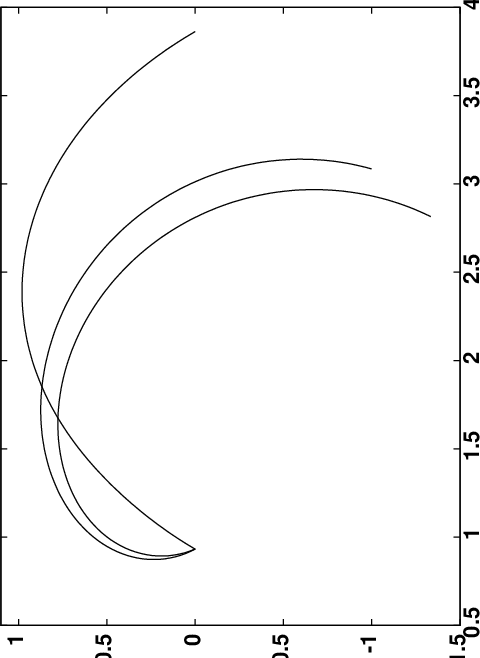} 
\includegraphics[angle=-90,width=0.25\textwidth]{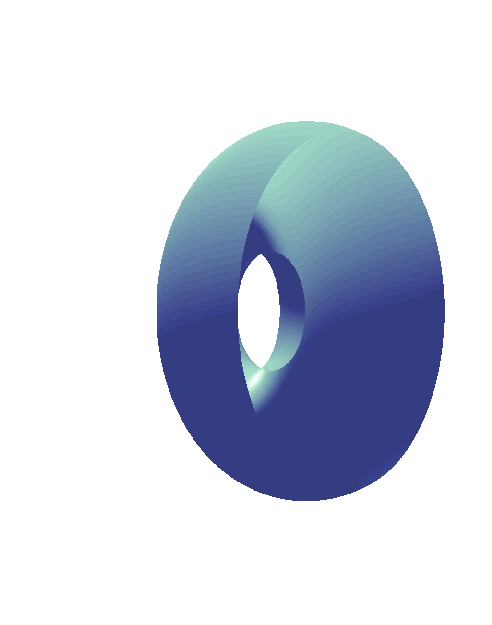} 
}
\caption{$(\BGNpwf^m)^h$
Surface area preserving 
generalized Willmore flow with mixed Navier and free boundary conditions 
for a torus-like cap. Left with $\spont=1$, right with $\spont=-1$.
The solutions are shown at times $t=0,1,5$ (left) and
at times $t=0,5,10$ (right), and we also}
visualize the axisymmetric surfaces generated by $\Gamma^M$.
\label{fig:kpwf_Navierfree_cons}
\end{figure}

\begin{figure}
\center
\includegraphics[angle=-90,width=0.15\textwidth]{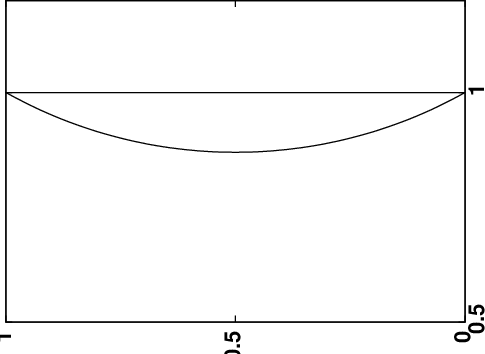} 
\includegraphics[angle=-90,width=0.25\textwidth]{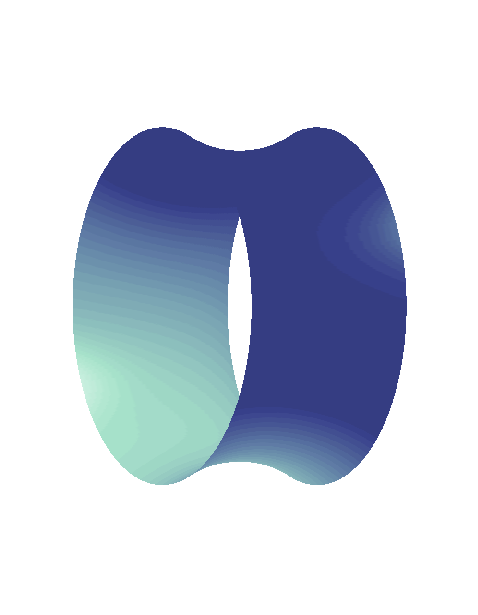} 
\includegraphics[angle=-90,width=0.15\textwidth]{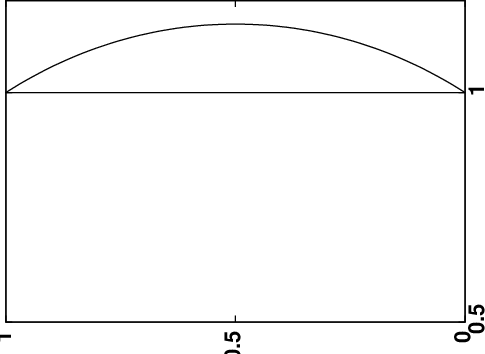} 
\includegraphics[angle=-90,width=0.25\textwidth]{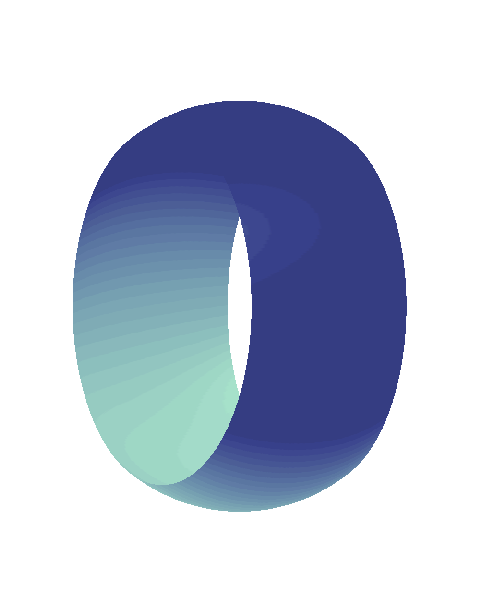} 
\caption{$(\BGNpwf^m)^h$
Generalized Willmore flow, for $\alpha_G = -1$, 
with Navier boundary conditions for an open cylinder.
Left with $\spont=1$, right with $\spont=-1$. 
The solutions are shown at times $t=0,0.5$, and we also visualize
the axisymmetric surfaces generated by $\Gamma^M$.
}
\label{fig:kpwf_Navier_cyl_Gauss}
\end{figure}

Repeating these experiments for Navier boundary conditions at the inner ring of
the torus cap, and enforcing surface area preservation, 
leads to the simulations shown in 
Figure~\ref{fig:kpwf_Navierfree_cons}. 
Here, for $\spont=1$ we observe pinch-off, and so we stop the simulation just
prior to the pinch-off happening.

Finally, we show two evolutions for Gaussian curvature effects. To this end,
we repeat the simulations in Figure~\ref{fig:kpwf_Navier_cylinder}, 
but now with $\alpha_G = -1$. In contrast to the earlier results,
the cylindrical initial data is no longer stationary for $\spont=-1$.
See Figure~\ref{fig:kpwf_Navier_cyl_Gauss} for a visualization of the results.

\begin{appendix}

\renewcommand{\theequation}{A.\arabic{equation}}
\setcounter{equation}{0}
\section{Consistency of weak formulations} \label{sec:A}

\subsection{Formulations based on $\varkappa$} \label{sec:A1}

Here we derive the strong form for (\ref{eq:weak3}), together with possible
boundary conditions for the $L^2$--gradient flow of (\ref{eq:EEax}).
We recall from (\ref{eq:kappay}), 
(\ref{eq:Axkappa}), (\ref{eq:ASxkappa}) and (\ref{eq:meanGaussS}) that
\begin{equation} \label{eq:A1:kappay}
\vec y\,.\,\vec\nu = 
2\,\pi\,\vec x\,.\,\vec\ek_1
\left[\alpha\,(\varkappa_{\mathcal{S}} - \spont) + \beta\,\AADE_{\mathcal{S}}
\right] ,
\quad \text{where}\quad
\varkappa_{\mathcal{S}} = \varkappa - 
\frac{\vec\nu\,.\,\vec\ek_1}{\vec x\,.\,\vec\ek_1}
\quad\text{and }
\AADE_{\mathcal{S}} = 2\,\pi\left(\vec x\,.\,\vec\ek_1\,\varkappa_{\mathcal{S}}, 
|\vec x_\rho| \right) - M_0\,.
\end{equation}

We begin by re-stating (\ref{eq:weak3a}), on noting (\ref{eq:A1:kappay}) and
(\ref{eq:xspace}), as
\begin{align}
& 2\,\pi
\left((\vec x\,.\,\vec\ek_1)\,\vec x_t\,.\,\vec\nu, \vec\chi\,.\,\vec\nu\,
|\vec x_\rho|\right) 
= \left(\vec y_\rho\,.\,\vec\nu, \vec\chi_\rho\,.\,\vec\nu
\, |\vec x_\rho|^{-1} \right) 
- \pi \left( \alpha
\left[\varkappa_{\mathcal{S}}-\spont \right]^2 + 2\,\lambda +
2\,\beta\,\AADE_{\mathcal{S}}\,\varkappa, 
\vec\chi\,.\,\vec\ek_1\,|\vec x_\rho| + (\vec x\,.\,\vec\ek_1)\,
\vec\chi_\rho\,.\,\vec\tau \right) 
\nonumber \\ & \qquad 
- 2\,\pi\,\alpha\left( \varkappa_{\mathcal{S}} - \spont, 
(\varkappa - \varkappa_{\mathcal{S}})
\,\vec\chi\,.\,\vec\ek_1 \,|\vec x_\rho| \right) 
- 2\,\pi\,\beta\,\AADE_{\mathcal{S}} \left(\vec\ek_2, \vec\chi_\rho \right) 
- 2\,\pi\,\alpha\left( \varkappa_{\mathcal{S}} - \spont,
(\vec\tau\,.\,\vec\ek_1)\,\vec\chi_\rho\,.\,\vec\nu \right) 
+ \left( \varkappa\, \vec y^\perp 
,\vec\chi_\rho \right) 
\nonumber \\ & \qquad 
- 2\,\pi\,\varsigma\,\sum_{p \in \partial_{SF} I \cup \partial_F I}
\vec\chi(p)\,.\,\vec\ek_1
= \sum_{i=1}^7 T_i(\vec\chi) 
 \qquad \forall\ \vec\chi \in \xspace \,.
\label{eq:A1:1} 
\end{align}
We now set
\begin{equation} \label{eq:matB}
\mathcal{B} = \partial_0 I \cup \partial_{SF}I \cup \partial_{F}I.
\end{equation}
On recalling (\ref{eq:tau}) and (\ref{eq:xspace}) 
throughout these $T_i(\vec\chi)$ calculations for any $\vec \chi \in \xspace$, 
we have that
\begin{align} \label{eq:A1:T1}
T_1(\vec\chi) = \left( (\vec y_s\,.\,\vec\nu)\,\vec\nu, \vec\chi_\rho \right)
= - \left( [(\vec y_s\,.\,\vec\nu)\,\vec\nu]_s, \vec\chi\,|\vec x_\rho|\right)
+ \sum_{p \in \mathcal{B}} (-1)^{p+1} 
\left[(\vec y_s\,.\,\vec\nu)\,\vec\chi\,.\,\vec\nu\right](p) 
= A_1(\vec\chi) + B_1(\vec\chi) \,.
\end{align}
Moreover,
\begin{align}
\sum_{i=2}^4 T_i(\vec\chi) & =
- \pi \left( \alpha
\left(\varkappa_{\mathcal{S}} -\spont \right) 
\left(2\,\varkappa - \varkappa_{\mathcal{S}} -\spont \right)
+ 2\,\lambda+2\,\beta\,\AADE_{\mathcal{S}}\,\varkappa,
\vec\chi\,.\,\vec\ek_1\,|\vec x_\rho| \right) 
\nonumber \\ & \qquad 
- \pi \left( \vec x\,.\,\vec\ek_1\left( \alpha
\left[\varkappa_{\mathcal{S}} -\spont \right]^2 
+ 2\,\lambda + 2\,\beta\,\AADE_{\mathcal{S}}\,\varkappa \right) \vec\tau 
+ 2\,\beta\,\AADE_{\mathcal{S}}\,\vec\ek_2, \vec\chi_\rho \right) 
= A_2(\vec\chi) + T_8(\vec\chi)\,,
\label{eq:A1:T23} 
\end{align}
where
\begin{align}
T_8(\vec\chi) & = \pi \left( \left[\vec x\,.\,\vec\ek_1\left( \alpha
\left[\varkappa_{\mathcal{S}}-\spont 
\right]^2 + 2\,\lambda + 2\,\beta\,\AADE_{\mathcal{S}}\,\varkappa \right) \vec\tau \right]_s, 
\vec\chi \,|\vec x_\rho| \right) 
\nonumber \\ & \ 
- \pi\, \sum_{p\in\mathcal{B}} (-1)^{p+1}
\left[\vec x\,.\,\vec\ek_1\left( \alpha
\left[\varkappa_{\mathcal{S}} -\spont 
\right]^2 + 2\,\lambda + 2\,\beta\,\AADE_{\mathcal{S}}\,\varkappa \right) 
\vec\chi\,.\,\vec\tau
+ 2\,\beta\,\AADE_{\mathcal{S}}\,\vec\chi\,.\,\vec\ek_2 \right](p)
= A_3(\vec\chi) + B_3(\vec\chi)\,.
\label{eq:A1:T6}
\end{align}
In addition,
\begin{equation} \label{eq:A1:T4}
T_5(\vec\chi) = 2\,\pi\,\alpha \left( \left[
\left(\varkappa_{\mathcal{S}} - \spont
\right) (\vec\tau\,.\,\vec\ek_1)\,\vec\nu \right]_s, \vec\chi\,|\vec x_\rho| 
\right) 
- 2\,\pi\,\alpha\,\sum_{p\in\mathcal{B}} (-1)^{p+1}
\left[ \left(\varkappa_{\mathcal{S}} - \spont \right) 
(\vec\tau\,.\,\vec\ek_1)\,\vec\chi\,.\,\vec\nu \right](p) 
= A_4(\vec\chi) + B_4(\vec\chi)\,.
\end{equation}
Finally,
\begin{equation} \label{eq:A1:T5}
T_6(\vec\chi) 
= - \left( (\varkappa\,\vec y^\perp)_s, \vec\chi\,|\vec x_\rho| \right)
+ \sum_{p\in \mathcal{B}} (-1)^{p+1}\left[\varkappa\,\vec\chi\,.\,\vec y^\perp
\right](p)
= A_5(\vec\chi) + B_5(\vec\chi)\,,
\end{equation}
and $B_2(\vec\chi) = T_7(\vec\chi)$. With the above definitions, we can write
(\ref{eq:A1:1}) as
\begin{equation} \label{eq:A1:1compact}
2\,\pi
\left((\vec x\,.\,\vec\ek_1)\,\vec x_t\,.\,\vec\nu, \vec\chi\,.\,\vec\nu\,
|\vec x_\rho|\right) 
= \sum_{i=1}^7 T_i(\vec\chi) 
= \sum_{i=1}^5 A_i(\vec\chi) + 
\sum_{i=1}^5 B_i(\vec\chi) \qquad \forall\ \vec\chi \in \xspace\,.
\end{equation}

It follows from (\ref{eq:ysnu}), (\ref{eq:kyperp}) 
and 
\begin{equation} \label{eq:nus}
\vec\nu_s = - \varkappa\,\vec\tau\,,
\end{equation}
that
\begin{align*}
( (\vec y_s\,.\,\vec\nu)\,\vec\nu + \varkappa\,\vec y^\perp)_s
& = ( (\vec y\,.\,\vec\nu)_s\,\vec\nu + \varkappa\,(\vec y\,.\,\vec\nu)\,
\vec\tau)_s 
=  (\vec y\,.\,\vec\nu)_{ss}\,\vec\nu 
+ (\vec y\,.\,\vec\nu)_s\,\vec\nu_s + \varkappa\,(\vec y\,.\,\vec\nu)_s\,
\vec\tau + (\vec y\,.\,\vec\nu)\,(\varkappa\,\vec\tau)_s
\nonumber \\ & 
= (\vec y\,.\,\vec\nu)_{ss}\,\vec\nu 
+ (\vec y\,.\,\vec\nu)\,(\varkappa\,\vec\tau)_s\,, 
\end{align*}
and so
\begin{equation}
A_1(\vec\chi) + A_5(\vec\chi) = 
- \left( \left[(\vec y\,.\,\vec\nu)_{ss} + \vec y\,.\,\vec\nu\,\varkappa^2
\right] \vec\nu + \vec y\,.\,\vec\nu\,\varkappa_s \,\vec\tau,
\vec\chi\,|\vec x_\rho| \right) . \label{eq:A1:A15}
\end{equation}

Choosing $\vec\chi = \chi\,\vec\tau$, for $\chi \in H^1_0(I)$, 
in (\ref{eq:A1:1compact}), 
and combining (\ref{eq:A1:T1}), (\ref{eq:A1:T23}), (\ref{eq:A1:T6}),
(\ref{eq:A1:T4}), (\ref{eq:A1:T5}) and (\ref{eq:A1:A15}), 
we obtain for the right hand side of (\ref{eq:A1:1compact}) 
\begin{align}
\sum_{i=1}^7 T_i(\chi\,\vec\tau) & =
\sum_{i=1}^5 A_i(\chi\,\vec\tau) 
= - \left( \vec y\,.\,\vec\nu\,\varkappa_s 
+ \pi \left[ \alpha
\left(\varkappa_{\mathcal{S}} -\spont \right) 
\left(2\,\varkappa - \varkappa_{\mathcal{S}} -\spont \right)
+ 2\,\lambda + 2\,\beta\,\AADE_{\mathcal{S}}\,\varkappa \right]
\vec\tau\,.\,\vec\ek_1 , \chi\,|\vec x_\rho| \right)
\nonumber \\ & \qquad \qquad \qquad \qquad
+ \pi  \left( \left[ \vec x\,.\,\vec\ek_1 \left( \alpha \left[\varkappa_{\mathcal{S}}-\spont 
\right]^2 + 2\,\lambda + 2\,\beta\,\AADE_{\mathcal{S}}\,\varkappa \right) \right]_s
- 2\,\alpha \left(\varkappa_{\mathcal{S}} - \spont
\right) \vec\tau\,.\,\vec\ek_1\,\varkappa, \chi \,|\vec x_\rho| \right) 
\nonumber \\ & 
= - \pi\,\alpha \left( 
\vec x\,.\,\vec\ek_1\left( 2\, 
(\varkappa_{\mathcal S}-\spont) \,\varkappa_s 
- \left[ \left[\varkappa_{\mathcal{S}}-\spont 
\right]^2 \right]_s \right), \chi \,|\vec x_\rho| \right) 
+ 2\,\pi\,\alpha \left( (\varkappa_{\mathcal{S}}-\spont)
(\varkappa_{\mathcal{S}}-2\,\varkappa) 
, \vec\tau\,.\,\vec\ek_1\,\chi \,|\vec x_\rho| \right) 
\nonumber \\ & 
= 2\,\pi\,\alpha \left(  
(\varkappa_{\mathcal{S}}-\spont)\left[ \vec x\,.\,\vec\ek_1\,
(\varkappa_{\mathcal{S}}-\varkappa)_s
+ \vec\tau\,.\,\vec\ek_1\,(\varkappa_{\mathcal{S}}-2\,\varkappa) \right]
, \chi \,|\vec x_\rho| \right), 
\label{eq:A1:tangential}
\end{align}
where we have noted that $\vec\tau_s\,.\,\vec\tau=0$, (\ref{eq:nus}), (\ref{eq:tau})
and (\ref{eq:A1:kappay}).
In addition, it holds, on noting (\ref{eq:A1:kappay}), (\ref{eq:tau}) and (\ref{eq:nus}), that
\[
 \vec x\,.\,\vec\ek_1\,(\varkappa_{\mathcal{S}}-\varkappa)_s
= - \vec x\,.\,\vec\ek_1\left[\frac{\vec\nu\,.\,\vec\ek_1}{\vec x\,.\,\vec\ek_1}
\right]_s = 
\frac{(\vec\nu \,.\,\vec\ek_1)\,\vec\tau \,.\,\vec\ek_1-
(\vec\nu_s \,.\,\vec\ek_1)\,\vec x \,.\,\vec\ek_1}{\vec x \,.\,\vec\ek_1}
=
\vec\tau\,.\,\vec\ek_1\,(2\,\varkappa - \varkappa_{\mathcal{S}})\,,
\]
and so substituting into (\ref{eq:A1:tangential}) yields that
$\sum_{i=1}^7 T_i(\chi\,\vec\tau) = 0$, 
as expected.

Choosing $\vec\chi = \chi\,\vec\nu$, for $\chi \in H^1_0(I)$, 
in (\ref{eq:A1:1compact}), 
and combining (\ref{eq:A1:T1}), (\ref{eq:A1:T23}), (\ref{eq:A1:T6}),
(\ref{eq:A1:T4}), (\ref{eq:A1:T5}) and (\ref{eq:A1:A15}), 
we obtain for the right hand side of (\ref{eq:A1:1compact}) 
\begin{align}
& \sum_{i=1}^7 T_i(\chi\,\vec\nu) =
\sum_{i=1}^5 A_i(\chi\,\vec\nu)
= - \left( (\vec y\,.\,\vec\nu)_{ss}
+ \vec y\,.\,\vec\nu\,\varkappa^2, \chi\,|\vec x_\rho| \right)
\nonumber \\ & \qquad
- \pi \left( \left[ \alpha
\left(\varkappa_{\mathcal{S}} -\spont \right) 
\left(2\,\varkappa - \varkappa_{\mathcal{S}} -\spont \right)
+ 2\,\lambda + 2\,\beta\,\AADE_{\mathcal{S}}\,\varkappa\right]\vec\nu\,.\,\vec\ek_1
, \chi\,|\vec x_\rho| \right)
\nonumber \\ & \qquad
+ \pi \left( \vec x\,.\,\vec\ek_1\left( \alpha
\left[\varkappa_{\mathcal{S}}-\spont 
\right]^2 + 2\,\lambda  + 2\,\beta\,\AADE_{\mathcal{S}}\,\varkappa\right) \varkappa , \chi 
\,|\vec x_\rho| \right) 
+ 2\,\pi\,\alpha \left( \left[
\left(\varkappa_{\mathcal{S}} - \spont
\right) \vec\tau\,.\,\vec\ek_1\right]_s, \chi\,|\vec x_\rho| \right) 
\nonumber \\ & \
= - 2\,\pi\left( \left(\vec x\,.\,\vec\ek_1\left[
\alpha\,(\varkappa_{\mathcal{S}} - \spont) +\beta\,\AADE_{\mathcal{S}} \right] \right)_{ss}
+ \vec x\,.\,\vec\ek_1\left[ \alpha\,(\varkappa_{\mathcal{S}} -
\spont) + \beta\,\AADE_{\mathcal{S}}\right] \varkappa^2, \chi\,|\vec x_\rho| \right)
\nonumber \\ & \qquad
+ \pi \left(
 \left[ \alpha \left[\varkappa_{\mathcal{S}}^2 -\spont^2 \right] 
-  2\,\lambda - 2\,\beta\,\AADE_{\mathcal{S}}\,\varkappa\right]\vec\nu\,.\,\vec\ek_1
, \chi\,|\vec x_\rho| \right)
\nonumber \\ & \qquad
+ \pi \left( \vec x\,.\,\vec\ek_1\left[ \alpha
\left[\varkappa_{\mathcal{S}}-\spont
\right]^2 + 2\,\lambda + 2\,\beta\,\AADE_{\mathcal{S}}\,\varkappa \right] \varkappa , \chi 
\,|\vec x_\rho| \right) 
+ 2\,\pi\,\alpha \left(
\left(\varkappa_{\mathcal{S}} \right)_s \vec\tau\,.\,\vec\ek_1
, \chi\,|\vec x_\rho| \right) 
\nonumber \\ & \
= - 2\,\pi\left( \alpha\,\vec x\,.\,\vec\ek_1\,(\varkappa_{\mathcal{S}})_{ss}
+ \alpha\,\vec\tau\,.\,\vec\ek_1\,(\varkappa_{\mathcal{S}})_s
+ \varkappa\,\vec\nu\,.\,\vec\ek_1\left[\alpha\,(\varkappa_{\mathcal{S}} -
\spont) + \beta\,\AADE_{\mathcal{S}}\right], \chi\,|\vec x_\rho| \right)
\nonumber \\ & \qquad
- 2\,\pi\left( \vec
x\,.\,\vec\ek_1\left[ \alpha\, (\varkappa_{\mathcal{S}} - \spont) +
\beta\,\AADE_{\mathcal{S}}\right] \varkappa^2, \chi\,|\vec x_\rho| \right)
+ \pi \left( \vec x\,.\,\vec\ek_1
 \left[ \alpha \left[\varkappa_{\mathcal{S}}^2 -\spont^2 \right] 
- 2\,\lambda - 2\,\beta\,\AADE_{\mathcal{S}}\,\varkappa \right] 
(\varkappa -\varkappa_{\mathcal{S}}) , \chi\,|\vec x_\rho| \right)
\nonumber \\ & \qquad
+ \pi \left( \vec x\,.\,\vec\ek_1\left[ \alpha
\left[\varkappa_{\mathcal{S}}-\spont 
\right]^2 + 2\,\lambda + 2\,\beta\,\AADE_{\mathcal{S}}\,\varkappa \right] \varkappa , \chi 
\,|\vec x_\rho| \right) 
\nonumber \\ & \ 
= - 2\,\pi\,\alpha\left( (\vec
x\,.\,\vec\ek_1\,(\varkappa_{\mathcal{S}})_{s})_s, \chi\,|\vec x_\rho| \right)
+ 2\,\pi\,\lambda \left( 
\vec x\,.\,\vec\ek_1\,\varkappa_{\mathcal{S}}, \chi\,|\vec x_\rho| \right)
+ \pi\, \alpha\left( \vec x\,.\,\vec\ek_1
 \left[  \varkappa_{\mathcal{S}}^2 -\spont^2 \right] 
(\varkappa -\varkappa_{\mathcal{S}})
, \chi\,|\vec x_\rho| \right)
\nonumber \\ &  \qquad
+ 4\,\pi\,\beta\,\AADE_{\mathcal{S}} \left( \vec x\,.\,\vec\ek_1,
\varkappa\,(\varkappa_{\mathcal{S}} - \varkappa), \chi\,|\vec x_\rho| \right)
+ \pi\, \alpha\left( \vec x\,.\,\vec\ek_1
 \left(  \varkappa_{\mathcal{S}} -\spont \right)
\left[ 2\,(\varkappa_{\mathcal{S}}-\varkappa)\,\varkappa
- 2\,\varkappa^2 + (\varkappa_{\mathcal{S}}-\spont)\,\varkappa \right]
, \chi\,|\vec x_\rho| \right)
\qquad \nonumber \\ & \
= - 2\,\pi\,\alpha\left( (\vec
x\,.\,\vec\ek_1\,(\varkappa_{\mathcal{S}})_{s})_s, \chi\,|\vec x_\rho| \right)
+ 2\,\pi\,\lambda \left( 
\vec x\,.\,\vec\ek_1\,\varkappa_{\mathcal{S}}, \chi\,|\vec x_\rho| \right)
- \pi\, \alpha\left( \vec x\,.\,\vec\ek_1
 \left[  \varkappa_{\mathcal{S}}^2 -\spont^2 \right]
\varkappa_{\mathcal{S}}
, \chi\,|\vec x_\rho| \right)
\nonumber \\ &  \qquad
+ 4\,\pi\,\beta\,\AADE_{\mathcal{S}} \left( \vec x\,.\,\vec\ek_1,
\varkappa\,(\varkappa_{\mathcal{S}} - \varkappa), \chi\,|\vec x_\rho| \right)
+ 4\,\pi\, \alpha\left( \vec x\,.\,\vec\ek_1
 \left(  \varkappa_{\mathcal{S}} -\spont \right) 
\left(\varkappa_{\mathcal{S}}-\varkappa\right)\varkappa
, \chi\,|\vec x_\rho| \right),
\label{eq:A1:normal1}
\end{align}
where we have noted that $\vec\nu_s\,.\,\vec\nu=0$, (\ref{eq:varkappa}),
(\ref{eq:A1:kappay}) and (\ref{eq:tau}).
Clearly, it follows from 
(\ref{eq:A1:1}), 
(\ref{eq:A1:normal1}) and (\ref{eq:meanGaussS})
that (\ref{eq:xtbgnlambda}) holds.

It remains to show that the weak formulation \eqref{eq:weak3} 
indeed enforces the desired boundary conditions.
Here, we recall that apart from the axisymmetric boundary conditions 
\eqref{eq:part0I} on $\partial_0 I$, we consider the following
four different types of boundary conditions on 
$\partial I \setminus \partial_0 I$.
\renewcommand{\theenumi}{(\roman{enumi})}
\begin{enumerate}
\item $\partial\mathcal{S}(t)$ is free, i.e.\ $\partial_F \mathcal{S}(t)$,
see (\ref{eq:free}) for $\partial\mathcal{S}(t)$ and
(\ref{eq:axifree}) in the axisymmetric case.
\item $\partial\mathcal{S}(t) \subset \partial\Domain$ is semifree, 
i.e.\ $\partial_{SF} \mathcal{S}(t)$, 
see (\ref{eq:sfree}) for $\partial\mathcal{S}(t)$ and
(\ref{eq:axisfree}) in the axisymmetric case. 
\item $\partial\mathcal{S}(t)$ 
clamped, 
i.e.\ $\partial_{C} \mathcal{S}(t)$, 
see (\ref{eq:clamped}) for $\partial\mathcal{S}(t)$ and
(\ref{eq:CL}) in the axisymmetric case.
\item $\partial\mathcal{S}(t)$ 
having Navier conditions,
i.e.\ $\partial_{N} \mathcal{S}(t)$, 
see (\ref{eq:Navier}) for $\partial\mathcal{S}(t)$ and
(\ref{eq:Nav}) in the axisymmetric case.
\end{enumerate}

Firstly, we recall from (\ref{eq:weak3}) that 
$\vec x(\cdot,t) \in \Vpartialzero$ with 
$\vec x_t(\cdot,t) \in \xspace$, for $t \in (0,T]$,
imposes strongly that
(\ref{eq:axisfreea}), (\ref{eq:fixedC}), (\ref{eq:fixedN}) and 
(\ref{eq:fixed0}) hold.
Furthermore, we note that (\ref{eq:weak3c}) weakly imposes
(\ref{eq:bcbc}) and (\ref{eq:clampedI}), 
recall the paragraph below (\ref{eq:varkappaweak}). This means that we still
need to show \eqref{eq:sdbca}, \eqref{eq:axifree}, \eqref{eq:axisfreeb},
\eqref{eq:axisfreec} and \eqref{eq:NavierI}.    

We begin with \eqref{eq:sdbca}. To this end, we choose test functions $\vec\chi
= \chi\,\vec\nu$ in \eqref{eq:A1:1}, where $\chi \in H^1(I)$ is zero away
from $\partial_0 I$, to obtain
\begin{align}
& 2\,\pi
\left((\vec x\,.\,\vec\ek_1)\,\vec x_t\,.\,\vec\nu, \chi\,
|\vec x_\rho|\right) 
= \left(\vec y_\rho\,.\,\vec\nu, \chi_\rho
\, |\vec x_\rho|^{-1} \right) 
 - \pi \left( \alpha
\left[\varkappa_{\mathcal{S}} -\spont 
\right]^2 + 2\,\lambda+ 2\,\beta\,\AADE_{\mathcal{S}}\,\varkappa, 
\chi\,\vec\nu\,.\,\vec\ek_1\,|\vec x_\rho| + (\vec x\,.\,\vec\ek_1)\,
\vec\tau\,.\,\vec\nu_\rho\,\chi\,\right) 
\nonumber \\ & 
- 2\,\pi\,\alpha
\left(\varkappa_{\mathcal{S}} - \spont, (\varkappa - \varkappa_{\mathcal{S}})
\,\chi\,\vec\nu\,.\,\vec\ek_1
\,|\vec x_\rho| \right) 
- 2\,\pi\,\alpha \left(
\varkappa_{\mathcal{S}} - \spont,
(\vec\tau\,.\,\vec\ek_1)\,\chi_\rho \right) 
+ \left( \varkappa\, \vec y^\perp - 2\,\pi\,\beta\,\AADE_{\mathcal{S}}
\,\vec\ek_2, \chi_\rho\,\vec\nu + \chi\,\vec\nu_\rho \right) .
\label{eq:A:weak3a}
\end{align}
We are interested in the boundary condition that is weakly enforced by 
\eqref{eq:A:weak3a} on $\partial_0 I$, and so the only relevant terms
are the ones involving $\chi_\rho$ on the right hand side of 
\eqref{eq:A:weak3a}. They simplify, on noting 
from (\ref{eq:tau}) that $\vec\nu\,.\,\vec\ek_2=\vec\tau\,.\,\vec\ek_1$, 
and on recalling \eqref{eq:kyperp} and \eqref{eq:A1:kappay}, to
\begin{align}
& \left(\vec y_\rho\,.\,\vec\nu, \chi_\rho \, |\vec x_\rho|^{-1} \right)
- 2\,\pi\,\alpha \left(
\varkappa_{\mathcal{S}} - \spont,
(\vec\tau\,.\,\vec\ek_1)\,\chi_\rho \right) 
+ \left( \varkappa\, \vec y^\perp 
- 2\,\pi\,\beta\,\AADE_{\mathcal{S}}\,\vec\ek_2,
\chi_\rho\,\vec\nu \right) \nonumber \\ & \qquad
 = \left(\vec y_s\,.\,\vec\nu
- 2\,\pi\,(\alpha\,(\varkappa_{\mathcal{S}} - \spont )
+ \beta\,\AADE_{\mathcal{S}})\,(\vec\tau\,.\,\vec\ek_1)
+ \varkappa\, \vec y^\perp\,.\,\vec\nu ,\chi_\rho \right)
\nonumber \\ & \qquad
 = \left( (\vec y\,.\,\vec\nu)_s
- 2\,\pi\,(\alpha\,(\varkappa_{\mathcal{S}} - \spont )
+ \beta\,\AADE_{\mathcal{S}})\,(\vec\tau\,.\,\vec\ek_1)
 ,\chi_\rho \right)
 =2\,\pi\,\alpha \left( \vec x\,.\,\vec\ek_1\,
(\varkappa_{\mathcal{S}})_\rho ,\chi_\rho\,|\vec x_\rho|^{-1} \right) .
\label{eq:A:weak3a2}
\end{align}
We can now argue as in \cite[Appendix~A.1]{axisd} to show 
that despite the degenerate weight in the last term in \eqref{eq:A:weak3a2},
the identity \eqref{eq:A:weak3a} gives rise to the boundary condition
$(\varkappa_{\mathcal{S}})_\rho = 0$ on $\partial_0 I$. In particular, 
we note that all the integrands in all the remaining terms converge to
zero for $\rho$ approaching $\partial_0 I$, on recalling \eqref{eq:bcnu}. 
This proves that \eqref{eq:weak3a} weakly enforces \eqref{eq:sdbca}.

Next we recall from (\ref{eq:partialM}) and (\ref{eq:my}) that
\begin{equation*} 
\vec y = 
2\,\pi\,\alpha_G\,\vec\ek_1 \qquad\text{ on } \quad \partial_M I
= \partial_N I \cup \partial_{SF} I \cup \partial_F I\,.
\end{equation*}
Combining with (\ref{eq:A1:kappay}) yields that
\begin{equation} \label{eq:A1:mye1}
\alpha\,(\varkappa_{\mathcal{S}} - \spont) + \beta\,\AADE_{\mathcal{S}}
- \alpha_G\,\frac{\vec\nu\,.\,\vec\ek_1}{\vec x\,.\,\vec\ek_1} = 0
\qquad\text{ on } \quad \partial_M I\,.
\end{equation}
Hence, we have that (\ref{eq:axifreec}), (\ref{eq:axisfreec}) and  
(\ref{eq:NavierI}) are imposed strongly.
Overall, we still need to show \eqref{eq:axifreea}, \eqref{eq:axifreeb} and
\eqref{eq:axisfreeb}. 
To this end, we derive conditions that make the second sum in
(\ref{eq:A1:1compact}) vanish for all test functions $\vec\chi \in \xspace$.
It follows from (\ref{eq:A1:T1}) and (\ref{eq:A1:T5}), 
on recalling (\ref{eq:kyperp}) and (\ref{eq:A1:kappay}), that
\begin{align*}
B_1(\vec\chi) + B_5(\vec\chi) & 
= \sum_{p \in \mathcal{B}} (-1)^{p+1} 
\left[\left((\vec y_s\,.\,\vec\nu)\,\vec\nu
+ \varkappa\,\vec y^\perp\right) .\,\vec\chi\right](p)
= \sum_{p \in \mathcal{B}} (-1)^{p+1} 
\left[\left((\vec y\,.\,\vec\nu)_s\,\vec\nu
+ \varkappa\,(\vec y\,.\,\vec\nu)\,\vec\tau\right) .\,\vec\chi\right](p) 
\nonumber \\ &
= 2\,\pi\,\sum_{p \in \mathcal{B}} 
(-1)^{p+1} \left[(\vec x\,.\,\vec\ek_1
\left[\alpha\,(\varkappa_{\mathcal{S}} - \spont) + \beta\,\AADE_{\mathcal{S}}
\right])_s\,\vec\chi\,.\,\vec\nu\right](p) 
\nonumber \\ & \quad
+ 2\,\pi\,\sum_{p \in \partial_{SF}I \cup \partial_{F}I} (-1)^{p+1} 
\left[ \varkappa\,(\vec x\,.\,\vec\ek_1
\left[\alpha\,(\varkappa_{\mathcal{S}} - \spont) + \beta\,\AADE_{\mathcal{S}}
\right])\,\vec\chi\,.\,\vec\tau\right](p) \,, 
\end{align*}
where in the last term we have noted that \eqref{eq:bcbc} implies that
$\vec\nu\,.\,\vec\ek_1 = 0$ on $\partial_0 I$, recall also \eqref{eq:bcnu}. 
Looking first at the normal boundary contributions, we compute
for a $\chi \in H^1(I)$, 
on noting 
from (\ref{eq:tau}) that $\vec\nu\,.\,\vec\ek_2=\vec\tau\,.\,\vec\ek_1$, that
\begin{align}
& \sum_{i=1}^5 B_i(\chi\,\vec\nu) 
= 2\,\pi\,\sum_{p \in \mathcal{B}} 
(-1)^{p+1} \left[(\vec x\,.\,\vec\ek_1
\left[\alpha\,(\varkappa_{\mathcal{S}} - \spont) + \beta\,\AADE_{\mathcal{S}}
\right])_s\,\chi\right](p) \nonumber \\ & \quad
- 2\,\pi\,\varsigma\,\sum_{p \in \partial_{SF} I \cup \partial_F I}
\left[\vec\nu\,.\,\vec\ek_1\,\chi\right](p)
- 2\,\pi\,\beta\,\AADE_{\mathcal{S}}\, 
\sum_{p\in \mathcal{B}} (-1)^{p+1}
\left[\vec\nu\,.\,\vec\ek_2\,\chi \right](p) 
- 2\,\pi\,\alpha\,
\sum_{p\in \mathcal{B}} (-1)^{p+1}
\left[ \left(\varkappa_{\mathcal{S}} - \spont \right) 
\vec\tau\,.\,\vec\ek_1\,\chi \right](p) \nonumber \\ &
= 2\,\pi\,\alpha\,
\sum_{p\in \mathcal{B}} (-1)^{p+1} 
\left[\vec x\,.\,\vec\ek_1\,
(\varkappa_{\mathcal{S}})_s\,\chi\right](p) 
- 2\,\pi\,\varsigma\,\sum_{p \in \partial_{SF} I \cup \partial_F I}
\left[\vec x\,.\,\vec\ek_1\,\frac{\vec\nu\,.\,\vec\ek_1}{\vec x\,.\,\vec\ek_1}
\,\chi\right](p) 
\nonumber \\ &
= 2\,\pi\,
\sum_{p\in \partial_{SF} I \cup \partial_F I} \left[ \vec x\,.\,\vec\ek_1
\left((-1)^{p+1} \,\alpha\,(\varkappa_{\mathcal{S}})_s
- \varsigma\,\frac{\vec\nu\,.\,\vec\ek_1}{\vec x\,.\,\vec\ek_1} \right)
\chi\right](p) 
\qquad \forall\ \chi \in H^1(I)\,. 
\label{eq:A1:Bnormal}
\end{align}
where in the last step we have observed \eqref{eq:sdbca} and \eqref{eq:matB}. 
This gives (\ref{eq:axifreea}) on $\partial_{F} I$, 
on recalling (\ref{eq:xpos}) and \eqref{eq:xspace}. 

Next we consider the tangential components. It holds, on noting 
(\ref{eq:A1:mye1}), (\ref{eq:meanGaussS}), 
$\vec\nu\,.\,\vec\ek_1=-\vec\tau\,.\,\vec\ek_2$, (\ref{eq:bclimit}),
(\ref{eq:mu}) and (\ref{eq:axibc}), that
\begin{align}
& \sum_{i=1}^5 B_i(\chi\,\vec\tau) 
= 2\,\pi\,\sum_{p \in \partial_{SF} I \cup \partial_F I} (-1)^{p+1} 
\left[ \varkappa\,(\vec x\,.\,\vec\ek_1
\left[\alpha\,(\varkappa_{\mathcal{S}} - \spont) + \beta\,\AADE_{\mathcal{S}}
\right])\,\chi\right](p) \nonumber \\ & \quad
- 2\,\pi\,\varsigma\,\sum_{p \in \partial_{SF} \cup \partial_F I}
\left[\vec\tau\,.\,\vec\ek_1\,\chi\right](p)  
- 2\,\pi\,\beta\,\AADE_{\mathcal{S}}\, 
\sum_{p\in \mathcal{B}} (-1)^{p+1}
\left[ \vec\tau\,.\,\vec\ek_2\,\chi \right](p)\nonumber \\ & \quad
- \pi\, \sum_{p\in\partial_{SF} I \cup \partial_F I} (-1)^{p+1}
\left[\vec x\,.\,\vec\ek_1\left( \alpha
\left[\varkappa_{\mathcal{S}} -\spont 
\right]^2 + 2\,\lambda + 2\,\beta\,\AADE_{\mathcal{S}}\,\varkappa \right) 
\chi \right](p) \nonumber \\ & 
= -2\,\pi\,\alpha_G\,\sum_{p \in \partial_{SF} I \cup \partial_F I} (-1)^{p+1} 
\left[ \vec x\,.\,\vec\ek_1\,\Gauss_{\mathcal{S}}\,\chi\right](p) 
- 2\,\pi\,\varsigma\,\sum_{p \in \partial_{SF} \cup \partial_F I}
(-1)^{p+1}
\left[\vec x\,.\,\vec\ek_1\,
\frac{\vec\mu\,.\,\vec\ek_1}{\vec x\,.\,\vec\ek_1}\,\chi\right](p)  
\nonumber \\ & \quad
+ 2\,\pi\,\beta\,\AADE_{\mathcal{S}}\, 
\sum_{p\in\partial_{SF} I \cup \partial_F I} (-1)^{p+1}
\left[\vec x\,.\,\vec\ek_1\,
 \frac{\vec\nu\,.\,\vec\ek_1}{\vec x\,.\,\vec\ek_1}\,\chi \right](p)
\nonumber \\ & \quad
- 2\,\pi\, \sum_{p\in\partial_{SF} I \cup \partial_F I} (-1)^{p+1}
\left[\vec x\,.\,\vec\ek_1\left( \tfrac12\,\alpha
\left[\varkappa_{\mathcal{S}} -\spont 
\right]^2 + \lambda + \beta\,\AADE_{\mathcal{S}}\,\varkappa \right) 
\chi \right](p) 
\nonumber \\ & 
= - 2\,\pi\,\varsigma\,\sum_{p \in \partial_{SF} \cup \partial_F I}
(-1)^{p+1}
\left[\vec x\,.\,\vec\ek_1\,
\frac{\vec\mu\,.\,\vec\ek_1}{\vec x\,.\,\vec\ek_1}\,\chi\right](p)  
\nonumber \\ & \quad
- 2\,\pi\, \sum_{p\in\partial_{SF} I \cup \partial_F I} (-1)^{p+1}
\left[\vec x\,.\,\vec\ek_1\left( 
\alpha_G\,\Gauss_{\mathcal{S}} + 
\tfrac12\,\alpha
\left[\varkappa_{\mathcal{S}} -\spont 
\right]^2 + \lambda + \beta\,\AADE_{\mathcal{S}}\,\varkappa_{\mathcal{S}}
 \right) \chi \right](p) 
\nonumber \\ & 
= 2\,\pi\,\sum_{p \in \partial_{SF} \cup \partial_F I}
(-1)^{p}
\left[\vec x\,.\,\vec\ek_1\left( 
\varsigma\,\frac{\vec\mu\,.\,\vec\ek_1}{\vec x\,.\,\vec\ek_1} +
\alpha_G\,\Gauss_{\mathcal{S}} + 
\tfrac12\,\alpha
\left[\varkappa_{\mathcal{S}} -\spont 
\right]^2 + \lambda + \beta\,\AADE_{\mathcal{S}}\,\varkappa_{\mathcal{S}}
 \right) \chi \right](p) 
\quad \forall\ \chi \in H^1(I)\,. 
\label{eq:A1:Btangential}
\end{align}
This yields (\ref{eq:axifreeb}) on $\partial_F I$, 
on recalling (\ref{eq:xpos}) and \eqref{eq:xspace}.
In order to prove \eqref{eq:axisfreeb} on $\partial_{SF} I$, we choose
a test function $\vec\chi \in \xspace$ and then combine
(\ref{eq:A1:Bnormal}) and \eqref{eq:A1:Btangential}.
For example, on $\partial_1 I$ we choose 
$\vec\chi = \chi\,\vec\ek_2 
= \chi\,(\vec\nu\,.\,\vec\ek_2)\,\vec\nu + \chi\,(\vec\tau\,.\,\vec\ek_2)\,
\vec\tau
= \chi\,(-1)^{p+1}\,(\vec\mu\,.\,\vec\ek_1)\,\vec\nu
+ \chi\,(- \vec\nu\,.\,\vec\ek_1)\,\vec\tau$, and hence we obtain the 
desired result. The case $\partial_2 I$ follows analogously.

\subsection{Formulations based on $\varkappa_{\mathcal{S}}$} \label{sec:A2}

Here we derive the strong form for (\ref{eq:PS}), together with possible
boundary conditions for the $L^2$--gradient flow of (\ref{eq:EEax}).
We recall from (\ref{eq:kappaSy}) and (\ref{eq:ASxkappa}) that
\begin{equation} \label{eq:A2:kappay}
\vec y_{\mathcal{S}}\,.\,\vec\nu =
2\,\pi\left[
\alpha\,(\varkappa_{\mathcal{S}} - \spont) + \beta\,\AADE_{\mathcal{S}}\right] ,
\quad\text{where }
\AADE_{\mathcal{S}} = 2\,\pi\left(\vec x\,.\,\vec\ek_1\,\varkappa_{\mathcal{S}}, 
|\vec x_\rho| \right) - M_0\,.
\end{equation}

We begin by re-stating (\ref{eq:PSa}), on noting (\ref{eq:tau}), as
\begin{align}
& 
2\,\pi\left((\vec x\,.\,\vec\ek_1)\,\vec x_t\,.\,\vec\nu, \vec\chi\,.\,\vec\nu\,
|\vec x_\rho|\right) 
=  \left((\vec x\,.\,\vec\ek_1)\,(\vec y_{\mathcal{S}})_\rho \,.\,\vec\nu , 
 \vec\chi_\rho\,.\,\vec\nu \,|\vec x_{\rho}|^{-1}\right)
\nonumber \\ & \qquad 
-\left(\pi\,\vec x\,.\,\vec\ek_1\left[\alpha\,
\left(\varkappa_{\mathcal{S}}-\spont \right)^2 
+ 2\,\lambda + 2\,\beta\,\AADE_{\mathcal{S}}\,\varkappa_{\mathcal{S}}\right]
-\vec y_{\mathcal{S}}\,.\,\vec\ek_1, \vec\chi_\rho\,.\,\vec\tau \right) 
\nonumber \\ & \qquad 
- \left( \pi\,[\,\alpha\left(\varkappa_{\mathcal{S}}-\spont \right)^2
+ 2\,\lambda + 2\,\beta\,\AADE_{\mathcal{S}}\,\varkappa_{\mathcal{S}}\,]
- \varkappa_{\mathcal{S}}\,\vec y_{\mathcal{S}}\,.\,\vec\nu 
- (\vec y_{\mathcal{S}})_s\,.\,\vec\tau, 
\vec\chi\,.\,\vec\ek_1\,|\vec x_\rho| \right) \nonumber \\ & \qquad 
+ \left(
\vec x\,.\,\vec\ek_1\,\varkappa_{\mathcal{S}}\,\vec y_{\mathcal{S}}^\perp,
\vec\chi_\rho \right) 
-\sum_{p \in \partial_{SF} I \cup \partial_F I}
\left[
\left[2\,\pi\,\varsigma + \vec{\rm m}\,.\,\vec y_{\mathcal{S}}
\right] \vec\chi\,.\,\vec\ek_1 \right](p) 
 = \sum_{i=1}^5 S_i(\vec\chi)
\qquad \forall\ \vec\chi \in \xspace\,.
\label{eq:A2:1}
\end{align}
On recalling (\ref{eq:tau}), (\ref{eq:xspace}) and (\ref{eq:matB})
 throughout these $S_i(\vec\chi)$ calculations, 
we have that
\begin{align}
S_1(\vec\chi) & = \left( \vec x\,.\,\vec\ek_1\,((\vec y_{\mathcal{S}})_s 
\,.\,\vec\nu)\,\vec\nu , \vec\chi_\rho \right) \nonumber \\ &
= - \left( \left[\vec x\,.\,\vec\ek_1\,((\vec y_{\mathcal{S}})_s 
\,.\,\vec\nu)\,\vec\nu \right]_s , \vec\chi\,|\vec x_\rho| \right)
+ \sum_{p\in \mathcal{B}} (-1)^{p+1} \left[
\vec x\,.\,\vec\ek_1\,((\vec y_{\mathcal{S}})_s 
\,.\,\vec\nu)\,\vec\chi\,.\,\vec\nu \right](p) 
= A_1(\vec\chi)  + B_1(\vec\chi)  \,.
\label{eq:A2:S1}
\end{align}
Moreover, it holds that
\begin{align}
& S_2(\vec\chi) = -\left( \left[\pi\,\vec x\,.\,\vec\ek_1\left[\alpha
\left(\varkappa_{\mathcal{S}}-\spont \right)^2
+ 2\,\lambda+ 2\,\beta\,\AADE_{\mathcal{S}}\,\varkappa_{\mathcal{S}}\right]
-\vec y_{\mathcal{S}} \,.\,\vec\ek_1\right] \vec\tau
, \vec\chi_\rho \right) \nonumber \\ &
= \left( \left [ \left(\pi\,\vec x\,.\,\vec\ek_1\left[\alpha
\left(\varkappa_{\mathcal{S}}-\spont \right)^2
+ 2\,\lambda+ 2\,\beta\,\AADE_{\mathcal{S}}\,\varkappa_{\mathcal{S}}\right]
-\vec y_{\mathcal{S}} \,.\,\vec\ek_1\right) \vec\tau \right]_s
, \vec\chi \,|\vec x_\rho| \right) \nonumber \\ & \quad
- \sum_{p\in \mathcal{B}} (-1)^{p+1} \left[
\left(\pi\,\vec x\,.\,\vec\ek_1\left[\alpha
\left(\varkappa_{\mathcal{S}}-\spont \right)^2
+ 2\,\lambda+2\,\beta\,\AADE_{\mathcal{S}}\,\varkappa_{\mathcal{S}}\right]
-\vec y_{\mathcal{S}} \,.\,\vec\ek_1\right) \vec\chi\,.\,\vec\tau \right](p)
= A_2(\vec\chi) + B_2(\vec\chi)\,.
\label{eq:A2:S2}
\end{align}
In addition,
\begin{align}
S_3(\vec\chi) & = 
- \left( \left[\pi\,[\alpha\left(\varkappa_{\mathcal{S}}-\spont \right)^2
+ 2\,\lambda+2\,\beta\,\AADE_{\mathcal{S}}\,\varkappa_{\mathcal{S}}]
- \varkappa_{\mathcal{S}}\,\vec y_{\mathcal{S}}\,.\,\vec\nu 
- (\vec y_{\mathcal{S}})_s\,.\,\vec\tau \right]
\vec\ek_1 , \vec\chi\,|\vec x_\rho| \right)
= A_3(\vec\chi)\,.
\label{eq:A2:S3}
\end{align}
Finally,
\begin{equation}
S_4(\vec\chi) = \left(
\vec x\,.\,\vec\ek_1\,\varkappa_{\mathcal{S}}\,\vec y_{\mathcal{S}}^\perp,
\vec\chi_\rho \right)
= - \left( \left[
\vec x\,.\,\vec\ek_1\,\varkappa_{\mathcal{S}}\,\vec y_{\mathcal{S}}^\perp 
\right]_s , \vec\chi\,|\vec x_\rho| \right)
+ \sum_{p\in \mathcal{B}} (-1)^{p+1} 
\left[ \vec x\,.\,\vec\ek_1\,\varkappa_{\mathcal{S}}\,
\vec\chi\,.\,\vec y_{\mathcal{S}}^\perp\right] (p)
= A_4(\vec\chi) + B_4(\vec\chi)\,,
\label{eq:A2:S4}
\end{equation}
and we also set $B_3(\vec\chi) = S_5(\vec\chi)$.
With the above definitions, we can write
(\ref{eq:A2:1}) as
\begin{equation} \label{eq:A2:1compact}
2\,\pi
\left((\vec x\,.\,\vec\ek_1)\,\vec x_t\,.\,\vec\nu, \vec\chi\,.\,\vec\nu\,
|\vec x_\rho|\right) 
= \sum_{i=1}^5 S_i(\vec\chi) 
= \sum_{i=1}^4 A_i(\vec\chi) + 
\sum_{i=1}^4 B_i(\vec\chi) \qquad \forall\ \vec\chi \in \xspace\,.
\end{equation}

On recalling (\ref{eq:aperp}), we observe that
\begin{equation}
A_4(\vec\chi) = - \left( \left[
\vec x\,.\,\vec\ek_1\,\varkappa_{\mathcal{S}}\,\vec y_{\mathcal{S}}^\perp 
\right]_s , \vec\chi\,|\vec x_\rho| \right)
= - \left( \left[
\vec x\,.\,\vec\ek_1\,\varkappa_{\mathcal{S}}\left(
(\vec y_{\mathcal{S}}\,.\,\vec\nu)\,\vec\tau - 
(\vec y_{\mathcal{S}}\,.\,\vec\tau)\,\vec\nu \right)
\right]_s , \vec\chi\,|\vec x_\rho| \right) . \label{eq:A2:A4}
\end{equation}

Choosing $\vec\chi = \chi\,\vec\tau$, for $\chi \in H^1_0(I)$, 
in (\ref{eq:A2:1compact}), 
and combining (\ref{eq:A2:S1}), (\ref{eq:A2:S2}), (\ref{eq:A2:S3}),
(\ref{eq:A2:S4}) and (\ref{eq:A2:A4}), 
and recalling $\vec\tau_s\,.\,\vec\tau=0$,
(\ref{eq:nus}), (\ref{eq:tau}), 
(\ref{eq:A2:kappay}) and (\ref{eq:meanGaussS}), 
we obtain for the right hand side of (\ref{eq:A2:1compact}) 
\begin{align*}
& \sum_{i=1}^5 S_i(\chi\,\vec\tau) =
\sum_{i=1}^4 A_i(\chi\,\vec\tau)  
\nonumber \\ & \,
= \left( \vec x\,.\,\vec\ek_1\,\varkappa\,(\vec y_{\mathcal{S}})_s 
\,.\,\vec\nu, \chi\,|\vec x_\rho| \right)
+ \left( \left[\pi\,\vec x\,.\,\vec\ek_1\left[\alpha
\left(\varkappa_{\mathcal{S}}-\spont \right)^2
+ 2\,\lambda +2\,\beta\,\AADE_{\mathcal{S}}\,\varkappa_{\mathcal{S}}\right]
-\vec y_{\mathcal{S}} \,.\,\vec\ek_1\right]_s, \chi \,|\vec x_\rho| \right) 
\nonumber \\ & \quad
- \left( \left[ \pi \left[\alpha\left(\varkappa_{\mathcal{S}}-\spont \right)^2
+ 2\,\lambda +2\,\beta\,\AADE_{\mathcal{S}}\,\varkappa_{\mathcal{S}}\right]
- \varkappa_{\mathcal{S}}\,\vec y_{\mathcal{S}}\,.\,\vec\nu 
- (\vec y_{\mathcal{S}})_s\,.\,\vec\tau \right] \vec\tau\,.\,
\vec\ek_1 , \chi\,|\vec x_\rho| \right) 
\nonumber \\ & \quad
- \left( \vec x\,.\,\vec\ek_1\,\varkappa_{\mathcal{S}}
\left[(\vec y_{\mathcal{S}}\,.\,\vec\nu)_s
+ \varkappa\,
\vec y_{\mathcal{S}}\,.\,\vec\tau\right]
+ (\vec x\,.\,\vec\ek_1\,\varkappa_{\mathcal{S}})_s\,
\vec y_{\mathcal{S}}\,.\,\vec\nu , \chi\,|\vec x_\rho| \right)
\nonumber \\ & \,
= \left( \vec x\,.\,\vec\ek_1\,(\varkappa-\varkappa_{\mathcal{S}})\,
(\vec y_{\mathcal{S}})_s \,.\,\vec\nu, \chi\,|\vec x_\rho| \right)
+ \left( \vec x\,.\,\vec\ek_1 \left(\pi \left[\alpha \,
\left(\varkappa_{\mathcal{S}}-\spont \right)^2 
+2\,\beta\,\AADE_{\mathcal{S}}\,\varkappa_{\mathcal{S}}\right]_s
-(\varkappa_{\mathcal{S}})_s\,
\vec y_{\mathcal{S}}\,.\,\vec\nu 
\right),\chi \,|\vec x_\rho| \right) 
\nonumber \\ & \quad
- \left( (\vec y_{\mathcal{S}})_s \,.\left(\vec\ek_1
-(\vec\tau\,.\,\vec\ek_1)\,\vec\tau\right), \chi\,|\vec x_\rho| \right) 
\nonumber \\ & \, = 0\,,  
\end{align*}
as expected.

Choosing $\vec\chi = \chi\,\vec\nu$, for $\chi \in H^1_0(I)$, 
in (\ref{eq:A2:1compact}), 
and combining (\ref{eq:A2:S1}), (\ref{eq:A2:S2}), (\ref{eq:A2:S3}),
(\ref{eq:A2:S4}) and (\ref{eq:A2:A4}), and on recalling 
(\ref{eq:varkappa}), $\vec\nu_s\,.\,\vec\nu=0$, 
(\ref{eq:A2:kappay}) and (\ref{eq:meanGaussS}),
we obtain for the right hand side of (\ref{eq:A2:1compact}) 
\begin{align}
& \sum_{i=1}^5 S_i(\chi\,\vec\nu) =
\sum_{i=1}^4 A_i(\chi\,\vec\nu)  
= - \left( \left[\vec x\,.\,\vec\ek_1\,((\vec y_{\mathcal{S}})_s 
\,.\,\vec\nu)\right]_s , \chi\,|\vec x_\rho| \right)
\nonumber \\ & \
+ \left(  \left[\pi\,\vec x\,.\,\vec\ek_1\left[\alpha
\left(\varkappa_{\mathcal{S}}-\spont \right)^2
+ 2\,\lambda+2\,\beta\,\AADE_{\mathcal{S}}\,\varkappa_{\mathcal{S}}\right]
-\vec y_{\mathcal{S}} \,.\,\vec\ek_1\right] \varkappa,
\chi \,|\vec x_\rho| \right) 
\nonumber \\ & \
- \left( \left[ \pi \left[\alpha\left(\varkappa_{\mathcal{S}}-\spont \right)^2
+ 2\,\lambda +2\,\beta\,\AADE_{\mathcal{S}}\,\varkappa_{\mathcal{S}} \right]
- \varkappa_{\mathcal{S}}\,\vec y_{\mathcal{S}}\,.\,\vec\nu 
- (\vec y_{\mathcal{S}})_s\,.\,\vec\tau \right] \vec\nu\,.\,
\vec\ek_1 , \chi\,|\vec x_\rho| \right) 
\nonumber \\ & \
- \left( \vec x\,.\,\vec\ek_1\,\varkappa\,\varkappa_{\mathcal{S}}\,
\vec y_{\mathcal{S}}\,.\,\vec\nu
- \left[\vec x\,.\,\vec\ek_1\,\varkappa_{\mathcal{S}}\,
\vec y_{\mathcal{S}}\,.\,\vec\tau \right]_s, \chi\,|\vec x_\rho| \right)
\nonumber \\ & 
= 
- \left( \left[\vec x\,.\,\vec\ek_1 \left((\vec y_{\mathcal{S}} 
\,.\,\vec\nu)_s + (\varkappa-\varkappa_{\mathcal{S}})\,\vec y_{\mathcal{S}}\,.\,\vec\tau
\right)\right]_s , \chi\,|\vec x_\rho| \right)
\nonumber \\ & \
+ \left(  \left[\pi\,\vec x\,.\,\vec\ek_1\left[\alpha
\left(\varkappa_{\mathcal{S}}-\spont \right)^2
+ 2\,\lambda+2\,\beta\,\AADE_{\mathcal{S}}\,\varkappa_{\mathcal{S}}\right]
-\vec y_{\mathcal{S}} \,.\,\vec\ek_1\right] \varkappa,
\chi \,|\vec x_\rho| \right) 
\nonumber \\ & \
- \left( \left[ \pi \left[\alpha\left(\varkappa_{\mathcal{S}}-\spont \right)^2
+  2\,\lambda+2\,\beta\,\AADE_{\mathcal{S}}\,\varkappa_{\mathcal{S}} \right]
+ (\varkappa- \varkappa_{\mathcal{S}})\,\vec y_{\mathcal{S}}\,.\,\vec\nu 
- (\vec y_{\mathcal{S}}\,.\,\vec\tau)_s \right] \vec\nu\,.\,
\vec\ek_1 , \chi\,|\vec x_\rho| \right) 
- \left( \vec x\,.\,\vec\ek_1\,\varkappa\,\varkappa_{\mathcal{S}}\,
\vec y_{\mathcal{S}}\,.\,\vec\nu, \chi\,|\vec x_\rho| \right)
\nonumber \\ & 
= - 2\,\pi\,\alpha \left( \left[\vec x\,.\,\vec\ek_1\,(\varkappa_{\mathcal{S}} 
)_s \right]_s , \chi\,|\vec x_\rho| \right)
- \left(\left[
\vec x\,.\,\vec\ek_1\,(\varkappa-\varkappa_{\mathcal{S}})\,\vec y_{\mathcal{S}}\,.\,\vec\tau
\right]_s , \chi\,|\vec x_\rho| \right)
\nonumber \\ & \
+ \left( \pi\,\vec x\,.\,\vec\ek_1\left[\alpha
\left(\varkappa_{\mathcal{S}}-\spont \right)^2
+ 2\,\lambda+2\,\beta\,\AADE_{\mathcal{S}}\,\varkappa_{\mathcal{S}}\right]
 \varkappa_{\mathcal{S}},
\chi \,|\vec x_\rho| \right) 
- \left( \vec x\,.\,\vec\ek_1\left[(\varkappa-\varkappa_{\mathcal{S}})^2
+ \varkappa\,\varkappa_{\mathcal{S}}\right]
\vec y_{\mathcal{S}}\,.\,\vec\nu, \chi\,|\vec x_\rho| \right)
\nonumber \\ & \
- \left( \vec y_{\mathcal{S}} \,.\,\vec\ek_1\,\varkappa 
- (\vec y_{\mathcal{S}}\,.\,\vec\tau)_s\, \vec x\,.\,
\vec\ek_1 \, (\varkappa- \varkappa_{\mathcal{S}}), \chi\,|\vec x_\rho| \right) 
\nonumber \\ & 
= - 2\,\pi\,\alpha \left( \left[\vec x\,.\,\vec\ek_1\,(\varkappa_{\mathcal{S}} 
)_s \right]_s , \chi\,|\vec x_\rho| \right)
- \left(\left[
\vec x\,.\,\vec\ek_1\,(\varkappa-\varkappa_{\mathcal{S}})\right]_s 
\vec y_{\mathcal{S}}\,.\,\vec\tau , \chi\,|\vec x_\rho| \right)
\nonumber \\ & \
- \pi\,\alpha \left( \vec x\,.\,\vec\ek_1 
\left[ (\varkappa_{\mathcal{S}})^2 + \spont\,\varkappa_{\mathcal{S}}
+ 2\,\varkappa^2 - 2\,\varkappa\,\varkappa_{\mathcal{S}} \right]
(\varkappa_{\mathcal{S}} - \spont), \chi\,|\vec x_\rho| \right) 
\nonumber \\ & \
- \left(  \vec y_{\mathcal{S}} \,.\,\vec\ek_1\, \varkappa,
\chi \,|\vec x_\rho| \right) 
+ 2\,\pi\,\lambda \left(\vec x\,.\,\vec\ek_1\, \varkappa_{\mathcal{S}},
\chi \,|\vec x_\rho| \right) 
- 2\,\pi\,\beta\,\AADE_{\mathcal{S}} \left( \varkappa\, \vec x\,.\,
\vec\ek_1 \, (\varkappa- \varkappa_{\mathcal{S}}), \chi\,|\vec x_\rho| \right) 
\nonumber \\ & 
= 
- 2\,\pi\,\alpha \left( \left[\vec x\,.\,\vec\ek_1\,(\varkappa_{\mathcal{S}} 
)_s \right]_s , \chi\,|\vec x_\rho| \right)
- \left(\left[
\vec x\,.\,\vec\ek_1\,(\varkappa-\varkappa_{\mathcal{S}})\right]_s
\vec y_{\mathcal{S}}\,.\,\vec\tau , \chi\,|\vec x_\rho| \right)
\nonumber \\ & \
- \left(  \pi\,\vec x\,.\,\vec\ek_1\left[\alpha
\left(\varkappa_{\mathcal{S}}^2-\spont^2 \right)
- 2\,\lambda\right] \varkappa_{\mathcal{S}}
+\vec y_{\mathcal{S}} \,.\,\vec\ek_1\, \varkappa, \chi \,|\vec x_\rho| \right) 
+ 2\,\pi \left( \vec x\,.\,\vec\ek_1 \left[\alpha\,
(\varkappa_{\mathcal{S}}-\spont) + \beta\,\AADE_{\mathcal{S}}\right]
\left(\varkappa_{\mathcal{S}}-\varkappa\right)\varkappa , \chi\,|\vec x_\rho| \right) 
\,.
\label{eq:A2:normal}
\end{align}
The desired result (\ref{eq:xtbgnlambda}) 
follows from (\ref{eq:A2:1}) 
and (\ref{eq:A2:normal}) on noting that
\begin{align*}
& -\left[
\vec x\,.\,\vec\ek_1\,(\varkappa-\varkappa_{\mathcal{S}})\right]_s
\vec y_{\mathcal{S}}\,.\,\vec\tau
-\vec y_{\mathcal{S}} \,.\,\vec\ek_1\, \varkappa
= 
-(\vec\nu_s \,.\,\vec\ek_1)\,
\vec y_{\mathcal{S}}\,.\,\vec\tau
-\vec y_{\mathcal{S}} \,.\,\vec\ek_1\, \varkappa
= -(\vec y_{\mathcal{S}}\,.\,\vec\nu)
\,\vec\nu \,.\,\vec\ek_1\, \varkappa
\nonumber \\ & \qquad
= 
2\,\pi\,\vec x\,.\,\vec\ek_1\left[\alpha\,(\varkappa_{\mathcal{S}}-\spont) 
+ \beta\,\AADE_{\mathcal{S}}\right]
(\varkappa_{\mathcal{S}}-\varkappa)\, \varkappa
= 
2\,\pi\,\vec x\,.\,\vec\ek_1\left[\alpha\,(\varkappa_{\mathcal{S}}-\spont)
+\beta\,\AADE_{\mathcal{S}}\right] 
\mathcal{K}_{\mathcal{S}}\,, 
\end{align*}
where we have recalled (\ref{eq:meanGaussS}), (\ref{eq:nus}) and 
(\ref{eq:A2:kappay}).

We now need to show that the weak formulation \eqref{eq:PS} enforces the
claimed boundary conditions. As in \S\ref{sec:A1}, the conditions
(\ref{eq:axisfreea}), (\ref{eq:fixedC}), (\ref{eq:fixedN}) and 
(\ref{eq:fixed0}) are enforces strongly. We also notice that \eqref{eq:PSc} 
enforces \eqref{eq:bcbc} and (\ref{eq:clampedI}), recall the discussion below 
\eqref{eq:varkappaSweak}. 
This means that we still
need to show \eqref{eq:sdbca}, \eqref{eq:axifree}, \eqref{eq:axisfreeb},
\eqref{eq:axisfreec} and \eqref{eq:NavierI}.

In order to show \eqref{eq:sdbca} we argue similarly
to \S\ref{sec:A1}. Choosing $\vec\chi=\chi\,\vec\nu$ in \eqref{eq:A2:1},
where $\chi \in H^1(I)$ is zero away from $\partial_0 I$, 
that the term $2\,\pi\,\alpha\,(\vec x\,.\,\vec\ek_1\,
(\varkappa_{\mathcal{S}})_\rho, \chi_\rho\,|\vec x_\rho|^{-1})$ equals a
sum of inner products that have integrands that converge to zero for $\rho$
approaching $\partial_0 I$. We can then again use \cite[Appendix~A.1]{axisd}
to show that \eqref{eq:PSa} implies the boundary condition \eqref{eq:sdbca}.

Next we recall from (\ref{eq:partialM}) and (\ref{eq:mSy}) 
that
\begin{equation} \label{eq:A2:my}
(\vec x\,.\, \vec\ek_1)\,\vec y_{\mathcal{S}} = 
2\,\pi\,\alpha_G\,\vec\ek_1 \qquad\text{ on } \quad \partial_M I
= \partial_N I \cup \partial_{SF} I \cup \partial_F I\,.
\end{equation}
Combining with (\ref{eq:A2:kappay}) and noting (\ref{eq:xpos}) yields 
(\ref{eq:A1:mye1}). 
Hence, we have that (\ref{eq:axifreec}), (\ref{eq:axisfreec}) and  
(\ref{eq:NavierI}) are imposed strongly.
Overall, we still need to show \eqref{eq:axifreea}, \eqref{eq:axifreeb} and
\eqref{eq:axisfreeb}. 
It follows from (\ref{eq:A2:S1}) and (\ref{eq:A2:S4}), 
on recalling (\ref{eq:nus}), (\ref{eq:aperp}), (\ref{eq:A2:kappay}), 
(\ref{eq:meanGaussS}) and (\ref{eq:axibc}), that
\begin{align*}
B_1(\vec\chi) + B_4(\vec\chi) &
= \sum_{p \in \mathcal{B}} (-1)^{p+1} \left[(\vec x \,.\,\vec\ek_1) 
\left( \left((\vec y_{\mathcal{S}})_s\,.\,\vec\nu\right)\vec\nu
+ \varkappa_{\mathcal{S}}\,\vec y_{\mathcal{S}}^\perp\right) 
.\,\vec\chi\right](p) \nonumber \\ &
=\sum_{p \in \mathcal{B}} (-1)^{p+1} \left[(\vec x \,.\,\vec\ek_1)
\left( \left[ (\vec y_{\mathcal{S}}\,.\,\vec\nu)_s +
(\varkappa-\varkappa_{\mathcal{S}})\,(\vec y_{\mathcal{S}}\,.\,\vec\tau)
\right] \vec\nu
+ \varkappa_{\mathcal{S}}\,(\vec y_{\mathcal{S}}\,.\,\vec\nu)
\,\vec\tau\right) .\,\vec\chi\right](p) 
\nonumber \\ &
= \sum_{p \in \mathcal{B}} (-1)^{p+1} \left[\,
 \left[ 2\,\pi\,\alpha\,(\vec x \,.\,\vec\ek_1)\,(\kappa_{\mathcal{S}})_s +
(\vec\nu\,.\,\vec\ek_1)\,(\vec y_{\mathcal{S}}\,.\,\vec\tau)
\right] \vec\nu\,.\,\vec\chi\right](p)
\nonumber \\ & \qquad
+ \sum_{p \in \partial_{SF}I \cup \partial_{F}I} (-1)^{p+1} 
\left[ (\vec x \,.\,\vec\ek_1)\,\varkappa_{\mathcal{S}}\,(\vec y_{\mathcal{S}}\,.\,\vec\nu)
\,\vec\tau\,.\,\vec\chi\right](p) \,. 
\end{align*}
Looking first at the normal boundary contributions, we compute, on noting
that $\vec{\rm m}= \vec\mu$ and (\ref{eq:mu}), 
\begin{align}
\sum_{i=1}^4 B_i(\chi\,\vec\nu) & 
= 
\sum_{p \in \mathcal{B}} (-1)^{p+1} \left[\,
 \left[ 2\,\pi\,\alpha\,(\vec x \,.\,\vec\ek_1)\,(\kappa_{\mathcal{S}})_s +
(\vec\nu\,.\,\vec\ek_1)\,(\vec y_{\mathcal{S}}\,.\,\vec\tau)
\right] \chi\right](p)
-\sum_{p \in \partial_{SF} I \cup \partial_F I}
\left[
\left[2\,\pi\,\varsigma + \vec{\rm m}\,.\,\vec y_{\mathcal{S}}
\right] (\vec\nu\,.\,\vec\ek_1)\,\chi \right](p) 
\nonumber \\ &
= 2\,\pi\,\alpha\,
\sum_{p\in \mathcal{B}} (-1)^{p+1} 
\left[\vec x\,.\,\vec\ek_1\,
(\varkappa_{\mathcal{S}})_s\,\chi\right](p) 
- 2\,\pi\,\varsigma\,\sum_{p \in \partial_{SF} I \cup \partial_F I}
\left[\vec x\,.\,\vec\ek_1\,\frac{\vec\nu\,.\,\vec\ek_1}{\vec x\,.\,\vec\ek_1}
\,\chi\right](p) \nonumber \\ & 
= 2\,\pi\,\sum_{p \in \partial_{SF} I \cup \partial_F I}
\left[ \vec x\,.\,\vec\ek_1\left(
(-1)^{p+1}\, \alpha\,
(\varkappa_{\mathcal{S}})_s
-\varsigma\,\frac{\vec\nu\,.\,\vec\ek_1}{\vec x\,.\,\vec\ek_1} \right)
\chi\right](p) 
\qquad \forall\ \chi \in H^1(I)\,.
\label{eq:A2:Bnormal}
\end{align}
where in the last step we have observed \eqref{eq:sdbca} and \eqref{eq:matB}.  
This gives (\ref{eq:axifreea}) on $\partial_{F} I$ on recalling (\ref{eq:xpos})
and \eqref{eq:xspace}. 

Next we consider the tangential components. It holds, on noting 
that $\vec {\rm m}= \vec\mu$, (\ref{eq:mu}), (\ref{eq:meanGaussS}) and (\ref{eq:A2:my}), that 
\begin{align}
\sum_{i=1}^4 B_i(\chi\,\vec\tau) &
=\sum_{p \in \partial_{SF}I \cup \partial_{F}I} (-1)^{p+1} 
\left[ (\vec x \,.\,\vec\ek_1)\,\varkappa_{\mathcal{S}}\,(\vec y_{\mathcal{S}}\,.\,\vec\nu)
\,\chi\right](p) 
 \nonumber \\ 
&  \qquad - \sum_{p \in \partial_{SF} I \cup \partial_F I} (-1)^{p+1} \left[
\left(\pi\,\vec x\,.\,\vec\ek_1\left[\alpha
\left(\varkappa_{\mathcal{S}}-\spont \right)^2
+ 2\,\lambda+2\,\beta\,\AADE_{\mathcal{S}}\,\varkappa_{\mathcal{S}}\right]
-\vec y_{\mathcal{S}} \,.\,\vec\ek_1\right) \chi \right](p)
\nonumber \\
& \qquad - 2\,\pi\,\varsigma \sum_{p \in \partial_{SF} I \cup \partial_F I}
\left[(\vec\tau\,.\,\vec\ek_1) \,\chi \right](p)
-\sum_{p \in \partial_{SF} I \cup \partial_F I}
(-1)^{p+1}
\left[ (\vec y_{\mathcal{S}}\,.\,\vec\tau)\,(\vec\tau\,.\,\vec\ek_1) \,\chi \right](p) 
\nonumber \\
& =\sum_{p \in \partial_{SF}I \cup \partial_{F}I} (-1)^{p+1} 
\left[ (\vec x \,.\,\vec\ek_1)\,\varkappa\,(\vec y_{\mathcal{S}}\,.\,\vec\nu)
\,\chi\right](p)
- 2\,\pi\,\varsigma \sum_{p \in \partial_{SF} I \cup \partial_F I}
\left[(\vec\tau\,.\,\vec\ek_1) 
\,\chi \right](p)  
 \nonumber \\ 
& \qquad - \sum_{p \in \partial_{SF} I \cup \partial_F I} (-1)^{p+1} \left[
\pi\,\vec x\,.\,\vec\ek_1\left[\alpha
\left(\varkappa_{\mathcal{S}}-\spont \right)^2
+ 2\,\lambda+2\,\beta\,\AADE_{\mathcal{S}}\,\varkappa_{\mathcal{S}}\right] \chi \right](p)
\nonumber \\
& =- 2\,\pi\,\alpha_G\sum_{p \in \partial_{SF}I \cup \partial_{F}I} (-1)^{p+1} 
\left[ (\vec x \,.\,\vec\ek_1)\,{\mathcal K}_{\mathcal{S}}
\,\chi\right](p)
  \nonumber \\ 
& \qquad  - \sum_{p \in \partial_{SF} I \cup \partial_F I} (-1)^{p+1} \left[
\pi\,\vec x\,.\,\vec\ek_1\left[\alpha
\left(\varkappa_{\mathcal{S}}-\spont \right)^2
+ 2\,\lambda+2\,\beta\,\AADE_{\mathcal{S}}\,\varkappa_{\mathcal{S}}\right] \chi \right](p)
\nonumber \\
& \qquad 
- 2\,\pi\,\varsigma \sum_{p \in \partial_{SF} I \cup \partial_F I} (-1)^{p+1}
\left[ \vec x \,.\,\vec\ek_1 \,
\frac{\vec\mu\,.\,\vec\ek_1}{\vec x\,.\,\vec\ek_1} 
\,\chi \right](p) 
\qquad \forall\ \chi \in H^1(I)\,.
\label{eq:A2:Btangential}
\end{align}
This yields (\ref{eq:axifreeb}) on $\partial_F I$ on recalling (\ref{eq:xpos}).
In order to prove \eqref{eq:axisfreeb} on $\partial_{SF} I$, we choose
a test function $\vec\chi \in \xspace$ and then combine
(\ref{eq:A2:Bnormal}) and \eqref{eq:A2:Btangential}, similarly to
the previous section, \S\ref{sec:A1}.

\renewcommand{\theequation}{B.\arabic{equation}}
\renewcommand{\thefigure}{B\arabic{figure}}
\setcounter{equation}{0}
\setcounter{figure}{0}
\section{Singularities for Willmore flow of genus-1 surfaces} 
\label{sec:B}

In the recent paper \cite{DallAcquaMSS20preprint} it was shown that a torus
of revolution $\mathcal{S}(0)$, 
with profile curve $\Gamma(0)$ such that its turning number
${\rm T}(\Gamma) = \frac1{2\,\pi}\, \int_I \varkappa\, |\vec x_\rho| \drho$
is zero, will develop a singularity under Willmore flow. 
Here we recall that for an immersed curve $\Gamma \subset \bR^2$, 
${\rm T}(\Gamma) \in \bZ$ is the winding number with respect to the origin of 
the tangent vector $\vec\tau$.
In particular, it is shown in \cite[Lemma~4.8]{DallAcquaMSS20preprint}
that either a singularity will develop in finite time, or as $t\to\infty$ 
one of the following quantities will grow unbounded:
$a(t) = \sup_I [ \varkappa^2 
+ ( \frac{\vec\nu\,.\,\vec\ek_1}{\vec x\,.\,\vec\ek_1} )^2 
]^\frac12$ or
$|\Gamma(t)|_{\bH^2} = ((\vec x\,.\,\vec\ek_1)^{-1},|\vec x_\rho|)$,
where $a(t)$ is a suitable $L^\infty$--norm of the second fundamental
form on $\mathcal{S}(t)$, and where
$|\Gamma(t)|_{\bH^2}$ denotes the length of the curve $\Gamma(t)$
in the hyperbolic plane $\bH^2$, see e.g.\ \cite[(2.6a)]{hypbol}. 

In this appendix, we will present numerical evidence that 
${\rm T}(\Gamma) = 0$ appears to be a necessary condition for a singularity to
develop, and that any such singularity is only attained as $t\to\infty$.
To this end, we consider the case $\partial I = \emptyset$ from now on and 
define 
\begin{equation} \label{eq:B:unboundedh}
|\Gamma^m|_{\bH^2,h} 
= \left( (\vec X^m\,.\,\vec\ek_1)^{-1},|\vec X^m_\rho| \right)^h 
\end{equation}
for $\Gamma^m = \vec X^m(I)$ and $\vec X^m \in \Vh$, as a discrete
analogue to $|\Gamma(t)|_{\bH^2}$. We remark that a corresponding discrete
analogue of $a(t)$ in practice behaved very similarly to 
\eqref{eq:B:unboundedh}, and so we omit its discussion here.
 
\begin{figure}
\center
\includegraphics[angle=-90,width=0.24\textwidth]{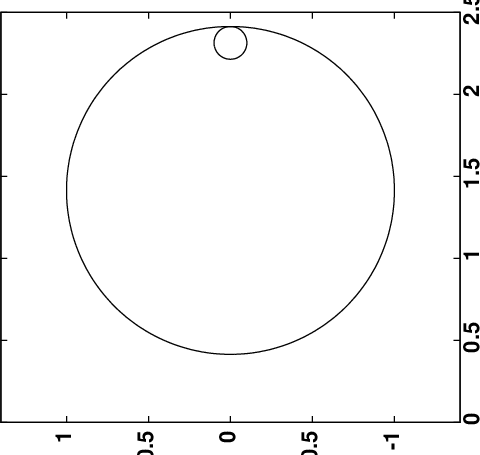}
\includegraphics[angle=-90,width=0.24\textwidth]{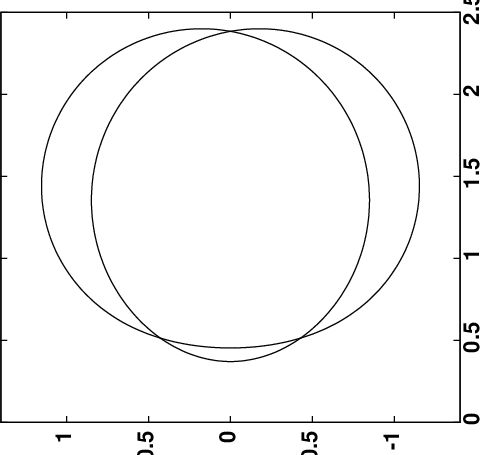}
\includegraphics[angle=-90,width=0.24\textwidth]{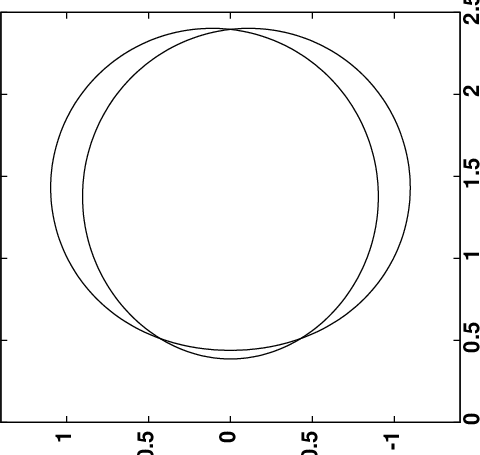}
\includegraphics[angle=-90,width=0.24\textwidth]{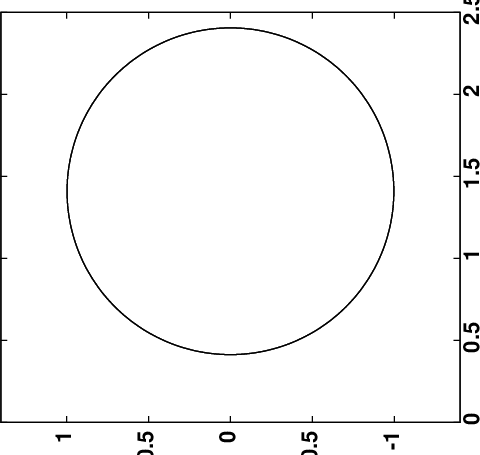}
\includegraphics[angle=-90,width=0.24\textwidth]{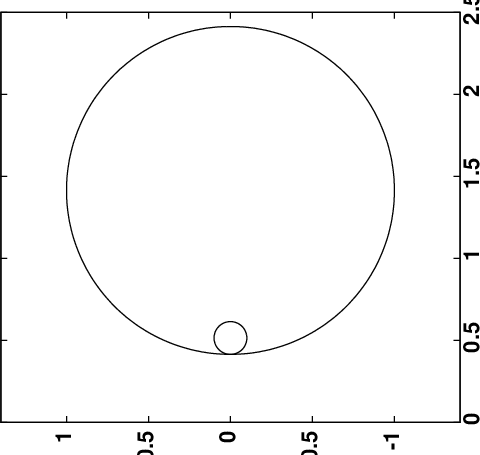}
\includegraphics[angle=-90,width=0.24\textwidth]{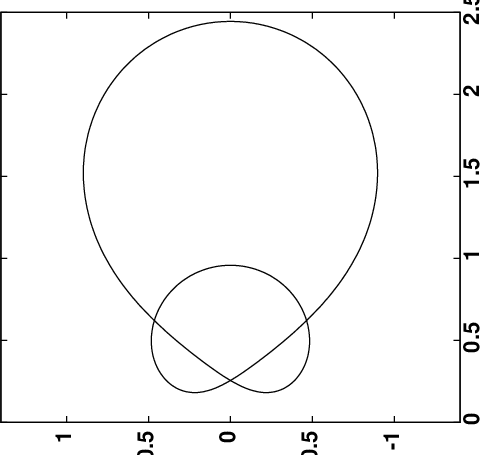}
\includegraphics[angle=-90,width=0.24\textwidth]{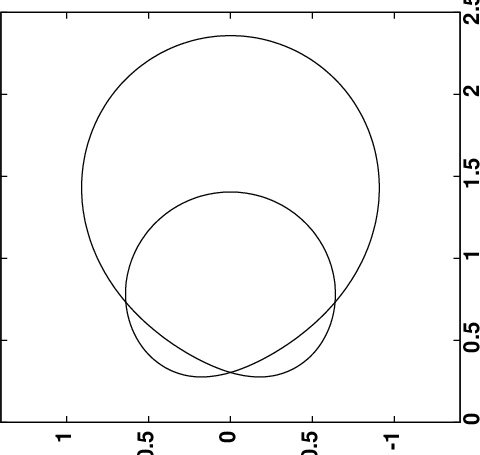}
\includegraphics[angle=-90,width=0.24\textwidth]{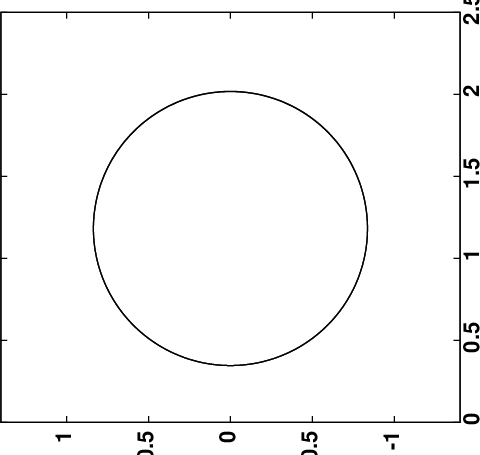}
\caption{
$(\BGNpwf^m)^h$
Willmore flow for ${\rm T}(\Gamma) = \pm2$.
We show the two evolutions at times $t=0,1,5,100$.
} 
\label{fig:B:clifford_limacon}
\end{figure}

We begin with two experiments for the case ${\rm T}(\Gamma) \not= 0$.
In particular, we choose as initial data a generating curve that is made up of 
a circle of radius 1 centered at $\sqrt{2}\,\vec\ek_1$ and a circle of radius
$r$ centered at $(\sqrt{2} \pm (1 - r))\,\vec\ek_1$, with $r=0.1$, so that
${\rm T}(\Gamma) = \pm 2$.
The simulations of Willmore
flow for the surfaces generated by these curves is shown in 
Figure~\ref{fig:B:clifford_limacon}. 
The discretization parameters for the scheme $(\BGNpwf^m)^{h}$
are $J=512$ and $\ttau = 10^{-4}$.
When the smaller circle is inscribed on the right, 
then the evolution immediately approaches a double-covering of the 
Clifford torus, recall also Figure~\ref{fig:pwfcigar41}.
If the smaller circle is inscribed on the left, a more complicated evolution 
ensues, but eventually a double covering of the Clifford torus is reached. 
At time $t=100$, both evolutions have reached a discrete Willmore energy of 
$78.96 \approx 8\,\pi^2$, i.e.\ about twice the Willmore energy of the 
Clifford torus. 
Overall the results in Figure~\ref{fig:B:clifford_limacon} indicate that 
${\rm T}(\Gamma) = 0$
is a necessary condition for a singularity under Willmore flow to occur.

\begin{figure}
\center
\includegraphics[angle=-90,width=0.3\textwidth]{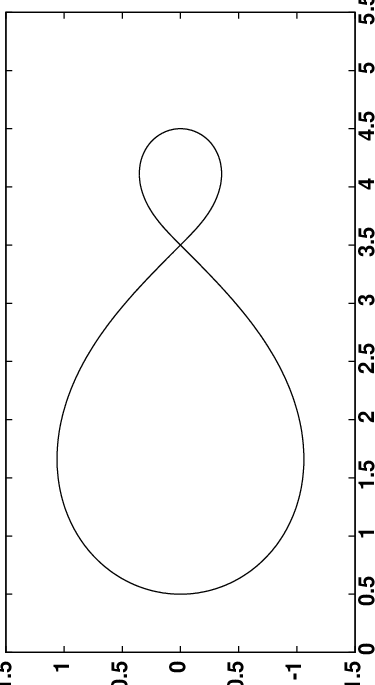}
\includegraphics[angle=-90,width=0.3\textwidth]{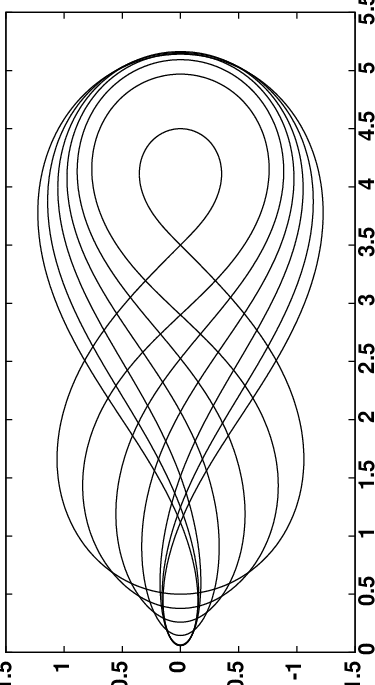}
\includegraphics[angle=-90,width=0.3\textwidth]{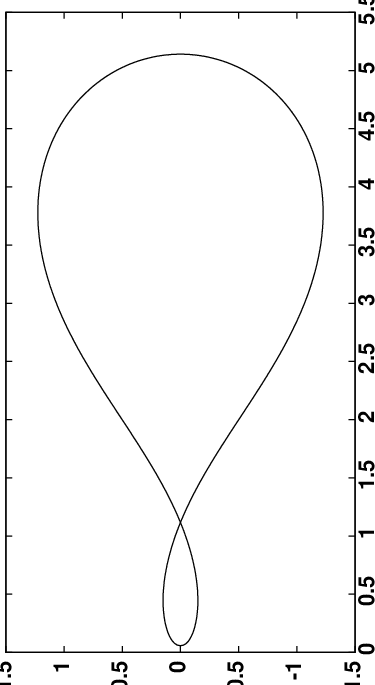}
\caption{$(\BGNpwf^m)^h$
Willmore flow for ${\rm T}(\Gamma) = 0$.
We show the solution at time $t=0$, at times
times $t=0,0.5,\ldots,3$ and at time $t=3$.
}
\label{fig:B:wf_lemniscate}
\end{figure}

\begin{figure}
\center
\includegraphics[angle=-90,width=0.3\textwidth]{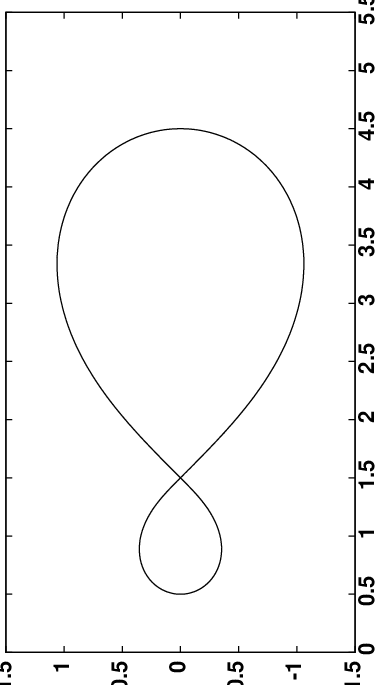}
\includegraphics[angle=-90,width=0.3\textwidth]{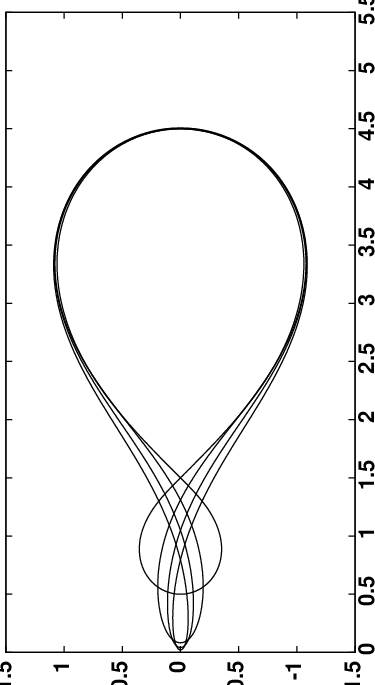}
\includegraphics[angle=-90,width=0.3\textwidth]{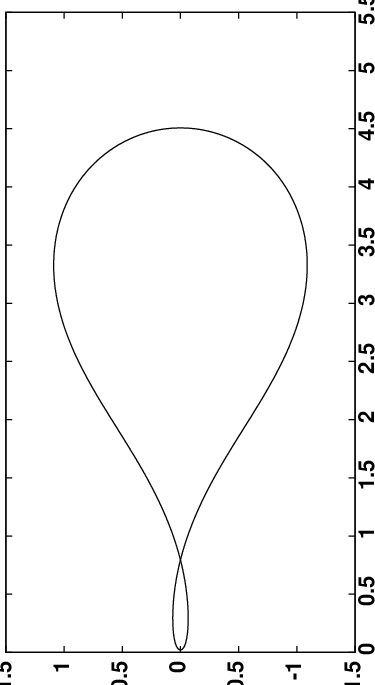}
\caption{$(\BGNpwf^m)^h$ 
Willmore flow for ${\rm T}(\Gamma) = 0$.
We show the solution at time $t=0$, at times
times $t=0,0.1,\ldots,0.3$ and at time $t=0.3$.
}
\label{fig:B:wf_leftlemniscate}
\end{figure}

We now concentrate on the possible onset of a singularity. To this end,
we show the evolution for two initial lemniscates in
Figures~\ref{fig:B:wf_lemniscate} and \ref{fig:B:wf_leftlemniscate}. 
The discretization parameters for the scheme $(\BGNpwf^m)^{h}$
are $J=1024$ and $\ttau = 10^{-6}$. 
The initial data satisfy ${\rm T}(\Gamma) = 0$, and so the result in 
\cite[Lemma~4.8]{DallAcquaMSS20preprint} yields that Willmore flow
will develop a singularity either in finite time or as $t\to\infty$.
In the first experiment the part of the lemniscate with larger radius is 
close to the $x_2$--axis. But during the evolution that part thins and
approaches the axis. In the second experiment, the initial lemniscate
is rotated by $180$ degrees, so that the part with the smaller radius is 
close to the $x_2$--axis. During the evolution that parts thins
even more and approaches the axis. 

\begin{figure}
\center
\includegraphics[angle=-90,width=0.3\textwidth]{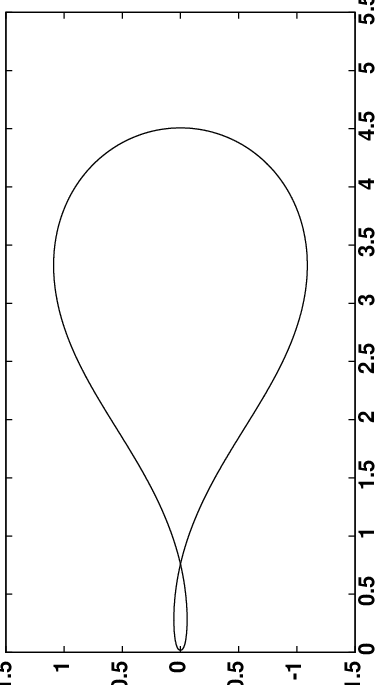}
\includegraphics[angle=-90,width=0.3\textwidth]{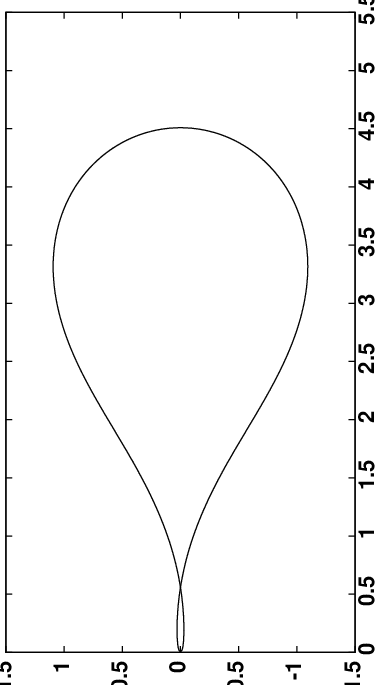}
\includegraphics[angle=-90,width=0.3\textwidth]{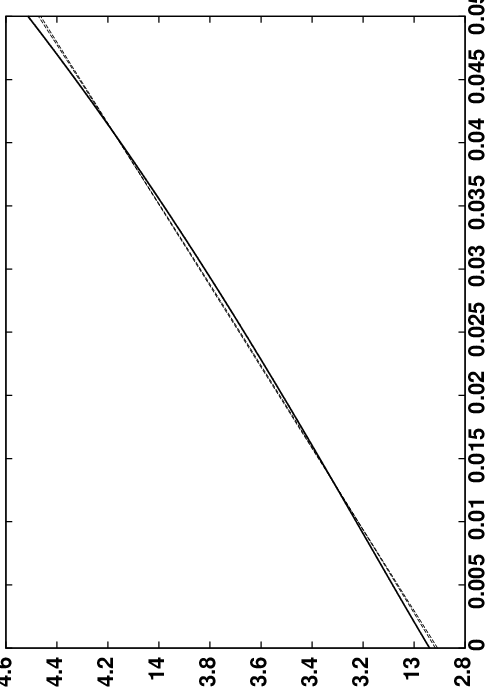}
\caption{$(\BGNpwf^m)^h$
Willmore flow for ${\rm T}(\Gamma) = 0$.
The solution at times $t=0,0.05$.
On the right a plot of $|\Gamma^m|_{\bH^2,h}$, recall \eqref{eq:B:unboundedh}, 
over time, together with the best fits $f(t) = a\,(1+t)^b$
and $g(t) = c + d\,t$, with $(a,b,c,d) = (12.9, 2.3, 12.9, 31.1)$.
}
\label{fig:B:wf_leftlemniscatefine}
\end{figure}

To investigate the behaviour close to the $x_2$--axis further, and to
help decide whether the evolution reaches a singularity in finite time, 
we start a refined computation from the solution at time $t=0.3$ in
Figure~\ref{fig:B:wf_leftlemniscate}. 
The results for $J = 4096$ and $\ttau = 10^{-8}$ are shown in 
Figure~\ref{fig:B:wf_leftlemniscatefine}.
The discrete quantity \eqref{eq:B:unboundedh}, which we
plot on the right of Figure~\ref{fig:B:wf_leftlemniscatefine}, appears to grow
polynomial in time. In particular, it seems to grow only slightly faster than
the best fits of the form $f(t) = a\,(1+t)^b$
and $g(t) = c + d\,t$, where for the best fits we observe
$(a,b,c,d) = (12.9, 2.3, 12.9, 31.1)$.
It is therefore difficult to draw definite conclusions. Hence we
conjecture that the singularity, where the toroidal surface closes up at the
origin, is reached only as $t\to\infty$.
This conjecture is in agreement with similar conclusion drawn by the authors
for the onset of a singularity for a genus-0 surface that converges to two
touching spheres. In particular, we stress that the evolution shown
in Figure~24 in \cite{axisd} closely matches the shape of the curves 
in Figure~\ref{fig:B:wf_leftlemniscatefine} near the $x_2$--axis.

\begin{figure}
\center
\includegraphics[angle=-90,width=0.3\textwidth]{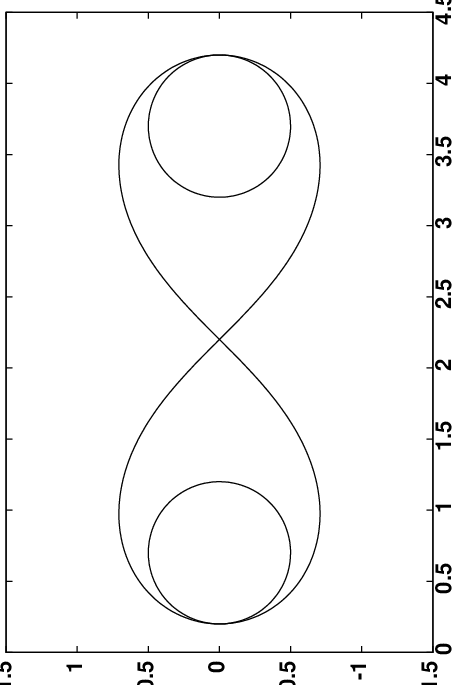}
\includegraphics[angle=-90,width=0.3\textwidth]{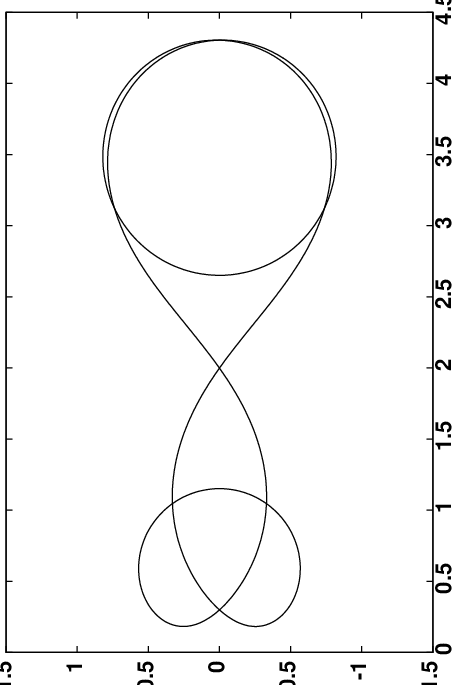}
\includegraphics[angle=-90,width=0.3\textwidth]{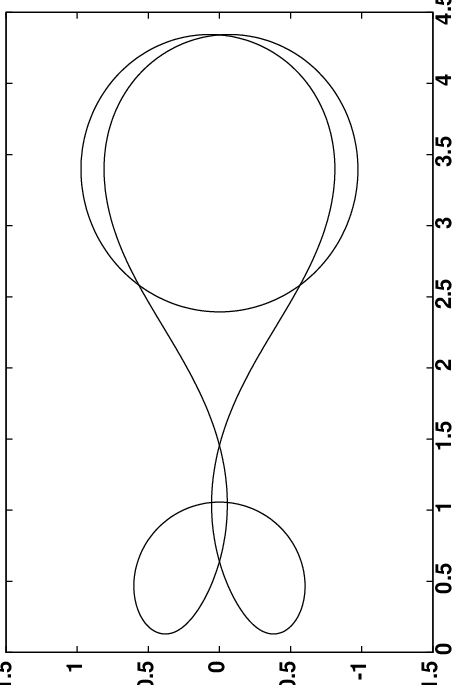}
\includegraphics[angle=-90,width=0.3\textwidth]{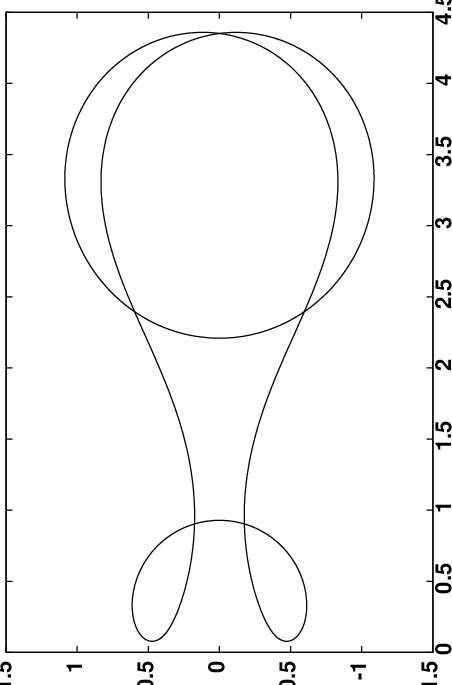}
\includegraphics[angle=-90,width=0.3\textwidth]{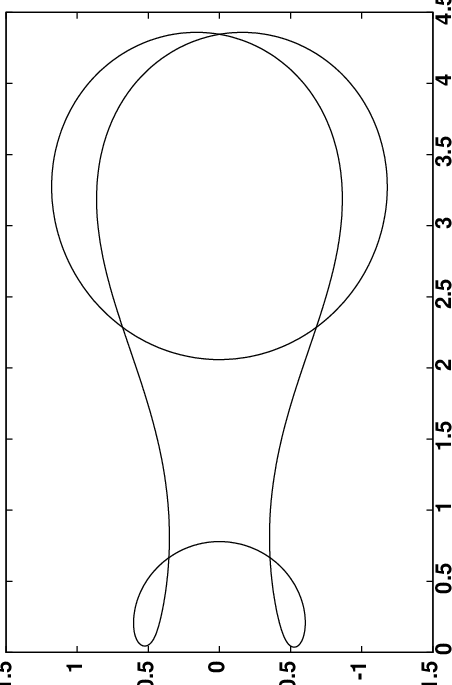}
\includegraphics[angle=-90,width=0.3\textwidth]{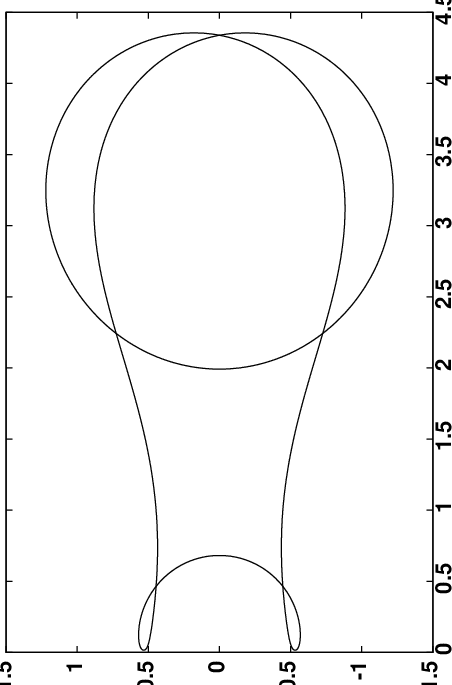}
\caption{$(\BGNpwf^m)^h$
Willmore flow for ${\rm T}(\Gamma) = 0$.
We show the solution at times $t=0,0.2,\ldots,0.8$ and $0.9$.
Here the evolution on $[0.8,0.9]$ was computed with the finer parameters
$J = 4096$ and $\ttau = 10^{-7}$.
}
\label{fig:B:wf_lemniscate2c}
\end{figure}

Finally, we also show a more interesting evolution for another initial data 
with ${\rm T}(\Gamma) = 0$. Here a symmetric lemniscate is inscribed with a
circle on each side. The Willmore flow for a torus with such a generating
curve can be seen in Figure~\ref{fig:B:wf_lemniscate2c}. 
The discretization parameters are $J = 1024$ and $\ttau = 10^{-5}$,
and for the final evolution on the time interval $[0.8,0.9]$ 
we use the finer parameters $J = 4096$ and $\ttau = 10^{-7}$.
We see that the circle
inscribed on the right of the lemniscate grows, while the circle on the left
untangles to create two loops close to the $x_2$--axis with large curvature.
The discrete quantity \eqref{eq:B:unboundedh} once again appears to grow
slightly faster than the best polynomial fits, with our best fits for the final
evolution in Figure~\ref{fig:B:wf_lemniscate2c} given by
$f(0.8+t) = a\,(1+t)^b$
and $g(0.8+t) = c + d\,t$, with
$(a,b,c,d) = (22, 1.5, 22, 33)$.

Overall, based on our numerical evidence, we conjecture the following:
If ${\rm T}(\Gamma(0)) = \alpha \not=0$, 
then Willmore flow will converge to a $|\alpha|$-covering of the Clifford 
torus. If ${\rm T}(\Gamma(0)) = 0$, then Willmore flow will develop a
singularity/blow-up as $t\to\infty$.

\end{appendix}

\noindent{\bf Acknowledgements}

\noindent The authors gratefully acknowledge the support 
of the Regensburger Universit\"atsstiftung Hans Vielberth.

\def\soft#1{\leavevmode\setbox0=\hbox{h}\dimen7=\ht0\advance \dimen7
  by-1ex\relax\if t#1\relax\rlap{\raise.6\dimen7
  \hbox{\kern.3ex\char'47}}#1\relax\else\if T#1\relax
  \rlap{\raise.5\dimen7\hbox{\kern1.3ex\char'47}}#1\relax \else\if
  d#1\relax\rlap{\raise.5\dimen7\hbox{\kern.9ex \char'47}}#1\relax\else\if
  D#1\relax\rlap{\raise.5\dimen7 \hbox{\kern1.4ex\char'47}}#1\relax\else\if
  l#1\relax \rlap{\raise.5\dimen7\hbox{\kern.4ex\char'47}}#1\relax \else\if
  L#1\relax\rlap{\raise.5\dimen7\hbox{\kern.7ex
  \char'47}}#1\relax\else\message{accent \string\soft \space #1 not
  defined!}#1\relax\fi\fi\fi\fi\fi\fi}


\begin{thebibliography}{10}

\bibitem{triplej}
J.~W. Barrett, H.~Garcke and R.~N\"urnberg, A parametric finite element method
  for fourth order geometric evolution equations. \textit{J. Comput. Phys.}
  \textbf{222} (2007) 441--462.

\bibitem{willmore}
J.~W. Barrett, H.~Garcke and R.~N\"urnberg, Parametric approximation of
  {W}illmore flow and related geometric evolution equations. \textit{SIAM J.
  Sci. Comput.} \textbf{31} (2008) 225--253.

\bibitem{pwftj}
J.~W. Barrett, H.~Garcke and R.~N\"urnberg, Elastic flow with junctions:
  Variational approximation and applications to nonlinear splines.
  \textit{Math. Models Methods Appl. Sci.} \textbf{22} (2012) 1250037.

\bibitem{pwf}
J.~W. Barrett, H.~Garcke and R.~N\"urnberg, Parametric approximation of
  isotropic and anisotropic elastic flow for closed and open curves.
  \textit{Numer. Math.} \textbf{120} (2012) 489--542.

\bibitem{pwfade}
J.~W. Barrett, H.~Garcke and R.~N\"urnberg, Computational parametric {W}illmore
  flow with spontaneous curvature and area difference elasticity effects.
  \textit{SIAM J. Numer. Anal.} \textbf{54} (2016) 1732--1762.

\bibitem{pwfopen}
J.~W. Barrett, H.~Garcke and R.~N\"urnberg, Stable variational approximations
  of boundary value problems for {W}illmore flow with {G}aussian curvature.
  \textit{IMA J. Numer. Anal.} \textbf{37} (2017) 1657--1709.

\bibitem{axisd}
J.~W. Barrett, H.~Garcke and R.~N\"urnberg, Finite element methods for fourth
  order axisymmetric geometric evolution equations. \textit{J. Comput. Phys.}
  \textbf{376} (2019) 733--766.

\bibitem{hypbolpwf}
J.~W. Barrett, H.~Garcke and R.~N\"urnberg, Stable discretizations of elastic
  flow in {R}iemannian manifolds. \textit{SIAM J. Numer. Anal.} \textbf{57}
  (2019) 1987--2018.

\bibitem{aximcf}
J.~W. Barrett, H.~Garcke and R.~N\"urnberg, Variational discretization of
  axisymmetric curvature flows. \textit{Numer. Math.} \textbf{141} (2019)
  791--837.

\bibitem{hypbol}
J.~W. Barrett, H.~Garcke and R.~N\"urnberg, Numerical approximation of curve
  evolutions in {R}iemannian manifolds. \textit{IMA J. Numer. Anal.}
  \textbf{40} (2020) 1601--1651.

\bibitem{bgnreview}
J.~W. Barrett, H.~Garcke and R.~N\"urnberg, Parametric finite element
  approximations of curvature driven interface evolutions. In A.~Bonito and
  R.~H. Nochetto (Eds.), \textit{Handb. Numer. Anal.}, Elsevier, Amsterdam,
  vol.~21, 2020 275--423.

\bibitem{BobenkoS05}
A.~I. Bobenko and P.~Schr\"{o}der, Discrete {W}illmore flow. In J.~Fujii (Ed.),
  \textit{ACM SIGGRAPH 2005 Courses}, ACM, New York, NY, SIGGRAPH '05, 5--es.

\bibitem{Canham70}
P.~B. Canham, The minimum energy of bending as a possible explanation of the
  biconcave shape of the human red blood cell. \textit{J. Theor. Biol.}
  \textbf{26} (1970) 61--81.

\bibitem{CapovillaGS02}
R.~Capovilla, J.~Guven and J.~A. Santiago, Lipid membranes with an edge.
  \textit{Phys. Rev. E} \textbf{66} (2002) 021607.

\bibitem{ClarenzDDRR04}
U.~Clarenz, U.~Diewald, G.~Dziuk, M.~Rumpf and R.~Rusu, A finite element method
  for surface restoration with smooth boundary conditions. \textit{Comput.
  Aided Geom. Design} \textbf{21} (2004) 427--445.

\bibitem{CoxL15}
G.~Cox and J.~Lowengrub, The effect of spontaneous curvature on a two-phase
  vesicle. \textit{Nonlinearity} \textbf{28} (2015) 773--793.

\bibitem{DallAcquaMSS20preprint}
A.~Dall'Acqua, M.~M\"uller, R.~Sch\"atzle and A.~Spener, The {W}illmore flow of
  tori of revolution. arXiv:2005.13500, 2020.

\bibitem{Davis04}
T.~A. Davis, Algorithm 832: {UMFPACK} {V}4.3---an unsymmetric-pattern
  multifrontal method. \textit{ACM Trans. Math. Software} \textbf{30} (2004)
  196--199.

\bibitem{DeckelnickD09}
K.~Deckelnick and G.~Dziuk, Error analysis for the elastic flow of parametrized
  curves. \textit{Math. Comp.} \textbf{78} (2009) 645--671.

\bibitem{DeckelnickS10}
K.~Deckelnick and F.~Schieweck, Error analysis for the approximation of
  axisymmetric {W}illmore flow by {$C^1$}-finite elements. \textit{Interfaces
  Free Bound.} \textbf{12} (2010) 551--574.

\bibitem{Dziuk08}
G.~Dziuk, Computational parametric {W}illmore flow. \textit{Numer. Math.}
  \textbf{111} (2008) 55--80.

\bibitem{Germain1821}
S.~Germain, \textit{Recherches sur la th{\'e}orie des surfaces {\'e}lastiques}.
  Veuve Courcier, Paris, 1821.

\bibitem{Helfrich73}
W.~Helfrich, Elastic properties of lipid bilayers: {T}heory and possible
  experiments. \textit{Z. Naturforsch. C} \textbf{28} (1973) 693--703.

\bibitem{JulicherL96}
F.~J{\"u}licher and R.~Lipowsky, Shape transformations of vesicles with
  intramembrane domains. \textit{Phys. Rev. E} \textbf{53} (1996) 2670--2683.

\bibitem{JulicherS94}
F.~J\"ulicher and U.~Seifert, Shape equations for axisymmetric vesicles: {A}
  clarification. \textit{Phys. Rev. E} \textbf{49} (1994) 4728--4731.

\bibitem{Kirchhoff1850}
G.~R. Kirchhoff, {\"U}ber das {G}leichgewicht und die {B}ewegung einer
  elastischen {S}cheibe. \textit{J. Reine Angew. Math.} \textbf{40} (1850)
  51--88.

\bibitem{Kuhnel15}
W.~K\"uhnel, \textit{Differential Geometry: {C}urves -- {S}urfaces --
  {M}anifolds}, vol.~77 of \textit{Student Mathematical Library}. Amer. Math.
  Soc., Providence, RI, 2015.

\bibitem{MarquesN14}
F.~C. Marques and A.~Neves, Min-max theory and the {W}illmore conjecture.
  \textit{Ann. of Math.} \textbf{179} (2014) 683--782.

\bibitem{MayerS02}
U.~F. Mayer and G.~Simonett, A numerical scheme for axisymmetric solutions of
  curvature-driven free boundary problems, with applications to the {W}illmore
  flow. \textit{Interfaces Free Bound.} \textbf{4} (2002) 89--109.

\bibitem{Nitsche93}
J.~C.~C. Nitsche, Boundary value problems for variational integrals involving
  surface curvatures. \textit{Quart. Appl. Math.} \textbf{51} (1993) 363--387.

\bibitem{Poisson1814}
S.~D. Poisson, M{\'e}moire sur les surfaces {\'e}lastiques.
  \textit{M{\'e}moires de l'Institut 1812} \textbf{9} (1814) 167--226.

\bibitem{Rusu05}
R.~E. Rusu, An algorithm for the elastic flow of surfaces. \textit{Interfaces
  Free Bound.} \textbf{7} (2005) 229--239.

\bibitem{Seifert97}
U.~Seifert, Configurations of fluid membranes and vesicles. \textit{Adv. Phys.}
  \textbf{46} (1997) 13--137.

\bibitem{TuO-Y03}
Z.~C. Tu and Z.~C. Ou-Yang, Lipid membranes with free edges. \textit{Phys. Rev.
  E} \textbf{68} (2003) 061915.

\bibitem{WangD08}
X.~Wang and Q.~Du, Modelling and simulations of multi-component lipid membranes
  and open membranes via diffuse interface approaches. \textit{J. Math. Biol.}
  \textbf{56} (2008) 347--371.

\bibitem{Willmore93}
T.~J. Willmore, \textit{Riemannian Geometry}. Oxford Science Publications, The
  Clarendon Press, Oxford University Press, New York, 1993.

\end{thebibliography}
\end{document}